%% file: local_law.tex
\documentclass[10pt]{extarticle}
\usepackage{psfrag}
\usepackage[parfill]{parskip}
\usepackage{float, graphicx}
\usepackage[]{epsfig}
\usepackage{amsmath, amsthm, amssymb,}
\usepackage{epsfig}
\usepackage{verbatim}
\usepackage{multicol}
\usepackage{url}
\usepackage{appendix}
\usepackage[colorlinks, bookmarks=true]{hyperref}
\usepackage{enumerate}
\usepackage[normalem]{ulem}
\usepackage{bm}
\usepackage{dsfont}
\usepackage{stmaryrd}
\usepackage{enumitem}
\usepackage{cite}
\usepackage{color}
\usepackage{xcolor}
\usepackage{authblk}
\usepackage{geometry}
\usepackage{imakeidx}
\makeindex
\newcommand{\xindex}[2]{\expandafter\index\expandafter{\pageref{#1}#2}}
\makeatletter

\def\@settitle{\begin{center}%
    \bfseries
 \normalfont\LARGE\@title
  \end{center}%
}
\def\@setauthors{\begin{center}%
 \normalsize\@author
  \end{center}%
}

\makeatother

\numberwithin{equation}{section}

\newcommand{\rd}{{\rm d}}

\newcommand{\pbb}[1]{\biggl({#1}\biggr)}

\def\cA{{\mathcal A}}
\def\cB{{\mathcal B}}

\def\cH{{\mathcal H}}
\def\cA{{\mathcal A}}

\def\sA{{\sf A}}
\def\sB{{\sf B}}

\def\sJ{{\sf U}}

\renewcommand{\frak}{\mathfrak}
\newcommand{\fa}{{\frak a}}
\newcommand{\fb}{{\frak b}}
\newcommand{\fc}{{\frak c}}
\newcommand{\fd}{{\frak q}}
\newcommand{\fe}{{\frak e}}

\newcommand{\fR}{{\frak R}}
\newcommand{\fr}{{\frak r}}

\newcommand{\fC}{{\frak C}}

\newcommand{\cT}{{\mathcal T}}

\newcommand{\al}{\alpha}

\newcommand{\be}{\begin{equation}}
\newcommand{\ee}{\end{equation}}

\newcommand{\e}{{\varepsilon}}

\newcommand{\la}{\lambda}

\newcommand{\td}{\tilde}

\newcommand{\cX}{{\cal X}}
\newcommand{\cG}{{\cal G}}

\newcommand{\cC}{\mathbb S}
\newcommand{\bI}{\mathbb I}

\newcommand{\bL}{{\mathbb L}}

\newcommand{\cY}{{\cal Y}}

\newcommand{\bW}{\mathbb W}
\newcommand{\bC}{{\mathbb C}}

\newcommand{\T}{\mathbb T}

\newcommand{\bT}{\T}

\newcommand{\bB}{{\mathbb B}}
\newcommand{\bE}{{\mathbb E}}
\newcommand{\bX}{{\mathbb X}}
\newcommand{\bG}{{\mathbb G}}

\newcommand{\bH}{{\mathbb H}}
\newcommand{\bV}{{\mathbb V}}

\newcommand{\bY}{\mathbb{Y}}
\newcommand{\bP}{\mathbb{P}}

\newcommand{\deq}{\mathrel{\mathop:}=}

\newcommand{\id}{\mspace{2mu}\mathrm{i}\mspace{-0.6mu}\mathrm{d}} 
\renewcommand{\epsilon}{\varepsilon}
\renewcommand{\leq}{\leqslant}
\renewcommand{\geq}{\geqslant}

\renewcommand{\le}{\leq}
\renewcommand{\ge}{\geq}
\renewcommand{\P}{\mathbb{P}}
\newcommand{\E}{\mathbb{E}}
\newcommand{\R}{\mathbb{R}}

\newcommand{\Z}{\mathbb{Z}}

\def\bC{{\mathbb C}}

\newcommand{\tQ}{{\tilde Q}}
\newcommand{\bfS}{{\bf S}}

\newcommand{\qq}[1]{[\![{#1}]\!]}

\newcommand{\absb}[1]{\bigl\lvert #1 \bigr\rvert}

\DeclareMathOperator{\im}{Im}

\DeclareMathOperator{\OO}{O}
\DeclareMathOperator{\oo}{o}

\DeclareMathOperator{\diam}{diam}

\newcommand{\del}{\partial}

\newcommand{\Tr}{\operatorname{Tr}}

\newcommand{\dist}{\mathrm{dist}}
\newcommand{\I}{Q}

\newcommand{\bS}{\mathbb{S}}
\newcommand{\cGT}{\cal G^{(\T)}}

\newcommand{\tcG}{\tilde{\cal G}}

\newcommand{\tcGT}{\tilde{\cal G}^{(\T)}}

\newcommand{\GT}{G^{(\T)}}
\newcommand{\tG}{\tilde{G}}
\newcommand{\tGT}{\tilde{G}^{(\T)}}

\newcommand{\tP}{\tilde{P}}
\newcommand{\ta}{\tilde{a}}

\newcommand{\sG}{{\sf G}}
\newcommand{\sS}{{\sf S}}
\newcommand{\GNd}{\sG_{N,d}}
\newcommand{\GNdp}{\tilde \sG_{N,d}}
\newcommand{\Pp}{\tilde \P}

\newcommand{\n}{\omega}
\newcommand{\msc}{m_{sc}}
\newcommand{\md}{m_d}

\renewcommand{\cal}{\mathcal}

\newcommand{\rn}[1]{%
       \lowercase\expandafter{\romannumeral#1}%
}

\newcommand{\RN}[1]{%
       \uppercase\expandafter{\romannumeral#1}%
}

\renewcommand{\Im}{{\mathrm{Im}}}
\renewcommand{\Re}{\mathrm{Re}}

\theoremstyle{plain} 

\makeatletter
\def\thm@space@setup{%
  \thm@preskip=4pt plus 2pt minus 2pt
}
\makeatother

\newtheorem{theorem}{Theorem}[section]
\newtheorem*{theorem*}{Theorem}
\newtheorem{lemma}[theorem]{Lemma}
\newtheorem*{lemma*}{Lemma}

\newtheorem*{corollary*}{Corollary}
\newtheorem{proposition}[theorem]{Proposition}
\newtheorem*{proposition*}{Proposition}

\newtheorem*{assumption*}{Assumption}
\newtheorem{claim}[theorem]{Claim}

\newtheorem{definition}[theorem]{Definition}
\newtheorem*{definition*}{Definition}

\newtheorem*{example*}{Example}
\newtheorem{remark}[theorem]{Remark}

\newtheorem*{remark*}{Remark}
\newtheorem*{remarks*}{Remarks}

\theoremstyle{definition}

\newcommand{\TE}{\mathrm{TE}}
\newcommand{\Ext}{\mathrm{Ext}}
\newcommand{\ri}{\mathrm{i}}

\newcommand{\Ba}{{\mathsf B}_1}
\newcommand{\Bb}{{\mathsf B}_2}
\newcommand{\Bc}{{\mathsf B}_3}
\newcommand{\cell}{{\mathsf{Cr}}}
\newcommand{\cella}{{\mathbb{S}}}
\newcommand{\As}{{\mathsf W}}

\title{Spectrum of Random $d$-regular Graphs Up to the Edge}
  \author[1]{Jiaoyang~Huang\thanks{huangjy@wharton.upenn.edu}}
  \author[2]{Horng-Tzer Yau\thanks{htyau@g.harvard.edu}}

\affil[1]{University of Pennsylvania}
\affil[2]{Harvard University}

 \date{}
\begin{document}

\maketitle



\begin{abstract}
 Consider  the normalized adjacency matrices of random $d$-regular graphs on $N$  vertices  with  fixed degree $d\geq3$.  We prove that,  with probability  $1-N^{-1+\varepsilon}$  for any $\varepsilon >0$,
 the following two properties hold
as $N \to \infty$ provided that $d\geq3$:
(i) The eigenvalues are close to the classical eigenvalue locations given by the Kesten-McKay distribution. In particular, the extremal eigenvalues are concentrated with polynomial error bound in $N$, i.e. $\lambda_2, |\lambda_N|\leq 2+N^{-c}$. (ii)  All  eigenvectors of  random $d$-regular graphs are completely delocalized.
\end{abstract}


{
  \hypersetup{linkcolor=black}
  \tableofcontents
}

\section{Introduction and Main Results}
\subsection{Introduction}
\label{sec:intro-intro} 
The random $d$-regular graph ensemble  is defined  to be the set     $\GNd$\index{$\mathsf G_{N,d}$} of  (simple) $d$-regular graphs on $N$ vertices equipped with  uniform probability.
It is a  fundamental
model of sparse random graphs 
and it  arises naturally in many different  contexts. %
The spectral properties of the adjacency matrices of  random $d$-regular graphs,  i.e. eigenvalues and eigenvectors,  are  of particular interest
in computer science,  combinatorics, and statistical physics. The relevant topics include 
the theory of expanders (see e.g.\ \cite{MR2072849}), quantum chaos (see e.g.\ \cite{MR3204183}),
and  graph $\zeta$-functions (see e.g.\ \cite{MR2768284}).

Throughout the paper, we fix $d\geq 3$, and denote $A=A(\cG)$  the adjacency matrix of a  random $d$-regular graph $\cG$ on $N$ vertices.
Thus $A$ is uniformly chosen among all symmetric $N\times N$ matrices with entries in $\{0,1\}$ with
$\sum_j A_{ij}=d$ and $A_{ii}=0$ for all $i$.
We normalize the adjacency matrix as $H=A/\sqrt{d-1}$, and denote its eigenvalues by $\la_1=d/\sqrt{d-1}\geq \la_2\geq \cdots\geq \la_N$.
The constant vector is a trivial eigenvector to $A$  with eigenvalue $d$. The adjacency matrices for random $d$-regular graphs are thus random matrices, 
albeit  the matrix entries  are not independent.

The spectral properties of random matrices with independent entries are well-understood. For example, the global and local spectral statistics of the generalized Wigner matrices
were analyzed in 
\cite{MR1810949,MR2810797,MR2639734,MR2662426,MR2981427,MR2784665,MR3372074,MR3541852,MR2964770, MR3699468,MR3034787,MR2930379,BY2016,MR3914908}.
These results were  extended to sparse matrices and  Erd\H{o}s-R\'enyi random graphs \cite{MR3098073,MR2964770,MR3429490,MR4058984,MR3800840, huang2020transition, he2019local,he2019bulk,he2021fluctuations}.  In addition to eigenvalue statistics,  delocalization of eigenvectors are often established provided the Green's function method was used. 

 Since matrix entries of random $d$-regular graphs  are  dependent, much less was known regarding  their spectral statistics.
If the degree $d$ is fixed, macroscopic eigenvalue statistics  were  studied
using the techniques of Poisson approximation of short cycles \cite{MR3078290,MR3315475} and    a (nonrigorous)  replica method \cite{PhysRevE.90.052109}.
These results showed that the macroscopic eigenvalue statistics for random $d$-regular graphs of fixed degree
are different from those of a Gaussian  ensemble. 
 The local eigenvalue statistics, however,  were conjectured to be the same as the Gaussian Orthogonal Ensemble statistics \cite{MR1691538,MR2433888,MR2072849}.
 It was known that 
the spectral density of   $H$  converges  to  the Kesten--McKay distribution\index{$\rho_{d}$}
with  density  given by
\begin{align}\label{kmlaw}
\rho_{d}(x)=\left(1+\frac{1}{d-1}-\frac{x^2}{d}\right)^{-1}\frac{\sqrt{[4-x^2]_+}}{2\pi},
\end{align}
 on spectral scales $(\log N)^{-c}$ \cite{MR2999215,MR3025715,MR3433288,MR3322309}. 
When the degree $d$ grows faster than  sufficiently  high  power of $\log N$, the eigenvalue rigidity down to the  scale $(\log N)^4/N$ (notice that, compared with previous works, the scale was greatly improved to essentially the optimal one) was established in \cite{MR3688032}. Using this result as an input, bulk universality and the normality of bulk eigenvectors in the same regime was proven in \cite{MR3690289,MR3729611}, and edge universality for degrees $d\in [N^{1/3+\varepsilon}, N^{2/3-\varepsilon}]$ was proven in \cite{bauerschmidt2020edge}.
In a joint paper with R. Bauerschmidt \cite{MR3962004}, we extended the  convergence of spectral density  to scale $(\log N)^{c}/N$ when   $d\geq d_0$ for some large $d_0$ fixed.
Spectral properties of  directed $d$-regular graphs have also been studied recently \cite{MR4029149,MR3602844,litvak2020circular,huang2021invertibility,AM}.

The second largest eigenvalue $\lambda_2$ of $d$-regular graphs is of particular interest in  theoretical computer science and combinatorics \cite{MR1421568,MR2247919}. The spectral gap, the gap between the first and second eigenvalues, measures the expanding property of the graph.
In \cite{MR963118}, Lubotzky, Phillips, and Sarnak defined \emph{Ramanujan graphs}, 
as $d$-regular graphs with all  non-trivial eigenvalues of $H$ bounded by $2$. 
 Since  $\lim_{N\rightarrow \infty} \la_2\geq 2$ for any deterministic family of $d$ regular graphs on $N$ vertices \cite{MR875835}, Ramanujan graphs asymptotically have the largest possible spectral gap among  $d$-regular graphs.   Ramanujan graphs have been a focus of substantial works  in theoretical computer science \cite{rosenthal2000constructions, charles2009cryptographic,kar2008topology}
and mathematics \cite{murty2003,terras2010zeta}.
Proving the existence of Ramanujan graphs with a
large number of vertices is however a difficult task,  which was only recently  solved for arbitrary $d\geq 3$ \cite{MR3374962, MR3374963}. 
On the other end, it was conjectured  by Alon \cite{MR875835}  and proven by Friedman \cite{MR2437174}  and  Bordenave  \cite{Bord15}  that with high probability   random
$d$-regular graphs are weakly Ramanujan, i.e. $\la_2, |\la_N|\leq 2+\oo(1)$. 
Their proofs are based on sophisticated moment methods, and the error term
in eigenvalue  is of order $\OO((\log \log N/ \log N)^2)$. 

 It was conjectured that the distribution of the second largest eigenvalue after normalizing by $N^{2/3}$ is the
same as that of the largest eigenvalue of the Gaussian Orthogonal Ensemble \cite{MR2072849,MR2433888};  i.e.,  there is a constant $c_{N, d}$ so that  $ N^{2/3} (\lambda_2- 2)- c_{N, d} $ has the Tracy-Widom distribution \cite{tracy1996orthogonal}. 
 This would imply that the fluctuations of extreme eigenvalues are of order $\OO(N^{-2/3})$ if $c_{N, d} $ is of order one. 
 If $c_{N, d}=0 $ (which seems to be the most probable scenario), then it would imply that  more than half of  all $d$-regular graphs are Ramanujan graphs,
namely $d$-regular graphs with $\lambda_2 \leq 2$. 
As a first step towards this conjecture, 
in this article we prove a polynomial error bound for extremal eigenvalues: with high  probability $\la_2, |\la_N|\leq 2+\OO(N^{-c})$ for some $c>0$. This gives a new proof of Alon's conjecture \cite{MR875835} with  an $\OO(N^{-c})$  error bound instead of $\OO((\log \log N/ \log N)^2)$ errors in  \cite{MR2437174, Bord15}.   In addition to this estimate on the second largest eigenvalue, our results imply estimates on all eigenvalues in the spectrum, see Theorem \ref{thm:evloc}.

 As discrete analogues of compact negatively curved surfaces \cite{MR3038543}, $d$-regular graphs provide a good ground to explore the phenomenon of  delocalization of eigenvectors. 
For $d$-regular graphs with locally tree-like structure,  it was known that 
the eigenvectors $v$ of $H$ are weakly delocalized   in the sense that
their entries are uniformly bounded by $(\log N)^{-c}\|v\|_2$ \cite{MR3025715,MR3433288,MR3038543} and their $\ell^2$-mass cannot concentrate on a small set \cite{MR3038543}.
If, in addition, the graphs are expanders,  
the eigenvectors of $H$ also satisfy the quantum ergodicity property \cite{MR3322309,MR3567266,MR3649482}. 
In  \cite{MR3962004}, we proved that the bulk eigenvectors of random $d$-regular graphs are completely delocalized with high probability  for $d\geq d_0$ for some large $d_0$ fixed.  In this paper, we extend  this result  to  any $d\geq 3$ and to all   eigenvectors including edge eigenvectors.  Our results  imply in particular that 
the eigenvectors of random $3$-regular graphs are completely delocalized with probability  $1-N^{-1+\varepsilon}$. 
We will explain later on in this section that the error probability  $1-N^{-1+\varepsilon}$ is optimal in the sense that there are events with probability  $N^{-1+\varepsilon}$ 
so that some eigenvectors are  localized. 
Our  proof will also  show that the existence of too many cycles in a small neighborhood is the main reason for the localization of eigenvectors.

All our results on both eigenvalues and eigenvectors  will be derived from a high probability estimate on the Green's function $G(z) = ( H - z)^{-1}$ of the normalized adjacency matrix
for random $d$-regular graphs. The  trace of the Green's function gives the Stieltjes transform of empirical eigenvalue distribution, which encodes the information of eigenvalue locations.  Delocalization of all eigenvectors  follows from  estimates on the entries of the Green's function provided that the imaginary part of the spectral parameter can be chosen to be nearly optimal, i.e.,  $\im z \gtrsim N^{-1 + \e}$. The recent self-consistent Green's  function method \cite{MR2871147,MR2481753,MR2537522, erdos2017dynamical,MR3800840},  invented for Wigner matrices,  was primarily applied to matrices with (nearly) independent entires. 
It is a powerful tool to obtain spectral information for the full range of spectrum. Our key  observation   is to construct  a  novel  functional of Green's function satisfying   a simple self-consistent equation.  This allows us to derive estimates on Green's function so that all results in this paper 
for random regular graphs with degree up to the optimal $d\geq 3$  are simple corollary of these estimates.    The correlations of entries in $H$ will cause major difficulties; these will be explained later in this  section. 

\subsection{Statements of main results: eigenvalue statistics}
\label{sec:intro}

In this paper we fix a parameter $0<\fc<1$\index{$\mathfrak c$} (e.g. $\fc=0.99$). Let $\fR=(\fc/4)\log_{d-1}N$\index{$\mathfrak R$}. 
For each $d\geq 3$, we associate  an integer $\omega_d>0$ (see Definition \ref{d:defwd} for the precise definition)  so that, roughly speaking, an infinite $d$-regular graph with excess (the smallest number of edges that must be removed to yield a graph with no cycles) less than or equal to  $\omega_d$ has good spectral properties. 

\begin{definition}\label{def:barOmega}
We define the  set  of 
 radius-$\fR$ tree like graphs,  $\bar\Omega\subset \GNd$\index{$\bar\Omega$},   to consist of graphs such that
\begin{itemize}
\item
the radius-$\fR$ neighborhood of any vertex has excess at most $\n_d$; 
\item
the number of vertices whose radius-$\fR$ neighborhood that contains a cycle is at most $N^{\fc}$.
\end{itemize}
\end{definition}

Our main results,  to be stated in the rest of this section, assert   that   on the set $\bar \Omega$ with probability $1-N^{-\fC}$ for any $\fC> 0$, 
  the eigenvalues of any $d$-regular graphs are near their  classical locations and the corresponding  eigenvectors are completely delocalized.

\begin{theorem}[Eigenvalue Rigidity]\label{thm:evloc}
Fix $d\geq 3$, $0<\fc<1$, $\fR=(\fc/4)\log_{d-1}N$ and recall the set of  radius-$\fR$ tree like graphs $\bar\Omega\subset \GNd$ from Definition \ref{def:barOmega}.  For any large $\fC>0$ and $N$ large enough, with probability $1-\OO(N^{-\fC})$ with respect to the uniform measure on $\bar\Omega$,  the eigenvalues of $\cG\in \bar\Omega$ satisfy:  
\begin{equation}\label{e:eigenloc}
|\lambda_i-\gamma_i|\leq N^{-\Omega(1)}, \quad 2\leq i\leq N,
\end{equation}
where $\gamma_i$ are classical eigenvalue locations given by the Kesten-McKay distribution \eqref{kmlaw}
\begin{align*}
\int_{\gamma_i}^2 \rho_d(x)\rd x=\frac{i}{N}, \quad 2\leq i\leq N.
\end{align*}
\end{theorem}

Thanks to Proposition \ref{prop:structure}, most $d$-regular graphs are locally tree-like
according to Definition \ref{def:barOmega} in the sense that   $\bP(\bar\Omega)\geq 1-N^{-(1-\fc)\omega_d}$. Theorem \ref{thm:evloc} thus implies that with high probability, the non-trivial 
extremal  eigenvalues of random $d$-regular graphs concentrate around $\pm2$.  More precisely, we have the following estimates on the extremal eigenvalues. 
We remark  that   the exponent  $\Omega(1)$ 
in the following   theorem can be chosen to   be $0.01$ by inspecting our proof.

\begin{theorem}[Extremal  eigenvalues]\label{thm:extremeev}
Fix $d\geq 3$ and $0<\fc<1$. There exists a positive integer $\omega_d$, defined in Definition \ref{d:defwd},  such that 
\begin{equation}\label{e:extremeev}
 \lambda_2, |\lambda_N|\leq 2+N^{-\Omega(1)},
\end{equation}
with probability $1-\OO(N^{-(1-\fc)\omega_d})$  with respect to the uniform measure on $\GNd$.   
\end{theorem}
%
%


\subsection{Delocalization of Eigenvectors}

For $d$-regular graphs with locally tree-like structure,
the eigenvectors $v$ of $H$ cannot concentrate on a small set, in the sense that any vertex set $\bV \subset \qq{N}$ with $\sum_{i\in \bV} |v_i|^2 \geq \varepsilon \|v\|_2$
must have at least $N^{c(\varepsilon)}$ elements \cite{MR3038543}.
Moreover, for deterministic locally tree-like  $d$-regular expander graphs,
it was proved that the eigenvectors $v$ satisfy a quantum ergodicity property:
for all $a \in \R^N$ with $\|a\|_\infty\leq 1$ and $\sum_i a_i=0$,
averages of $|\sum_i a_i v_i^2|^2$ over macroscopically  many eigenvectors $v$
are close to $0$ \cite{MR3322309,MR3649482,MR3567266}. 
Our main results  regarding eigenvectors of sparse $d$-regular graphs  are summarized in the following theorem:


\begin{theorem}[Eigenvector delocalization]\label{thm:delocalizationev}
Fix $d\geq 3$, $0<\fc<1$, $\fR=(\fc/4)\log_{d-1}N$ and recall the set of  radius-$\fR$ tree-like graphs $\bar\Omega\subset \GNd$ from Definition \ref{def:barOmega}.  For any large $\fC>0$ and $N$ large enough, with probability $1-\OO(N^{-\fC})$ with respect to the uniform measure on $\bar\Omega$,  the eigenvectors of $\cG\in \bar\Omega$ satisfy 
\begin{equation}\label{evest}
  \|v\|_\infty \leq 
\frac{(\log N)^{\Omega(1)}}{\sqrt{N}} \|v\|_2.
\end{equation}
\end{theorem}

\begin{figure}
\begin{center}
 \includegraphics[scale=0.4,trim={4cm 10.5cm 1cm 10.5cm},clip]{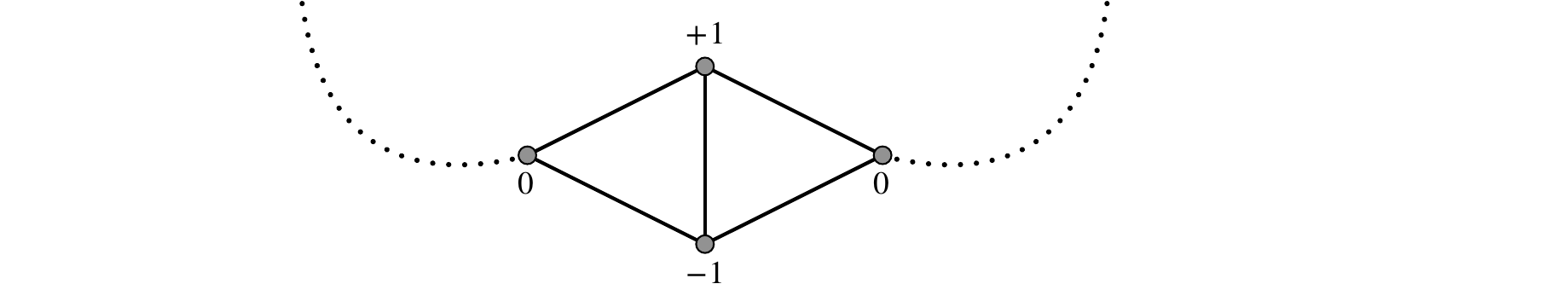}
\caption{
Theorem~\ref{thm:delocalizationev} shows
that a random $d$-regular graph has only completely delocalized eigenvectors with probability $1-\OO(N^{-(1-\oo(1))\omega_d})$.
On the other hand, it is not difficult to show that a random $d$-regular graph has localized eigenvectors with probability $\Omega(N^{-d+2})$. For example,
a random $3$-regular graph contains the subgraph shown above with probability $\Omega(N^{-1})$. Our results  imply that 
the eigenvectors of random $3$-regular graphs are completely delocalized with probability  $1-N^{-1+\oo(1)}$, which is optimal.
\label{fig:fork}}
\end{center}
\end{figure}

Thanks to Proposition \ref{prop:structure}, 
$\bP(\bar\Omega)\geq 1-N^{(-1+\fc)\omega_d}$. Theorem \ref{thm:evloc} thus implies that, with probability $1-\OO(N^{-(1-\fc)\omega_d})$  with respect to the uniform measure on $\GNd$, \eqref{evest} holds.
In Theorem \ref{thm:delocalizationev}, $\bar\Omega$  is the set of $d$-regular graphs, which do not have many cycles in a small neighborhood. 
The existence of too many cycles in a small neighborhood may result in localized eigenvectors.  For example, a random $3$-regular graph containing the subgraph in Figure \ref{fig:fork},  which happens with probability $\Omega(N^{-1})$,  has a localized eigenvector. Using $\omega_3=1$ from Proposition \ref{p:wd} and taking $\fc=\oo(1)$, we conclude from Theorem \ref{thm:delocalizationev} that
the eigenvectors of random $3$-regular graphs are completely delocalized with probability $1-N^{-1+\oo(1)}$, which is optimal. However, there do exist infinite families of locally tree-like $d$-regular graphs with localized eigenvectors, as constructed in \cite{alon2021high,ganguly2021non}. Finally, we remark that for random $d$-regular graphs of fixed degree,
a Gaussian wave correlation structure for the eigenvectors was predicted in \cite{0907.5065}
and partially confirmed in \cite{MR3945757}.

\subsection{Comments on key ideas}
\begin{figure}
\begin{center}
 \includegraphics[scale=0.3,trim={0cm 8cm 0 6cm},clip]{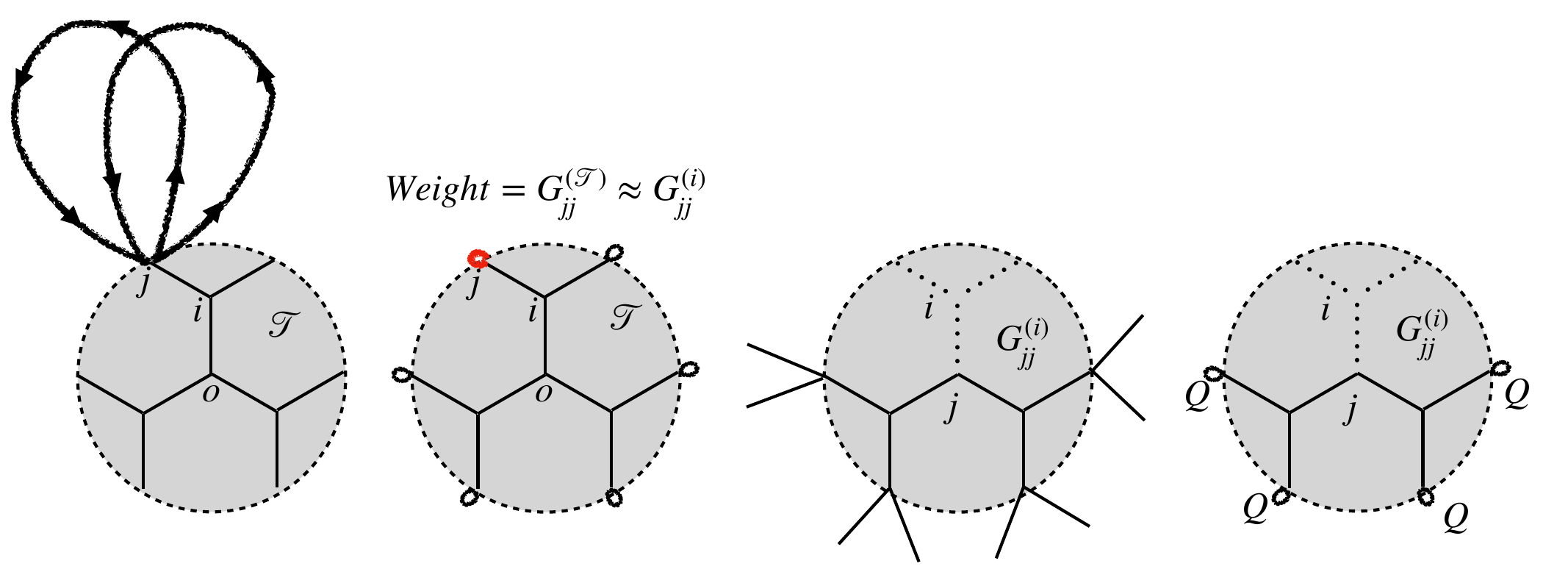}
 \caption{The sum of weighted paths from $j$ to itself staying outside $\cT$ is approximately given by $G_{jj}^{(i)}$. In average, $G_{jj}^{(i)}$ can be approximated by the Green's function at root vertex of a depth $\ell$ truncated $(d-1)$-ary tree with boundary weights $Q(\cG,z)$.}
 \label{f:sumpath}
 \end{center}
 \end{figure}

The Green's function of a graph encodes its spectral information. Eigenvalue rigidity and delocalization of eigenvectors follow from estimates of the Green's function down to optimal scale.
The main difficulty to  estimate Green's function near   the spectral edge is due to the square root singularity of the Kesten-McKay law near the spectral edge, which is a general phenomenon in random matrix theory. To overcome it,  instead of approximating Green's function by tree extensions as in \cite{MR3962004}, we approximate it by extensions  of graphs with certain boundary weights. Informally, the entries of the Green's function can be interpreted as a sum over weighted paths:
\begin{align}\label{e:weightsum}
G(z)=(H-z)^{-1}=-\frac{1}{z}-\frac{H}{z^2}-\frac{H^2}{z^3}-\frac{H^3}{z^4}\cdots,\quad z\in \bC^+.
\end{align}
The $i j$-th entry of $G$  is a sum over weighted paths from vertex $i$ to vertex $j$. For any vertex $o$, to compute $G_{oo}$, we need to sum over weighted paths from $o$ to itself. To do that we fix a radius-$\ell$ (which is of order $\log \log N$, see Definition \ref{cop}) neighborhood of $o$, and denote it by $\cT$ (which may not be a tree) and its vertex set $\bT$. Then paths from $o$ to itself either stay completely inside $\cT$, or leave at some boundary point of $\cT$ and come back. For random $d$-regular graphs, most vertices have large tree neighborhoods (which contain no cycles). The main contribution comes from those paths which leave $\cT$ at some boundary point, say vertex $j$ and come back at $j$ as in Figure \ref{f:sumpath}. The sum of such weighted paths is the sum  (over $j$)  of weighted paths from $j$ to itself staying outside $\cT$.  Denote by  $\cGT$ the graph $\cG$ with vertices in $\bT$ removed and $G^{(\bT)}$ its Green's function.  With these notations, the weighted paths from $j$ to itself staying outside $\cT$ is just  $G^{(\bT)}_{jj}$.    We view $\cT$ as a graph rooted at $o$, and denote  the parent vertex of $j$ by $i$ as in Figure \ref{f:sumpath}.
Since  most vertices in random $d$-regular graphs have large tree neighborhoods, the paths from $j$ to itself not containing $i$ are likely to stay outside $\cT$. We can further approximate $G^{(\bT)}_{jj}$ by $G_{jj}^{(i)}$, where $G^{(i)}$
is the Green's function of $\cG^{(i)}$,  the graph $\cG$ with  vertex $i$ removed. Therefore we can approximate $G_{oo}$ by the Green's function of $\cT$ with  weights $G_{jj}^{(i)}$ at 
the boundary vertex $j$. 


The quantity we will use to derive a self-consistent equation is the average of $G_{jj}^{(i)}$ for all pairs of adjacent vertices $i\sim j$:
\begin{align}\label{e:Qsum}
Q(\cG,z)=\frac{1}{Nd}\sum_{i\sim j}G_{jj}^{(i)}.
\end{align}
For $d$-regular graphs, the quantity $Q(\cG,z)$, although complicated, turns out to be a more fundamental object than the Stieltjes transform of the empirical measure 
of the eigenvalues of the matrix $H$. 
To compute $G_{jj}^{(i)}$, we can approximate it by the Green's function of a radius $\ell$  neighborhood of $j$ with vertex $i$ removed  and with a suitable weight at each boundary vertex.
Since most vertices in  random $d$-regular graphs have large tree neighborhoods,  most vertices $j$ in the summation of \eqref{e:Qsum} have large tree neighborhoods. For those vertices $j$, we can future replace the boundary weights in the approximation by the  weight $Q(\cG,z)$ (Figure \ref{f:sumpath}). Moreover, the radius $\ell$ neighborhoods of those vertices $j$ with vertex $i$ removed are truncated $(d-1)$-ary trees of depth $\ell$. 
Let $Y_\ell(Q(\cG,z), z)$ be the Green's function at root vertex of a depth $\ell$ truncated $(d-1)$-ary tree with boundary weights $Q(\cG,z)$. We expect to have the following self-consistent equation of $Q(\cG,z)$
\begin{align}\label{e:selfeq}
Q(\cG,z)=\frac{1}{Nd}\sum_{i\sim j}G_{jj}^{(i)}\approx Y_\ell(Q(\cG,z), z).
\end{align}

The above heuristics based on expressing Green's function in terms of a sum  over paths, i.e.,  \eqref{e:weightsum},   are not rigorous, because  \eqref{e:weightsum}  does not converge for $z$ close to the spectrum. Instead of using this expansion, we will analyze the self-consistent equations and use heavily the idea of local resampling (see  Section \ref{sec:switch}). 
The local resampling in random graphs was known for a long time, but the first implementation to derive the spectral statistics appeared only recently in  \cite[Section 3]{MR3688032}.
The local resampling was modified  in \cite[Section 7]{MR3962004}, which  randomizes  the boundary of a neighborhood $\cT$ (instead of randomizing edges near a vertex \cite{MR3688032})  by switching the edge boundary of $\cT$ with randomly chosen edges in the remaining of the graph. The local resampling used in this paper performs a simple switching on a  boundary edge of $\cT$ only when this switching pair is far away from other switching pairs. While one can put any restrictions on  switchings, it is critical in our application that  this restriction preserves the reversibility of the switching, i.e.,   the law for  the graphs and their switched graphs  is exchangeable.

More precisely, we will use the exchangeability in the following way. 
Since the local resampling preserves uniform distribution on random $d$-regular graphs, 
instead of proving \eqref{e:selfeq} for $\cG$, we prove it for $\tcG$, which is the graph 
obtained from $\cG$ by a  local resampling.
Since  the local resampling is equivalent to taking out the subgraph $\cT$ and replanting it back to $\cG$ at random locations, the new boundary vertices of $\cT$ after local resampling are typically far away from each other in the sense of graph distance (Section \ref{sec:dist}) and Green's function distance (Section \ref{sec:distGF}). Therefore we can rigorously approximate the Green's function $\tilde G_{oo}$ of $\tcG$ by the Green's function of $\cT$ with boundary weights $\tilde G_{jj}^{(i)}$ at 
the boundary vertex $j$. On the other hand, since the edge boundary of $\cT$ contains more than $ \log N$ edges, the average of the weights at boundary vertices of $\cT$ concentrates around $Q(\tcG,z)$ (Section \ref{sec:concentration}). This allows us to replace boundary weights by the same weight $Q(\tcG,z)$ at each boundary vertex, and leads to the self-consistent equation \eqref{e:selfeq}.

While the delocalization of eigenvectors is a direct consequence of \eqref{e:selfeq},  the  error term in the self-consistent equation \eqref{e:selfeq} we will derive    is not small enough to detect the locations of extremal eigenvalues, i.e. $\la_2,\la_N$.  This is a standard problem  in random matrix theory, namely,  the self-consistent equation yields a ``weak law" while detecting  extremal  eigenvalue locations requires some form of  a ``strong law". 
Without getting into terminology,  our goal is to  show that   the high moments of the self-consistent equation
\begin{align}\label{e:hmm}
\bE[|Q(\cG,z)-Y_\ell(Q(\cG,z), z)|^{2p}],
\end{align}
 are sufficiently small for large integer $p>0$.  This step is typically called the fluctuation averaging estimate \cite{MR2871147,MR3800840}. In our setting, the matrix elements of  $d$-regular graphs are not independent 
 and the traditional method clearly fails.   Our  key observation is that the independence property  can be replaced by the reversibility of 
 the local resampling, i.e.,  the pair  $(\cG, \tcG) $ with  $\tcG$ denoting the graph obtained from a random $d$-regular graph $\cG$ by a  local resampling forms an exchangeable pair, i.e.
$
(\cG, \tcG)\stackrel{law}{=}(\tcG, \cG).
$ 
Besides the lack of independence,    there is an outlier eigenvalue at $d/\sqrt{d-1}$ for $d$-regular graphs. We need to remove the contribution of the outlier eigenvalue from the self-consistent equation \eqref{e:hmm}. With these ideas, we  improve  
 the error bound  in the self-consistent equation   \eqref{e:selfeq} to its square in the high moments sense.  
With this improved error bound, we conclude that non-trivial extremal eigenvalues concentrate around $\pm2$.

Finally, we remark on the difficulties  in extending the results  \cite{MR3962004} for large $d$ all the way to  $d\geq 3$ in this paper. 
We first note that  the spectral property of the graphs  for small $d$ are very sensitive to  small perturbations;   one can easily break the nice spectral properties, e.g. delocalized eigenvectors and large spectral gap, by changing a small number of edges.  To overcome this, we design the  new local resampling to perform simple switchings only when the switching pair has tree neighborhood and is far away from other switching pairs. In this way,  all the intermediate graphs from local resampling will have good local geometry, i.e. they don't have too many cycles in small neighborhoods. Then we can show all the intermediate graphs have desired spectral properties for any $d\geq 3$. 
Another major difficulty for small $d$ is that  typical estimates of the Green's function do not hold with very high probability. For example when $d=3$, the  delocalization of eigenvectors hold at most with probability $1-\OO(N^{-1})$ (See Figure \ref{fig:fork}).  In \cite{MR3962004}, we need to take a union bound for various  estimates of the Green's function over the upper half plane. 
A similar procedure would have been  impossible for small $d$.  To overcome this difficulty,  we  restrict our estimates to  the set $\bar\Omega$ of  ``radius-$\fR$ tree-like graphs", which are graphs with benign local structures. 
In  the high moments estimates of the self-consistent equation \eqref{e:hmm}, we also need to  restrict the expectation to the set $\bar\Omega$. Here we used a key property of 
the local resampling, namely, it preserves the radius-$\fR$ tree like property with high probability. In this way, we  are able to show that for any $d\geq 3$ on the set  $\bar\Omega$ 
key  spectral properties hold with probability  $1-\OO(N^{-\fC})$ for any $\fC>0$. 

Random $d$-regular graphs can also be constructed from $d$ copies of random perfect matchings, or random lifts of a base graph containing two vertices and $d$ edges between them. This class of random graphs obtained from random lifts and in  particular  their extremal  eigenvalues have been extensively studied \cite{amit2002random, amit2006random, friedman2003relative, linial2010word,friedman2014relativized, puder2015expansion,lubetzky2011spectra, bordenave2019eigenvalues}. It is interesting to see if the approach in this paper can be used to analyze   extremal eigenvalues in this setting.

In summary, the method in this paper,  based on the self-consistent equation  \eqref{e:selfeq} and a new resampling mechanism in Section \ref{sec:switch}, provides a new powerful method 
to analyze spectral properties of random $d$-regular graphs up to the optimal $d=3$. This method not only yields  estimates on both eigenvalues and eigenvectors 
over the entire range of spectrum, but its  estimates on the location of the second eigenvalue (in fact, on all extremal eigenvalues)  
are  stronger than those obtained by the 
techniques of irreducible traces or non-backtracking walks (which are sophisticated moment methods).
 In addition, we have isolated a class of events beyond which all our estimates hold with very high probability. We believe that this method can be further tighten to improve the error bound on the second largest eigenvalue
and we hope that it will eventually  lead to the optimal  bound $\la_2, |\la_N|\leq 2+N^{-\fc}$ for all $\fc < 2/3$.


\subsection{Some notations}
We reserve letters in mathfrak mode, e.g. $\fa, \fb, \fc$, to represent universal constants, and use $\fC,\fC_1, \fC_2,\cdots$, to represent large universal constants, which
may be different from line by line. We use letters in mathcal mode, e.g. $\cB, \cG, \cT$, to represent graphs, or subgraphs, and letters in mathbb mode, e.g. $\bB, \bG, \bT$, to represent set of vertices. We use letters in mathsf mode, e.g. $ \mathsf A, \mathsf B, \mathsf U, \mathsf W$ to represent index sets.
For two quantities $X$ and $Y$ depending on $N$, 
we write that $X = \OO(Y )$ if there exists some universal constant such
that $|X| \leq \fC Y$ . We write $X = \oo(Y )$, or $X \ll Y$ if the ratio $|X|/Y\rightarrow \infty$ as $N$ goes to infinity. We write
$X\asymp Y$ or $X=\Omega(Y)$ if there exists a universal constant $\fC>0$ such that $ Y/\fC \leq |X| \leq  \fC Y$. We remark that the implicit constants may depend on $d$.
 We write $a\vee b=\max\{a,b\}$ and $a\wedge b=\min\{a,b\}$. We denote $\qq{a,b} = [a,b]\cap\Z$ and $\qq{N} = \qq{1,N}$. We denote the complex upper half plane as $\bC^+$. We write a function $g\equiv 0$, if it is identically zero, and $g\not\equiv 0$, if it is not identically zero.

\noindent\textbf{Acknowledgements }
The research of J.H. is supported by the Simons Foundation as a Junior Fellow at
the Simons Society of Fellows, and NSF grant DMS-2054835.
The research of H.-T. Y. is partially supported by NSF grant DMS-1855509 and  DMS-2153335.

\section{Tree-like Graphs}

In Section \ref{sec:graph}, we collect some basic definitions and terminology from graph theory, and recall some basic structure properties of random $d$-regular graphs. In Section \ref{sec:TE} and \ref{sc:ext}, we introduce the tree extension of graphs and extension of graphs with general weights, which includes tree extension as a special case. We also obtain estimates for the Green's functions for those extended graphs.
\subsection{Graphs}\label{sec:graph}

In this section we collect some definitions and terminology about graphs, and basic structure properties of random $d$-regular graphs.

\paragraph{Graphs, adjacency matrices, Green's functions}
Throughout this paper, graphs $\cG$ are always simple (i.e., have no self-loops or multiple edges) and have vertex degrees at most $d$
(non-regular graphs are also used).
The degree of a vertex (the number of its adjacent vertices) in the graph $\cG$ is denoted by $\deg_{\cG}(\cdot)$.
The distance (length of the shortest path between two vertices) in the graph $\cG$ is denoted by $\dist_{\cG}(\cdot, \cdot)$. 
The diameter $\diam(\cG)$ of the graph $\cG$ is the greatest distance between any pair of vertices. If two vertices $i,j\in \cG$ are adjacent, i.e. $\dist_{\cG}(i,j)=1$, we write $i\sim j$.
For any graph $\cG$, the adjacency matrix is the (possibly infinite) symmetric matrix $A$ indexed by the vertices of the graph,
with $A_{ij}= A_{ji} = 1$ if there is an edge between $i$ and $j$, and $A_{ij}=0$  otherwise.
Throughout the paper, we denote the normalized adjacency matrix by $H=A/\sqrt{d-1}$.
Moreover, we denote the (unnormalized) adjacency matrix of a directed edge $(i,j)$ by $e_{ij}$, i.e.\ $(e_{ij})_{kl} = \delta_{ik}\delta_{jl}$.
The Green's function of a graph $\cG$ is the unique matrix $G = G(\cG;z)$ defined by $G(H-z)=I$ for $z \in \bC^+$, where $\bC^+$ is the upper half plane. The Stieltjes transform of the empirical eigenvalue distribution of $\cG$ is $m_N(z)=\Tr G(z)/N$.\index{$m_N$}
In Appendix~\ref{app:Green}, several well-known properties of Green's function are summarized;
they will be used throughout the paper.
The Stieltjes transform $m_N(z)$ encodes information of eigenvalues of $H$ (and thus of $A$), and the Green's function $G(z)$ encodes both information of eigenvalues and eigenvectors.
In particular, the spectral resolution is given by $\eta = \im [z]$:
the macroscopic behavior corresponds to $\eta$ of order $1$,
the mesoscopic behavior corresponds to $1/N \ll \eta \ll 1$,
and the microscopic behavior of individual eigenvalues corresponds
to $\eta\asymp 1/N$.

\paragraph{Subsets and Subgraphs}
Let $\cG$ be a graph, and denote the set of its edges by the same symbol $\cG$ and its vertices by $\bG$.
For any subset $\bT\subset \bG$, we define the graph $\cG^{(\bT)}$ from removing the vertices $\bT$ and edges adjacent to $\bT$ from $\cG$. Then the adjacency matrix of $\cG^{(\bT)}$ is the restriction of that of $\cG$ to $\bG \setminus \bT$.
We write $G^{(\bT)}$ for the Green's function of $\cG^{(\bT)}$.
For any subgraph $\cal T \subset \cG$, we denote by $\del \cal T=\{v\in \bG: \dist_{\cG}(v,\cal T)=1\}$
the vertex boundary of $\cal T$ in $\cG$,
and by $\partial_E \cal T = \{e\in \cG: \text{$e$ is adjacent to $\cal T$ but $e\not\in \cal T$}\}$ the  edge boundary of $\cal T$ in $\cG$.

\paragraph{Neighborhoods}
Given a subset $\bT$ of the vertex set of a graph $\cG$ and an integer $r>0$,
we denote the radius-$r$ neighborhood of $\bT$ in $\cG$ by $\cB_r(\bT,\cG)$,
i.e., it is the subgraph induced by $\cG$ on the set $\{j \in \bG: \dist_{\cal G}(\bT,j) \leq r\}$.
In particular, $\cB_r(i,\cG)$ is the radius-$r$ neighborhood of the vertex $i$.

{

\paragraph{Trees}
The infinite $d$-regular tree is the unique (up to isometry) infinite connected $d$-regular graph without cycles, and we denote it by $\cX$.
The infinite $(d-1)$-ary tree is the unique (up to isometry) infinite connected tree that is $d$-regular
at every vertex except for a distinguished root vertex $o$, which has degree $d-1$, we denote it by $\cY$.

}

\paragraph{Excess of random $d$-regular graphs}
For any graph $\cG$, we define its excess to be the smallest number of edges that must be removed to yield a graph with no cycles
(a forest). It is given by
\begin{align}\label{def:excess}
\text{excess}(\cG)\deq \#\text{edges}(\cG) -\#\text{vertices}(\cG) +\#\text{connected components}(\cG).
\end{align} 
There are different conventions for the normalization of the excess.
Our normalization is such that the excess of a tree or forest is $0$.
Note that if $\cH \subset \cG$ is a subgraph, then $\text{excess}(\cH) \leq \text{excess}(\cG)$.

We will use the following well-known estimates for the excess in random $d$-regular graphs. We recall the constants $0<\fc<1$ and $\fR=(\fc/4)\log_{d-1}N$ from the beginning of Section \ref{sec:intro}. The following proposition implies that with high probability, random $d$-regular graphs are radius-$\fR$ tree like
as in Definition \ref{def:barOmega}. 
\begin{proposition} \label{prop:structure}
Fix $d\geq 3$, let $\n \geq 1$ be an integer. If $R\leq\fR$,
then the following holds for a uniformly chosen random $d$-regular graph $\cG$ on $\qq{N}$,
with probability at least $1-\OO(N^{-(1-\fc)\n})$
for $N$ large enough we have 
\begin{itemize}
\item All radius-$R$ neighborhoods have excess at most $\omega$:
\begin{equation}
\label{e:structure1}
\text{for all $i \in \qq{N}$, the subgraph $\cB_{R}(i,\cG)$ has excess at most $\n$.}
\end{equation}
\item Most radius-$R$ neighborhoods are trees:
\begin{equation}
\label{e:structure2}
|\{i \in \qq{N}: \text{the subgraph $\cB_{R}(i,\cG)$ is not a tree}\}| \leq N^{\fc}.
\end{equation}
\end{itemize}
\end{proposition}

\begin{proof}
We have that 
$
\bP(\eqref{e:structure1} \text{ holds})=1-\OO(N^{-(1-\fc)\omega}),
$
and
$
\bP(\eqref{e:structure2} \text{ holds})=1-\OO(e^{-N^{2\fc}}) 
$
from \cite[Proposition 4.1]{MR3962004}.
Proposition \ref{prop:structure} follows from combining them.

\end{proof}

\subsection{Trees and tree extension}
\label{sec:TE}
In this section, we calculate the Green's functions of the infinite $d$-regular tree and the infinite $(d-1)$-ary tree, and introduce the concepts of deficit function and tree extension.

We recall that $\rho_{d}(x)$ is the density of the Kesten--McKay distribution \eqref{kmlaw}
and let $\rho_{sc}(x):=\sqrt{[4-x^2]_+}/2\pi$\index{$\rho_{sc}$} be the semicircle distribution. We denote their Stieltjes transforms as $\md(z)$ and $\msc(z)$,\index{$m_d(z)$, $\msc(z)$}
\begin{equation*}
  m_d(z) = \int \frac{ \rho_d(\lambda)}{\lambda-z} \rd\lambda, \qquad
  \msc(z) = \int \frac{\rho_{sc}(\la)}{\lambda-z}  \rd\lambda.
\end{equation*}
For any $z\in \bC^+$, $m_{sc}(z)$ satisfies the quadratic equation $m_{sc}^2(z)+zm_{sc}(z)+1=0$  and
$\md(z)$ is explicitly related to $\msc(z)$ by the equation 
\begin{align} \label{e:mdstieltjes}
\md(z)=\frac{1}{-z-d(d-1)^{-1}\msc(z)}=\frac{\msc(z)}{1-({d-1})^{-1}\msc^2(z)}.
\end{align}

For the infinite $d$-regular tree and the infinite $(d-1)$-ary tree,
the following proposition computes their Green's function explicitly.
\begin{proposition}\label{greentree}
Let $\cX$ be the infinite $d$-regular tree.
For all $z \in \bC^+$, its Green's function is
\begin{equation} \label{e:Gtreemkm}
  G_{ij}(z)=m_{d}(z)\left(-\frac{\msc(z)}{\sqrt{d-1}}\right)^{\dist_{\cX}(i,j)}.
\end{equation}
Let $\cY$ be the infinite $(d-1)$-ary tree with root vertex $o$.
Its Green's function is
\begin{equation} \label{e:Gtreemsc}
  G_{ij}(z)=m_{d}(z)\left(1-\left(-\frac{\msc(z)}{\sqrt{d-1}}\right)^{2{\rm anc}(i,j)+2}\right)\left(-\frac{\msc(z)}{\sqrt{d-1}}\right)^{\dist_{\cY}(i,j)},
\end{equation}
where ${\rm anc}(i,j)$ is the distance from the common ancestor of the vertices $i,j$ to the root $o$. 
In particular,
\begin{align}\label{e:Gtreemsc2}
G_{oi}(z)=\msc(z)\left(-\frac{\msc(z)}{\sqrt{d-1}}\right)^{\dist_{\cY}(o,i)}.
\end{align}
\end{proposition}

\begin{proof}
See \cite[Proposition 5.1]{MR3962004}.
\end{proof}

\begin{definition}[deficit function]
Given a graph $\cG$ with vertex set $\bG$ and degrees bounded by $d$,
a deficit function for $\cG$ is a function
$g: \bG \to \qq{0,d}$ satisfying $\deg_{\cG} (v) \leq d-g(v)$ for all vertices $v\in \bG$.
We call a vertex $v\in \bG$ extensible if $\deg_\cG(v)<d-g(v)$.
\end{definition}

\begin{definition}[tree extenstion]
Let $\cG$ be a  graph with deficit function $g$.
The tree extension of $\cG$ is the (possibly infinite) graph $\TE(\cG)$ defined
by attaching to any extensible vertex $v$ in $\cG$, $d-g(v)-\deg_{\cG}(v)$ copies of $(d-1)$-ary trees. In this way, each vertex $v\in \bG$ in the extended graph has degree $d-g(v)$.
\end{definition}

\begin{figure}[t]
\centering
\scalebox{.6}{\input{TE.pspdftex}}
\caption{The left figure illustrates a finite graph $\cG$; its extensible vertices are shown as grey circles.
The right figure shows the tree extension $\TE(\cG)$, in which a rooted tree (darkly shaded) is attached to every extensible vertex.
\label{fig:TE}}
\end{figure}
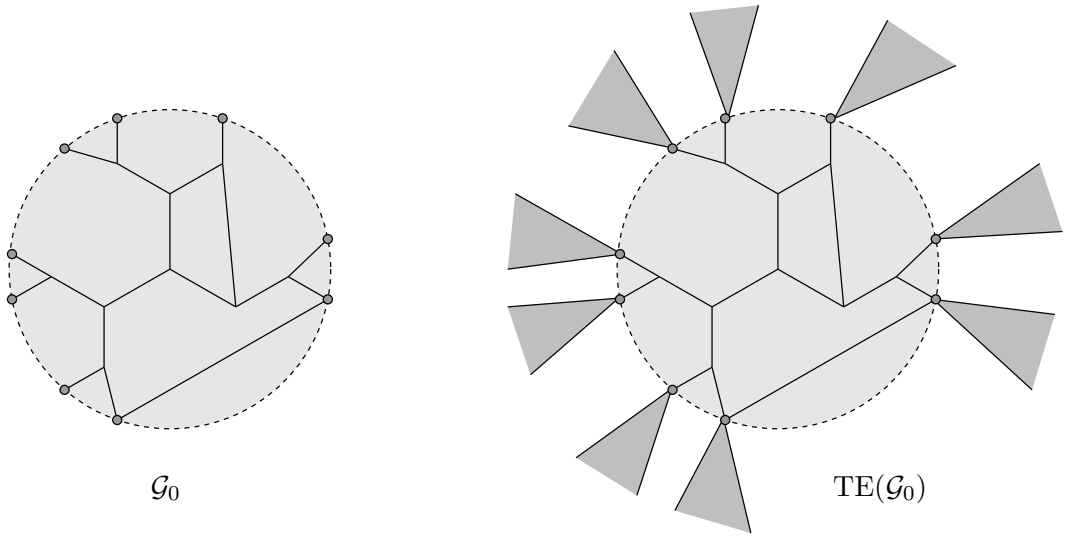

\begin{remark}\label{r:convention}(Conventions for deficit functions)\
Throughout this paper, all graphs $\cG$ are equipped with a deficit function $g$.
The interpretation of the deficit function $g(v)$ is that it measures the difference to the desired degree of the vertex $v$.
We use the following conventions for deficit functions.
\begin{itemize}
\item
For $d$-regular graphs $\cG\in \GNd$ and the infinite $d$-regular graph $\cX$, their deficit function is identically zero, i.e. $g\equiv 0$. For the infinite $(d-1)$-ary tree $\cY$, its deficit function $g(o)=1$ at the root vertex, and $g=0$ elsewhere. 
\item
If $\bT$ is a subset of the vertices of $\cG$, and $g$ is the deficit function of $\cG$,
then the deficit function $g'$ of $\cG^{(\bT)}$ is given by 
\begin{align}\label{e:removing}
g'(v) = g+\deg_{\cG}(v)-\deg_{\cG^{(\bT)}}(v).
\end{align} 
Thus $\TE(\cG^{(\bT)})$ can be obtained from $\TE(\cG)$ by removing the edges incident to $\bT$.
\item
If $\cH \subset \cG$ is a subgraph (which was not obtained as $\cG^{(\bT)}$),
then the deficit function of $\cH$ is given by the restriction of the deficit function of $\cG$ on $\cH$,
unless specified explicitly.
\end{itemize}
\end{remark}

For each $d\geq3$, we associate  an integer $\omega_d>0$ so that an infinite $d$-regular graph with excess less or equal than $\omega_d$ has good spectral properties: 

\begin{definition}[$\n_d$\index{$\omega_d$} definition]\label{d:defwd}
Fix $d\geq 3$, we define $\n_d$ to be the largest $\n$ such that the following holds:
For any connected graph (possibly infinite) $\cG$ with vertex set $\bG$ and deficit function $g$, if the sum of deficit function $\sum_{v\in \bG}g(v)$  plus the excess of $\cG$ is at most   $\n$ i.e., $\sum_{v\in \bG}g(v)+\text{excess}(\cG)\le \n$, then for all $z \in \bC^+$ and all $i,j\in \bG$, the Green's function of the tree extension of $\cG$,  $G_{ij}(\TE(\cG),z)$, satisfies
\begin{equation} \label{e:boundPij} 
\left|G_{ij}(\TE(\cG),z)\right|
\lesssim \left(\frac{|m_{sc}(z)|}{\sqrt{d-1}}\right)^{\dist_{\cG}(i,j)},\quad  |\Im[G_{ij}(\TE(\cG),z)]|\lesssim \Im[\md(z)],
\end{equation}
and the diagonal terms satisfy the estimate
\begin{equation} \label{e:boundPii}
|G_{ii}(\TE(\cG),z)|\asymp 1.
\end{equation}

\begin{remark}
For $d=3$ we have $\omega_3=1$, see Proposition \ref{p:wd}.
In general, $\omega_d\geq 1$ is an increasing function of $d$, and $\omega_d$ goes to infinity as $d$ increases. 
\end{remark}
\begin{remark}\label{r:wdext}
Thanks to the resolvent identity \eqref{e:Schurixj}, if a graph $\cG$ satisfies \eqref{e:boundPij} and \eqref{e:boundPii}, it still satisfies \eqref{e:boundPij} after removing one vertex. But \eqref{e:boundPii} might fail.  In fact, after removing vertex $j$ from $\cG$, the deficit function of the new graph $\cG^{(j)}$ is as given in \eqref{e:removing}. And  the Green's function of its tree extension satisfies
\begin{align*}
G_{ii}(\TE(\cG^{(j)}),z)=G^{(j)}_{ii}(\TE(\cG,z)=G_{ii}(\TE(\cG),z)-\frac{G_{ij}(\TE(\cG),z)G_{ji}(\TE(\cG),z)}{G_{jj}(\TE(\cG),z)},
\end{align*}
which is still at most $\OO(1)$, but we may not have the lower bound.
Therefore,
if a connected graph $\cG$ with vertex set $\bG$ and deficit function $g$ satisfying
that 
\begin{itemize}
\item either the sum of deficit function $\sum_{v\in \bG}g(v)$  plus the excess is at most $\n_d$: $\sum_{v\in \bG}g(v)+\text{excess}(\cG)\le \n_d$,
\item or $g\not\equiv 0$ and the sum of deficit function $\sum_{v\in \bG}g(v)$  plus the excess is at most $\n_d+1$: $\sum_{v\in \bG}g(v)+\text{excess}(\cG)\le \n_d+1$,
\end{itemize}
then \eqref{e:boundPij} still holds for its Green's function. 
\end{remark}

\end{definition}

\subsection{Green's function for extension of graphs with weight \texorpdfstring{$\Delta(z)$}{D(z)}}\label{sc:ext}

Fix $z=E+\ri\eta\in \bC^+$, and $\kappa\deq \dist(E, \{-2,2\})=||E|-2|$. In this section we introduce the extension of graphs with general weight $\Delta(z)$. The tree extension introduced in last section is a special case, which corresponds to the extension of graphs with the weight given by the Stieltjes transform of the semicircle distribution $\msc(z)$. By a perturbative argument, then we obtain estimates for the Green's function for extensions of graphs with general weight $\Delta(z)$ close to $\msc(z)$.

\begin{definition}[extension with weight $\Delta(z)$] \label{def:ext}
Let $\cG$ be a finite graph with vertex set $\bG$ and deficit function $g: \bG\mapsto \qq{1,d}$.
\begin{itemize}
\item
The extension of $\cG$ with weight $\Delta(z)$ is the graph $\Ext(\cG, \Delta(z))$ defined
by assigning any extensible vertex $v$ in $\cG$ a weight of $-(d-g(v)-\deg_{\cG}(v))\Delta(z)/\sqrt{d-1}$.
\item
We denote $G(\Ext(\cG, \Delta(z)),z)$ the Green's function of the graph $\Ext(\cG,\Delta(z))$.
\end{itemize}
\end{definition}
Let $H$ be the normalized adjacency matrix of $\cG$ with vertex set $\bG$. The normalized adjacency matrix of the extension of $\cG$ with weight $\Delta(z)$ is given by:
\begin{align*}
H-\sum_{v\in \bG}\frac{d-g(v)-\deg_{\cG}(v)}{d-1}\Delta(z)e_{vv},
\end{align*}
where the matrix $(e_{vv})_{ij}=\delta_{vi}\delta_{vj}$. Its Green's function is 
\begin{align}\label{e:formula}
G(\Ext(\cG,\Delta(z)),z)
=\left(H-z-\sum_{v\in \bG}\frac{d-g(v)-\deg_{\cG}(v)}{d-1}\Delta(z)e_{vv}\right)^{-1}.
\end{align}
From \eqref{e:formula}, one can see that the Green's function of $\TE(\cG)$ restricted on the vertex set $\bG$ of $\cG$ is the same as that of the extension of $\cG$ with weight $m_{sc}(z)$.  Recall that in our convention all graphs have degree $d$ up to the associated  deficit function. 
Therefore, if $\cG$ has no cycles and zero deficit function, i.e. $g\equiv 0$, the Green's function of the graph $\Ext(\cG, m_{sc}(z))$ is the same as the Green's function of the infinite $d$-regular tree,
\begin{align}\label{e:Gtreemsccopy}
G_{ij}(\Ext(\cG,m_{sc}(z)), z)=m_d(z)\left(-\frac{m_{sc}(z)}{\sqrt{d-1}}\right)^{\dist_{\cG}(i,j)},\quad i,j\in \bG.
\end{align}
If $\cG$ has no cycles and its deficit function equals $1$ at vertex $o$ and zero elsewhere, i.e. $g(v)=\bm1(v=o)$, the Green's function of the graph $\Ext(\cG, m_{sc}(z))$ is the same as the Green's function of the infinite $(d-1)$-ary tree with root vertex $o$,
\begin{align}\label{e:Gtreemsccopy2}
G_{oi}(\Ext(\cG,m_{sc}(z)), z)=m_{sc}(z)\left(-\frac{m_{sc}(z)}{\sqrt{d-1}}\right)^{\dist_{\cG}(o,i)},\quad i\in \bG.
\end{align}

In the rest of this section we study the Green's function $G(\Ext(\cG,\Delta(z)), z)$ for extension of graphs with weight $\Delta(z)$,  when the weight $\Delta(z)$ is close to $m_{sc}(z)$, i.e.
\begin{align*}
|\Delta(z)-m_{sc}(z)| \ll1.
\end{align*}
Since the proofs are straightforward computations, we postpone them to Appendix \ref{a:ext}.
For any integer $\ell\geq 1$, we define the functions $X_\ell(\Delta(z),z), Y_\ell(\Delta(z),z)$ as
\begin{align}\begin{split}\label{def:Y}
&X_\ell(\Delta(z),z)=G_{oo}(\Ext(\cB_\ell(o,\cX),\Delta(z)), z),\\
&Y_\ell(\Delta(z),z)=G_{oo}(\Ext(\cB_\ell(o,\cY),\Delta(z)), z),
\end{split}\end{align}
where $\cX$ is the infinite $d$-regular tree, and $\cY$ is the infinite $(d-1)$-ary tree with root vertex $o$. Then $m_{sc}(z)$ is a fix point of the function $Y_\ell$, i.e. $Y_\ell(m_{sc}(z),z)=m_{sc}(z)$. And $X_\ell(\msc(z),z)=\md(z)$, The following proposition gives the stability estimates of the function $Y_\ell$. If $\Delta(z)$ is sufficiently close to $m_{sc}(z)$, then $Y_\ell(\Delta(z),z)$ is close to $m_{sc}(z)$ and we have explicit estimates.
\begin{proposition}\label{p:recurbound}
If $\ell|\Delta(z)-m_{sc}(z)|\ll 1$, then the functions $X_\ell(\Delta(z),z), Y_\ell(\Delta(z),z)$  defined in \eqref{def:Y}
satisfy
\begin{align}\begin{split}\label{e:Xrecurbound}
&\phantom{{}={}}X_\ell(\Delta(z),z)-\md(z)
=\frac{d}{d-1}(\md(z))^2(m_{sc}(z))^{2\ell}(\Delta(z)-m_{sc}(z))
+\OO\left(\ell|\Delta(z)-m_{sc}(z)|^2\right),
\end{split}\end{align}%
and
\begin{align}\begin{split}\label{e:recurbound}
&\phantom{{}={}}Y_\ell(\Delta(z),z)-m_{sc}(z)
=(m_{sc}(z))^{2\ell+2}(\Delta(z)-m_{sc}(z))\\
&+\msc^{2\ell+2}(z)\md(z)\left(\frac{1-\msc^{2\ell+2}(z)}{d-1}+\frac{d-2}{d-1}\frac{1-\msc^{2\ell+2}(z)}{1-\msc^2(z)}\right)(\Delta-\msc(z))^2
+\OO\left(\ell^2|\Delta(z)-m_{sc}(z)|^3\right).
\end{split}\end{align}
\end{proposition}

More generally if $\cG$ is a finite rooted graph $\cG$ with root $o$ and vertex set $\bG$. We assume that $\cG$ has no cycles, and its deficit function takes value $1$ at root $o$ and $0$ elsewhere, i.e., $g(v)={\bf1}_{v=o}$, then the tree extension of $\cG$ is the infinite $(d-1)$-ary tree with root vertex $o$. Therefore, $m_{sc}(z)$ is a fixed point, $G_{oo}(\Ext(\cG, m_{sc}(z)),z)=m_{sc}(z)$. The following Proposition gives its stability estimates, i.e. if $\Delta(z)$ is sufficiently close to $m_{sc}(z)$ then $G_{oo}(\Ext(\cG,\Delta(z)), z)$ is sufficiently close to $\Delta(z)$.
\begin{proposition}\label{p:fixpoint}
Fix a finite rooted graph $\cG$ with root $o$ and vertex set $\bG$. We assume that $\cG$ has no cycles, and its deficit function takes value $1$ at root $o$ and $0$ elsewhere, i.e.,  $g(v)={\bm 1}_{v=o}$. Take a weight $\Delta(z)$ satisfying $\diam(\cG)|\Delta(z)-m_{sc}|\ll1$.  Then the Green's function of the graph $\Ext(\cG, \Delta(z))$ satisfies
\begin{align*}
G_{oo}(\Ext(\cG,\Delta(z)),z)-\Delta(z)
=
\OO\left((1+\diam(\cG)\left(\sqrt{\kappa+\eta}|\Delta(z)-m_{sc}(z)|+|\Delta(z)-m_{sc}(z)|^2\right)\right),
\end{align*}%
where $\eta=\Im[z]$ and $\kappa=\min\{|\Re[z]-2|, |\Re[z]+2|\}$.
\end{proposition}

By a perturbation argument, if $\Delta(z)$ is sufficiently close to $m_{sc}(z)$, the Green's function $G(\Ext(\cG,\Delta(z)), z)$ of the graph $\Ext(\cG, \Delta(z))$ is close to the Green's function of its tree extension $\TE(\cG)$. As a consequence, the Green's function $G(\Ext(\cG,\Delta(z)), z)$ also satisfies estimates \eqref{e:boundPij}, \eqref{e:boundPii}.

\begin{proposition}\label{p:subgraph}
Fix a finite connected graph $\cG$ with vertex set $\bG$ and deficit function $g$, and recall $\n_d$ from Definition \ref{d:defwd}. If the sum of deficit function $\sum_{v\in \bG}g(v)$  plus the excess of $\cG$ is at most  $\n_d$ or $\n_d+1$ if $g\not\equiv 0$, and the weight $\Delta(z)$ satisfies $\diam(\cG)|\Delta(z)-m_{sc}(z)|\ll1$, then for any two vertices $i,j\in \bG$ the Green's function of $\Ext(\cG, \Delta(z))$ satisfies
{\begin{align}\begin{split} \label{e:smalldiff}
&\phantom{{}={}}\left|G_{ij}(\Ext(\cG, \Delta(z)), z)-G_{ij}(\Ext(\cG, m_{sc}),z)\right|\\
&\lesssim (1+\dist_{\cG}(i,j))\left(\frac{|m_{sc}|}{\sqrt{d-1}}\right)^{\dist_\cG(i,j)}|\Delta(z)-m_{sc}(z)|.
\end{split}\end{align}}%
As a consequence, we have
\begin{align}\label{e:boundPijcopy}
&|G_{ij}(\Ext(\cG,\Delta(z)), z)|
\lesssim \left(\frac{|m_{sc}|}{\sqrt{d-1}}\right)^{\dist_\cG(i,j)},\\
\label{e:sumPixcopy}
&\sum_{x\in \bG: \deg_\cG(x)<d-g(x)}|G_{ix}(\Ext(\cG,\Delta(z)), z)|
\lesssim(d-1)^{\diam(\cG)/2},\\
\label{e:sumPix2copy}
&\sum_{x\in \bG: \deg_\cG(x)<d-g(x)}|G_{ix}(\Ext(\cG,\Delta(z)), z)||G_{jx}(\Ext(\cG,\Delta(z)), z)|
\lesssim
(1+\dist_{\cG}(i,j))\left(\frac{|m_{sc}(z)|}{\sqrt{d-1}}\right)^{\dist_{\cG}(i,j)}.
\end{align}%
Moreover, if the sum of deficit function $\sum_{v\in \bG}g(v)$  plus the excess of $\cG$ is at most $\n_d$, then 
\begin{align}\label{e:boundPiicopy} 
|G_{ii}(\Ext(\cG,\Delta(z)), z)|\asymp 1.
\end{align}

\end{proposition}

\begin{remark}\label{r:subgraph}
Under the same assumptions as in  Propostion \ref{p:subgraph}, if two weights $\Delta_1(z)$ and $\Delta_2(z)$ satisfy
$\diam(\cG)|\Delta_1(z)-m_{sc}(z)|, \diam(\cG)|\Delta_2(z)-m_{sc}(z)|\ll1$, then by the same argument as in Propostion \ref{p:subgraph} we have 
{\begin{align*}
&\phantom{{}={}}\left|G_{ij}(\Ext(\cG, \Delta_1(z)),z)-G_{ij}(\Ext(\cG,\Delta_2(z)),z)\right|\\
&\lesssim (1+\dist_{\cG}(i,j))\left(\frac{|m_{sc}|}{\sqrt{d-1}}\right)^{\dist_\cG(i,j)}|\Delta_1(z)-\Delta_2(z)|.
\end{align*}}%
\end{remark}

The Green's function $G(\Ext(\cG,\Delta(z)), z)$ depends only weakly on $\cG$, i.e. if we replace $\cG$ by a sufficiently large subgraph $\cH\subset \cG$, then the difference between $G(\Ext(\cG,\Delta(z)), z)$ and $G(\Ext(\cH,\Delta(z)), z)$ is small. The following proposition quantifies this phenomenon. We call it the localization principle. This localization principle will be used repeatedly throughout Sections \ref{s:switchstable} and \ref{s:improveG}.
\begin{proposition}[Localization principle]\label{p:localization}
Fix a finite connected graph $\cG$ with vertex set $\bG$ and deficit function $g$, and recall $\n_d$ from \ref{d:defwd}.  Assume that the sum of deficit function $\sum_{v\in \bG}g(v)$  plus the excess of $\cG$ is at most  $\n_d$ or $\n_d+1$ if $g\not\equiv 0$, and $\diam(\cG)|\Delta(z)-m_{sc}(z)|\ll1$. Then for any two vertices $i,j$ in $\cG$, and a subgraph $\cH$ of $\cG$ which contains    $\cB_r(i,j,\cG)$, the Green's functions of $\Ext(\cG, \Delta(z))$ and $\Ext(\cH, \Delta(z))$ satisfy
\begin{align} \label{e:compatibility}\begin{split}
&\phantom{{}={}}\left|G_{ij}(\Ext(\cG, \Delta(z)),z)-G_{ij}(\Ext(\cH,\Delta(z)),z)\right|\\
&\lesssim (1+\diam(\cG))\left(\sqrt{\kappa+\eta}|\Delta(z)-m_{sc}(z)|+|\Delta(z)-m_{sc}(z)|^2\right)+\left(\frac{|m_{sc}|}{\sqrt{d-1}}\right)^{2r}.
\end{split}\end{align}
\end{proposition}
The proofs of Propositions \ref{p:recurbound}, \ref{p:fixpoint},  \ref{p:subgraph} and \ref{p:localization} are proven by explicit computation, and perturbation arguments when $|\Delta(z)-m_{sc}(z)|\ll1$. We postpone them to Appendix \ref{a:ext}.

\section{Local Resampling}
\label{sec:switch}

In this section, we introduce the local resampling of a random $d$-regular graph, which is a modification of the local resampling in \cite[Section 7]{MR3962004}.
The procedure effectively resamples the edges on the boundary of balls of radius $\ell$,
by switching them with random edges from the remainder of the graph.
This resampling generalizes the local resampling introduced in \cite{MR3688032},
where switchings were used to resample the neighbors of a vertex (corresponding to $\ell=0$).
The local resampling provides an effective access to the randomness of the random $d$-regular graph,
which is fundamental for the proof of Theorem \ref{thm:mrmsc}.

To introduce the local resampling, we require some definitions.
We consider simple $d$-regular graphs on vertex set $\qq{N}$
and identify such graphs with their sets of edges throughout this section.
(Deficit functions do not play a role in this section.)
For any graph $\cG$, we denote the set of unoriented edges by $E(\cG)$,
and the set of oriented edges by $\vec{E}(\cG):=\{(u,v),(v,u):\{u,v\}\in E(\cG)\}$.
For a subset $\vec{S}\subset \vec{E}(\cG)$, we denote by $S$ the set of corresponding non-oriented edges.
For a subset $S\subset E(\cG)$ of edges we denote by $[S] \subset \qq{N}$ the set of vertices incident to any edge in $S$.
Moreover, for a subset $\bV\subset\qq{N}$ of vertices, we define $E(\cG)|_{\bV}$ to be the subgraph of $\cG$ induced on $\bV$.

\begin{figure}[h]
\centering
\input{switch.pspdftex}
\caption{
The switching encoded by the two directed edges $\vec S=\{(v_1, v_2), (v_3, v_4)\}$
replaces the unoriented edges $\{v_1,v_2\}, \{v_3,v_4\}$ by $\{v_1,v_4\},\{v_2,v_3\}$.
\label{fig:switching}}
\end{figure}
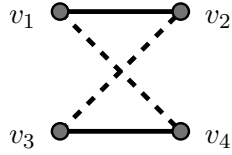

\begin{definition}[Switchings]
A (simple) switching is encoded by two oriented edges $\vec S=\{(v_1, v_2), (v_3, v_4)\} \subset \vec{E}$.
We assume that the two edges are disjoint, i.e.\ that $|\{v_1,v_2,v_3,v_4\}|=4$.
Then the switching consists of
replacing the edges $\{v_1,v_2\}, \{v_3,v_4\}$ by the edges $\{v_1,v_4\},\{v_2,v_3\}$,
as illustrated in Figure~\ref{fig:switching}.
We denote the graph after the switching $\vec S$ by $T_{\vec S}(\cG)$,
and the new edges $\vec S' = \{(v_1,v_4), (v_2,v_3)\}$ by
$
  T(\vec S) = \vec S'
$.
\end{definition}

Our local resampling involves a fixed center vertex, we now assume to be vertex $o$,
and a radius $\ell$.
Given a $d$-regular graph $\cG$, we abbreviate $\cT=\cB_{\ell}(o,\cG)$ (which may not be a tree) and its vertex set $\bT$.\index{$\cT, \bT$}
The edge boundary $\del_E \cT$ of $\cT$ consists of the edges in $\cG$ with one vertex in $\T$ and the other vertex in $\qq{N}\setminus\T$.
We enumerate $\del_E \cT$ as $ \del_E \cT = \{e_1,e_2,\dots, e_\mu\}$, where $e_\al=\{l_\al, a_\al\}$ with $l_\al \in \T$ and $a_\al \in \qq{N} \setminus \T$. We orient the edges $e_i$ by defining $\vec{e}_i=(l_\al, a_\al)$.
We notice that $\mu$ and the edges $e_1,e_2, \dots, e_\mu$ depend on $\cG$. The edges $e_\al$ are distinct, but
the vertices $a_\al$ are not necessarily distinct and neither are the vertices $l_\al$. Our local resampling switches the edge boundary of $\cT$ with randomly chosen edges in $\cGT$
if the switching is admissible (see below), and leaves them in place otherwise.
To perform our local resampling, we choose $(b_1,c_1), \dots, (b_\mu,c_\mu)$ to be independent, uniformly chosen oriented edges from the graph $\cGT$, i.e.,
the oriented edges of $\cG$ that are not incident to $\T$,
and define 
\begin{equation}\label{e:defSa}
  \vec{S}_\al= \{\vec{e}_\al, (b_\al,c_\al)\},
  \qquad
  {\bf S}=(\vec S_1, \vec S_2,\dots, \vec S_\mu).
\end{equation}
The sets $\bf S$ will be called the \emph{resampling data} for $\cG$. We remark that repetitions are allowed in the resampling data $(b_1, c_1), (b_2, c_2),\cdots, (b_\mu, c_\mu)$.

For $\al\in\qq{1,\mu}$,
we define the indicator functions
$I_\al \equiv I_\al(\cG,{\bf S})=1$\index{$I_\alpha$}
if the subgraph $\cB_{\fR/4}(\{a_\al, b_\al, c_\al\}, \cG^{(\bT)})$ after adding the edge $\{a_\al, b_\al\}$ is a tree, otherwise $I_\al \equiv I_\al(\cG,{\bf S})=0$; the indicator functions $J_\al \equiv J_\al(\cG,{\bf S})=1$\index{$J_\alpha$} if $\dist_{\cG^{(\bT)}}(\{a_\al,b_\al,c_\al\}, \{a_\beta,b_\beta,c_\beta\})> {\fR/4}$ for all $\beta\in \qq{1,\mu}\setminus \{\al\}$, otherwise $J_\al \equiv J_\al(\cG,{\bf S})=0$. We define the \emph{admissible set}
\begin{align}\label{Wdef}
\As_{\bf S}:=\{\al\in \qq{1,\mu}: I_\al(\cG,{\bf S}) J_\al(\cG,{\bf S})=1\}.\index{$\mathsf W_{\bf S}$}
\end{align}
We say that the index $\al \in \qq{1,\mu}$ is \emph{switchable} if $\al\in \As_{\bf S}$. We denote the set $\bW_{\bf S}=\{b_\al:\al\in \As_{\bf S}\}$\index{$\bW_{\bf S}$}. Let $\nu:=|\As_{\bf S}|$ be the number of admissible switchings and $\al_1,\al_2,\dots, \al_{\nu}$
be an arbitrary enumeration of $\As_{\bf S}$.
Then we define the switched graph by
\begin{equation} \label{e:Tdef1}
T_{\bf S}(\cG) := \left(T_{\vec S_{\al_1}}\circ \cdots \circ T_{\vec S_{\al_\nu}}\right)(\cG),
\end{equation}
and the switching data by
\begin{equation} \label{e:Tdef2}
  T({\bf S}) := (T_1(\vec{S_1}), \dots, T_\mu(\vec{S_\mu})), 
  \quad
  T_\al(\vec{S}_\al) \deq
  \begin{cases}
    T(\vec{S}_\al) & (\al \in \As_{\bf S}),\\
    \vec{S}_\al & (\al \not\in \As_{\bf S}).
  \end{cases}
\end{equation}

 We remark, the indicator functions $I_\al$ and $J_\al$ are different from those in \cite[Section 7]{MR3962004}. 
 Our indicator function $I_\alpha$ imposes a ``tree" condition, which ensures that 
 $a_\al$ and $\{b_\al, c_\al\}$ are far away from each other, and their neighborhoods are trees. 
And our indicator function $J_\al$ imposes an ``isolation" condition, which ensures that we only perform simple switching when the switching pair is far away from other switching pairs. In this way, we do not need to keep track of the interaction between different simple switchings. Notice that all conditions related to  $I_\al$ and $J_\al$ are imposed on  balls of radius $\fR/4$. 
 

To make the structure more clear, we introduce an enlarged probability space.
Equivalently to the definition above, the sets $\vec{S}_\al$ as defined in \eqref{e:defSa} are uniformly distributed over 
\begin{align*}
{\sf S}_{\al}(\cG)=\{\vec S\subset \vec{E}: \vec S=\{\vec e_\al, \vec e\}, \text{$\vec{e}$ is not incident to $\T$}\},
\end{align*}
i.e., the set of pairs of oriented edges in $\vec{E}$ containing $\vec{e}_\al$ and another oriented edge in $\cGT$.
Therefore ${\bf S}=(\vec S_1,\vec S_2,\dots, \vec S_\mu)$ is uniformly distributed over the set
${\sf S}(\cG)=\sf S_1(\cG)\times \cdots \times \sf S_\mu(\cG)$.

\begin{definition}\label{def:enlarge}
For any graph $\cG\in \GNd$,  denote by $\iota(\cG) = \{\cG\} \times \sf S(\cG)$ the fibre
of local resamplings of $\cG$ (with respect to vertex $o$),
and define the enlarged probability space
\begin{align*}
\GNdp = \iota(\GNd) = \bigsqcup_{\cG\in \GNd}\iota(\cG),
\end{align*}
with the probability measure $\Pp(\cG, {\bf S}):= \P(\cG)\P_{\cG}({\bf S}) = (1/|\GNd|)(1/|\sS(\cG)|)$
for any $(\cG, {\bf S})\in  \GNdp$.
Here $\P(\cG)=1/|\GNd|$ is the uniform probability measure on $\GNd$,
and $\P_{\cG}$ is the uniform probability measure on $\sS(\cG)$.\index{$\tilde \bP$}\index{$\bP_\cG$}
\end{definition}

Let $\pi: \GNdp \to \GNd$, $(\cG,{\bf S}) \mapsto \cG$ be the canonical projection onto the first component.
It is easy to see that 
$\pi$ is measure preserving: $\P = \Pp \circ \pi^{-1}$.

On the enlarged probability space, we define the maps
\begin{alignat}{2}
\label{e:Ttildedef}
\tilde T &: \GNdp \to \GNdp, &\quad
\tilde T(\cal G, {\bf S}) &:= (T_{\bf S}(\cal G), T({\bf S})),
\\
\label{e:Tdef}
T &: \GNdp \to \GNd, &\quad
T(\cal G, {\bf S}) &:= \pi(\tilde T(\cG,{\bf S})) = T_{\bf S}(\cal G).
\end{alignat}
For the statement of the next proposition,
recall that $\GNd$ denotes the set of simple $d$-regular graph on $\qq{N}$.
For any finite graph $\cT$ on a subset of $\qq{N}$,
we define $\GNd(\cT):=\{\cG\in \GNd: \cB_\ell(o,\cG)=\cT\}$
to be the set of $d$-regular graphs whose radius-$\ell$ neighborhood of the vertex $o$ in $\cG$ is $\cT$.

\begin{proposition} \label{prop:reverse}
For any graph $\cT$, we have
\begin{equation}\label{e:tildeTT}
  \tilde T(\iota(\GNd(\cT))) = \iota(\GNd(\cT)),
\end{equation}
and $\tilde T$ is an involution: $\tilde T \circ \tilde T = \id$.
\end{proposition}

\begin{proof}
Since our local resampling does not change the radius $\ell$ neighborhood of $o$, from our construction, $\tilde T(\iota(\GNd(\cT))) \subset \iota(\GNd(\cT))$. Next we show that $\tilde T$ is an involution, then it follows $ \tilde T(\iota(\GNd(\cT))) = \iota(\GNd(\cT))$.

To verify that $\tilde T$ is an involution,
let $(\cG, {\bf S}) \in \GNdp$ and abbreviate $(\tcG, \td {\bf S}) = \tilde T(\cG, {\bf S})$.
Then, thanks to \eqref{e:tildeTT}, the edge boundaries of the radius-$\ell$ neighborhoods of $o$ have the same number of edges $\mu$
in $\tcG$ and $\cG$.
Moreover, we can choose the  enumeration of the boundary of the $\ell$-ball in $\tcG$ such that,
for any $\al \in \qq{1,\mu}$, we have $T_\al(\vec S_\al) \in \sS_\al(\tcG)$.
Define
\begin{align*}
\td \As_{\tilde {\bf S}}:=\{\al\in \qq{1,\mu}: I_\al(\tcG,\td{\bf S}) J_\al(\tcG,\td{\bf S})=1\}.
\end{align*}
We claim that $\td \As_{\tilde {\bf S}} = \As_{\bf S}$.

First, by the definition of switchings, it is easy to see that $J_\al(\tcG, \td {\bf S})=J_\al(\cG,{\bf S})$. It suffices to verify that $I_\al(\tcG, \td{\bf S}) = I_\al(\cG,{\bf S})$ also holds for all $\al \in \qq{1,\mu}$.
If $J_\al(\tcG, \td {\bf S})=J_\al(\cG,{\bf S})=1$ and $I_\al(\cG,{\bf S})=1$, then $\{a_\beta,b_\beta,c_\beta\}\not\in \cB_{\fR/4}(\{a_\al,b_\al,c_\al\},\cGT)$ for all $\beta\neq \al$. Moreover, $\al$ is switchable and $\cB_{\fR/4}(\{a_\al,b_\al,c_\al\},\cG^{(\bT)})\cup \{b_\al,c_\al\}=\cB_{\fR/4}(\{a_\al,b_\al,c_\al\},\tcG^{(\bT)})\cup \{a_\al,b_\al\}$ have zero excess. This implies $I_\al(\tcG,{\tilde{\bf S}})=1$. 
Otherwise if $J_\al(\tcG, \td {\bf S})=J_\al(\cG,{\bf S})=0$, or $I_\al(\cG,{\bf S})=0$ then index $\al$ is not switchable. Therefore we have that the $\fR/4-$neighborhoods of $\{a_\al, b_\al, c_\al\}$ never change, i.e. $\cB_{\fR/4}(\{a_\al,b_\al,c_\al\}, \cG^{(\bT)})=\cB_{\fR/4}(\{a_\al,b_\al,c_\al\}, \tcG^{(\bT)})$. In this case, we also have $I_\al(\tcG, \td{\bf S}) = I_\al(\cG,{\bf S})$.
In summary, we have verified the claim $\td \As_{\tilde {\bf S}}=\As_{\bf S}$. By the definition of our switchings,
it follows that $T(\td {\bf S})=\bf S$ and $T_{\td{\bf S}}(\tcG)=\cG$.
Therefore $\tilde T$ is an involution.
\end{proof}

\begin{remark}\label{r:radiusR}
We remark that from the proof of Proposition \ref{prop:reverse}, the radius $\fR/4$ neighborhoods of the vertex $o$ in the original graph $\cG$ and in the switched graph $\tcG$ are isomorphic to each other.
\end{remark}

\begin{proposition}
$\tilde T$ and $T$ are measure preserving: $\Pp \circ \tilde T^{-1} = \Pp$ and $\Pp \circ T^{-1} = \P$.
\end{proposition}

In other words, that $T$ is measure preserving means that
if $\cG$ is uniform over $\GNd$, and given $\cG$, we choose $\bf S$ uniform over $\sS(\cG)$,
then $T_{\bf S}(\cal G)$ is uniform over $\GNd$.

\begin{proof}
We decompose the enlarged probability space according to the radius-$\ell$ neighborhood of $o$ as
\begin{equation} \label{e:GNdpdecomp}
  \GNdp = \bigcup_{\cT} \GNdp(\cT), \quad \text{where } \GNdp(\cT) = \iota(\GNd(\cT)).
\end{equation}
Notice that, given any $\cT$, the size of the set $\sS(\cG)$ is (by construction) independent of the graph $\cG \in \GNd(\cT)$. 
Therefore, given any $\cT$, the restricted measure $\Pp|{\GNdp(\cT)}$ is uniform, i.e.,
proportional to the counting measure on the finite set $\GNdp(\cT)$.
Since, by Proposition~\ref{prop:reverse},
the map $\tilde T$ is an involution on $\GNdp(\cT)$, it is in particular a bijection
and preserves the uniform measure $\Pp|{\GNdp(\cT)}$.
Since $\tilde T$ acts diagonally in the decomposition \eqref{e:GNdpdecomp},
this implies that the map $\tilde T$ preserves the measure $\Pp$.
Since $\P = \Pp \circ \pi^{-1}$ and $T = \pi \circ \tilde T$,
it immediately follows that also $T$ is measure preserving:
\begin{equation*}
\Pp \circ T^{-1} = \Pp \circ \tilde T^{-1} \circ \pi^{-1} = \Pp \circ \pi^{-1} = \P,
\end{equation*}
as claimed.
\end{proof}

\section{Weak Local Law and Delocalization of Eigenvectors}
\label{sec:weak}

In this section we prove the following theorem. It states that with high probability the Green's function of a random $d$-regular graph $\cG$ at vertices $i,j$, can be approximated by the Green's function of the tree extension of a radius-$r$ neighborhood $\cB_r(i,j,\cG)$ of vertices $i,j$. The delocalization of eigenvectors Theorem \ref{thm:delocalizationev} follows from Theorem \ref{thm:mrmsc} by a standard argument  \cite[Section 18.5]{MR3699468}.

Before stating Theorem \ref{thm:mrmsc}, we need to introduce some parameters

\begin{definition}[Choices of parameters]\label{cop}
We fix parameters $0<\fc<1$, $\fa\geq 12$,  $\fb\geq 25\fa$,  $0 < \fr\leq \fc/32$ and a large constant $\fd$. Let $\fR=(\fc/4)\log_{d-1}N$,  $r=\fr\log_{d-1}N$ and $\ell\in\qq{\fa \log_{d-1}\log N, 2\fa\log_{d-1}\log N}$.  We restrict ourself to the spectral domain $z\in \bC^+$ with $\Im[z]\geq (\log N)^\fb/N$. 
With this  choice of parameters, 
\begin{align}\label{e:relation1}
\fR/8\geq r=\OO(\log N)\gg \ell=\OO(\log\log N), \quad 
(\log N)^\fa\leq (d-1)^\ell\leq (\log N)^{2\fa}.
\end{align}
\end{definition} 
For $z=E(z)+\ri \eta(z)\in \bC^+$, we define $\kappa(z)=\min\{|E(z)-2|, |E(z)+2|\}$, the distance from $z$ to the spectral edges $\pm 2$. 
From the quadratic equation of $m_{sc}(z)$, we have
\begin{align*}
\Im[m_{sc}(z)]\asymp 
\left\{
\begin{array}{cc}
\sqrt{\kappa(z)+\eta(z)} & \text{ if } |E(z)|\leq 2,\\
\eta(z)/\sqrt{\kappa(z)+\eta(z)} & \text{ if } |E(z)|\geq 2.
\end{array}
\right.
\end{align*}
For the rest of this  paper, we fix some error parameters: for $z\in\bC^+$,  
\begin{align}\label{e:defeps0}
\varepsilon_0(z):=
(\log N)^{8\fa}\left(\frac{1}{(d-1)^r}+\sqrt{\frac{\Im[\md(z)]}{N\eta(z)}}+\frac{1}{(N\eta(z))^{2/3}}\right),
\end{align}
and
\begin{align}\begin{split}
\label{e:defeps}
&\varepsilon(z)=\varepsilon_0(z) \text{ if } \varepsilon_0(z)\leq \frac{\kappa(z)+\eta(z)}{\log N},\\
&\varepsilon(z)=(\log N)^4\varepsilon_0(z) \text{ if } \varepsilon_0(z)> \frac{\kappa(z)+\eta(z)}{\log N}. 
\end{split}\end{align}\index{$\varepsilon_0(z)$}\index{$\varepsilon(z)$}

We remark that  the error $\sqrt{\Im[\md(z)]/(N\eta(z))}$ in \eqref{e:defeps0} is a common error appearing in random matrix theory;
the error $1/(d-1)^r$ is to take into account of the approximation of the Green's function.
From the definition \eqref{e:defeps} of $\varepsilon(z)$, there is a dichotomy, i.e., either $\varepsilon(z)\leq (\kappa(z)+\eta(z))/\log N$ or $\varepsilon(z)  \geq  (\log N)^3(\kappa(z)+\eta(z))$.
This fact will be crucial in the proof of Proposition \ref{prop:bootstrap}.
We also define two additional parameters 
\begin{align}\label{e:defphi}
\begin{split}
&\varepsilon'(z)=(\log N)^3\varepsilon(z),\\
&\varphi(z)=(\log N)^{ 2\fa}\sqrt{\frac{\Im[\md(z)]+\varepsilon'(z)+\varepsilon(z)/\sqrt{\kappa(z)+\eta(z)+\varepsilon(z)}}{N\eta(z)}}.
\end{split}
\end{align}
where the choice of $\varphi(z)$ will be clear in \eqref{Vcondition}.\index{$\varepsilon'(z)$}\index{$\varphi(z)$} Notice that from our choice of parameters, on the spectral domain $z\in \bC^+$ with $\Im[z]\geq (\log N)^\fb/N$, we have 
\begin{align}\label{e:relation2}
\varepsilon(z)\leq \frac{1}{(\log N)^{4\fa}},\quad \varphi (z) \leq \frac{\varepsilon(z)}{(\log N)^{4\fa}}.
\end{align}

\begin{theorem}\label{thm:mrmsc}
Fix $d\geq 3$, $0<\fc<1$, $\fR=(\fc/4)\log_{d-1}N$  and recall the set of  radius-$\fR$ tree like graphs $\bar\Omega\subset \GNd$ from Definition \eqref{def:barOmega}.  For any large $\fC>0$ and $N$ large enough, with probability $1-\OO(N^{-\fC})$ with respect to the uniform measure on $\bar\Omega$,  the Green's function of $\cG\in \bar\Omega$ satisfies
\begin{align}\label{e:locallaw}
  \left|G_{ij}(\cG,z)-G_{ij}(\Ext(\cB_{r}(i,j,\cG),m_{sc}(z)),z)\right|\lesssim  \frac{\varepsilon(z)}{\sqrt{\kappa(z)+\eta(z)+\varepsilon(z)}},  
\end{align}
for any vertices $i,j\in\qq{N}$, and 
the Stieltjes transform of its empirical eigenvalues satisfies
\begin{align} \label{e:localmckay}
  \left|m_N(z) - \md(z)\right| \lesssim \frac{\varepsilon(z)}{\sqrt{\kappa(z)+\eta(z)+\varepsilon(z)}},
 \end{align}
 uniformly for $z\in \bC^+$, with $\Im[z]\geq (\log N)^\fb/N$.
 
\end{theorem}

We remark that while Theorem~\ref{thm:mrmsc} shows that the spectral density (or its Stieltjes transform, which is the trace of the Green's function) does concentrate around the Kesten-Mckay law. 
The individual entries of the Green's function of the random $d$-regular graph with bounded degree is approximated  by the Green's function of a neighborhood of these entries. Recall that the off-diagonal entries of the Green's function of  a Wigner matrix is uniformly small. This property clearly fails for 
the Green's function of the random $d$-regular graph with a fixed degree.

In Section \ref{sec:outline}, we give the proof of Theorem \ref{thm:mrmsc}, which uses Propositions \ref{p:Omega+} and \ref{p:Omega-} as input. The proofs of Propositions \ref{p:Omega+} and \ref{p:Omega-} are given in Sections \ref{s:switchstable} and \ref{s:improveG} respectively.

\subsection{Notations and Definitions}
\label{sec:outline-parameters}

To study the change of graphs before and after local resampling, we define the following graphs (which are not $d$-regular).
\begin{itemize}
\item $\cG$ is the original unswitched graph;
\item $\tcG$ is the switched graph $T_{\bf S}(\cG)$\index{$\tilde\cG, \cGT,\tcGT,\cG^{(\bT\bW_{\bf S})},\tilde\cG^{(\bT\bW_{\bf S})}$} ;
\item $\cGT$ is the unswitched graph obtained from $\cG$ with vertices $\T$ removed;
\item $\tcGT$ is the switched graph obtained from $\tcG$ with vertices $\bT$ removed;
\item $\cG^{(\bT\bW_{\bf S})}=\tilde\cG^{(\bT\bW_{\bf S})}$ is the intermediate graph obtained from $\cGT$ with vertices $\bW_{\bf S}=\{b_\al: \al\in \As_{\bf S}\}$ removed, or equivalently obtained from $\tcGT$ with vertices $\bW_{\bf S}$ removed.
\end{itemize}
Following the conventions of Remark \ref{r:convention},
the deficit functions of these graphs are given by \eqref{e:removing}. More explicitly, the deficit functions are simply $g(v)=d-\deg(v)$. We abbreviate their Green's functions by $G$, $\GT$, $\tG^{(\bT\bW_{\bf S})}=G^{(\bT\bW_{\bf S})}$, $\tGT$, and $\tG$ respectively\index{$G, \GT, \tG^{(\bT\bW_{\bf S})},G^{(\bT\bW_{\bf S})},\tGT$}.
The local resampling as defined in Section \ref{sec:switch}
has a smaller effect in $\cGT$ than they do in $\cG$. Indeed, in the original graph $\cG$,
simple switchings have the effect of removing two edges and adding two edges, while in $\cGT$ simple switchings only remove
the edges $\{b_\al,c_\al\}_{\al\in \As_{\bf S}}$ and add the edges $\{a_\al,b_\al\}_{\al\in \As_{\bf S}}$.

The small distance behavior is captured in terms of cycles in neighborhoods of radius $\OO(\fR)$.
For any graph, we recall that excess is  the smallest number of edges that must be removed to yield a graph with no cycles
(a forest). We also recall the set of  radius-$\fR$ tree like graphs $\bar\Omega\subset \GNd$ from Definition \eqref{def:barOmega}. 
We will also need the following set $\bar\Omega^+$. Roughly speaking, for any fixed large constant $\fd$, 
 from a $d$-regular graph in $\cG\in\bar \Omega$, with probability $1-\OO(N^{-\fd})$ with respect to the randomness of the resampling data $\bfS$, the switched graph $T_{\bfS}(\cG)\in \bar \Omega^+$. The $d$-reguar graphs in  $\bar \Omega^+$ are also tree-like at small distances. However, compared with $d$-regular graphs in  $\bar \Omega$, they are a bit more deviated from trees locally, depending on $\fd$.

\begin{definition}\label{def:barOmega1}
We define the  set $\bar\Omega^+\subset \GNd$\index{$\bar\Omega^+$}  to consist  of graphs such that
\begin{itemize}
\item
the radius-$ \fR/8$ neighborhood of any vertex has excess at most $\n_d$, the radius-$\fR/2$ neighborhood of any vertex has excess at most $C_\fd$, where $C_\fd$ is a constant depending on $\fd$ and explicitly specified in the proof of Proposition \ref{p:newG}.

\item
the number of vertices that have a radius-$\fR$ neighborhood that contains a cycle is at most $2N^{\fc}$.
\end{itemize} 
\end{definition}

Notice that the sets $\bar\Omega$ and  $ \bar\Omega^+$ are subsets of $\GNd$ defined by deterministic properties. Moreover,  it is  clearly that 
 $\bar\Omega \subset  \bar\Omega^+$.  
For any graph in  $\bar\Omega,\bar\Omega^+$, radius-$ \fR/8$ neighborhood of any vertex has excess at most $\n_d$.  The Green's functions of their extensions with weight close to $m_{sc}(z)$ are stable, and satisfy Proposition \ref{p:subgraph}.
However, the small distance  graphic   behavior captured by sets $\bar\Omega,\bar\Omega^+$ does no guarantee nice spectral properties, i.e. stability of their Green's function.
We need the following notion of  spectral regular graphs.
For any finite $d$-regular graph $\cG$, we introduce the following quantity
\begin{equation} \label{e:IG}
Q(\cG,z)=\frac{1}{Nd}\sum_{i\in \qq N}\sum_{j: i\sim j} G_{ii}^{(j)}(\cG,z), 
\end{equation}
where $G^{(j)}(\cG,z)$ is the Green's function of the graph obtained from $\cG$ by removing the vertex $j$.\index{$Q(\cG,z)$}
We define the following set of $d$-regular graphs, whose Green's function has good estimates.
\begin{definition}\label{def:Omega}
For $z\in\bC^+$, we define the 
set of  spectral regular graphs,  $\Omega(z) \subset \bar\Omega$ \index{$\Omega(z)$}, to be the set of graphs $\cG$
such that 
\begin{align}\begin{split}\label{e:defOmega}
 &|Q(\cG,z)-m_{sc}(z)|\leq \frac{\varepsilon(z)}{\sqrt{\kappa(z)+\eta(z)+\varepsilon(z)}}, \\
  &\left|G_{ij}(z)-G_{ij}(\Ext(\cB_r(i,j,\cG),Q(\cG,z)),z)\right|
  \leq \varepsilon(z),
\end{split}\end{align}
for any two vertices $i,j\in\qq{N}$.
\end{definition}

Notice that the last inequality in particular holds for $i=j$ and thus the error for the diagonal terms of  Green's  function is $\e(z) \ll \frac{\varepsilon(z)}{\sqrt{\kappa(z)+\eta(z)+\varepsilon(z)}}$. In other words, we are able to approximate  Green's  function of the original graph with that of a  tree like graph 
provide that extension uses the weight $Q(\cG,z)$. A key novelty of this paper is to derive a self-consistent equation for the quantity $Q(\cG,z)$: $Q(\cG,z)\approx Y_\ell(Q(\cG,z),z)$, where the function $Y_\ell(\cdot, z)$ is defined in \eqref{def:Y}.
The self-consistent equation captures the square root behavior at the spectral edge of the empirical density. This is reminiscent of the quadratic equation for the Stieltjes transform of semicircle distribution. This enables us to improve the results in \cite{MR3962004} to the spectral edge.

 We have so far defined three subsets of  $\GNd$:    $\Omega(z) \subset  \bar\Omega  \subset \bar \Omega^+$.  
The set  $ \bar\Omega^+   $ consists of  graphs with regularity conditions weaker than those of  $  \bar\Omega  $.  And  the 
 spectral regular graph set   $\Omega(z)$,  consists of  $d$-regular graphs satisfying  certain  resolvent  estimates.

\subsection{Two key Propositions}
\label{sec:outline}

The key inputs to prove 
Theorem \ref{thm:mrmsc} are the following two Propositions:  Proposition \ref{p:Omega+} states for $\cG\in \Omega(z)$, with high probability with respect to the resampling data $\bf S$, we have good estimates for the Green's function of the switched graphs $\tcG=T_{\bfS}(\cG)$ and $\tcGT$. We denote the set of such graphs by $\Omega^+_o(z)$. Proposition \ref{p:Omega-} states for $\cG\in \Omega^+_o(z)$, with high probability with respect to the resampling data $\bf S$, the Green's function of the switched graph $\tcG=T_{\bf S}(\cG)$ has an improved estimate near vertex $o$.  For the rest of this section, all statements 
about local resampling refer to  a neighborhood centered at a vertex labelled by $o$. 
The proofs of Propositions \ref{p:Omega+} and \ref{p:Omega-} will occupy Sections \ref{s:switchstable} and \ref{s:improveG}.

\begin{definition}\label{def:Omegao+}
For $z\in\bC^+$ with $\Im[z]\geq (\log N)^\fb/N$ and vertex $o\in \qq{N}$,  we define the 
set of  $d$-regular graphs,  $\Omega_o^+(z) \subset \bar\Omega^+$ \index{$\Omega_o^+(z)$}, to be the set of graphs $\cG$
such that its
Green's function  satisfies the following estimates:
\begin{align}
 &  \quad |Q(\cG, z)-m_{sc}(z)|\leq \frac{\fC_1\varepsilon(z)}{\sqrt{\kappa(z)+\eta(z)+\varepsilon(z)}}, \label{e:weakQm} \\
&  \quad 
|G_{ij}(z)-G_{ij}(\Ext(\cB_r(i,j,\cG), Q(\cG,z)),z)|\leq \frac{\fC_1}{\log N}, \quad i,j\in \qq{N}, \label{e:weakrigid1} \\
&
\quad |\GT_{ij}(z)-G_{ij}(\Ext(\cB_r(i,j,\cGT), Q(\cG,z),z))|\leq \fC_1(\log N)^3\varepsilon(z), \quad i,j\in \qq{N}\setminus\bT,  \label{e:weakrigid2}
\end{align}%
where $\bT=\bB_\ell(o,\cG)$ (which may not be a tree) and the large constant $\fC_1$ will be chosen in Proposition \ref{p:Omega+}.
\end{definition}

The following proposition  says from any graph $\cG\in \Omega(z)$, with high probability with respect to the switching data, the switched graph is in $ \Omega_o^+(z)$. We  notice that  \eqref{e:weakQm} is weaker than the statement in \eqref{e:defOmega}, while 
\eqref{e:weakrigid1} is much weaker than the corresponding bound in \eqref{e:defOmega} before the switching.  We have that $\Omega(z)\subset  \Omega_o^+(z)$.
The statement  \eqref{e:weakrigid2} is weaker than  the corresponding bound in \eqref{e:defOmega}  if we ignore the removal of $\bT$. 

\begin{proposition}\label{p:Omega+}
Let $d\geq 3$, $\omega_d$ as in Definition \ref{d:defwd} and $z\in \bC^+$ with $\Im[z]\geq (\log N)^\fb/N$. The sets $\Omega(z)$ and $\Omega_o^+(z)$ are defined in Definitions \ref{def:Omega} and \ref{def:Omegao+}.
Then for any $\cG\in \Omega(z)$, there exist constants $\fd, \fC_1>0$ (large and depending  on $d$) and 
events $F_1(\cG), F_2(\cG) \subset \sS(\cG)$ with $\P_{\cG}(F_1(\cG)\cap F_2(\cG))= 1-\OO(N^{-\fd})$, 
explicitly defined in Proposition \ref{p:newG} and Section \ref{sec:stability} below, such that for any switching data $\bfS\in F_1(\cG)\cap F_2(\cG)$, 
the switched graph $\tcG=T_{\bf S}(\cG)\in \Omega_o^+(z)$.

\end{proposition}

\begin{definition}\label{def:Omegao-}
For $z\in\bC^+$ with $\Im[z]\geq (\log N)^\fb/N$ and vertex $o\in \qq{N}$,  we define the set of  $d$-regular graphs,  $\Omega_o^-(z) \subset \GNd$ \index{$\Omega_o^-(z)$}, to be the set of graphs $\cG$
such that its
Green's function  satisfies the following estimates:
There exists  a large constant $\fC_2$ which will be chosen in Proposition \ref{p:Omega-},
\begin{align}\label{e:weakone}
&|Q(\cG,z)-m_{sc}(z)|\leq \frac{\fC_2\varepsilon(z)}{\sqrt{\kappa(z)+\eta(z)+\varepsilon(z)}},\\
\begin{split}\label{e:improverigid1}
  &|G_{oi}(z)-G(\Ext(\cB_r(o,i,\cG),  Q(\cG, z)),z)|
  \leq \fC_2\left(
  \left(\frac{\log N }{(d-1)^{\ell/2}}+\frac{(\log N)^2}{(d-1)^{\ell-\ell_i/2}}\right)\varepsilon'(z)\right.\\
&\left.+\log N (\sqrt{\kappa(z)+\eta(z)}|Q(\cG,z)-m_{sc}(z)|+|Q(\cG,z)-m_{sc}(z)|^2)\right), \quad \forall  i\in \bT.
\end{split}
\\
\begin{split}
\label{e:improverigid1.5}
  &|G_{oi}(z)-G(\Ext(\cB_r(o,i,\cG),  Q(\cG, z)),z)|
  \leq\fC_2\left(
  \frac{(\log N)^2 \varepsilon'(z)}{(d-1)^{\ell/2}}\right.\\
&\left.+\log N (\sqrt{\kappa(z)+\eta(z)}|Q(\cG,z)-m_{sc}(z)|+|Q(\cG,z)-m_{sc}(z)|^2)\right), \quad \forall  i\in \qq{N}\setminus\bT,
\end{split}\end{align}%
where $\ell_i=\dist_{\cG}(o,i)$, $\bT=\bB_\ell(o, \cG)$. If the vertex $o$ has radius-$\fR$ tree neighborhood in the graph $\cG$, the following holds
\begin{align}
\label{e:improverigid2}
  \frac{1}{d}\sum_{i: o\sim i}|G_{oo}^{(i)}(z)-Y_\ell(Q(\cG, z),z)| \leq \fC_2\frac{\log N\varepsilon'(z)}{(d-1)^{\ell/2}}.
\end{align}

\end{definition}

The following proposition  says from any graph $\cG\in \Omega_o^+(z)$, with high probability with respect to the switching data, the switched graph is in $ \Omega_o^-(z)$. We  notice that although \eqref{e:weakone} is weaker than the statement in \eqref{e:defOmega}, the estimates \eqref{e:improverigid1} and \eqref{e:improverigid1.5} for Green's function centered at vertex $o$ is better than \eqref{e:defOmega} before the switching, thanks to the $1/(d-1)^{\ell/2}$ factor.
In other words, local resampling centered around vertex $o$ improves the Green's function estimates centered at vertex $o$. 

\begin{proposition}\label{p:Omega-}
Let $d\geq 3$,  $\omega_d$ as  in Definition \ref{d:defwd} and $z\in \bC^+$ with $\Im[z]\geq (\log N)^\fb/N$. The sets $\Omega_o^+(z)\subset \bar\Omega^+$ are defined in Definitions \ref{def:barOmega1} and \ref{def:Omegao+}.
Then for any $\cG\in \Omega_o^+(z)$, 
there exist constants $\fd,  \fC_2>0$ (large and  depending  on $d$) and 
an event $F_3(\cG) \subset \sS(\cG)$ with $\P_{\cG}(F_3(\cG))= 1-\OO(N^{-\fd})$, 
explicitly defined in Proposition \ref{tGconcentration} below, such that,  for any switching data $\bfS\in F_3(\cG)$, 
 the switched graph $\tcG=T_{\bf S}(\cG)\in \Omega_o^-(z)$. 

\end{proposition}

%

\subsection{Proof of Theorem  \ref{thm:mrmsc}}

Assuming the propositions in the previous subsection, Theorem \ref{thm:mrmsc} is an easy consequence of the following proposition. Instead of approximating the Green's function $G_{ij}$ by the Green's function of the tree extension of a neighborhood $\cB_r(i,j,\cG)$ of vertices $i,j$, we approximate $G_{ij}$ by the extension of $\cB_r(i,j,\cG)$ with the weight $Q(\cG,z)$. The approximation with $Q(\cG,z)$ as the boundary weight, leads naturally to the quadratic self-consistent equation of $Q(\cG,z)$, which can be used to obtain improved estimates of $Q(\cG,z)$ for both the bulk and edge regions. 

\begin{proposition}\label{thm:mrQ}
Fix $d\geq 3$, $0<\fc<1$, $\fR=(\fc/4)\log_{d-1}N$  and recall the set of  radius-$\fR$ tree like graphs $\bar\Omega\subset \GNd$ from Definition \eqref{def:barOmega}. 
For any large $\fC>0$ and $N$ large enough, with probability $1-\OO(N^{-\fC})$ with respect to the uniform measure on $\bar\Omega$,  
the quantity $Q(\cG,z)$ as defined in \eqref{e:IG} satisfies
\begin{align}\label{e:Q-msc}
 |Q(\cG,z)-m_{sc}(z)|\leq \frac{\varepsilon(z)}{\sqrt{\kappa(z)+\eta(z)+\varepsilon(z)}},
\end{align}
and the Green's function satisfies
\begin{align}\label{e:G-PQ}
  \left|G_{ij}(z)-G_{ij}(\Ext(\cB_{r}(i,j,\cG),Q(\cG,z)),z)\right|\leq \varepsilon(z), 
\end{align}
uniformly in $i,j\in\qq{N}$, and $z=E(z)+\ri\eta(z)\in \bC^+$, with $\Im[z]\geq (\log N)^\fb/N$.
\end{proposition}

\begin{proof}[Proof of Theorem \ref{thm:mrmsc}]
We prove that if $\cG\in \bar\Omega$, then \eqref{e:Q-msc} and \eqref{e:G-PQ} together imply Theorem \ref{thm:mrmsc}.
If $\cG\in \bar\Omega$, then $\cB_r(i,j,\cG)$ has excess at most $\omega_d$, Proposition \ref{p:subgraph} and \eqref{e:Q-msc} imply
\begin{align*}
&\phantom{{}={}}|G_{ij}(\Ext(\cB_{r}(i,j,\cG),Q(\cG,z)),z)-G_{ij}(\Ext(\cB_{r}(i,j,\cG), m_{sc}(z)),z)|\\
&\lesssim |Q(\cG,z)-m_{sc}(z)|\lesssim \frac{\varepsilon}{\sqrt{\kappa+\eta+\varepsilon}}.
\end{align*}
Combining with  \eqref{e:G-PQ}, we have with probability $1-\OO(N^{-\fC})$ with respect to the uniform measure on $\bar\Omega$, it holds
\begin{align*}
 \left|G_{ij}(z)-G_{ij}(\Ext(\cB_{r}(i,j,\cG), m_{sc}(z)),z)\right|
 \lesssim \varepsilon+\frac{\varepsilon}{\sqrt{\kappa+\eta+\varepsilon}}\lesssim \frac{\varepsilon}{\sqrt{\kappa+\eta+\varepsilon}},
\end{align*}
uniformly in $i,j\in\qq{N}$, and $z=E+\ri\eta\in \bC^+$, with $\eta\geq (\log N)^\fb/N$. If vertex $i$ has radius-$r$ neighborhood in $\cG$, then $G_{ii}(\Ext(\cB_{r}(i,\cG),m_{sc}(z)),z)=\md(z)$. By averaging over all vertices, and recalling for $\cG\in \bar\Omega$, all vertices have radius $\fR$ tree neighborhood except for at most $N^\fc$ of them, we get
\begin{align*}
|m_N(z)-\md(z)|\leq \frac{1}{N}\sum_{i\in \qq{N}}\left|G_{ii}(z)-\md(z)\right|\lesssim\frac{\varepsilon}{\sqrt{\kappa+\eta+\varepsilon}}.
\end{align*}
Theorem \ref{thm:mrmsc} follows.
\end{proof}

In the rest of this section, 
we prove Proposition \ref{thm:mrQ}. We recall the set of  spectral regular graphs  $\Omega(z) \subset \bar\Omega$ from Definition \ref{e:defOmega}, and we can use it to reformulate Proposition \ref{thm:mrQ} as that 
\begin{align}\label{e:eqmrQ}
\P\left(\bar\Omega\setminus\bigcap_{\Im[z]\geq (\log N)^\fb/N}\Omega(z)\right)\lesssim N^{-\fC}.
\end{align} 
We  prove \eqref{e:eqmrQ} by an iteration scheme. To explain the iteration scheme, we need to introduce the following set $\Omega^-(z)$. The defining relations \eqref{e:defOmega-} are similar to the defining relations \eqref{e:defOmega} of $\Omega(z)$, except that the righthand sides in \eqref{e:defOmega-} is smaller by a factor $1/2$. Thus it holds that $\Omega^-(z)\subset\Omega(z)$.
\begin{definition}\label{def:Omega-}
For $z\in\bC^+$, we define the 
set $\Omega^-(z) \subset \bar\Omega$\index{$\Omega^-(z)$} be the set of graphs $\cG$
such that 
\begin{align}\begin{split}\label{e:defOmega-}
 &|Q(\cG,z)-m_{sc}(z)|\leq \frac{\varepsilon(z)}{2\sqrt{\kappa(z)+\eta(z)+\varepsilon(z)}},\\
 & \left|G_{ij}(z)-G_{ij}(\Ext(\cB_r(i,j,\cG),Q(\cG,z)),z)\right|
  \leq \frac{\varepsilon(z)}{2}.
\end{split}\end{align}
for any two vertices $i,j\in\qq{N}$.
\end{definition}
In Proposition \ref{prop:betacomp}, we prove that for $|z| \geq 2d$ the event $\Omega^-(z)$ holds deterministically. Since the Green's functions $G_{ij}(z)$ are Lipschitz in $z$, conditioning on the event $\Omega^-(w)$, $\Omega(z)$ holds deterministically for all $\{z\in \bC^+: |z-w|\lesssim 1/N, \Im[z]\geq (\log N)^\fb/N\}$. In Proposition \ref{prop:bootstrap}, using Propositions \ref{p:Omega+} and \ref{p:Omega-}, we show that the difference of the two sets $\Omega(z), \Omega^{-}(z)$ is negligible, i.e. $\bP\left(\Omega(z)\setminus\Omega^-(z)\right)=\OO(N^{-\fd+1})$. As a consequence, on the event $\Omega^{-}(z)$, by losing a small probability, the event $\Omega^-(z-\ri/N)$ holds. Starting from $|z| \geq 2d$, where $\Omega^-(z)$ holds deterministically, it follows by an iterative scheme, with high probability, the events $\Omega^-(z)$ holds for all $z$ in the form
$\left(-2d +\alpha/N\right)+\ri\left(2d-\beta/N\right)$ with 
$\al\in \qq{0, 4dN}$ and $\beta\in \qq{0, 2dN-(\log N)^\fb}$. Using again the Green's functions $G_{ij}(z)$ are Lipschitz in $z$, this implies \eqref{e:eqmrQ}, and Proposition \ref{thm:mrQ} follows.

\begin{proposition} \label{prop:betacomp}
Let $d\geq 3$ and $\omega_d$ as in Definition \ref{d:defwd}. 
Let $\cG\in \bar \Omega$. Then for any $z \in \bC^+$ with $|z|\geq 2d$, the quantity $Q(\cG,z)$ as defined in \eqref{e:IG} satisfies
\begin{align*}
|Q(\cG,z)-m_{sc}(z)|\leq \frac{\varepsilon(z)}{2\sqrt{\kappa(z)+\eta(z)+\varepsilon(z)}},
\end{align*}
and  the Green's function satisfies
\begin{align*}
  \left|G_{ij}(z)-G_{ij}(\Ext(\cB_{r}(i,j,\cG),Q(\cG,z)),z)\right|\leq \frac{\varepsilon(z)}{2}, 
\end{align*}
for any $i,j\in\qq{N}$.
\end{proposition}

\begin{proof}[Proof of Proposition \ref{prop:betacomp}]
It follows by the Combes--Thomas method \cite{MR0391792}. For any finite simple graph $\cal G$ with degree bounded by $d$, 
and any $z$ with $|z|\geq 2d$,
\begin{equation}\label{betabound}
 |G_{ij}(z)| \leq 1/d, \qquad
 |G_{ij}(z)|\leq (d-1)^{ -\dist_\cal G(i,j)/2}.
\end{equation}
Let $\cG$ be a $d$-regular graph on $N$ vertices, with excess at most $\omega_d$ in any radius-$\fR$ neighborhood. Using \eqref{betabound} as input, by the same argument as in \cite[Propsition 6.1]{MR3962004}, we have for any $z \in \bC^+$ with $|z|\geq 2d$, and any $i,j \in \qq{N}$, the Green's function of $\cG$ satisfies
 \begin{align}\begin{split}\label{e:mscapprox}
  &\phantom{{}={}}\left|G_{ij}(z)-G_{ij}(\Ext(\cB_r(i,j,\cal G),m_{sc}(z)),z)\right|\\
  &\leq \left(\frac{1}{(d-1)^{r}}+\sqrt{\frac{\Im[\md(z)]}{N\eta}}+\frac{1}{(N\eta)^{2/3}}\right)
  \lesssim \frac{\varepsilon(z)}{\log N}.
\end{split} \end{align}

If the radius-$\fR$ neighborhood of vertex $i$ is a truncated $d$-regular tree, then the neighbhorhood $\cB_r(i,\cG)$ is a truncated $d$-regular tree with root vertex $i$. \eqref{e:Gtreemsccopy} gives  $G_{ii}(\Ext(\cB_{r}(i,\cG), m_{sc}(z)),z)=\md(z)$, and 
\begin{align}\label{e:itreecasemd}
|G_{ii}(z)-G_{ii}(\Ext(\cB_{r}(i,\cG), m_{sc}(z)),z)|
=|G_{ii}(z)-\md(z)|
\lesssim \varepsilon(z)/\log N.
\end{align}
If the radius-$\fR$ neighborhood of vertex $i$ is a truncated $d$-regular tree, for any $j\sim i$, $\cB_r(i,\cG^{(j)})$ is a truncated $(d-1)$-ary tree with root vertex $i$.  \eqref{e:Gtreemsccopy2} gives  $G_{ii}(\Ext(\cB_{r}(i,\cG^{(j)}), m_{sc}(z)),z)=m_{sc}(z)$.  Using the Schur complement formula \eqref{e:Schurixj}  and \eqref{e:mscapprox} we get
\begin{align}\label{e:itreecase}
|G_{ii}^{(j)}(z)-G_{ii}(\Ext(\cB_{r}(i,\cG^{(j)}),  m_{sc}(z)),z)|
=|G_{ii}^{(j)}(z)-m_{sc}(z)|
\lesssim \varepsilon(z)/\log N.
\end{align}
For $\cG\in \bar\Omega$, the number of vertices whose radius-$\fR$ neighborhoods is not a tree is at most $N^\fc$. We can average  \eqref{e:itreecase} over all $i\in \qq{N}$,
\begin{align}\label{e:Q-mbound0}
&|Q(\cG,z)-m_{sc}(z)|
= \frac{1}{Nd}\sum_{i\in\qq{N}}\sum_{j:i\sim j}|G_{ii}^{(j)}-m_{sc}(z)|
\lesssim \varepsilon(z)/\log N, 
\end{align}
and thanks to Proposition \ref{p:subgraph}
 \begin{align}\label{largedgreen}
  \left|G_{ij}(\Ext(\cB_r(i,j,\cal G),m_{sc}(z)),z)-G_{ij}(\Ext(\cB_r(i,j,\cal G),Q(\cG, z)),z)\right|\lesssim \varepsilon(z)/\log N.
 \end{align}
This finishes the proof of Proposition \ref{prop:betacomp}.

\end{proof}

\begin{proposition} \label{prop:bootstrap}
Let $d\geq 3$, $\omega_d$ as in Definition \ref{d:defwd} and the sets $\Omega(z), \Omega^-(z)$ as in Definitions \ref{def:Omega} and \ref{def:Omega-} respectively. For any  $z\in \bC^+$ with $\Im[z]\geq (\log N)^\fb/N$, we have
\begin{align*}
\mathbb P\left(\Omega(z)\setminus\Omega^-(z)\right)=\OO(N^{-\fd+1}).
\end{align*}
\end{proposition}

\begin{proof}[Proof of Proposition \ref{prop:bootstrap}]

We recall the sets $\Omega_o^+(z), \Omega_o^-(z)$ from  Definitions \ref{def:Omegao+} and \ref{def:Omegao-}.
For simplicity of notations, we write $\Omega=\Omega(z)$, $\Omega^-=\Omega^-(z)$, $\Omega_o^+=\Omega_o^+(z)$ and $\Omega_o^-=\Omega_o^-(z)$. Then we have $\Omega \subset \Omega_o^+\subset \GNd$ and $\Omega^-_o\subset \GNd$. 
We take $\ell\in\qq{\fa \log_{d-1}\log N, 2\fa\log_{d-1}\log N}$ such that $|1+(\msc(z))^2+(\msc(z))^4+\cdots+(\msc(z))^{2\ell}|\gtrsim 1$.
In the following we first show that
\begin{align}\label{e:smalllose}
  \P(\Omega \setminus (\Omega \cap \Omega^-_o)) \lesssim N^{-\fd}.
\end{align}
We recall the maps $\tilde T$ and $T$ from \eqref{e:Tdef}. We define $\tilde\Omega = T^{-1}(\Omega)$, $\tilde\Omega_o^-= T^{-1}(\Omega_o^-)$ and $\tilde \Omega_o^+ = T^{-1}(\Omega_o^+)$.
Since $T$ is measure preserving, we have
\begin{equation}\label{e:lose0}
 \P(\Omega\setminus \Omega_o^-)=\tilde \P( T^{-1}(\Omega \setminus \Omega_o^-))
  =\tilde \P( T^{-1}(\Omega) \setminus T^{-1}(\Omega_o^-) )=\tilde \P( \tilde\Omega \setminus \tilde \Omega^-_o).
\end{equation}
We can reformulate Proposition \ref{p:Omega+} as
$\P_{\cG}(T_{\bf S}(\cal G) \in  \GNd\setminus \Omega_o^+) \lesssim N^{-\fd}$ for all $\cal G \in \Omega$. It implies
\begin{align}\label{e:lose1}
\tilde\P(\tilde T (\tilde\Omega)\setminus \tilde\Omega_o^+)=\tilde\P(\tilde \Omega\setminus \tilde T (\tilde \Omega_o^+))\lesssim N^{-\fd},
\end{align}
since $\tilde T$ is a measure preserving involution. Similarly, Proposition \ref{p:Omega-} can be reformulated as
$\P_{\cG}(T_{\bf S}(\cal G) \in  \GNd\setminus \Omega^-_o) \lesssim N^{-\fd}$ for all $\cal G \in \Omega_o^+$, which implies 
\begin{align}\label{e:lose2}
\tilde \P(\tilde T (\tilde \Omega_o^+)\setminus \tilde \Omega_o^-)\lesssim N^{-\fd}.
\end{align}
The claim \eqref{e:smalllose} follows from combining \eqref{e:lose0}, \eqref{e:lose1} and \eqref{e:lose2},
\begin{align*}
\P(\Omega\setminus \Omega_o^-)=
\tilde \P( \tilde\Omega \setminus \tilde \Omega^-_o)\leq
\tilde\P(\tilde\Omega\setminus \tilde T(\tilde \Omega_o^+))+\tilde \P(\tilde T (\tilde \Omega_o^+)\setminus \tilde \Omega_o^-)\lesssim N^{-\fd}.
\end{align*}

By a union bound over all the indices $o\in\qq{N}$, \eqref{e:smalllose} gives that
\begin{align*}
\P\left(\Omega\setminus \Omega\bigcap_{o\in \qq{N}}\Omega_o^-\right)\lesssim N^{-\fd+1}.
\end{align*}
In the following we prove that $\Omega\cap_{o\in\qq{N}}\Omega^-_o\subset \Omega^-$, thus Proposition \ref{prop:bootstrap} follows. For $\cG\in\Omega\cap_{o\in\qq{N}}\Omega^-_o$, since $\Omega\in \bar\Omega$, the number of vertices whose radius-$\fR$ neighborhoods is not a tree is at most $N^{\fc}$.  We can average \eqref{e:improverigid2} over all $o\in \qq{N}$,
\begin{align*}
Q(\tcG,z)-Y_\ell(Q(\tcG,z),z)=\OO\left(\frac{\log N \varepsilon'}{(d-1)^{\ell/2}}\right).
\end{align*}
Thanks to Proposition \ref{p:recurbound}, we can expand $Y_\ell(Q(\tcG,z),z)$ around $m_{sc}(z)$ and get
\begin{align}\begin{split}\label{e:expand}
&\OO\left((\log N)^2|Q(\tcG,z)-m_{sc}(z)|^3+\frac{\log N\varepsilon'}{(d-1)^{\ell/2}}\right)=(1-(m_{sc}(z))^{2\ell+2})(Q(\tcG,z)-m_{sc}(z))\\
&-\msc^{2\ell+2}(z)\md(z)\left(\frac{1-\msc^{2\ell+2}(z)}{d-1}+\frac{d-2}{d-1}\frac{1-\msc^{2\ell+2}(z)}{1-\msc^2(z)}\right)(Q(\tcG,z)-m_{sc}(z))^2.
\end{split}\end{align}%
We consider the quadratic equation $ax^2+bx+c=0$, with
\begin{align}\begin{split}\label{e:coeff}
&a=-\msc^{2\ell+2}(z)\md(z)\left(\frac{1-\msc^{2\ell+2}(z)}{d-1}+\frac{d-2}{d-1}\frac{1-\msc^{2\ell+2}(z)}{1-\msc^2(z)}\right)(\Delta-\msc(z))^2,\\
&b=(1-(m_{sc}(z))^{2\ell+2}),\quad c=\OO\left((\log N)^2|Q(\tcG,z)-m_{sc}(z)|^3+\frac{\log N\varepsilon'}{(d-1)^{\ell/2}}\right).
\end{split}\end{align}
We recall that by our choice of $\ell\in\qq{\fa \log_{d-1}\log N, 2\fa\log_{d-1}\log N}$, it holds $|1+(\msc(z))^2+(\msc(z))^4+\cdots+(\msc(z))^{2\ell}|\gtrsim 1$.
It follows that $\sqrt{\kappa+\eta}\lesssim |1-(m_{sc}(z))^2|\lesssim |b|\lesssim \ell|1-(m_{sc}(z))^2|\lesssim \ell \sqrt{\kappa+\eta}$ and we have $|b/a|=|1-m_{sc}^2(z)|/|m_{sc}(z)|^{2\ell+3}\gtrsim \sqrt{\kappa+\eta}$. One can directly check that
\begin{align}\label{e:cbound}
c\lesssim (\log N)^2\left(\frac{\varepsilon}{\sqrt{\kappa+\eta+\varepsilon}}\right)^{3}+\frac{\varepsilon}{(\log N)^{\fa/2-4}}\lesssim \frac{\varepsilon}{(\log N)^2},
\end{align}
where we used the relations \eqref{e:relation2} for the choices of parameters.
From our choice of $\varepsilon(z)$ as in \eqref{e:defeps}, there are two cases: $\varepsilon\lesssim (\kappa+\eta)/\log N$ or $\varepsilon\gtrsim (\log N)^3(\kappa+\eta)$.
If  $\varepsilon\lesssim (\kappa+\eta)/\log N$ then $\varepsilon/\sqrt{\kappa+\eta+\varepsilon}\asymp \varepsilon/\sqrt{\kappa+\eta}$.
Moreover in this case, thanks to \eqref{e:cbound}, we have $|ac/b^2|\lesssim |c|/(\kappa+\eta)\lesssim \varepsilon/((\log N)^2(\kappa+\eta))\ll1$. The two roots of $ax^2+bx+c=0$ satisfies:$|x_1|=|b/a|+\OO(|c/b|)\asymp |b/a|\gtrsim \sqrt{\kappa+\eta}$ and 
\begin{align*}
|x_2|=\OO(|c/b|)\lesssim \frac{|c|}{\sqrt{\kappa+\eta}}\lesssim \frac{1}{(\log N)^2}\frac{\varepsilon}{\sqrt{\kappa+\eta+\varepsilon}}.
\end{align*}
Moreover we also have the estimate for the difference $|Q(\tcG,z)-m_{sc}(z)|\lesssim \varepsilon/\sqrt{\kappa+\eta}$ from \eqref{e:weakone}. Using that $\varepsilon\lesssim (\kappa+\eta)/\log N$,  $|Q(\tcG,z)-m_{sc}(z)|\lesssim \varepsilon/\sqrt{\kappa+\eta}\lesssim \sqrt{\kappa+\eta}/\log N\ll |x_1|$. Thus we conclude that
\begin{align}\label{e:upcase}
|Q(\tcG,z)-m_{sc}(z)|=|x_2|\lesssim \frac{|c|}{\sqrt{\kappa+\eta}}\lesssim \frac{1}{(\log N)^2}\frac{\varepsilon}{\sqrt{\kappa+\eta+\varepsilon}}
.
\end{align}
In the other case, we have $\varepsilon\gtrsim (\log N)^{3}(\kappa+\eta)$ and $\varepsilon/\sqrt{\kappa+\eta+\varepsilon}\asymp \sqrt{\varepsilon}$, by the quadratic formula, we have 
\begin{align}\begin{split}\label{e:lowcase}
&\phantom{{}={}}|Q(\tcG,z)-m_{sc}(z)|
\lesssim \frac{|b|+\sqrt{|b|^2+4|a||c|}}{2|a|}
\lesssim  \frac{|b|}{|a|}+\sqrt{\frac{|c|}{|a|}}\\
&\lesssim \sqrt{\kappa+\eta}+\sqrt{|c|}\lesssim\sqrt{\kappa+\eta}+\frac{1}{\log N}\frac{\varepsilon}{\sqrt{\kappa+\eta+\varepsilon}}\lesssim \frac{\sqrt{\varepsilon}}{\log N}.
\end{split}\end{align}
In the first case, using  \eqref{e:upcase}, we have $\log N(\sqrt{\kappa+\eta}|Q(\tcG, z)-\msc(z)|+|Q(\tcG, z)-\msc(z)|^2)\lesssim\varepsilon/\log N$. In the second case, using \eqref{e:lowcase}, we also have
\begin{align*}
&\phantom{{}={}}\log N(\sqrt{\kappa+\eta}|Q(\tcG, z)-\msc(z)|+|Q(\tcG, z)-\msc(z)|^2)
\lesssim\log N \left(\frac{\sqrt \epsilon}{\log N}\right)^2\lesssim \frac{\varepsilon}{\log N},
\end{align*}
where we used that $\varepsilon\gtrsim (\log N)^{3}(\kappa+\eta)$.

Finally, thanks to \eqref{e:improverigid1},\eqref{e:improverigid1.5}, and the discussion above
\begin{align*}
&\phantom{{}={}}\left|G_{ij}(z)-G_{ij}(\Ext(\cB_r(i,j,\cG),Q(\cG,z)),z)\right|\lesssim \frac{1}{(d-1)^r}
+\frac{\varepsilon}{(\log N)^{\fa/2-5}}\\
&+\log N (\sqrt{\kappa+\eta}|Q(\tcG,z)-m_{sc}(z)|+|Q(\tcG, z)-\msc(z)|^2)
\lesssim \frac{\varepsilon}{\log N},
\end{align*}
where we used the relations \eqref{e:relation1} and \eqref{e:relation2} for the choices of parameters. We conclude that $\Omega\cap_{o\in\qq{N}}\Omega^-_o\subset \Omega^-$, this finishes the proof of Proposition \ref{prop:bootstrap}.

\end{proof}

\begin{proof}[Proof of Proposition \ref{thm:mrQ}]
We first show that the Green's functions are Lipschitz in $z$, and the Lipschitz constant is at most $\OO(1/\Im[z])$.
\begin{claim}\label{l:Lip}
Fix $z\in \bC^+$. If $\max_{1\leq i\leq N}|G_{ii}(z)|\lesssim 1$, then for any $1\leq i,j\leq N$ and $|z'-z|\leq \Im[z]/2$, we have
$|G_{ij}(z')-G_{ij}(z)|\lesssim |z'-z|/\Im[z]$.
\end{claim}
\begin{proof}
Let $z_t=z+t(z'-z)$ for $0\leq t\leq 1$. Then $z_0=z$ and $z_1=z'$. Moreover, since $|z'-z|\leq \Im[z]/2$, we have that $\Im[z_t]\geq\Im[z]/2$ for $0\leq t\leq 1$. By taking derivative of $G_{ij}(z_t)$ with respect to $t$, we have
\begin{align}\begin{split}\label{e:tder}
|\del_tG_{ij}(z_t)|
&\leq |z'-z|\sum_{k=1}^N |G_{ik}(z_t)| |G_{kj}(z_t)|\leq |z'-z|\sum_{k=1}^N (|G_{ik}(z_t)|^2+ |G_{kj}(z_t)|^2)/2\\
&=\frac{|z'-z|}{\Im[z_t]}\frac{\Im[G_{ii}(z_t)]+\Im[G_{jj}(z_t)]}{2}\leq \frac{|z'-z|}{\Im[z]}(|G_{ii}(z_t)|+|G_{jj}(z_t)|).
\end{split}\end{align}
Let $\Gamma_t=\max_{1\leq k\leq N}|G_{kk}(z_t)|$. Then $\Gamma_0\lesssim 1$ and \eqref{e:tder} implies $|\del_t\Gamma_t|\lesssim4\Gamma_t$. It follows that $\Gamma_t\lesssim 1$ for all $0\leq t\leq 1$, and $|\del_tG_{ij}(z_t)|\lesssim |z'-z|/\Im[z]$. The Claim \ref{l:Lip} follows: $|G_{ij}(z')-G_{ij}(z)|=|\int_0^1\del_tG_{ij}(z_t)\rd t|\lesssim |z'-z|/\Im[z]$.
\end{proof}

We recall the graph sets $\bar\Omega, \Omega(z), \Omega^-(z)$ from  Definitions \ref{def:barOmega}, \ref{def:Omega} and \ref{def:Omega-}. We take a lattice grid $\bL$:
\begin{align*}
z_{\al\beta}=\left(-2d +\frac{\al}{N}\right)+\ri\left(2d-\frac{\beta}{N}\right).
\end{align*}
for $\al\in \qq{0, 4dN}$ and $\beta\in \qq{0, 2dN-(\log N)^\fb}$.
And for any $\beta\in \qq{0, 2dN-(\log N)^\fb}$, let $\bL_\beta=\{z_{\al \beta}: \al\in \qq{0, 4dN}\}$. In the following we first prove that \begin{align}\label{e:eqmrQ1}
\P\left(\bar\Omega \setminus \bigcap_{z\in \bL}\Omega^-(z)\right)\lesssim N^{-\fd+3}.
\end{align}
We prove \eqref{e:eqmrQ1} by induction, i.e. we inductively prove for any parameter $\beta\in \qq{0, 2dN-(\log N)^\fb}$
\begin{align}\label{e:betacase}
\P\left(\bar\Omega\setminus \bigcap_{z\in \bL, \Im[z]\geq (2d-\beta/N)}\Omega^-(z)\right)
\lesssim \beta N^{-\fd+2}.
\end{align}
The base case for $\beta=0$ follows from Proposition \ref{prop:betacomp},
\begin{align*}
\P\left(\bar\Omega\setminus \bigcap_{z\in \bL, \Im[z]\geq 2d}\Omega^-(z)\right)
=0.
\end{align*}
Thanks to Claim \ref{l:Lip}, the Green's functions are Lipchitz in $z$, which implies
\begin{align}\label{e:Om-O}
\bigcap_{z\in \bL\atop \Im[z]\geq (2d-\beta/N)}\Omega^-(z)
\subset\bigcap_{|z|\leq 2d\atop \Im[z]\geq (2d-(\beta+1)/N)}\Omega(z),
\end{align}
and especially,
\begin{align}\label{e:Om-O2}
\bigcap_{z\in \bL\atop \Im[z]\geq (2d-\beta/N)}\Omega^-(z)
\subset\bigcap_{z\in \bL, \Im[z]\geq (2d-\beta/N)}\Omega^-(z)
\bigcap_{z\in \bL_{\beta+1}}\Omega(z)
,
\end{align}
If the statement \eqref{e:betacase} holds for $\beta$, then using \eqref{e:Om-O2}, we get
\begin{align}\begin{split}\label{e:lowerz}
&\phantom{{}={}}\P\left(\bar\Omega\setminus\bigcap_{z\in \bL, \Im[z]\geq (2d-\beta/N)}\Omega^-(z)
\bigcap_{z\in \bL_{\beta+1}}\Omega(z)
\right)\\
&=\P\left(\bar\Omega\setminus\bigcap_{z\in \bL, \Im[z]\geq (2d-\beta/N)}\Omega^-(z)\right)
\lesssim \beta N^{-\fd+2}.
\end{split}\end{align}
Then we can use Proposition \ref{prop:bootstrap} to replace $\Omega(z_{\al\beta+1})$ by $\Omega^-(z_{\al\beta+1})$ in \eqref{e:lowerz}. And for each replacement, we lose $\OO(N^{-\fd+1})$ in probability,
\begin{align*}\begin{split}
\P\left(\bar\Omega\setminus\bigcap_{z\in \bL, \Im[z]\geq (2d-(\beta+1)/N)}\Omega^-(z)\right)
&\lesssim \beta N^{-\fd+2}+\OO(4dN^{-\fd+2})\lesssim \OO((\beta+1) N^{-\fd+2}).
\end{split}\end{align*}
This gives the statement \eqref{e:betacase} for $\beta+1$. Therefore statement \eqref{e:betacase} holds for any $\beta\in \qq{0, 2dN-(\log N)^\fb}$.

Similar to \eqref{e:Om-O}, using Claim \ref{l:Lip} again, we have
\begin{align}\label{e:Om-O3}
\bigcap_{z\in \bL}\Omega^-(z)
\subset\bigcap_{|z|\leq 2d\atop \Im[z]\geq (\log N)^\fb/N}\Omega(z).
\end{align}
The claim \eqref{e:eqmrQ1} and \eqref{e:Om-O3} together imply \eqref{e:eqmrQ}, where we can take $\fC=\fd-3$, and Proposition \ref{thm:mrQ} follows.
\end{proof}

\section{Proof of Proposition \ref{p:Omega+}: Stability Under Local Resampling}\label{s:switchstable}

In this Section,  mostly  by perturbation arguments and Schur complement formulas,  we prove Proposition \ref{p:Omega+}: if the Green's functions of $\cG$ or $\cGT$ satisfy certain estimates, the Green's functions of the switched graphs, i.e. $\tcG$ or $\tcGT$, satisfy similar (slightly worse)  estimates.
In Section \ref{sec:stabilityT}, we obtain the estimates for the Green's function of $\cGT$. 
In Sections \ref{sec:dist}--\ref{sec:stability}, we obtain the estimates for the Green's functions of $\cG^{(\bT\bW_{\bfS})}$ and $\tcGT$.
In Section \ref{sec:weakstab}, we obtain the estimates for the Green's functions of $\tcG$. In Section \ref{sec:changeQG}, we show that the difference between  $Q(\cG,z)$ by $Q(\tcG,z)$ is small. The proof of Proposition \ref{p:Omega+} is given in Section \ref{sec:provekey1}. We will outline  the major steps and ideas of the  proofs in Remarks \ref{r1} and \ref{r2}.


\subsection{Stability under removal of a neighborhood}
\label{sec:stabilityT}


We recall the notations and definitions  from  Section \ref{sec:switch} and Section \ref{sec:weak}. Our local resampling involves a fixed center vertex $o$,
and a radius $\ell$.
Given a $d$-regular graph $\cG$, we abbreviate $\cT=\cB_{\ell}(o,\cG)$ (which may not be a tree) and its vertex set $\bT$.
We enumerate edge boundary $\del_E \cT$ as $ \del_E \cT = \{e_1,e_2,\dots, e_\mu\}$, where $e_\al=\{l_\al, a_\al\}$ with $l_\al \in \T$ and $a_\al \in \qq{N} \setminus \T$. 
The following deterministic estimate shows that
removing the neighborhood $\T$, which is the vertex set of $\cB_\ell(o,\cG)$, from the graph $\cG$
has a small effect on the Green's function in the complement of $\T$.

\begin{proposition} \label{prop:stabilityGT}
Fix $d\geq 3$,  $\omega_d$ as   in Definition \ref{d:defwd}, and $z\in \bC^+$ with $\Im[z]\geq (\log N)^{\fb}/N$.
Let $\cG\in \Omega(z)$ be a spectral regular graph, as defined in Definition \ref{def:Omega}.
Then, for all vertices $i,j\in \qq{N} \setminus \T$, we have
\begin{equation}\label{e:stabilityGT}
|\GT_{ij}(z)-G_{ij}(\Ext(\cB_r(i,j, \cGT), Q(\cG,z)),z)|
\lesssim (\log N)^3\varepsilon(z).
\end{equation}
Especially, we have $\Omega(z)\subset \Omega_o^+(z)$, where  $\Omega_o^+(z)$ is defined in Definition \ref{def:Omegao+}.
\end{proposition}

\subsubsection{Boundary of a neighborhood}

\begin{proposition}\label{GTstructure}
Let $\cal G$ be a $d$-regular graph on $\qq{N}$ and fix $R\geq 2\ell$. Assume
that $\cB=\cB_{R}(o,\cG)$ has excess at most $\omega$. 
Then the following hold.
\begin{itemize}
\item
After removing $\T$, most boundary vertices of $\cT$ are far away from the other boundary vertices:
\begin{align}\label{fewclose}
|\{\al\in \qq{1,\mu}: \exists \beta\in \qq{1,\mu} \setminus \{\al\}, \dist_{\cGT}(a_\al, a_\beta)\leq R/2\}|\leq 2\omega.
\end{align}
\item
After removing $\T$, any vertex $x\in \qq{N} \setminus \T$ can only be close to at most $\omega+1$  boundary vertices of $\bT$:
\begin{align}
\label{distfinitedega}
|\{\al\in\qq{1,\mu}: \dist_{\cGT}(x,a_\al)\leq {R}/2\}|\leq \omega+1.
\end{align}
\end{itemize}
\end{proposition}


Next, for $R\geq 2\ell$, we have the following deterministic bound on the deficit functions for the
connected components of the subgraph obtained from $\cB=\cB_{R}(o,\cG)$ by removing  $\bT$.

\begin{proposition} \label{p:deficitbound}
Under the same assumption as in Proposition \ref{GTstructure}, the following holds:
Let $\cA$ be the annulus obtained by removing $\T$ from $\cB$, i.e. $\cA=\cB^{(\bT)}$,
then for any connected component of $\cA$, its deficit function is nonzero and the {sum of its} deficit function plus its excess is at most $\omega+1$. 
\end{proposition}

For the above statements,
recall that we view $\cal B$ as subgraphs of $\cG$ (which has zero deficit function), the deficit function of $\cal B$ is also zero.
By our conventions in Remark \ref{r:convention}, the deficit function of $\cA=\cB^{(\bT)}$ is given by \eqref{e:removing}.
In the remainder of this section, we prove Propositions~\ref{GTstructure} and \ref{p:deficitbound}.

By assumption the neighborhood $\cB$ has excess at most $\omega$. 
We partition $\qq{1,\mu}$ into sets  $\{\sA_1, \sA_2, \cdots, \sA_\kappa\}$,
such that $\al$ and $\beta$ are in the same set $\sA_s$ iff  $a_\al$ and $a_\beta$ are in the same connected component of $\cA$.

\begin{lemma} \label{lem:beta2n}
Under the same assumption as in Proposition \ref{GTstructure}, we have $\sum_{s=1}^\kappa {\bf1}_{|\sA_s|\geq 2}|\sA_s|\leq 2\omega$, and the connected component of $\cA$ containing $\{a_\al: \al\in \sA_s\}$ has excess at most $\omega-|\sA_s|+1$, for any $1\leq s\leq \kappa$.
\end{lemma}

\begin{proof}
We recall the definition  of $\text{excess}(\cG)$ from \eqref{def:excess},
\begin{align}\label{c-c}
\text{excess}(\cG)
  =\#\text{connected components}(\cG)-
  \#\text{vertices}(\cG)+\#\text{edges}(\cG).
\end{align}
The annulus $\cA$ is obtained from $\cB$ by removing the edges $e_1, e_2,\cdots, e_\mu$. Thanks to \eqref{c-c}, the number of excess decreases by $|\sA_s|-1$ if we remove edges $\{e_\al: \al\in \sA_s\}$  from $\cB$. Therefore, the excess of $\cA$ satisfies $\text{excess}(\cA)\leq \omega-\sum_{s=1}^{\kappa}(|\sA_s|-1)$. By rearranging it, we have in particular that 
\begin{align*}
\sum_{s=1}^\kappa {\bf 1}_{|\sA_s|\geq 2}|\sA_s|\leq 
2\sum_{s=1}^{\kappa}(|\sA_s|-1)\leq 2\omega.
\end{align*}
Since the excess of a graph equals the sum of excesses of its connected components, the connected component of $\cA$ containing $\{a_\al: \al\in \sA_s\}$ has excess at most $\text{excess}(\cA)\leq \omega-|\sA_s|+1$, for any $1\leq s\leq \kappa$.
\end{proof}

\begin{proof}[Proof of Proposition \ref{GTstructure}]
We use the same notation as in Lemma \ref{lem:beta2n}. For \eqref{fewclose}, since we have $R\geq 2\ell$, if $\dist_{\cGT}(a_\al, a_\beta)\leq {R/2}$ then $a_\al$ and $a_\beta$ are in the same connected component of $\cA$, which implies that there exists some $1\leq s\leq \kappa$, with $\al \in \sA_s$ and $|\sA_s|\geq 2$. Therefore,  
\begin{align*}
|\{\al\in \qq{1,\mu}: \exists \beta\in \qq{1,\mu} \setminus \{\al\}, \dist_{\cGT}(a_\al, a_\beta)\leq {R/2}\}|\leq
\sum_{s: |\sA_s|\geq 2}|\sA_s|\leq 2\omega,
\end{align*}
where we used Lemma \ref{lem:beta2n} in the last inequality.
For \eqref{distfinitedega}, if $\dist_{\cGT}(x,a_\al)\leq {R}/2$, then $x$ is in the annulus $\cA$. Say $x$ belongs to connected component of $\cA$ containing $\{a_\al: \al\in \sA_s\}$. Then thanks to Lemma \ref{lem:beta2n} we have
\begin{align*}
|\{\al\in\qq{1,\mu}: \dist_{\cGT}(x,a_\al)\leq {R}/2\}|\leq |\sA_s|\leq \omega+1.
\end{align*}
\end{proof}

\begin{proof}[Proof of Proposition~\ref{p:deficitbound}]
We use the same notations as in Lemma \ref{lem:beta2n}.
For the connected component containing $\{a_\al: \al\in\sA_s\}$,  its excess is at most $\omega-|\sA_s|+1$. 
The deficit function of $\cB$ is zero, and the deficit function of $\cA=\cB^{(\bT)}$ is given by 
\eqref{e:removing}:
\begin{align*}
g(v)=\deg_{\cB}(v)-\deg_{\cA}(v).
\end{align*}
Therefore, the deficit function satisfies 
\begin{align*}
\sum_{v \in \text{ component of $\{a_\al: \al\in\sA_s\}$}} g(v)= \sum_{\alpha: \alpha \in \sA_s}1= |\sA_s|.
\end{align*}
The statement follows. 

\end{proof}

\subsubsection{Proof of Proposition \ref{prop:stabilityGT}}

Let $\cG_o=\cB_{r}(\{\bT,i,j\},\cG)$ with zero deficit function, and the deficit function of $\cG_o^{(\bT)}$ is given by \eqref{e:removing}, according to remark \ref{r:convention}.
We abbreviate
\begin{equation*}
P=G(\Ext(\cG_o, Q(\cG,z))),\qquad
\quad P^{(\bT)}=G(\Ext(\cG_o^{(\bT)}, Q(\cG,z))).
\end{equation*}
Thanks to Proposition \ref{p:deficitbound}, each connected component of $\cG_o$ and $\cGT_o$ satisfies the assumptions in Proposition \ref{p:localization}. Thus we have
\begin{align}\begin{split}\label{replaceEr0}
&\absb{G_{ij}(\Ext(\cB_{r}(i,j,\cG),  Q(\cG,z)))-P_{ij}},\\
&\absb{G_{ij}(\Ext(\cB_{r}(i,j,\cG^{(\bT)}),  Q(\cG,z)))-P_{ij}^{(\bT)}}
\lesssim (\log N)\varepsilon,
\end{split}\end{align}
where the $\log N$ factor is from the term $1+\diam(\cG)$ in \eqref{e:compatibility}.
The normalized adjacency matrices of $\cG$ and $\Ext(\cG_o, Q(\cG,z))$ have the block  matrix form
\begin{align}\label{decomp}
 \left[\begin{array}{cc}
        H & B'\\
        B & D
       \end{array}
\right],\quad 
\left[\begin{array}{cc}
        H& {B}_o'\\
        {B}_o & D_o
       \end{array}
\right],
\end{align}
where $H$ is the normalized adjacency matrix of $\cT$.
The nonzero entries of $B$ and ${B}_o$ occur for the indices
$(i,j) \in \{(a_1,l_1), (a_2,l_2),\cdots, (a_\mu,l_\mu)\}$ and take values $1/\sqrt{d-1}$.
Notice that for nonzero entries $ B_{ij}=( {B}_o)_{ij}$. In the rest of this section,  by abuse of notations, we will  not distinguish
$B$ and ${B}_o$.

By the Schur complement formula \eqref{e:Schur1}, we have
\begin{align}\begin{split}\label{GTGPTPS}
 &\GT=(D-z)^{-1}
 =G|_{\bT^c}-G|_{\bT^c \bT}(G|_{\T})^{-1}G|_{\bT \bT^c},\\
 & P^{(\bT)}=(D_o-z)^{-1}
 =P-P(P|_{\T})^{-1}P,
\end{split}\end{align}
and also
\begin{align}\label{invGP}
(G|_{\T})^{-1}=H-z-B'(D-z)^{-1}B,\quad
(P|_{\T})^{-1}=H-z-B'(D_o-z)^{-1}B.
\end{align}

\begin{claim}\label{c:Piysum}
Under the same assumption as in Proposition \ref{prop:stabilityGT},
for any vertex $i$ in $\cG_o$,
\begin{align*}
\sum_{y\in \bT}|P_{iy}|\lesssim \log N,\quad
\sum_{y\in \bT}|G_{iy}|\lesssim \log N.
\end{align*}
\end{claim}
\begin{proof}
Since $\cG_o$ satisfies the assumptions in Proposition \ref{p:subgraph}, the first estimate follows from \eqref{e:boundPijcopy} and a counting argument. For the second estimate, since $\cG\in \Omega(z)$, \eqref{e:defOmega} and \eqref{replaceEr0} imply
\begin{align*}
\sum_{y\in \bT}|G_{iy}|\lesssim \sum_{y\in \bT}|P_{iy}|+\varepsilon\log N\lesssim \log N+ |\bT|\varepsilon \log N
\lesssim \log N,
\end{align*}
where we used that $|\bT|\varepsilon\lesssim (d-1)^{\ell} \varepsilon\lesssim (\log N)^{2\fa} \varepsilon\ll1$ from \eqref{e:relation2}. 
\end{proof}

\begin{claim}\label{c:sumPTSbd}
Under the same assumption as in Proposition \ref{prop:stabilityGT}, for any $x \in \T$,
\begin{equation} \label{e:sumPTSbd}
\sum_{y\in \T} |(P|_{\T})^{-1}_{xy}|
\lesssim 1.
\end{equation}
\end{claim}

\begin{proof}
For any $x\in \T$, by \eqref{invGP} we have
\begin{align}\begin{split} \label{express01}
&\phantom{{}={}}\sum_{y\in \T} |(P|_{\T})^{-1}_{xy}|
  =\sum_{y\in \T}|(H-z-B'(D_o-z)^{-1}B)_{xy}|\\
 &\leq |z|+\sum_{y\in \bT } H_{xy}   +\sum_{y\in \T} |(B'(D_o-z)^{-1}B)_{xy}|\leq \OO(1)+\sum_{y\in \T} |(B'(D_o-z)^{-1}B)_{xy}|.
\end{split}\end{align} 
Let $\bT_\ell=\{v\in \bT: \dist_{\cT}(o,v)=\ell\}=\{l_1, l_2,\cdots, l_\mu\}$.
The terms in the last sum in \eqref{express01} vanish unless $x,y\in \T_\ell$, equivalently $x=l_\al, y=l_\beta$ for some $\al,\beta\in \qq{1,\mu}$.
If $x\not\in \bT_\ell$, we have $\sum_{y\in \T} |(P|_{\T})^{-1}_{xy}|\lesssim 1$. If $x\in \bT_\ell$, the sum is bounded by
\begin{align}\begin{split}\label{e:D-z}
\sum_{\al,\beta\in \qq{1,\mu}} B'_{x a_\al }|(D_o-z)^{-1}_{a_\al a_\beta}|B_{a_\beta l_\beta}\leq \frac{1}{d-1}\sum_{\al, \beta\in\qq{1,\mu}}1_{l_\al=x} |(D_o-z)^{-1}_{a_\al a_\beta}|
.
\end{split}\end{align}
Thanks Proposition \ref{p:deficitbound},
the  graph $\cGT_o$ satisfies the assumptions in Proposition \ref{p:subgraph},  and \eqref{e:boundPijcopy} gives
\begin{align}\label{e:bb1}
|(D_o-z)^{-1}_{a_\al a_\beta}|=|G_{a_\al a_\beta}(\Ext(\cG_o^{(\T)}, Q(\cG,z)),z) |\lesssim 1.
\end{align}
The number of $\al$ such that $l_\al=x$ is at most $d-1$. For each such $\al$,
lemma \ref{lem:beta2n} implies the number of $\beta$ such that $a_\al, a_\beta$ are in the same connected component of $\cG_o^{(\bT)}$ is at most $\OO(\omega_d)$. Thus, there are at most $(d-1)\OO(\omega_d)$ nonzero terms in the sum of \eqref{e:D-z}. Combining with \eqref{e:bb1}, we conclude that
\begin{equation} \label{express02}
\sum_{y\in \T} |(B'(D_o-z)^{-1}B)_{xy}|\lesssim 1.
\end{equation}
The claim \eqref{e:sumPTSbd} follows from combining \eqref{express01} and \eqref{express02}.
\end{proof}

\begin{claim}\label{c:GTinv}
Under the same assumption as in Proposition \ref{prop:stabilityGT}, 
for any $x,y\in \T$,
\begin{align}
 \label{GTinvbd} &|(G|_{\T})_{xy}^{-1}-(P|_{\T})_{xy}^{-1}|  \lesssim (\log N)\varepsilon,
\end{align}
\end{claim}

\begin{proof}
Define matrices $ W$ and $ E$ by
\begin{align}
 \label{GTdiff} &G|_{\T}=P|_{\T}+ W,\\
 \label{GTinv}&(G|_{\T})^{-1}=(P|_{\T})^{-1}+ E.
\end{align}
From \eqref{e:defOmega} and \eqref{replaceEr0}, for any $x,y\in \T$, we have $| W_{xy}|\leq (\log N)\varepsilon$. 
We claim the same estimate holds for the entries of the matrix $ E$. 
Notice from \eqref{invGP} that $ E_{xy}\neq 0$ only for $x,y\in \T_\ell=\{v\in \bT: \dist_{\cT}(o,v)=\ell\}$.
Let $\Gamma:=\max_{x,y\in \T}| E_{xy}|=\max_{x,y\in \T_\ell}| E_{xy}|$.
By taking the product of \eqref{GTdiff} and \eqref{GTinv},
\begin{align}\label{selfcontrol1}
  E+(P|_{\T})^{-1} W E+(P|_{\T})^{-1} W (P|_{\T})^{-1}=0.
\end{align}
For any $x,y\in \T$, taking the $(x,y)$-th entry
\begin{align*}
| E_{xy}|
&\leq \sum_{i,j\in \T} |(P|_{\T})^{-1}_{xi}|| W_{ij}|| E_{jy}|+\sum_{i,j\in \T} |(P|_{\T})^{-1}_{xi}|| W_{ij}| |(P|_{\T})^{-1}_{jy}|\\
&\leq |\T|\OO(\varepsilon\log N)\Gamma \sum_{i\in \T} |(P|_{\T})^{-1}_{xi}|
+\OO(\varepsilon\log N)\sum_{i\in \T} |(P|_{\T})^{-1}_{xi}|\sum_{j\in \T}|(P|_{\T})^{-1}_{jy}|\\
&\leq \oo(1)\Gamma+\OO(\varepsilon\log N).
\end{align*}
For the second inequality, we used $| W_{xy}|=\OO(\varepsilon\log N)$;
for the third inequality we used $ |\T|\lesssim(d-1)^{\ell}$, \eqref{e:sumPTSbd}, and the fact $(d-1)^{\ell}\varepsilon\log N\ll1$ from \eqref{e:relation1}, \eqref{e:relation2}.
Taking the maximum on the right-hand side of the above inequality, and rearranging, we get
$
 \Gamma \lesssim \varepsilon \log N,
$
as claimed.
\end{proof}

\begin{proof}[Proof of Proposition \ref{prop:stabilityGT}]
We take difference of the two expressions in \eqref{GTGPTPS},\begin{align}\begin{split}\label{e:GTbbbbb}
 &\phantom{{}={}}|G^{(\bT)}_{ij}-P^{(\bT)}_{ij}|
 \lesssim |G_{ij}-P_{ij}|
 +\sum_{x,y\in \bT}|G_{ix}-P_{ix}||(G|_\bT)^{-1}_{xy}||G_{yj}|\\
&+\sum_{x,y\in \bT}|P_{ix}||(G|_\bT)^{-1}_{xy}-(P|_\bT)^{-1}_{xy}||G_{yj}|+
\sum_{x,y\in \bT}|P_{ix}||(P|_\bT)^{-1}_{xy}||G_{yj}-P_{yj}|\\
&\lesssim \varepsilon\log N\left(1+
 \sum_{x,y\in \bT} |((P|_\bT)^{-1}_{xy}|+\varepsilon\log N)|G_{yj}|+
|P_{ix}||G_{yj}|+
+|P_{ix}||(P|_\bT)^{-1}_{xy}|\right)\\
&\lesssim (\log N)^3 \varepsilon.
\end{split}
\end{align}
where in the second inequality, we used
$|G_{ix}-P_{ix}|\lesssim \varepsilon\log N$ from \eqref{e:defOmega} and \eqref{replaceEr0},  and $|(G|_\bT)^{-1}_{xy}-(P|_\bT)^{-1}_{xy}|\lesssim \varepsilon \log N$ from Claim \ref{c:GTinv}; in the last inequality we used  Claims \ref{c:Piysum} and \ref{c:sumPTSbd}.
Proposition \ref{prop:stabilityGT} follows from combining \eqref{e:GTbbbbb} and \eqref{replaceEr0}.
\end{proof}

\begin{remark}[Remark on Methods] 
The decomposition \eqref{decomp} and the Schur complement formulas in this case \eqref{GTGPTPS} and \eqref{invGP} are the main tool in this paper. 
These formulas provide  estimates  on $\GT$ and $P^{(\bT)}$  based on informations of  $G$ and $P$. In order to use the Schur complement formulas \eqref{GTGPTPS} and \eqref{invGP}, 
we need to estimate $G|_{\T}^{-1}$. Since $G|_{\T}$ can be approximate locally by $P|_{\T}$, taking the inverse can be easily achieved by the standard resolvent identity since 
$\T$ does not contain too many vertices. 
\label{r1}
\end{remark}

\subsection{Graph distance between switched vertices}
\label{sec:dist}

We recall the notations and definitions  from  Section \ref{sec:switch} and Section \ref{sec:weak}. 
In this section, we give some basic estimates for the local resampling and graph distance between switched vertices.
\begin{proposition} \label{lem:switchable}
For any vertex $x \in \qq{N}\setminus \T$, and $\al\in\qq{1,\mu}$
\begin{equation} \label{e:approxunifbd}
  \P_{\cG}(b_\al= x)=  \P_\cG(c_\al = x)
  \leq 2/N,
\end{equation}
where $\P_{\cG}(\cdot)$ is defined in Definition \ref{def:enlarge}.
\end{proposition}

\begin{proof}
We recall that, for any $\al\in \qq{1,\mu}$,
the oriented edge $(b_\al,c_\al)$ is uniformly chosen from the oriented edges of $\cGT$.
By the definition of $\T$, there are $Nd-\OO((d-1)^{\ell})$ oriented edges in $\cGT$,
and since for any vertex $x\in \cGT$, the degree obeys $\deg_{\cGT}(x)\leq d$,
\begin{align*}
  \P_{\cG}(b_\al = x)= \P_{\cG}(c_\al = x)\leq \frac{d}{Nd-\OO((d-1)^{\ell})}
  \leq \frac{2}{N}.
\end{align*}
\end{proof}

We remark that the edges $\{b_\al, c_\al\}$ are randomly chosen in the graph $\cGT$. In the following proposition, we show that, by paying a small probability, we can make sure that 
most $\{b_\al, c_\al\}$ are at least distance $\fR$ away from any $x\in \qq{N}\setminus \bT$, and their radius $\fR$ neighborhoods  are trees.
\begin{proposition} \label{prop:distdeg}
For any graph $\cG \in \bar\Omega$ or $\bar\Omega^+$ (as in Definitions \ref{def:barOmega} and \ref{def:barOmega1}), and any large number $\fd>0$, 
the following holds with $\P_{\cG}$-probability at least $1-\OO(N^{-\fd})$:
\begin{itemize}
\item
Any vertex $x\in \qq{N}\setminus\T$ is far away from most vertices in the set $\{a_1, b_1, c_1,\dots, a_\mu, b_\mu, c_\mu\}$:
\begin{align}
{\label{distfinitedegb}} |\{\al\in\qq{1,\mu}:\dist_{\cGT}(x,\{a_\al, b_\al, c_\al\})\leq {\fR/4}\}|= \OO_\fd(1).
\end{align}
\item
Most indices $\al \in \qq{1,\mu}$ are good in the following sense: 
\begin{align}
\label{distdega}
&|\Ba| = \OO_\fd(1),
\text{with } \Ba =\{\al\in \qq{1,\mu}: {\dist_{\cGT}(a_\al, \{b_\al,c_\al\})\leq {\fR}}\},
\\
\label{distdegb}
&|\Bb| = \OO_\fd(1), 
\text{with }\Bb =\{\al\in \qq{1,\mu}: \exists \beta\neq \al, {\dist_{\cGT}(\{a_\al, b_\al,c_\al\}, \{a_\beta, b_\beta, c_\beta\})\leq {\fR}/4}\},
\\
\label{treeneighbor}
&|\Bc| = \OO_\fd(1),
\text{with } \Bc =\{\al\in \qq{1, \mu}: \text{$\cB_{\fR}(\{b_\al,c_\al\}, \cGT)$ is not a tree}\}.
\end{align}
\end{itemize}
\end{proposition}

\begin{proof}

For \eqref{distfinitedegb}, by taking $R=\fR/2$ in \eqref{distfinitedega}, we have 
\begin{equation} \label{e:fana}
|\{\al\in \qq{1,\mu}: {\dist_{\cGT}(x, a_\al)\leq {\fR/4}}\}|= \OO_\fd(1).
\end{equation}
We notice that 
in any graph with degree bounded by $d$,
the number of vertices at distance at most ${\fR}$ from vertex $x$ is bounded by $\OO((d-1)^{{\fR}})$.
By \eqref{e:approxunifbd} we have
\begin{equation*} 
\P_{\cG}(\dist_{\cGT}(x,\{b_\al,c_\al\})\leq {\fR})\lesssim \frac{(d-1)^{{\fR}}}{N}.
\end{equation*}
Since the $\{b_1,c_1\}, \{b_2,c_2\},\dots, \{b_\mu, c_\mu\}$ are independent, it follows from a union bound
\begin{align}\begin{split}\label{e:distxbi}
&\phantom{{}={}}\P_\cG\left(|\{\al\in \qq{1,\mu}:\dist_{\cGT}(x,\{b_\al, c_\al\})\leq {\fR}\}|
\geq C\fd\right)\\
&\leq \binom{\mu}{C\fd} \pbb{\frac{(d-1)^{{\fR}}}{N}}^{C\fd}\ll N^{-\fd},
\end{split}\end{align}
provided $C$ is large enough. \eqref{e:distxbi} trivially implies
\begin{align}\label{e:distxbi2}
\P_\cG\left(|\{\al\in \qq{1,\mu}:\dist_{\cGT}(x,\{b_\al, c_\al\})\leq {\fR/4}\}|\geq C\fd\right)\ll N^{-\fd}; 
\end{align}
\eqref{e:fana} and \eqref{e:distxbi2} imply \eqref{distfinitedegb}.
The same argument used in the proof of  \eqref{e:distxbi} gives \eqref{distdega}.

For \eqref{distdegb}, using \eqref{fewclose} with $R=\fR/2$, we have
\begin{align}\label{fewclosecopy}
\{\al\in \qq{1,\mu}: \exists \beta\neq \al, {\dist_{\cGT}(a_\al, a_\beta)\leq {\fR/4}}\}= \OO_\fd(1),
\end{align}
For those $\al$ which are not in the set \eqref{fewclosecopy}, thanks to \eqref{e:approxunifbd}, we have
\begin{align*}
\P_\cG( {\dist_{\cGT}(\{a_\al, b_\al,c_\al\}, \{a_\beta, b_\beta, c_\beta\})\leq {\fR/4}})\lesssim \frac{(d-1)^{{\fR/4}}}{N}.
\end{align*}
The claim \eqref{distdegb} follows from a union bound with \eqref{fewclosecopy}.

For \eqref{treeneighbor}, by the assumption $\cG\in \bar\Omega$ or $\bar\Omega^+$, all except at most $2N^{\fc}$ many vertices have radius-$\fR$ tree neighborhoods.
In particular, the same holds for $\cGT$.
By \eqref{e:approxunifbd}, it follows that
\begin{align*}
\P_{\cG}\left(\text{the radius-$\fR$ neighborhood of $b_\al, c_\al$ contains cycles}\right)\lesssim N^{-1+\fc}.
\end{align*}
The claim \eqref{treeneighbor} follows from a union bound.
\end{proof}

The estimates  \eqref{distfinitedegb}  and \eqref{distdegb} also imply the following
estimates for the switched graph $\tcGT$.

\begin{proposition} \label{structuretGT}
For any graph $\cG \in \bar\Omega$ or $\bar\Omega^+$ (as in Definitions \ref{def:barOmega} and \ref{def:barOmega1}),  assume \eqref{distfinitedegb} and \eqref{distdegb}, then
\begin{itemize}
\item
For any index $\al\in \qq{1,\mu}\setminus \Bb$ (as defined in \eqref{distdegb}),
\begin{align}\label{tcGTdist}
\dist_{\tcGT}(\{a_\al,b_\al, c_\al\}, \{a_\beta,b_\beta,c_\beta: \beta\in \qq{1,\mu}\setminus \{\al\}\})>\fR/4. 
\end{align}
\item
For any vertex $x\in \qq{N}\setminus \T$, 
\begin{align}\label{e:lessshortdist}
|\{\al\in\qq{1,\mu}: \dist_{\tcGT}(x,\{a_\al,b_\al,c_\al\})\leq \fR/4\}|= \OO_\fd(1).
\end{align}
\end{itemize}
\end{proposition}
\begin{proof}
For our local resampling, it holds that if $\dist_{\cGT}(\{a_\al,b_\al, c_\al\}, \{a_\beta,b_\beta,c_\beta: \beta\in \qq{1,\mu}\setminus \{\al\}\})>\fR/4$, then $\dist_{\tcGT}(\{a_\al,b_\al, c_\al\}, \{a_\beta,b_\beta,c_\beta: \beta\in \qq{1,\mu}\setminus \{\al\}\})>\fR/4$. 
The first claim \eqref{tcGTdist} follows from the definition of the set $\Bb$,
as in \eqref{distdegb}.
If $\dist_{\tcGT}(x,\{a_\al,b_\al,c_\al\})\leq \fR/4$, then either $\al\in \Bb$ or
$\dist_{\cGT}(x,\{a_\al,b_\al,c_\al\})\leq \fR/4$. Thus \eqref{e:lessshortdist} follows from \eqref{distfinitedegb} and \eqref{distdegb}.

\end{proof}

\begin{proposition}\label{p:newG}
For any graph $\cG\in \bar\Omega$ (as in Definition \ref{def:barOmega}), 
there exists
event $F_1(\cG) \subset \sS(\cG)$ with probability $\bP_{\cG}(F_1(\cG))=1-\OO(N^{-\fd})$, it holds that 
for any switching data $\bfS\in \sS(\cG)$\index{$F_1(\cG)$},  $\tcG=T_{\bfS}(\cG)\in \bar\Omega^+$ as defined in Definition \ref{def:barOmega1}. 
\end{proposition}
\begin{proof}
We denote $F_1(\cG)\subset \sS(\cG)$ the set of switching data $\bfS\in \sS(\cG)$ such that the following statements hold
\begin{align}\label{e:exp2}
&|\{\al\in\qq{1,\mu}:\dist_{\cGT}(x,\{a_\al, b_\al, c_\al\})\leq \fR/2\}|\lesssim \fd, \quad x\in \qq{N}\setminus \bT,
\end{align}
and
\begin{align}\begin{split}\label{e:exp1}
&|\{\al\in \qq{1,\mu}: \exists \beta\neq \al, {\dist_{\cGT}(a_\al, \{a_\beta, b_\beta, c_\beta\})\leq {\fR/2}}\}|\lesssim \fd,\\
&|\{\al\in \qq{1,\mu}: \exists \beta\neq \al, {\dist_{\cGT}(\{ b_\al,c_\al\}, \{a_\beta, b_\beta, c_\beta\})\leq {\fR}}\}|\lesssim \fd,\\
&|\{\al\in \qq{1, \mu}: \text{$\cB_{\fR}(\{b_\al,c_\al\}, \cGT)$ is not a tree}\}|\lesssim \fd.
\end{split}\end{align}
The same argument as for \eqref{distfinitedegb}, \eqref{distdegb} and \eqref{treeneighbor}, by using Proposition \ref{GTstructure} with $R=\fR$, implies that $\bP_{\cG}(F_1(\cG))=1-\OO(N^{-\fd})$. In the following, we prove that for any ${\bf S}\in F_1(\cG)$, $\tcG=T_{\bfS}(\cG)\in \bar\Omega^+$.

We recall the local resampling from Section \ref{sec:switch}, that $\tcG$ is obtained from $\cG$ by removing those edges $\{a_\al, l_\al\}, \{b_\al, c_\al\}$, and adding edges $\{l_\al, c_\al\}, \{a_\al, b_\al\}$ for $\al\in \As_{\bf S}$. If a vertex $x$ is far away from those vertices involving in local resampling, i.e. $\dist_{\cG}(x,\{l_\al, a_\al, b_\al, c_\al:\al\in \As_{\bf S} \})\geq \fR$, then its radius $\fR$ neighborhood does not change $\cB_{\fR}(x, \cG)=\cB_{\fR}(x, \tcG)$. Especially if $\cB_{\fR}(x, \cG)$ is a truncated $d$-regular tree, so is $\cB_{\fR}(x, \tcG)$. Therefore the number of vertices that have a radius-$\fR$ neighborhood that contains a cycle in $\tcG$ is at most $N^\fc+4|\As_{\bf S}|(d-1)^\fR\leq 2N^{\fc}$.

If $x\not\in \bT$ or $x\not\in \cB_{\fR/8}(\{a_\al: \al\in \qq{1,\mu}\setminus \As_{\bfS}\}\cup \{a_\al, b_\al, c_\al:\al\in \As_{\bfS}\}, \tcGT)$, the radius $\fR/8$ neighborhood of $x$ in $\cG$ and $\tcG$ are the same. Thus the excess of $\cB_{\fR/8}(x, \cG)$ is the same as that of $\cB_{\fR/8}(x, \tcG)$, which is at most $\omega_d$. By the definition of switchable set $\As_{\bfS}$ as in \eqref{Wdef}, if $\al\in \As_{\bfS}$, then $\dist_{\cGT}(\{a_\al, b_\al, c_\al\}, \{a_\beta, b_\beta, c_\beta\})\geq \fR/4$ for all $\beta\in \qq{1,\mu}\setminus\{\al\}$ and the subgraph $\cB_{ \fR/4}(\{a_\al, b_\al, c_\al\}, \cG^{(\bT)})$ after adding the edge $\{a_\al, b_\al\}$ is a tree. If $x\in \bT$ or $x\in \cB_{\fR/8}(\{a_\al: \al\in \qq{1,\mu}\setminus \As_{\bfS}\}\cup \{a_\al, b_\al, c_\al:\al\in \As_{\bfS}\}, \tcGT)$, the neighborhood $\cB_{\fR/8}(x, \tcG)$ is either a subgraph of $\cB_{\fR/4}(\{a_\al, b_\al\})$ for some $\al\in  \As_{\bfS}$, or a subgraph of $\cB_{\fR/4}(\{a_\al: \al\in \qq{1,\mu}\setminus \As_{\bfS}\}, \cG\setminus\{\{l_\al, a_\al\}: \al\in \As_{\bfS}\})$ attaching some subtrees at $\{c_\al: \al\in \As_{\bfS}\}$. In both cases the excess is at most $\omega_d$.

If $x\in \bT$, let $\sB\subset \As_{\bf S}$ be the set of indices in \eqref{e:exp1}; If $x\in \qq{N}\setminus\bT$, let $\sB\subset \As_{\bf S}$ be the set of indices in \eqref{e:exp1}, union with 
$
\{\al\in \As_{\bf S}: \dist_{\tcGT}(x, \{a_\al, b_\al, c_\al\})\leq \fR/2\}.
$
From our choice of $F_1(\cG)$, it holds that $|\sB|\lesssim \fd$. The neighborhood $\cB_{\fR/2}(x, \tcG)$ is a subgraph of
\begin{align}\label{e:ggf}
\cB_{\fR/2}(\{x\}\cup\{a_\al, b_\al, c_\al:\al\in {\sB}\}\cup\{a_\al:\al\in \qq{1,\mu}\setminus \As_\bfS\},  \cG\setminus\{\{l_\al, a_\al\}:\al\in \As_\bfS\})
\end{align}
adding the edges $\{\{a_\al, b_\al\}: \al\in {\sB}\}$, removing the edges $\{\{b_\al, c_\al\}: \al\in {\sB}\}$, and attaching subtrees at $\{c_\al, \al \in \As_{\bfS}\setminus {\sB}\}$. Therefore the excess of $\cB_{\fR/2}(x, \tcG)$ is upper bounded by $|{\sB}|=\OO(\fd)$ plus the excess of the graph \eqref{e:ggf}. By our assumption, $\cG\in \bar \Omega$, each radius $\fR$ neighborhood of $\cG$ contains at most $\omega_d$ cycles. So is the graph $\cG\setminus\{\{l_\al, a_\al\}:\al\in \As_\bfS\}$. The graph \eqref{e:ggf} is an union of $1+3|{\sB}|+(\mu-|\As_\bfS|)=\OO(\fd)$ radius $\fR/2$ neighborhoods. Its excess is at most $\OO_\fd(1)$. We conclude that the excess of $\cB_{\fR/2}(x, \tcG)$ is at most $\OO_\fd(1)$. This finishes the proof of Proposition \ref{p:newG}.

\end{proof}

\subsection{The Green's function distance and switching cells}
\label{sec:distGF}

The bounds proved in Section~\ref{sec:dist} provide accurate control for distances at most $\OO(\fR)$.
However, random vertices are typically much further away from each other,
and as mentioned in Section~\ref{sec:outline},
we require stronger upper bounds on the Green's function for such large distances.
These bounds are in fact a general consequence of the Ward identity,
\begin{equation} \label{e:Ward-bis}
  \sum_{j} |G_{ij}(z)|^2 = \frac{\im [G_{ii}(z)]}{\im [z]},
\end{equation}
which holds for the Green's function of any symmetric matrix (see \eqref{e:Ward}).
To make use of it, we introduce a much coarser measure of distance in terms of the size
of the Green's function as follows.

\subsubsection{Definition}\label{sec:defcells}
We recall the control parameter $\varphi(z)$ from \eqref{e:defphi}, and the admissible set $\As_{\bf S}\in \qq{1,\mu}$ from \eqref{Wdef}.
 We define a relation $\simeq$ on $\qq{N} \setminus \T$ by
setting $x \simeq y$ if and only if $\dist_{\cGT}(x,y)\leq  \fR/4$ or  there exist vertices $u,v$ with $\dist_{\cGT}(u,x)\leq1, \dist_{\cGT}(y,v)\leq 1$, and 
\begin{equation} \label{e:simdef}
|\GT_{uv}(z)|\geq
\varphi(z).
\end{equation}
We will say $x, y$ are Green's function correlated at the  threshold $\varphi(z)$ if  $x \simeq y$. 
For two sets $\bX,\bY\subset \qq{N}\setminus \bT$, we write $\bX\simeq \bY$\index{$\simeq$} if there exist $x\in \bX$ and $y\in \bY$ such that $x\simeq y$. Otherwise, we write $\bX\not\simeq \bY$.

The concept of Green's function correlated is an important tool introduced in \cite{MR3962004}.  The rationale is that the Ward identity \eqref{e:Ward-bis} provides only an $L_2$ bound on the Green's function. 
This bound is not strong enough since some  entires $G_{ij}$ can be large. These special entries needed to be treated with additional arguments which explore in particular the fact that the resampling  data are chosen uniformly and independently. We recall the index set $\As_{\bf S}$ from \eqref{Wdef}, for $\al\neq \beta\in \As_{\bf S}$ we have $\dist_{\cGT}(\{a_\al, b_\al, c_\al\}, \{a_\beta, b_\beta, c_\beta\})> \fR/4$.

\begin{definition}[Cells]\label{def:cells}
We denote $\cell_0$ the set 
\begin{align*}
\cell_0=\{\al \in \As_{\bf S}: \exists \beta\in \qq{1,\mu}: \{b_\al, c_\al\}\simeq a_\beta\}.
\end{align*}
The relation $\simeq$ and graph distance induce a relation on the set $\As_{\bf S}\setminus \cell_0$. We write $\al\approx \beta$ if and only if 
either $\{b_\al,c_\al\}\simeq \{b_\beta, c_\beta\}$. The relation $\approx$ induces a partition of the set
$\As_{\bf S}\setminus \cell_0=\cell_1\cup \cell_2\cup\cdots\cell_\kappa$. 
We define sets $\cella_0, \cella_1,\dots, \cella_\kappa\subset \qq{N}\setminus \bT$ called \emph{cells} by
\begin{align}\label{e:defC_i}
&\cella_0=\{a_1,a_2,\cdots, a_\mu\}\cup_{\al \in \cell_0}\{b_\al,c_\al\},\quad \cella_s= \cup_{\al\in \cell_s}\{b_\al, c_\al\}, \quad s=1,2,\cdots, \kappa.
\end{align}
\end{definition}
For later use, we note the following elementary properties of cells:
\begin{itemize}
\item
For any $x\in \cB_1(\cC_s, \cGT)$ and $y\in \cB_1(\cC_t,\cGT)$ such that $s \neq t$ we have $|\GT_{xy}(z)|< \varphi(z)$.
For any vertex $x\in \qq{N}\setminus\T$, if $x\not\simeq \cC_s$, then for any $y\in \cB_1(\cC_s,\cGT)$, $|\GT_{xy}(z)|<\varphi(z)$.
\item The cells are far away from each other, for any $s\neq t$ we have $\dist_{\cGT}(\cella_s, \cella_t)\geq \fR/4$. 
\end{itemize}

\subsubsection{Estimates}
The next proposition shows that the cells do not cluster. 
\begin{proposition} \label{prop:Rdist2}
Let $\cG \in \Omega(z)$ or $\Omega_o^+(z)$ as defined in Definitions \ref{def:Omega} and \ref{def:Omegao+} respectively. Then for any large number $\fd$,
with probabililty at least $1-\OO(N^{-\fd})$ under $P_{\cG}(\cdot)$ (as in Definition \ref{def:enlarge}),
the following estimates hold:
\begin{itemize}
\item
Any $x \in \qq{N} \setminus \T$ is $\simeq$-connected to fewer than $\OO(\log N)$ of vertices in $\{b_1, c_1,b_2,c_2,\dots, b_\mu, c_\mu\}$,
\begin{align}\label{finitedeg}
|\{\al\in \qq{1,\mu}: x\simeq \{b_\al,c_\al\}\}|\lesssim \log N.
\end{align}
In particular, $x$ is $\simeq$-connected to  $\OO(\log N)$ of the cells.
\item 
Most cells are in the form $\{b_\al, c_\al\}$
\begin{align}\begin{split}\label{e:Spcellbd}
|\{\al\in \As_{\bf S}: \al\in \cell_0 \text{ or } \al \in \cell_s, |\cell_s|>1\}|\lesssim \log N.
\end{split}\end{align}
In particular, each cell contains at most $\OO(\log N)$ of $\{b_\al, c_\al: \al\in \As_{\bf S}\}$.  
\end{itemize}
\end{proposition}

In the remainder of this section, we prove the above proposition.
It is essentially a straightforward
consequence of the definitions, combined with union bounds.

\begin{lemma}\label{lem:Rdist1}
Let $\cG \in \Omega(z)$ or $\Omega_o^+(z)$ as defined in Definitions \ref{def:Omega} and Proposition \ref{def:Omegao+} respectively. Then for any vertex $x\in \qq{N}\setminus \bT$, 
$
\Im[G_{xx}^{(\bT)}(z)]\lesssim  \Im[\md(z)]+\varepsilon'(z)+\varepsilon(z)/\sqrt{\kappa(z)+\eta(z)+\varepsilon(z)}
$ and
\begin{equation} \label{e:bsimx}
\P_{\cG}(\{b_\al,c_\al\}\simeq x)\lesssim  1/(\log N)^{ 4\fa}.
\end{equation}
\end{lemma}

\begin{proof}
We recall that for $\cG \in \Omega(z)$ or $\Omega_o^+(z)$ then 
\begin{align}\label{asumpGT}
 |Q(\cG,z)-m_{sc}(z)|\lesssim \frac{\varepsilon}{\sqrt{\kappa+\eta+\varepsilon}},\quad \left|\GT_{ij}-P_{ij}(\cB_r(i,j,\cGT),z, Q(\cG,z))\right|\lesssim \varepsilon',
\end{align}
for any vertices $i,j\in \qq{N}\setminus\bT$.
We can estimate $\Im[G_{xx}^{(\bT)}]$, 
\begin{align}\begin{split}\label{e:Imbound}
&\phantom{{}={}}\Im[G_{xx}^{(\bT)}]\lesssim |\Im [P_{xx}(\cB_r(x,\cGT), z,Q(\cG,z))]| +\varepsilon'\\
&\lesssim |\Im [P_{xx}(\cB_r(x,\cGT), z,\msc(z))]| +|Q(\cG,z)-\msc(z)|+\varepsilon'\\
&\lesssim \Im[\md(z)]+\frac{\varepsilon}{\sqrt{\kappa+\eta+\varepsilon}}+\varepsilon',
\end{split}\end{align}
where we used \eqref{e:smalldiff} in the second inequality, and \eqref{e:boundPij} for the last inequality.
The Ward identity \eqref{e:Ward} implies
\begin{align}\label{e:Wardpf}
\sum_i |\GT_{xi}|^2=\Im [\GT_{xx}]/\eta\lesssim (\Im[\md(z)]+\varepsilon'+\varepsilon/\sqrt{\kappa+\eta+\varepsilon})/\eta.
\end{align}
For any vertex $x\in \qq{N} \setminus \T$, set
\begin{align*}
\bV_x &:=\left\{u\in \qq{N} \setminus \T: |\GT_{xu}|\geq \varphi\right\},
\end{align*}
The inequality \eqref{e:Wardpf} and the definition of $\varphi$ in \eqref{e:defphi} imply 
\begin{align}\label{Vcondition}
|\bV_x|\lesssim N/(\log N)^{ 4\fa}.
\end{align}
Thus
\begin{align*}
\P_{\cG}(\{b_\al,c_\al\} \simeq x) 
&\leq P_{\cG}(\dist_{\cGT}(x, \{b_\al, c_\al\})\leq {\fR/4})+ \P_{\cG}(b_\al\text{ or }c_\al\in \cB_1( \cup_{u\in \cB_1(x, \cGT)}\bV_u, \cGT))\\
&\lesssim (d-1)^{\fR/4}/N+d^2|\bV_x|/N\lesssim (d-1)^{\fR/4}/N +d^2/(\log N)^{ 4\fa}
 \lesssim 1/(\log N)^{4\fa},
\end{align*}%
where the second inequality holds because $b_\al, c_\al$ are approximately uniform \eqref{e:approxunifbd}.
\end{proof}

\begin{proof}[Proof of Proposition \ref{prop:Rdist2}]
The proof of \eqref{finitedeg} is similar to that of \eqref{distfinitedegb}.
 From our construction of local resampling, $ \{b_\al,c_\al\}\simeq x$ is an independent event for different $\alpha$.  By the union bound and \eqref{e:bsimx}, we have
\begin{align*}
&\phantom{{}={}}\P( |\{\al\in \qq{1,\mu }: \{b_\al,c_\al\}\simeq x\}|\geq C\log N)
\lesssim {\mu \choose C\log N} \left(\frac{\OO(1)}{(\log N)^{4\fa}}\right)^{C\log N}\\
&\leq \left(\frac{\mu}{(\log N)^{4\fa}}\right)^{C\log N}\frac{\OO(1)^{C\log N}}{(C\log N)!}
\ll \left(\frac{\OO(1)}{(\log N)^{2\fa}}\right)^{C\log N} \lesssim  N^{-2\fa C \log \log N}
\ll N^{-\fd}
\end{align*}
where we used that $\mu\lesssim (d-1)^\ell\leq (\log N)^{2\fa}$.

For \eqref{e:Spcellbd},
we recall the index sets $\Ba, \Bb \subset \qq{1,\mu}$ from \eqref{distdega}, \eqref{distdegb}:
\begin{align*}
\Ba &=\{\al\in \qq{1,\mu}: {\dist_{\cGT}(a_\al, \{b_\al,c_\al\})\leq {\fR}}\},\\
\Bb &=\{\al\in \qq{1,\mu}: \exists \beta\neq \al, {\dist_{\cGT}(\{a_\al, b_\al,c_\al\}, \{a_\beta, b_\beta, c_\beta\})\leq {\fR/4}}\}.
\end{align*}
We notice that if $\al\in\cell_0$ or $\al\in \cell_s, |\cell_s|>1$, then either $\al\in \Ba \cup \Bb$ or $\{b_\al, c_\al\}\simeq \{a_\beta: \beta\in\qq{1,\mu}\}\cup\{b_\beta, c_\beta: \beta\in \As_{\bf S}\setminus\{\al\}\}$. For any $\al\in \qq{1,\mu}$, we apply  \eqref{e:bsimx} with $x=a_\beta$ or $ \{b_\beta, c_\beta,\beta\in \As_{\bf S}\setminus\{\al\}\}$, 
\begin{align}\begin{split}\label{e:cross0}
&\phantom{{}={}}\P_{\cG}\left( \{b_\al, c_\al\}\simeq \{a_\beta: \beta\in\qq{1,\mu}\}\cup \big \{b_\beta, c_\beta:  \beta\in \As_{\bf S}
\setminus\{\al\} \big \}\right)\\
&\leq \frac{\mu+2|\As_{\bf S}|}{(\log N)^{4\fa}}\lesssim\frac{1}{(\log N)^{2\fa}},
\end{split}\end{align}
where we used that $|\As_{\bf S}|\leq \mu\leq (\log N)^{2\fa}$.
Then by a union bound, we get
\begin{align}\begin{split}\label{e:cross}
&\phantom{{}={}}\P_{\cG}\Big ( \Big | \Big \{\al\in \qq{1,\mu}: \{b_\al, c_\al\}\simeq \{a_\beta: \beta\in\qq{1,\mu}\}\cup \big \{b_\beta, c_\beta:  \beta\in \As_{\bf S}
\setminus\{\al\} \big \}  \Big \} \Big |\geq C\log N       \Big )\\
&\leq {\mu \choose C\log N}\left(\frac{\OO(1)}{(\log N)^{2\fa}}\right)^{C\log N}
\leq \left(\frac{\mu}{\log N)^{2\fa}}\right)^{C\log N}\frac{\OO(1)^{C\log N}}{(C\log N)!}
\lesssim \left(\frac{e\OO(1)}{C\log N}\right)^{C\log N}\ll N^{-\fd},
\end{split}\end{align}%
provided $C$ is large enough.

From the discussion above,  only $\al\in \Ba \cup \Bb$ or $\{b_\al, c_\al\}\simeq \{a_\beta: \beta\in\qq{1,\mu}\}\cup\{b_\beta, c_\beta, \beta\in \As_{\bf S}\setminus\{\al\}\setminus\{\al\}\}$ contribute to the statement \eqref{e:Spcellbd}. 
Thus, combining \eqref{e:cross} with the estimate $|\Ba \cup \Bb|=\OO_\fd(1)$ from \eqref{distdega} and  \eqref{distdegb}, the estimate \eqref{e:Spcellbd} follows.
\end{proof}

\subsection{Stability under switching}
\label{sec:stability}

We recall the cells from Definition \ref{def:cells}, and the set of switching data $\sS(\cG)$ from Section \ref{sec:switch}.
In this section, we derive estimates for the Green's function $\tGT$ of the graph $\tcGT$, which is the switched graph obtained from $\tcG$ with vertices $\bT$ removed. Before stating Proposition~\ref{prop:stabilitytGT} , we need to introduce some sets. 

\begin{definition}\label{def:F2}
We define the set $F_2(\cG) \subset \sS(\cG)$\index{$F_2(\cG)$} of switching data $\bf S$, such that for the switching data $\bf S\in\sS(\cG)$, the following holds:
\begin{enumerate}
\item
All except for $\OO_{\fd}(1)$ of the vertices $\{c_1,c_2,\dots, c_\mu\}$ have radius-$\fR$ tree neighborhoods in $\cGT$,
i.e.\ \eqref{treeneighbor} holds.
\item  
The vertices $\{a_1,\dots, a_\mu,b_1,\dots, b_\mu, c_1,\cdots,c_\mu\}$ do not cluster in the sense of graph distance,
i.e. \eqref{distfinitedegb}--\eqref{distdegb} hold. \item  
The vertices $\{b_1,b_2,\dots, b_\mu, c_1, c_2,\cdots,c_\mu\}$ do not cluster in the sense of the Green's function distance,
i.e.\ \eqref{finitedeg} and \eqref{e:Spcellbd} hold.
\end{enumerate}
\end{definition}

Then, for any $\cG \in \Omega_o^+(z)$, we have 
\begin{align}\label{FGmeasure}
\P_{\cG}(F_2(\cG))= 1-\OO(N^{-\fd}).
\end{align}
Indeed,
(i) and (ii) follow from Proposition \ref{prop:distdeg},
(iii) follows from Proposition \ref{prop:Rdist2}.

\begin{definition}\label{def:sfJ}
We recall the admissible set $\As_{\bf S}\subset \qq{1,\mu}$ from \eqref{Wdef}, and define $\sJ\subset \As_{\bf S}$ to be the set  that for any $\al \in \sJ$ the following holds:

\begin{itemize}
\item In $\cGT$, $c_\al$ has a radius-${\fR}$ tree neighborhood;
\item There exists some $1\leq s\leq \kappa$, such that $\cella_{s}=\{b_\al, c_\al\}$.
In other words, the cell containing $\{b_\al, c_\al\}$ is $\cella_{s}=\{b_\al, c_\al\}$.
\end{itemize}

\end{definition}

For ${\bf S}\in F_2(\cG)$, at least $\mu-\OO_{\fd}(1)$ edges are switchable:
\begin{align} \label{switchnum} 
\nu=|\As_{\bf S}|\geq \mu-\OO_{\fd}(1).
\end{align}
Especially, for $\al\neq \beta\in \As_{\bf S}$, we have 
\begin{align}\label{e:ffar}
\dist_{\cGT}(\{a_\al, b_\al, c_\al\}, \{a_\beta, b_\beta,c_\beta\}), \dist_{\tcGT}(\{a_\al, b_\al, c_\al\}, \{a_\beta, b_\beta,c_\beta\})\geq \fR/4.
\end{align}
Moreover, Proposition \ref{structuretGT} implies that  \eqref{tcGTdist} and \eqref{e:lessshortdist} hold.
Finally thanks to Proposition \ref{prop:Rdist2}, we have that the size of the set $\sJ$ is at least $\mu-\OO(\log N)$. The results of this section are the following stability estimates.

\begin{proposition} \label{prop:stabilitytGT}
Let $d\geq 3$, $\omega_d$ as in Definition \ref{d:defwd} and $z\in \bC^+$ with $\Im[z]\geq (\log N)^\fb/N$.
 Let $\cG \in \Omega(z)$ or $\Omega_o^+(z)$ as defined in Definitions \ref{def:Omega} and \ref{def:Omegao+} respectively. 
Then there exists an event $F_2(\cG) \subset \sS(\cG)$ as in Definition \ref{def:F2} with $\P_{\cG}(F_2(\cG))= 1-\OO(N^{-\fd})$,
such that for any ${\bf S}\in F_2(\cG)$ with $\tcG  = T_{\bf S}(\cG)$ the following hold:
\begin{itemize}
\item
For $i,j \in \qq{N}\setminus\T\bW_{\bf S}$ (we write $\T\bW_{\bf S} = \T \cup \bW_{\bf S}$), 
\begin{align}\label{boundhGT}
 |G^{(\T\bW_{\bf S})}_{ij}(z)-G_{ij}(\Ext(\cB_{r}(i,j,\cG^{(\bT \bW_{\bf S})}),Q(\cG, z)),z)|\lesssim \varepsilon'(z).
\end{align}

\item 
For  vertices $i, j \in \qq{N} \setminus \T\bW_{\bf S}$ and indices $0\leq s\neq t\leq \kappa$, with either {\rm (i)}  $i\not\simeq \cC_t$ and $j\in \cB_1(\cC_t, \cGT)$, or {\rm (ii)} $i\in \cB_1(\cella_s, \cGT), j\in \cB_1(\cella_t, \cGT)$, the following holds
\begin{equation} \label{e:diffcellest}
|G^{(\bT \bW_{\bf S})}_{ij}(z)|\lesssim \varphi(z).
\end{equation}
We notice that {\rm (ii)} is almost but not exactly a special case of {\rm (i)}. It's possible  that  there is a vertex $i\in \cB_1(\cella_s, \cGT) \setminus \cella_s$ 
with  $i\simeq \cC_t$.

\item
For $i,j \in \qq{N}\setminus\T$,
\begin{equation} \label{e:stabilitytGT}
|\tGT_{ij}(z)-G_{ij}(\Ext(\cB_{r}(i,j,\tcGT),Q(\cG, z)),z)|
\lesssim \varepsilon'(z).
\end{equation}

\item
For the set $\sJ\subset \As_{\bf S}$ as in Definition \ref{def:sfJ}, it holds $|\As_{\bf S}\setminus \sJ|\lesssim \log N$ and
for any $\al\in \sJ$, 
\begin{alignat}{2}
\label{e:greendistIJ}
|\tGT_{ic_\al}(z)| &\lesssim (\varepsilon'(z))^2 +\varphi(z), &&\quad\text{if $i=c_\beta$ for some $\beta\in \qq{1,\mu}\setminus \sJ$},
\\
\label{e:greendistJJ}
|\tGT_{ic_\al}(z)| &\lesssim (\varepsilon'(z))^3 +\varphi(z), &&\quad \text{if $i=c_\beta$ for some $\beta\in \sJ\setminus \{\al\}$},
\\
\label{e:greendistNJ}
|\tGT_{ic_\al}(z)| &\lesssim (\varepsilon'(z))^2 +\varphi(z), &&\quad  \text{if $i\not\simeq \{b_\al,c_\al\}$ and $\dist_{\tcGT}(i,a_\al)\geq { \fR/4}$}.
\end{alignat}
\end{itemize}
\end{proposition}

\subsubsection{Proof of \eqref{boundhGT}}\label{proofofboundhGT}

We recall that $\bW_{\bf S}=\{b_\al: \al\in \As_{\bf S}\}$, and  let $\cG_o=\cB_{r}(\{\bW_{\bf S}, i,j\},\cGT)$. The graph $\cG_o^{(\bW_{\bf S})}$ is obtained from $\cG_o$ by removing the vertex set $\bW_{\bf S}$, and its deficit function is given by \eqref{e:removing}.
We abbreviate
\begin{equation*}
P=G(\Ext(\cG_o, Q(\cG,z))),\qquad
\quad P^{(\bW_{\bf S})}=G(\Ext(\cG_o^{(\bW_{\bf S})}, Q(\cG,z))).
\end{equation*}
Thanks to Proposition \ref{p:deficitbound} and our assumption that $\cG\in \Omega(z)$ or $\Omega_o^+(z)$, each connected component of $\cG_o$ and $\cG^{(\bW_{\bf S})}_o$ satisfies the assumptions in Proposition \ref{p:localization}. Thus we have
\begin{align}\begin{split}\label{replaceEr2}
&\absb{G_{ij}(\Ext(\cB_{r}(i,j,\cG^{(\bT)}),  Q(\cG,z)))-P_{ij}}
\lesssim (\log N)\varepsilon\lesssim \varepsilon',\\
&\absb{G_{ij}(\Ext(\cB_{r}(i,j,\cG^{(\bT\bW_{\bf S})}),  Q(\cG,z))-P_{ij}^{(\bW_{\bf S})}}
\lesssim (\log N)\varepsilon\lesssim \varepsilon',
\end{split}\end{align}
where the $\log N$ factor is from the term $1+\diam(\cG)$ in \eqref{e:compatibility}.
Thanks to  Proposition \ref{prop:stabilityGT} and \eqref{replaceEr2}, 
 we have
\begin{align}\begin{split}\label{e:bTb}
\left|\GT_{ij}-P_{ij}\right|
&\leq
\left|\GT_{ij}-G_{ij}(\Ext(\cB_r(i,j,\cGT), Q(\cG,z)))\right|\\
&+\left|G_{ij}(\Ext(\cB_r(i,j,\cGT),Q(\cG,z)))-P_{ij}\right|\lesssim \varepsilon'.
\end{split}\end{align}
To remove $\bW_{\bf S}$, we apply the Schur complement formula \eqref{e:Schur1}:
for any $i,j\in \qq{N}\setminus \bT \bW_{\bf S}$,
\begin{align}\label{GSschur}\begin{split}
 G_{ij}^{(\bT)}-G_{ij}^{(\bT \bW_{\bf S})}=&\sum_{x,y\in \bW_{\bf S}}G_{ix}^{(\bT)}(G^{(\bT)}|_{\bW_{\bf S}})^{-1}_{xy}G^{(\bT)}_{yj},\\
 P_{ij}-P_{ij}^{(\bW_{\bf S})}=&\sum_{x,y\in \bW_{\bf S}}P_{ix}(P|_{\bW_{\bf S}})^{-1}_{xy}P_{yj}.
 \end{split}
\end{align}
From \eqref{e:ffar}, any two points $x\neq y\in \bW_{\bf S}$ are far away from each other, i.e. $\dist_{\cG^{(\bT)}}(x,y)\geq {\fR/4}$. Therefore, $P|_{\bW_{\bf S}}$ is a diagonal matrix  of order one, for $x=b_\al,y=b_\beta\in \bW_{\bf S}$,
\begin{align}\label{3125}
|G^{(\bT)}_{xy}-P_{xy}|
\lesssim 
 \varepsilon' {\bf1}_{x\simeq y}+ {\bf1}_{x \not \simeq y} \varphi .
\end{align}
where we used the definition of Green's function correlated $\simeq$ in \eqref{e:simdef}.
Moreover, on $F_2(\cG)$, for any $ b_\al\in \bW_{\bf S}$ we have $|\{b_\beta\in \bW_{\bf S}: b_\beta\simeq b_\al\}|\lesssim \log N$. 
Therefore by the resolvent identity \eqref{e:resolv},  the fact  that $P|_{\bW_{\bf S}}$ is a diagonal matrix  of order one and the norm bound 
\[
\|(\GT-P)|_{\bW_\bfS}\|_{L_1\rightarrow L_\infty} \lesssim \e' \log N + |\bW_\bfS| \varphi \ll 1,
\]
we have 
\begin{align}\label{e:Ginv}
\Big | (G^{(\bT)}|_{\bW_{\bf S}})^{-1}_{xy} -(P|_{\bW_{\bf S}})^{-1}_{xy} \Big |   \lesssim   \varepsilon'{\bf1}_{x\simeq y}+ {\bf1}_{x \not \simeq y} \varphi. 
\end{align}
Notice that  $(G^{(\bT)}|_{\bW_{\bf S}})^{-1}_{xy} -(P|_{\bW_{\bf S}})^{-1}_{xy}$ and $G^{(\bT)}_{xy}-P_{xy}$ satisfy the same bound.

Using \eqref{e:bTb} and \eqref{e:Ginv} as input,  by taking difference of \eqref{GSschur}, we get
\begin{align}\begin{split}\label{e:GTbbb}
&\phantom{{}={}}|G_{ij}^{(\bT \bW_{\bf S})}-P_{ij}^{(\bW_{\bf S})}|
 \lesssim |\GT_{ij}-P_{ij}|
 +\sum_{x,y\in {\bW_{\bf S}}}|\GT_{ix}-P_{ix}||(G|_{\bW_{\bf S}})^{-1}_{xy}||\GT_{yj}|\\
&+\sum_{x,y\in {\bW_{\bf S}}}|P_{ix}||(G|_{\bW_{\bf S}})^{-1}_{xy}-(P|_{\bW_{\bf S}})^{-1}_{xy}||\GT_{yj}|+
\sum_{x,y\in {\bW_{\bf S}}}|P_{ix}||(P|_{\bW_{\bf S}})^{-1}_{xy}||\GT_{yj}-P_{yj}|
\lesssim \varepsilon'.
\end{split}
\end{align}
This together with \eqref{replaceEr2} imply the claim \eqref{boundhGT}.

\subsubsection{Proof of \eqref{e:diffcellest}}

Under both conditions given for \eqref{e:diffcellest},
we have $|\GT_{ij}| \leq \varphi$ by the definition of $\simeq$ as in \eqref{e:simdef}.
By the resolvent identity \eqref{GSschur}, we therefore have
\begin{equation}
\label{hGTdiffcell}
|G^{(\bT\bW_{\bf S})}_{ij}|\leq \varphi+\left|\sum_{x,y\in \bW_{\bf S}}G_{ix}^{(\bT)}(G^{(\bT)}|_{\bW_{\bf S}})^{-1}_{xy}G^{(\bT)}_{yj}\right|.
\end{equation}

We first prove \eqref{e:diffcellest} for case (i) that  $i\not\simeq \cC_t$ and $j\in \cB_1(\cC_t, \cGT)$. 
 From \eqref{e:ffar}, any two vertices $x, y\in \bW_{\bf S}$ with $x\neq y$, they are far away from each other, i.e. $\dist_{\cGT}(x,y)\geq {\fR/4}$. 
Moreover, on $F_2(\cG)$, thanks to \eqref{distfinitedegb},  
\[
|\{x\in \bW_{\bf S}: \dist_{\cGT}(i,x)\leq \fR/4\}|  \vee    |\{y\in \bW_{\bf S}: \dist_{\cGT}(y,j)\leq \fR/4\}|=\OO_\fd(1).
\]
There are several cases for the term $|G_{ix}^{(\bT)}(G^{(\bT)}|_{\bW_{\bf S}})^{-1}_{xy}G^{(\bT)}_{yj}|$ for   $x,  y\in \bW_{\bf S}$.
\begin{itemize}
\item If $\dist_{\cGT}(i,x)  \vee     \dist_{\cGT}(y,j)\leq \fR/4$, then $x\not\in \bS_t$, $y\in \bS_t$ and we have $|(G^{(\bT)}|_{\bW_{\bf S}})^{-1}_{xy}|\lesssim \varphi$ from \eqref{e:Ginv}. The number of such terms is at most $\OO_\fd(1)$, and total contribution is $\OO(\varphi)$.

\item If $\dist_{\cGT}(i,x)> \fR/4$, $\dist_{\cGT}(y,j)\leq \fR/4$ or $\dist_{\cGT}(i,x)\leq \fR/4$, $\dist_{\cGT}(y,j)> \fR/4$. For the former case, we have $y\in \bS_t$.  (i) If $x\not\in \bS_t$ then $x\not \simeq y$.  Hence   $|G_{ix}^{(\bT)}|\lesssim \varepsilon'$ (by \eqref{asumpGT}) and $|(G^{(\bT)}|_{\bW_{\bf S}})^{-1}_{xy}|\lesssim \varphi$ (by \eqref{e:Ginv}).  The number of such terms is at most $\OO(|\bW_{\bf S}|)$, 
and total contribution in   the last term in \eqref{hGTdiffcell}   is $\OO(\varphi \varepsilon' |\bW_{\bf S}| )$ where we have bounded $G^{(\bT)}_{yj}$ by a constant of order one. 
(ii) If $x, y\in \bS_t$ with $x\neq y$,  then  $|(G^{(\bT)}|_{\bW_{\bf S}})^{-1}_{xy}|\lesssim \varepsilon'$, $i\not\simeq x$ and $|G_{ix}^{(\bT)}|\lesssim \varphi$, by the definition of $\simeq$ as in \eqref{e:simdef}. The number of such terms is at most $\OO(|\bS_t|)=\OO(\log N)$, and total contribution in the last term in \eqref{hGTdiffcell} is $\OO(\log N \varepsilon'\varphi)$.
If $x= y\in \bS_t$ then $i\not\simeq x$ and $|G_{ix}^{(\bT)}|\lesssim \varphi$. The number of such terms is at most $\OO_\fd(1)$, and total contribution is $\OO(\varphi)$.
We have the same estimates for the latter case.

\item If $\dist_{\cGT}(i,x)  \wedge   \dist_{\cGT}(y,j)> \fR/4$, then by \eqref{asumpGT},  $|G_{ix}^{(\bT)}|  \vee  |G_{yj}^{(\bT)}|\lesssim \varepsilon'$.  (i) If vertices $x,y$ are in different cells, then  $|(G^{(\bT)}|_{\bW_{\bf S}})^{-1}_{xy}|\lesssim \varphi$ from \eqref{e:Ginv}, the number of such terms is at most 
$\OO(|\bW_{\bf S}|^2)$, and total contribution in   the last term in \eqref{hGTdiffcell}  is $\OO(\varphi(\varepsilon' |\bW_{\bf S}|)^2)$.
(ii) If $x,y$ are in the same cell $\bS_t$, then $|G_{ix}^{(\bT)}|\lesssim \varphi$.  The number of such terms is at most 
$\OO(|\bS_t|^2)=\OO((\log N)^2)$, and total contribution in the last term in \eqref{hGTdiffcell} is $\OO((\log N)^2 \varepsilon'\varphi)$.
(iii) If $x\neq y$ are in the same cell $\bS_s$ with some $s\neq t$, then $|(G^{(\bT)}|_{\bW_{\bf S}})^{-1}_{xy}|\lesssim \varepsilon'$ and $|G_{yj}^{(\bT)}|\lesssim \varphi$, the number of such terms is at most 
$\OO(\log N|\bW_{\bf S}|)$, and total contribution is $\OO(\log N\varphi\varepsilon'^2 |\bW_{\bf S}|)$.
(iv) If $x=y$ are in the same cell $\bS_s$ with some $s\neq t$, then $|G_{yj}^{(\bT)}|\lesssim \varphi$, the number of such terms is at most 
$\OO(|\bW_{\bf S}|)$, and total contribution is $\OO(\varphi \varepsilon' |\bW_{\bf S}|)$.
\end{itemize}
 Therefore, from the discussion above, we have
\begin{align*}
\sum_{x,y\in \bW_{\bf S}}\left|G_{ix}^{(\bT)}(G^{(\bT)}|_{\bW_{\bf S}})^{-1}_{xy}G^{(\bT)}_{yj}\right|
\lesssim \varphi+\varphi \varepsilon'  ( |\bW_{\bf S}|+ \log N )  +\varphi\varepsilon'  ( |\bW_{\bf S}|+ (\log N ) ^2)\lesssim \varphi,
\end{align*}
where we used that $  ( |\bW_{\bf S}|+ (\log N )^2) \varepsilon'\ll1$ from our choice of parameters \eqref{e:relation1} and \eqref{e:relation2} in the last inequality. 

For case (ii) that  $i\in \cB_1(\cella_s, \cGT), j\in \cB_1(\cella_t, \cGT)$ and $s\neq t$, although it is possible that $i\simeq \cC_t$, by our definition of cells, we have that $|\GT_{ij}|\leq \varphi$. It follows from essentially the same argument as in case (i). 
This finishes the proof of \eqref{e:diffcellest}.

\subsubsection{Proof of \eqref{e:stabilitytGT}}
The graph $\tcG^{(\bT)}$ is obtained from $\cG^{(\bT\bW_{\bf S})}$ by adding the vertices $\bW_{\bf S}$ and  for each $b_\al\in \bW_{\bf S}$, 
the edges between  $b_\al\in \bW_{\bf S}$ and the set  $(\{x: \dist_{\cG^{(\bT)}}(x, b_\al)=1\} \cup \{a_\alpha\})\setminus \{c_\alpha\}$. 
Let $\cG_o=\cB_{r}(\{\bW_{\bf S}, i,j\},\tcGT)$. The graph $\cG_o^{(\bW_{\bf S})}$ is obtained from $\cG_o$ by removing the vertex set $\bW_{\bf S}$, and its deficit function is given by \eqref{e:removing}.
We abbreviate
\begin{equation*}
P=G(\Ext(\cG_o, Q(\cG,z))),\qquad
\quad P^{(\bW_{\bf S})}=G(\Ext(\cG_o^{(\bW_{\bf S})}, Q(\cG,z))).
\end{equation*}
Thanks to Proposition \ref{p:deficitbound} and our assumption that $\cG\in \Omega(z)$ or $\Omega_o^+(z)$, each connected component of $\cG_o$ and $\cG^{(\bW_{\bf S})}_o$ satisfies the assumptions in Proposition \ref{p:localization}. Then we have  for $i,j \in \qq{N}\setminus\T$,
\begin{align}\begin{split}\label{replaceEr3}
&\absb{G_{ij}(\Ext(\cB_{r}(i,j,\tcG^{(\bT)}),  Q(\cG,z)))-P_{ij}}
\lesssim \varepsilon',\\
&\absb{G_{ij}(\Ext(\cB_{r}(i,j,\cG^{(\bT\bW_{\bf S})}), Q(\cG,z)))-P_{ij}^{(\bW_{\bf S})}}
\lesssim \varepsilon'.
\end{split}\end{align}
Thanks to \eqref{boundhGT} and \eqref{replaceEr3}, we have for all $i,j \in \qq{N}\setminus\T \bW_{\bf S} $
\begin{align}\begin{split}\label{e:TWTW}
\left|G_{ij}^{(\bT\bW_{\bf S})}-P^{(\bW_{\bf S})}_{ij}\right|
&\leq
\left|G_{ij}^{(\bT\bW_{\bf S})}-G_{ij}(\Ext(\cB_r(i,j,\cG^{(\bT\bW_{\bf S})}), Q(\cG,z)))\right|\\
&+\left|G_{ij}(\Ext(\cB_r(i,j,\cG^{(\bT\bW_{\bf S})}),Q(\cG,z)))-P^{(\bW_{\bf S})}_{ij}\right|\lesssim \varepsilon'.
\end{split}\end{align}

The normalized adjacency matrices of $\tcGT$ and $\cG_o$ respectively have the block form
\begin{align*}
\left[
\begin{array}{cc}
H & {B'}\\
{B} & D
\end{array}
\right],\quad
\left[
\begin{array}{cc}
H & {B}_o'\\
{B}_o & D_o
\end{array}
\right],
\end{align*}
where $H$ is the normalized adjacency matrix for $\tcGT|_{\bW_{\bf S}}$ (which is a zero matrix),
and $ B$ (respectively $ B_o$) corresponds to the edges from $\bW_{\bf S}$ to $\qq{N}\setminus \bT \bW_{\bf S}$.
Notice that for nonzero entries $ B_{ij}=( B_o)_{ij}$, in the rest of this section we will therefore not distinguish
$B$ and ${B}_o$.

To  estimate  $\tGT_{ij}$ for  $i,j\in \bW_{\bf S}$, by the Schur complement formula \eqref{e:Schur}, we have
\begin{align}
  \label{G1xSchur1d}
 \tGT|_{\bW_{\bf S}}&=(H-z-{B}'G^{(\bT\bW_{\bf S})}{B})^{-1},\\
 \label{G1xSchur2d}
 P|_{\bW_{\bf S}}&=(H-z-{B}'P^{(\bW_{\bf S})}{B})^{-1}.
\end{align}
From \eqref{e:ffar}, any two vertices $x\neq y \in \bW_{\bf S}$, they are far away from each other, i.e. $\dist_{\tcG^{(\bT)}}(x, y)\geq \fR/4$. 
Thus $P|_{\bW_{\bf S}}$ is a diagonal matrix, with diagonal entries $|P_{xx}|\asymp 1$. Therefore, for any $i,j\in \bW_{\bf S}$, by taking difference of \eqref{G1xSchur1d} and \eqref{G1xSchur2d}, and using \eqref{e:TWTW},
\begin{align}\label{e:ijAA}
 |\tG^{(\bT)}_{ij}- P_{ij}|
 =|(( P|_{\bW_{\bf S}})^{-1}+B'(P^{(\bW_{\bf S})}-G^{(\bT\bW_{\bf S})})B)^{-1}_{ij}-P_{ij}|\lesssim \varepsilon'.
\end{align}

For the estimates of $\tGT_{ij}$ where $i\in \bW_{\bf S}, j\in \qq{N}\setminus \bT \bW_{\bf S}$, we have by the Schur complement formula \eqref{e:Schur1}:
\begin{align}\label{e:ss1}
 \tilde{G}^{(\bT)}=-\tilde{G}^{(\bT)}{B}'G^{(\bT\bW_{\bf S})},\quad
  P=-P{B}'P^{( \bW_{\bf S})}.
\end{align}
Therefore, by taking the difference of these two equations,
\begin{align}\begin{split}\label{e:iAj}
|\tilde{G}^{(\bT)}_{ij}-P_{ij}|
&\leq \sum_{x\in \bW_{\bf S}}|\tGT_{ix}||(B
'(G^{(\bT\bW_{\bf S})}-P^{(\bW_{\bf S})}))_{xj}| \\
&+  \sum_{x\in \bW_{\bf S}}  |(\tGT-P)_{ix}||(B'P^{(\bW_{\bf S})})_{xj}|
\lesssim \varepsilon'+(\varepsilon')^2|\bW_{\bf S}|\lesssim \varepsilon'.
\end{split}\end{align}
Here to bound  the first summation, we divide it into two cases, $x=i$ and $x \neq i$. For $x\neq i$, we have $\dist_{\tcGT}(i,x)\geq \fR/4$ thus $|\tGT_{ix}|\lesssim \varepsilon'$ by \eqref{e:ijAA}, and 
 $|(B'(G^{(\bT\bW_{\bf S})}-P^{(\bW_{\bf S})}))_{xj}|\lesssim \varepsilon'$ 
by \eqref{e:TWTW}. This gives the contribution  $\OO((\varepsilon')^2|\bW_{\bf S}|)$. For $x=i$ with similar argument, it contributes $\OO(\varepsilon')$.
For the second summation in \eqref{e:iAj}, we divide the summation according to  $ \dist_{\cGT}(x,j)\leq \fR/4$ or $ \dist_{\cGT}(x,j)>  \fR/4$. We notice that   
 on $F_2(\cG)$,  $|\{x\in \bW_{\bf S}: \dist_{\cGT}(x,j)\leq \fR/4\}|=\OO_\fd(1)$ by \eqref{distfinitedegb}. The rest of the argument is similar to the one for bounding the first summation.

%
%
%
%
%
%
%
%
%
%

To  estimate  $\tGT_{ij}$ for  $i,j\in \qq{N}\setminus \bT$, 
we have the Schur complement formula \eqref{e:Schur}:
\begin{align}\label{e:ss3}
 \tGT=\tG^{(\bT\bW_{\bf S})}+\tG^{(\bT\bW_{\bf S})} B\tGT B' \tG^{(\bT\bW_{\bf S})},\quad
{P}=P^{(\bW_{\bf S})}+P^{(\bW_{\bf S})} B P B'P^{(\bW_{\bf S})}.
\end{align}
By taking the difference, 
\begin{align}\begin{split}\label{e:inoAj}
 \tGT_{ij}-P_{ij} =&\tG^{(\bT\bW_{\bf S})}_{ij}-P^{(\bW_{\bf S})}_{ij}
 +\sum_{x,y\in \bW_{\bf S}}\left(((\tG^{(\bT\bW_{\bf S})}-P^{(\bW_{\bf S})}) B)_{ix}\tGT_{xy} (B' \tG^{(\bT\bW_{\bf S})})_{yj}\right.\\
 +&\left.({P}^{(\bW_{\bf S})} B)_{ix}(\tGT-P)_{xy} (B' \tG^{(\bT\bW_{\bf S})})_{yj}+(P^{(\bW_{\bf S})} B)_{ix}P_{xy} (B' (\tG^{(\bT\bW_{\bf S})}-\P^{(\bW_{\bf S})}))_{yj}\right).
\end{split}\end{align}

 We claim the following bounds hold: 
\begin{align}\label{311}
|((\tG^{((\bT\bW_{\bf S})}-P^{(\bW_{\bf S})})B)_{ix}|\lesssim  {\bf1}_{i \simeq x}\varepsilon' + {\bf1}_{i \not \simeq x} \varphi,\quad |\{x\in \bW_{\bf S}: i\simeq x\}|\lesssim \log N,
\end{align}
and
\begin{align}\begin{split}\label{312}
&|(B'\tG^{(\bT\bW_{\bf S})})_{yj}|\lesssim  {\bf1}_{\{\dist_{\cG^{(\bT)}}(j,y)\leq \fR/4\}} +{\bf1}_{j \simeq y}\varepsilon' + {\bf1}_{j \not \simeq y} \varphi,\\
&|\{y\in \bW_{\bf S}: \dist_{\cG^{(\bT)}}(j,y)\leq \fR/4\}|=\OO_{\fd}(1),\quad |\{y\in \bW_{\bf S}: j\simeq y\}|\lesssim \log N. 
\end{split}\end{align}
In fact, by definition, 
 \begin{align}\label{e:difGexp}
((\tG^{(\bT\bW_{\bf S})}-P^{(\bW_{\bf S})})B)_{ix} = \sum_{z} (\tG^{(\bT\bW_{\bf S})}-P^{(\bW_{\bf S})})_{iz} B_{zx}= \frac{1}{\sqrt{d-1}}\sum_{z\sim x} (\tG^{(\bT\bW_{\bf S})}-P^{(\bW_{\bf S})})_{iz}.
\end{align}%
We remind that  from our construction $\tG^{(\bT\bW_{\bf S})} =  G^{(\bT\bW_{\bf S})}$.
By  \eqref{e:TWTW}, $|\tG^{(\bT\bW_{\bf S})}_{iz}-P^{(\bW_{\bf S})}_{iz}|\lesssim \varepsilon'$.
Since the number of $z$ in the summation is bounded (by the degree $d$), we have proved that $
|(\tG^{((\bT\bW_{\bf S})}-P^{(\bW_{\bf S})})B)_{ix}|\lesssim  \varepsilon'$. Suppose now that $ i\not\simeq x $. By  \eqref{e:diffcellest} and \eqref{e:difGexp}, we get
$|(\tG^{((\bT\bW_{\bf S})}-P^{(\bW_{\bf S})})B)_{ix}|\lesssim \varphi$. The first relation of  \eqref{311} follows.
On $F_2(\cG)$, thanks to  \eqref{finitedeg},  most vertices in $\bW_{\bf S}$ are far away from $i$ in terms of the Green's function distance; more precisely, we have $|\{x\in \bW_{\bf S}: i\simeq x\}|\lesssim \log N$. 

In the following, we prove  \eqref{312}. 
For $y\in \bW_{\bf S}$ with $\dist_{\cG^{(\bT)}}(j,y)\leq \fR/4$, we have $|(B'\tG^{(\bT\bW_{\bf S})})_{yj}|\lesssim 1$. 
For $\dist_{\cG^{(\bT)}}(j,y)> \fR/4$,  we notice that $|(B' \tG^{(\bT\bW_{\bf S})})_{yj}  |= |(B' (\tG^{(\bT\bW_{\bf S})} - P^{(\bW_{\bf S})})_{yj}  |$ satisfies a bound similar to \eqref{311}. Therefore, the first relation in claim  \eqref{312} follows.
On $F_2(\cG)$, thanks to \eqref{distfinitedegb} and \eqref{finitedeg},  most vertices in $\bW_{\bf S}$ are far away from $j$ in terms of graph distance and Green's function distance; more precisely, 
we have $|\{y\in \bW_{\bf S}: \dist_{\cG^{(\bT)}}(j,y)\leq \fR/4\}|=\OO_{\fd}(1)$ and $|\{y\in \bW_{\bf S}: j\simeq y\}|\lesssim \log N$.

We can use \eqref{311} and \eqref{312} to estimate the second term on the righthand side of \eqref{e:inoAj}, when $x\neq y$, 
\begin{align*}
&\phantom{{}={}}\sum_{x\neq y\in \bW_{\bf S}}\left|((\tG^{(\bT\bW_{\bf S})}-P^{(\bW_{\bf S})}) B)_{ix}\tGT_{xy} (B' \tG^{(\bT\bW_{\bf S})})_{yj}\right|\\
&\lesssim (\log N\varepsilon'+|\bW_\bfS|\varphi)\varepsilon'(\OO_\fd(1)+\log N \varepsilon'+|\bW_\bfS|\varphi),
\end{align*}
and when $x=y$, we simply bound $|\tGT_{xy}|$ by $\OO(1)$,
\begin{align*}
\sum_{x= y\in \bW_{\bf S}}\left|((\tG^{(\bT\bW_{\bf S})}-P^{(\bW_{\bf S})}) B)_{ix}\tGT_{xy} (B' \tG^{(\bT\bW_{\bf S})})_{yj}\right|
\lesssim \varepsilon'+\log N(\varepsilon')^2+|\bW_\bfS|\varphi^2.
\end{align*}
We have similar estimates for the rest terms on the righthand side of \eqref{e:inoAj}.
Therefore, we have the following estimate for \eqref{e:inoAj}
\begin{align}\label{e:ijnoA}
| \tGT_{ij}-P_{ij}|\lesssim \varepsilon' + \log N (\varepsilon')^2+|\bW_{\bf S}|\varepsilon' \varphi\lesssim \varepsilon',
\end{align}
where we used that $(\log N+|\bW_{\bfS}|) \varepsilon'\ll1$ from our choice of parameters \eqref{e:relation1} and \eqref{e:relation2}.  
The estimates \eqref{e:ijAA}, \eqref{e:iAj} and \eqref{e:ijnoA} together give the claim \eqref{e:stabilitytGT}.

\subsubsection{Proof of \eqref{e:greendistIJ}-- \eqref{e:greendistNJ}}

Under extra assumptions of the locations of the vertices $i,j$ we have better estimates of $\tGT_{ij}$. In the following we prove claims \eqref{e:greendistIJ}--\eqref{e:greendistNJ}.
By the first identity in \eqref{e:ss3}, we have 
\begin{align}\begin{split}\label{e:ss3copy}
|\tGT_{ij}|&=|\tG^{(\bT \bW_{\bf S})}_{ij}|+\sum_{x,y\in \bW_{\bf S}}|(\tG^{(\bT\bW_{\bf S})} B)_{ix}\tGT_{xy} (B' \tG^{(\bT\bW_{\bf S})})_{yj}|.
\end{split}\end{align}
In all three cases in \eqref{e:greendistIJ}--\eqref{e:greendistNJ}, \eqref{e:diffcellest} implies $|\tG^{(\bT \bW_{\bf S})}_{ij}|\lesssim \varphi$.

To prove  \eqref{e:greendistIJ}, we choose  $j=c_\al$ for some $\al\in \sJ$ in the previous equation.  By Definition  \ref{def:sfJ} of the set $\sJ$, the same as \eqref{312}, we have
\begin{align}\label{e:BGB}
|(B' \tG^{(\bT\bW_{\bf S})})_{yj}|\lesssim 
\left\{
\begin{array}{cc}
\varepsilon', &y=b_\al,\\
\varphi, & y\in \bW_{\bf S}\setminus\{ b_\al\}.
\end{array}
\right.
\end{align}
We divide  \eqref{e:ss3copy} into two cases $y\neq b_\al$ or $y=b_\al$ so that 
\begin{align}\begin{split}\label{e:divide}
|\tGT_{ij}|
&\lesssim \varphi+\sum_{x,y\in \bW_{\bf S},y\neq b_\al}|(\cdots)|+\sum_{x\in \bW_{\bf S},y= b_\al}|(\cdots)|.
\end{split}\end{align}

Since $i=c_\beta$ for some $\beta\in \qq{1, \mu}\setminus \sJ$, 
by \eqref{e:stabilitytGT} 
we have 
\[
\max_{x\in \bW_{\bf S}} \sum_{y\in \bW_{\bf S}} |\tGT_{xy} |\lesssim \max_{x\in \bW_{\bf S}} \sum_{y\in \bW_{\bf S}} | P_{xy}| +|\bW_{\bf S}|\varepsilon' \lesssim   1+|\bW_{\bf S}|\varepsilon'.
\]
For  $y\neq b_\al$,  we use \eqref{e:BGB},    sum over $y\in \bW_{\bf S}$ and then $x\in \bW_{\bf S}$ to get 
\begin{align}\begin{split}\label{e:sumbb0}
\sum_{x,y\in \bW_{\bf S},y\neq b_\al}|(\cdots)|
&\lesssim\sum_{x,y\in \bW_{\bf S}}|(\tG^{(\bT\bW_{\bf S})} B)_{ix}\tGT_{xy}|\varphi
\lesssim\sum_{x\in \bW_{\bf S}}|(\tG^{(\bT\bW_{\bf S})} B)_{ix}|(1+|\bW_{\bf S}|\varepsilon')\varphi\\
&\lesssim(\log N \varepsilon'+|\bW_{\bf S}|\varphi)(1+|\bW_{\bf S}|\varepsilon')\varphi\lesssim \varphi, 
\end{split}\end{align}
where we have used the relations \eqref{e:relation1} and \eqref{e:relation2} for the choices of parameters. 
For the case $y=b_\al$, we further divide it into two cases $x=b_\al$ or $x\neq b_\al$ and use \eqref{e:BGB} so that 
\begin{align}\begin{split}\label{e:sumbb1}
\sum_{x\in \bW_{\bf S}}|(\cdots)| 
&\lesssim\sum_{x\in \bW_{\bf S}}|(\tG^{(\bT\bW_{\bf S})} B)_{ix}\tGT_{xb_\al}|\varepsilon'
\lesssim (\varepsilon')^2+\sum_{b_\al\neq x\in \bW_{\bf S}}|(\tG^{(\bT\bW_{\bf S})} B)_{ix}|(\varepsilon')^2\\
&\lesssim(1+\log N \varepsilon'+|\bW_{\bf S}|\varphi)(\varepsilon')^2\lesssim (\varepsilon')^2.
\end{split}\end{align}
The Claim \eqref{e:greendistIJ} follows from plugging \eqref{e:sumbb0} and \eqref{e:sumbb1} into \eqref{e:divide}.

For \eqref{e:greendistJJ}, $i=c_\beta$ for some $\beta\in \sJ\setminus \{\alpha\}$, we have the same estimate \eqref{e:BGB} for $(\tG^{(\bT\bW_{\bf S})} B)_{ix}$. Similarly to \eqref{e:divide}, we divide the sum into four cases depending on $x=b_\beta$ or $x\neq b_\beta$, and $y=b_\al$ or $y\neq b_\al$,
\begin{align*}\begin{split}
|\tGT_{i   c_\alpha}|
&\lesssim \varphi+\sum_{x,y\in \bW_{\bf S}\atop x\neq b_\beta, y\neq b_\al}|(\cdots)|+\sum_{x\in \bW_{\bf S}\setminus \{b_\beta\},y= b_\al}|(\cdots)|+\sum_{x=b_\beta,y\in \bW_{\bf S}\setminus \{b_\al\}}|(\cdots)|+\sum_{x=b_\beta,y=b_\al}|(\cdots)|.
\end{split}\end{align*}%
The same argument as in \eqref{e:sumbb0} and \eqref{e:sumbb1} gives
\begin{align*}
|\tGT_{ij}|\lesssim \varphi+(\varepsilon')^3,
\end{align*}
where the leading order term $\OO((\varepsilon')^3)$ is from the case that $x=b_\beta$ and $y=b_\al$.
This proves  \eqref{e:greendistJJ}.

For \eqref{e:greendistNJ}, if $i=b_\beta\in \bW_{\bf S}$, $i\not\simeq\{b_\al, c_\al\}$, and {$\dist_{\tcGT}(i,a_\al)\geq { \fR/4}$}, then $|\tGT_{ib_\al}|\lesssim \varepsilon'$  by \eqref{e:stabilitytGT}. 
On $F_2(\cG)$, thanks to \eqref{e:lessshortdist}, the number of $x\in \bW_{\bf S}$ such that $|\tGT_{ix}|\lesssim 1
$ is $\OO_\fd(1)$, and $ |\tGT_{ix}| \lesssim \varepsilon'$ for the rest of $x\in \bW_{\bf S}$. 
Therefore \eqref{e:ss1} gives that
\begin{align}\begin{split}\label{e:iAj2}
|\tGT_{ij}|
&\lesssim\sum_{x\in \bW_{\bf S}}|\tGT_{ix}||(B
'G^{(\bT\bW_{\bf S})})_{xj}|=
\sum_{x\in \bW_{\bf S},x\neq b_\al}|(\cdots)|
+\sum_{x= b_\al}|(\cdots)|\\
&\lesssim (1+|\bW_{\bf S}|\varepsilon')\varphi+(\varepsilon')^2\lesssim \varphi+(\varepsilon')^2.
\end{split}\end{align}
 Finally for the case $i\in \qq{N}\setminus \bT\bW_{\bf S}$, $i\not\simeq \{b_\al, c_\al\}$ and $\dist_{\tcGT}(i,a_\al)\geq \fR/4$, similar argument as for \eqref{e:inoAj} implies 
\begin{align*}\begin{split}
|\tGT_{ij}|&=|\tG^{(\bT\bW_{\bf S})}_{ij}|+\sum_{x,y\in \bW_{\bf S}}|(\tG^{(\bT\bW_{\bf S})} B)_{ix}\tGT_{xy} (B' \tG^{(\bT\bW_{\bf S})})_{yj}|\\
&\lesssim \varphi+\sum_{x,y\in \bW_{\bf S},y\neq b_\al}|(\cdots)|+\sum_{x\in \bW_{\bf S}\setminus\{b_\al\},y= b_\al}|(\cdots)|+\sum_{x=y=b_\al}|(\cdots)|
\lesssim \varphi+(\varepsilon')^2,
\end{split}\end{align*}
where the leading order term $\OO((\varepsilon')^2)$ is from the case $x=y=b_\al$.
This together with \eqref{e:iAj2} proves  \eqref{e:greendistNJ}.

\begin{remark}[Remark on Methods] The Remark \ref{r1} can be applied to estimate $G^{(\bT\bW_{\bf S})}$ from $\GT$. We  had also  used the Ward identity \eqref{e:Ward-bis} and  Green's function correlated notion. Once 
 $G^{(\bT\bW_{\bf S})}$ is estimated, we note that   by definition $G^{(\bT\bW_{\bf S})} = \tG^{(\bT\bW_{\bf S})}$. Thus our  next task is to estimate $\tG^{(\bT)}$ from  $G^{(\bT\bW_{\bf S})}$. For  $i,j\in \bW_{\bf S}$, we  use the  Schur complement formulas \eqref{G1xSchur1d} and \eqref{G1xSchur2d}. For $i\in \bW_{\bf S}, j\in \qq{N}\setminus \bT \bW_{\bf S}$, we have  the Schur  complement formula \eqref{e:ss1}. Finally, for  $i,j\in \qq{N}\setminus \bT$, 
we use  the Schur complement  formula \eqref{e:ss3}. Similar ideas will be applied  to estimate $\tG$ from bounds on $\tG^{(\bT)}$. We have thus bounded $\tG$ from bounds on $G$ through the path 
\begin{align*}
G \longrightarrow \GT  \longrightarrow G^{(\bT\bW_{\bf S})} = \tG^{(\bT\bW_{\bf S})}  \longrightarrow \tG^{(\bT)}  \longrightarrow \tG.
\end{align*}
Notice that  in each step in this path most estimates deteriorate (in the sense that most estimates get worse by constant factors)  and hence this procedure only provides  stability bounds on the Green's functions. In order to use continuity arguments, 
we will need to improve the final estimates by a small constant factor in certain key terms (otherwise after a few iterations the constants will increase exponentially). The key improvement will be achieved by a concentration argument exploring the fact that the switching data are essentially independent variables. This will be done in Section \ref{s:improveG}. 
\label{r2}
\end{remark}

\subsection{Stability estimate for the switched graph}
\label{sec:weakstab}
In this section we prove the following proposition, which gives
the estimates on the Green's function of the switched graph $\tcG$.
\begin{proposition}\label{prop:tGweakstab}
Let $d\geq 3$,  $\omega_d$ as in Definition \ref{d:defwd} and $z\in \bC^+$ with $\Im[z]\geq (\log N)^\fb/N$.
Recall that  $\cG \in \Omega(z)$ or $\Omega_o^+(z)$ were  defined in Definitions \ref{def:Omega} and \ref{def:Omegao+} respectively. For any vertex $i\in \bT$, let $\ell_i=\dist_{\cT}(o,i)$ the distance from $i$ to $o$.
Then for ${\bf S} \in F_2(\cG)$ (defined  in Definition \ref{def:F2}) and with notation 
$\tcG=T_{\bf S}(\cG)$, 
the Green's function of the switched graph satisfies the stability estimates
\begin{align}
\label{e:Gaxbound}
\begin{split}
&
|\tG_{ij}(z)-G_{ij}(\Ext(\cB_r(\{\bT,i,j\},\tcG), Q(\cG, z)),z)|\\
&\lesssim\left(\frac{(1+|\ell_i-\ell_j|)}{(d-1)^{|\ell_i-\ell_j|/2}}+\frac{(\log N)^2}{(d-1)^{(2\ell-\ell_i-\ell_j)/2}}\right)\varepsilon'(z),   \; \; \forall \ i,j \in \bT; 
\end{split}
\\
\begin{split}
\label{e:GTNbound} 
 |\tG_{ij}(z)-G_{ij}(\Ext(\cB_r(\{\bT, i,j\},\tcG),  Q(\cG, z)),z)|  
  \lesssim\frac{(\log N)^2\varepsilon'(z)}{(d-1)^{(\ell-\ell_i)/2}}, \quad \forall \ i\in \bT, j\in\qq{N}\setminus \bT;  %
\end{split}
\\
&|\tG_{ij}(z)-G_{ij}(\Ext(\cB_r(\{\bT, i,j\},\tcG),  Q(\cG, z)),z)|\lesssim \log N \varepsilon'(z), \quad \forall\  i,j\in \qq{N}\setminus\bT.\label{e:GNNbound}
\end{align}%
\end{proposition}

We enumerate new outer boundary vertices in the switched graph $\tcG$,
 by $\bI=\{\ta_1,\ta_2,\dots,\ta_\mu\}$, where $\ta_\al=a_\al$ if $\al\in \qq{1,\mu}\setminus \As_{\bf S}$ and $\ta_\al=c_\al$ if $\al\in \As_{\bf S}$.
Let $\sJ$ be the index set defined in Definition \ref{def:sfJ}.
Let $\cG_o=\cB_{r}(\{\bT, i,j\},\tcG)$. The graph $\cG_o^{(\bT)}$ is obtained from $\cG_o$ by removing the vertex set $\bT$, and its deficit function is given by \eqref{e:removing}. We abbreviate
\begin{align*}
P=G(\Ext(\cG_o,  Q(\cG, z)),z),\qquad
\quad P^{(\bT)}=G(\Ext(\cG_o^{(\bT)}, Q(\cG, z)),z).
\end{align*}
The normalized adjacency matrices of $\tcG$ and $\Ext(\cG_o, Q(\cG, z))$ respectively have the block form
\begin{align*}
\left[
\begin{array}{cc}
H & {B'}\\
{B} & D
\end{array}
\right],\quad
\left[
\begin{array}{cc}
H & {B}_o'\\
{B}_o & D_o
\end{array}
\right],
\end{align*}
where $H$ is the normalized adjacency matrix for $\cT$,
and $B$ (respectively $B_o$) corresponds to the edges from $\bI$ to $\T_\ell=\{v\in \bT: \dist_{\cT}(o,v)=\ell\}$,
where $\bI$ is the set of outer boundary vertices of $\cT$ in the switched graph $\tcG$, and $\bT_\ell$ is the inner vertex boundary of $\cT$.
The nonzero entries of $B$ and $B_o$ occur for the indices
$(i,j) \in \bI\times \T_\ell$ and take values $1/\sqrt{d-1}$.
Notice that $B_{ij}=(B_o)_{ij}$. In the rest of this section we will therefore not distinguish
$B$ and $B_o$.

By the Schur Complement formula \eqref{e:Schur}, we have
\begin{align}
  \label{G1xSchur1}
 \tG|_\T&=(H-z-{B}'\tGT{B})^{-1},\\
 \label{G1xSchur2}
  P|_\T&=(H-z-{B}'P^{(\T)}{B})^{-1},
\end{align}
and, by the resolvent identity \eqref{e:resolv}, the difference of \eqref{G1xSchur1} and \eqref{G1xSchur2} is
\begin{equation}\label{Tterm}
 \tG|_\T-P|_{\T}=P{B}'(\tGT-P^{(\T)}){B}P + (\tG-P){B}'(\tGT-P^{(\T)}){B}P
\end{equation}
We remark that the last term is an higher order term in the sense that it is of order  $(\tG-P)^2$ instead of linear in $(\tG-P)$ .
\subsubsection{Boundary estimates}
\begin{proposition}\label{l:smallsum}
  Under the assumptions of Proposition~\ref{prop:tGweakstab}, 
  for any vertex $j\in \qq{N} \setminus \T$ we have
  \begin{align}\label{sumtGtPi}
    \sum_{\al\in \qq{1,\mu}}|\tGT_{\ta_\al j}(z)-P^{(\T)}_{\ta_\al j}(z)|\lesssim (\log N)\varepsilon'(z),
  \end{align}
  and
   \begin{align}\label{sumtGtPij}
    \sum_{\al\neq \beta\in \qq{1,\mu}}|\tGT_{\ta_\al\ta_\beta}(z)-P^{(\T)}_{\ta_\al \ta_\beta}(z)|\lesssim  (\log N)^2\varepsilon'(z).
  \end{align}
\end{proposition}

\begin{proof}[Proof of Proposition \ref{l:smallsum}]
Recall the index set $\sJ\subset \As_{\bf S}$ as in Definition \ref{def:sfJ}.
To prove \eqref{sumtGtPi}, we decompose $\qq{1,\mu}$ according to the relations between $\{a_\al, b_\al, c_\al\}$ and  vertex $j$ as
$\qq{1,\mu}=\sJ_1\cup \sJ_2$, where 
 \begin{align*}
&\sJ_1=\{\al\in \sJ: j\not\simeq \{b_\al, c_\al\},\text{ }\dist_{\tcGT}(j,a_\al)\geq \fR/4, \text{ and } \dist_{\tcGT}(j,\ta_\al)\geq \fR/4\},\\
&\sJ_2=\qq{1,\mu}\setminus \sJ_1.
\end{align*}%
By the defining relations of $F_2(\cG)$ as in Definition \ref{def:F2}, we have $|\sJ|\geq \mu-\OO(\log N)$.
Combining with \eqref{e:lessshortdist} and \eqref{finitedeg}, we get 
\begin{align*}
|\sJ_1|\geq \mu-\OO(\log N),\quad |\sJ_2|\lesssim \log N.
\end{align*}
Now, for $\al\in \sJ_1$, we have $\ta_\al=c_\al$ and the conditions for \eqref{e:greendistNJ} are satisfied.
Moreover, since $\dist_{\tcGT}(\ta_\al,j)\geq \fR/4$, 
the vertices $\ta_\al$ and $j$ are in different connected components of $\cGT_o$;
it follows that $P^{(\T)}_{\ta_\al j}=0$.
Therefore, by \eqref{e:greendistNJ},
\begin{align}\label{sumtGtPiexp1}
 \sum_{\al\in \sJ_1}|\tGT_{\ta_\al j}-P^{(\T)}_{\ta_\al j}|
 =\sum_{\al\in \sJ_1}|\tGT_{\ta_\al j}|
 \lesssim |\sJ_1|\left((\varepsilon')^2+\varphi\right) \lesssim \varepsilon',
\end{align}
where we used that  $|\sJ_1|\leq \mu\lesssim (d-1)^{\ell}$,  $(d-1)^\ell \varepsilon'\ll1$ and $(d-1)^\ell \varphi\ll \varepsilon'$ from \eqref{e:relation2}.
For $\al\in \sJ_2$, by \eqref{e:stabilitytGT}, we have
\begin{align}\label{sumtGtPiexp2}
\sum_{\al\in \sJ_2}|\tGT_{\ta_\al j}-P^{(\T)}_{\ta_\al j}|\lesssim |\sJ_2|\varepsilon'\lesssim\log N \varepsilon'.
\end{align}
Then \eqref{sumtGtPi} follows by combining \eqref{sumtGtPiexp1} and \eqref{sumtGtPiexp2}.

For \eqref{sumtGtPij}, we split the sum over 
\begin{align*}
\{\al\neq \beta\in \qq{1,\mu}\}=&\{\al\neq \beta \in \qq{1,\mu}\setminus \sJ\}\cup \{\al\in \qq{1,\mu}\setminus \sJ,  \beta \in \sJ\}\\
\cup& \{\al\in \sJ,  \beta \in \qq{1,\mu}\setminus \sJ\}\cup \{\al\neq \beta\in \sJ\}.
\end{align*}
For $\al\neq \beta\in \qq{1,\mu}\setminus \sJ$, by \eqref{e:stabilitytGT} and $|\qq{1,\mu}\setminus \sJ|\lesssim \log N$, we have
\begin{align*}
\sum_{\al\neq \beta\in \qq{1,\mu}\setminus \sJ}|\tGT_{\ta_\al\ta_\beta}-P^{(\T)}_{\ta_\al\ta_\beta}|\lesssim (\log N)^2 \varepsilon'.
\end{align*}
For $\al\in \qq{1,\mu}\setminus \sJ,  \beta \in \sJ$,  $\ta_\al$ and $\ta_\beta$ are in different connected components of $\tcGT_o$, and thus 
$|P_{\ta_\al\ta_\beta}^{(\T)}|=0$. By \eqref{e:greendistIJ}, 
\begin{align*}
\sum_{\al\in \qq{1,\mu}\setminus \sJ,  m \in \sJ}|\tGT_{\ta_\al\ta_\beta}-P^{(\T)}_{\ta_\al\ta_\beta}|
&=\sum_{\al\in \qq{1,\mu}\setminus \sJ,  \beta \in \sJ}|\tGT_{\ta_\al\ta_\beta}|
 \lesssim (\log N)|\sJ|\left((\varepsilon')^2+\varphi\right)\lesssim (\log N)\varepsilon'.
\end{align*}%
The same estimate holds for $\al\in \sJ,  \beta \in \qq{1,\mu}\setminus \sJ$.
For $\al\neq \beta\in \sJ$, the same reasoning as above gives $P_{\ta_\al\ta_\beta}^{(\T)}=0$, and using \eqref{e:greendistJJ}
\begin{align*}
\sum_{\al\neq \beta\in \sJ}|\tGT_{\ta_\al\ta_\beta}-P^{(\T)}_{\ta_\al\ta_\beta}|=\sum_{\al\neq \beta\in \sJ}|\tGT_{\ta_\al\ta_\beta}|\lesssim |\sJ|^2\left((\varepsilon')^3+\varphi\right)\lesssim \varepsilon'.
\end{align*}
Now \eqref{sumtGtPij} follows by combining the above four cases.
\end{proof}

\subsubsection{Proof of \eqref{e:Gaxbound}}

\begin{proof}
For any vertices $i,j\in\bT$, by 
 \eqref{e:boundPijcopy}, we have 
 \begin{align}\label{e:Pijcc}
|P_{ij}|\lesssim \left(\frac{|m_{sc}|}{\sqrt{d-1}}\right)^{\dist_{\cT}(i,j)}\leq \left(\frac{1}{\sqrt{d-1}}\right)^{\dist_{\cT}(i,j)}.
\end{align}
We recall $\T_\ell=\{v\in \bT: \dist_{\cT}(o,v)=\ell\}$, and $\ell_i=\dist_{\cT}(o,i)$ the distance from $i$ to $o$. Let $\Gamma_i=\max_{v\in \bT_\ell}|\tilde{G}_{iv}-P_{iv}|$.
For $j\in \bT$ with $\dist_{\cT}(o,j)=\ell_j$,  we claim the following bounds
\begin{align}\label{e:sumPpj}
\sum_{\alpha\in\qq{1,\mu}}|P_{il_\alpha}|
\lesssim \sum_{\alpha\in\qq{1,\mu}}\left(\frac{1}{\sqrt{d-1}}\right)^{\dist_{\cT}(i, l_\alpha)}\lesssim (1+\ell_i)(d-1)^{(\ell-\ell_i)/2},
\end{align}
and
 \begin{align}\label{e:sPPcc}
\sum_{\al \in\qq{1,\mu}}|P_{il_\al}||P_{l_\al j}|
\lesssim \sum_{\alpha\in\qq{1,\mu}}\left(\frac{1}{\sqrt{d-1}}\right)^{\dist_{\cT}(i, l_\alpha)}\left(\frac{1}{\sqrt{d-1}}\right)^{\dist_{\cT}(l_\alpha,j)}
\lesssim
\frac{(1+|\ell_i-\ell_j|)}{(d-1)^{|\ell_i-\ell_j|/2}}.
\end{align}%
The first inequality in \eqref{e:sumPpj} and \eqref{e:sPPcc} follows from \eqref{e:Pijcc}. To prove the second inequality in \eqref{e:sumPpj}, we construct a new graph $\cal C$ from $\cT$ by repeating the following pruning procedure. For any leaf vertex $x$, if $x\neq i,o$, we remove $x$. The each vertex in the graph $\cal C$ has degree at least two, except for $i,o$. 
The original graph $\cT$ can be obtained from $\cal C$ by attaching certain truncated $d$-regular trees. Therefore, for any two vertices $x,y\in \cal C$, we have $\dist_{\cal C}(x,y)=\dist_{\cT}(x,y)$. For any $\al \in \qq{1,\mu}$, the shortest path from $i$ to $l_\al$ first stays in $\cal C$, leaves it at some vertex $x(l_\al)$ and enters the subtree in $\cT\setminus \cal C$ containing $l_\al$. The depth of the subtree is $\dist_\cT(x(l_\al), l_\al)=\ell-\dist_{\cal C}(o, x(l_\al))$. Therefore, for the sum over $\al\in \qq{1,\mu}$ in \eqref{e:sumPpj}, we can first sum over $\al$ such that $x(l_\al)=x\in \cal C$, the total number of such indices $\al$ is $\OO((d-1)^{\ell-\dist_{\cal C}(o, x(l_\al))})$
\begin{align}\begin{split}\label{e:dto}
 &\phantom{{}={}}\sum_{\alpha\in\qq{1,\mu}}\left(\frac{1}{\sqrt{d-1}}\right)^{\dist_{\cT}(i, l_\alpha)}
 =\sum_{x\in \cal C}  \sum_{\alpha\in\qq{1,\mu}, x=x(l_\al)}\left(\frac{1}{\sqrt{d-1}}\right)^{\dist_{\cT}(i, l_\alpha)}\\
 &= \sum_{x\in \cal C}  \sum_{\alpha\in\qq{1,\mu}, x=x(l_\al)}\left(\frac{1}{\sqrt{d-1}}\right)^{\dist_{\cal C}(i, x)+\dist_{\cT}(x,l_\al)}
\lesssim \sum_{x\in \cal C} \sqrt{d-1}^{\ell-\dist_{\cal C}(o,x)}\left(\frac{1}{\sqrt{d-1}}\right)^{\dist_{\cal C}(i, x)}\\
&=(d-1)^{\ell/2}\sum_{x\in \cal C} \left(\frac{1}{\sqrt{d-1}}\right)^{\dist_{\cal C}(i, x)+\dist_{\cal C}(o,x)}\leq (d-1)^{\ell/2}\sum_{q\geq 0}\sum_{x\in \cal C, \dist_{\cal C}(o, x)=q} \left(\frac{1}{\sqrt{d-1}}\right)^{q+|\ell_i-q|},
\end{split}\end{align}%
where for the last inequality we used that for $x\in \cal C$ with $\dist_{\cal C}(o,x)=q$, $\dist_{\cal C}(i,x)\geq |\dist_{\cal C}(o,i)-\dist_{\cal C}(i,x)|=|\ell_i-q|$.

Since $\cG\in \bar \Omega$, $\cT$ has excess at most $\omega_d$. So does the graph $\cal C$. Next we show the number of vertices $x\in \cal C$ with $\dist_{\cal C}(o,x)=q$ is at most $2\omega_d+2$. If $q=0$, it is necessary that $x=o$. If $q\geq 1$, we enumerate such vertices as $x_1,x_2,\cdots,x_k$. For any $x_s$, the shortest path in $\cal C$ from $o$ to $x_s$ contains the edge $\{y_s, x_s\}$ (which may not be unique). We denote $\cA$ the graph from $\cal C$ by removing the edges $\{y_s, x_s\}$ for $1\leq s\leq k$. 
We partition $\qq{1,k}$ into sets  $\{\sA_1, \sA_2, \cdots, \sA_\kappa\}$,
such that $s,t$ are in the same set $\sA_\tau$ iff  $x_s$ and $x_t$ are in the same connected component of $\cA$. We recall the definition  of $\text{excess}(\cG)$ from \eqref{def:excess},
\begin{align}\label{c-c2}
\text{excess}(\cG)
  =\#\text{connected components}(\cG)-
  \#\text{vertices}(\cG)+\#\text{edges}(\cG).
\end{align}
The graph $\cA$ is obtained from $\cal C$ by removing the edges $\{y_s, x_s\}$ for $1\leq s\leq k$. If $\sA_\tau$ is a singleton, i.e. $\sA_\tau=\{s\}=1$, then the connected component of $\cA$ containing $x_s$ contains vertex $i$ or has excess at least $1$. In the second case,  the number of excess decreases by at least one if we remove the edge $\{y_s, x_s\}$ from $\cal C$. If $|\sA_\tau|\geq 2$, thanks to \eqref{c-c2}, the number of excess decreases by $|\sA_\tau|-1$ if we remove edges $\{\{x_s, y_s\}: s\in \sA_\tau\}$  from $\cal C$. Therefore, the excess of $\cal C$ satisfies $\omega_d\geq \text{excess}(\cal C)\geq -1+ \sum_{\tau=1}^{\kappa}\max\{1, |\sA_\tau|-1\}$. By rearranging it, we have in particular that 
\begin{align}\label{e:size}
k=\sum_{\tau=1}^\kappa |\sA_\tau|\leq 
2\sum_{\tau=1}^{\kappa}\max\{1, |\sA_\tau|-1\}\leq 2\omega_d+2.
\end{align}
By plugging \eqref{e:size} into \eqref{e:dto}, we get 
\begin{align}\begin{split}\label{e:countq}
 (d-1)^{\ell/2}\sum_{q\geq 0}\sum_{x\in \cal C, \dist_{\cal C}(o, x)=q} \left(\frac{1}{\sqrt{d-1}}\right)^{q+|\ell_i-q|}
&\leq  (d-1)^{\ell/2}\sum_{q\geq 0} \frac{2\omega_d+2}{(d-1)^{(q+|\ell_i-q|)/2}}\\
&\lesssim (1+\ell_i)(d-1)^{(\ell-\ell_i)/2}.
\end{split}\end{align}
This finishes the  second inequality in \eqref{e:sumPpj}.

For the second inequality in \eqref{e:sPPcc}, similarly to the proof of \eqref{e:sumPpj}, we construct a new graph $\cal C$ from $\cT$ by repeating the pruning procedure: if $x$ is a leaf vertex and $x\neq i,j,o$, we remove $x$. 
For any $\al \in \qq{1,\mu}$, the shortest paths from $i,j$ to $l_\al$ first stay in $\cal C$, leave it at some vertex $x(l_\al)$ and enter the subtree in $\cT\setminus \cal C$ containing $l_\al$. The depth of the subtree is $\dist_\cT(x(l_\al), l_\al)=\ell-\dist_{\cal C}(o, x(l_\al))$. Therefore, for the sum over $\al\in \qq{1,\mu}$ in \eqref{e:sPPcc}, we can first sum over $\al$ such that $x(l_\al)=x\in \cal C$, the total number of such indices $\al$ is $\OO((d-1)^{\ell-\dist_{\cal C}(o, x(l_\al))})$
\begin{align}\begin{split}\label{e:sPPcc2}
 &\phantom{{}={}}
 \sum_{\alpha\in\qq{1,\mu}}\left(\frac{1}{\sqrt{d-1}}\right)^{\dist_{\cT}(i, l_\alpha)+\dist_{\cT}(l_\alpha,j)}
 = \sum_{x\in \cal C}  \sum_{\alpha\in\qq{1,\mu}, x=x(l_\al)}\left(\frac{1}{\sqrt{d-1}}\right)^{\dist_{\cal C}(i, x)+2\dist_{\cT}(x,l_\al)+\dist_{\cal C}(j, x)}\\
&\lesssim\sum_{x\in \cal C} \left(\frac{1}{\sqrt{d-1}}\right)^{\dist_{\cal C}(i, x)+\dist_{\cal C}(x,j)}\leq \sum_{q\geq 0}\sum_{x\in \cal C, \dist_{\cal C}(i, x)=q} \left(\frac{1}{\sqrt{d-1}}\right)^{q+|\dist_\cT(i,j)-q|}\\
&\lesssim \frac{1+\dist_\cT(i,j)}{(d-1)^{\dist_{\cT}(i,j)/2}}\lesssim
\frac{(1+|\ell_i-\ell_j|)}{(d-1)^{|\ell_i-\ell_j|/2}}
\end{split}\end{align}%
where for the second line we used that for $x\in \cal C$ with $\dist_{\cal C}(i,x)=q$, $\dist_{\cal C}(j,x)\geq |\dist_{\cal C}(i,j)-\dist_{\cal C}(i,x)|=|\dist_{\cal T}(i,j)-q|$. For the last two inequalities we used $\dist_{\cal T}(i,j)\geq|\ell_i-\ell_j|$ and \eqref{e:countq}.

%

We can bound the second term on the right-hand side of \eqref{Tterm}  by
\begin{align}\begin{split}\label{firstterm}
&\phantom{{}={}} |((\tG-P){B}'(\tGT-P^{(\bT)}){B}P)_{ij}|
  \lesssim \Gamma_i\sum_{\al,\beta\in \qq{1,\mu}}|\tGT_{\ta_\al\ta_\beta}-P^{(\bT)}_{\ta_\al\ta_\beta}||P_{l_\beta j}|\\
 &\lesssim \log N \varepsilon'\Gamma_i \sum_{\beta\in\qq{1,\mu}}|P_{l_\beta j}|
 \lesssim(\log N) (1+\ell_j)(d-1)^{(\ell-\ell_j)/2}\varepsilon' \Gamma_i=\oo(1)\Gamma_i,
 \end{split}
\end{align}
where we used \eqref{sumtGtPi} and \eqref{e:sumPpj}.

The first term on the right-hand side of \eqref{Tterm} can be bounded by 
\begin{align}\begin{split}\label{sumtwo}
|(P{B}'(\tGT-P^{(\bT)}){B}P)_{ij}|
 &\lesssim\sum_{\al\in\qq{1,\mu}}|P_{il_\al}||\tGT_{\ta_\al\ta_\al}-P^{(\bT)}_{\ta_\al\ta_\al}||P_{l_\al j}|\\
 &+\sum_{\al\neq \beta\in \qq{1,\mu}}|P_{il_\al}||\tGT_{\ta_\al\ta_\beta}-P^{(\bT)}_{\ta_\al\ta_\beta}||P_{l_\beta j}|.
\end{split}\end{align}
We can estimate the first term on the righthand side of \eqref{sumtwo} by 
\begin{align}\label{e:bbt}
\sum_{\al\in\qq{1,\mu}}|P_{il_\al}||\tGT_{\ta_\al\ta_\al}-P^{(\bT)}_{\ta_\al\ta_\al}||P_{l_\al j}|
\lesssim \varepsilon'\sum_{\al \in\qq{1,\mu}}|P_{il_\al}||P_{l_\al j}|
\lesssim \frac{\varepsilon'(1+|\ell_i-\ell_j|)}{(d-1)^{|\ell_i-\ell_j|/2}}
,
\end{align}
where we used \eqref{e:stabilitytGT} in the first inequality and \eqref{e:sPPcc} in the second inequality. 

For the second term on the righthand side of \eqref{sumtwo},  by  \eqref{e:greendistIJ}, \eqref{e:greendistJJ}  and \eqref{e:greendistNJ}, we have,
\begin{align}\begin{split}\label{e:off-d0}
&\phantom{{}={}}\sum_{\al\neq \beta\in \qq{1,\mu}}|P_{il_\al}||\tGT_{\ta_\al\ta_\beta}-P^{(\bT)}_{\ta_\al\ta_\beta}||P_{l_\beta j}|\lesssim
\sum_{\al\neq \beta \in \qq{1,\mu}\setminus \sJ}|P_{il_\alpha}|\varepsilon'|P_{l_\beta j}|\\
&+\left (\sum_{\beta \in \qq{1,\mu}\setminus \sJ\atop \al\in \sJ}+\sum_{\al\in \qq{1,\mu}\setminus \sJ\atop \beta \in \sJ}\right )|P_{il_\al}|\left((\varepsilon')^2+\varphi\right) |P_{l_\beta j}|
+ \sum_{\al\neq \beta \in \sJ}|P_{il_\al}|\left((\varepsilon')^3+\varphi\right)|P_{l_\beta j}|.
\end{split}\end{align}
Using \eqref{e:Pijcc} and $\mu - |\sJ| \lesssim \log N$, the first term on the righthand side of \eqref{e:off-d0} is bounded as 
\begin{align*}
\sum_{\al\neq \beta \in \qq{1,\mu}\setminus \sJ}|P_{il_\alpha}|\varepsilon'|P_{l_\beta j}|\lesssim \frac{(\log N)^2\varepsilon'}{(d-1)^{(2\ell-\ell_i-\ell_j)/2}}.
\end{align*}
By  \eqref{e:Pijcc},  \eqref{e:sumPpj}  and $\mu - |\sJ| \lesssim \log N$,  we can bound the second term on the righthand side of \eqref{e:off-d0} as
\begin{align*}\begin{split}
\sum_{\beta \in \qq{1,\mu}\setminus \sJ\atop \al\in \sJ}|P_{il_\al}|\left((\varepsilon')^2+\varphi\right) |P_{l_\beta j}|
&\lesssim \log N (1+\ell_i)(d-1)^{(\ell-\ell_i)/2}(d-1)^{-(\ell-\ell_j)/2}
\left((\varepsilon')^2+\varphi\right)\\
&\lesssim \log N (1+\ell_i)(d-1)^{(\ell_j-\ell_i)/2}
\left((\varepsilon')^2+\varphi\right).
\end{split}\end{align*}
Using  \eqref{e:sumPpj}, the last term on the righthand side of \eqref{e:off-d0} is bounded by
\begin{align*}
\sum_{\al\neq \beta \in \sJ}|P_{il_\al}|\left((\varepsilon')^3+\varphi\right)|P_{l_\beta j}|\lesssim
(1+\ell_i)(1+\ell_j)\left(d-1\right)^{(2\ell-\ell_i-\ell_j)/2}\left((\varepsilon')^3+\varphi\right).
\end{align*}
Collecting the above estimates together, we get
\begin{align}  \begin{split}
&\eqref{e:off-d0}
\lesssim \frac{(\log N)^2\varepsilon'}{(d-1)^{(2\ell-\ell_i-\ell_j)/2}}+(\log N)(1+\min\{\ell_i,\ell_j\})(d-1)^{ |\ell_i-\ell_j|/2}\left((\varepsilon')^2+\varphi\right)\\
&+(1+\ell_i)(1+\ell_j)\left(d-1\right)^{(2\ell-\ell_i-\ell_j)/2}\left((\varepsilon')^3+\varphi\right)\lesssim \frac{(\log N)^2\varepsilon'}{(d-1)^{(2\ell-\ell_i-\ell_j)/2}}
,
\end{split} \label{e:off-d}
\end{align}
where we used \eqref{e:relation2}.
It follows from plugging \eqref{firstterm}, \eqref{sumtwo}, \eqref{e:bbt} and  \eqref{e:off-d} into \eqref{Tterm},
\begin{align}\begin{split}\label{tG1xminusP1x}
&\phantom{{}={}}|\tG_{ij}-P_{ij}|\leq  
\oo(1)\Gamma_i+\left(\frac{(1+|\ell_i-\ell_j|)}{(d-1)^{|\ell_i-\ell_j|/2}}+\frac{(\log N)^2}{(d-1)^{(2\ell-\ell_i-\ell_j)/2}}\right)\varepsilon'.
\end{split}\end{align}
By taking the maximum over $j\in \bT_\ell$ and rearranging it,
we get 
\begin{align}\label{e:Gammaxb}
\Gamma_i
&\lesssim
\frac{(\log N)^2}{(d-1)^{|\ell_i-\ell|/2}}\varepsilon'.
\end{align}
Plugging \eqref{e:Gammaxb} into \eqref{firstterm} and combining with \eqref{tG1xminusP1x}, we get 
\begin{align*}\begin{split}
&\phantom{{}={}}|\tG_{ij}-P_{ij}|\lesssim
\left(\frac{(1+|\ell_i-\ell_j|)}{(d-1)^{|\ell_i-\ell_j|/2}}+\frac{(\log N)^2}{(d-1)^{(2\ell-\ell_i-\ell_j)/2}}\right)\varepsilon'.
\end{split}\end{align*}
This finishes the proof of \eqref{e:Gaxbound}.

\subsubsection{Proof of \eqref{e:GTNbound}}

For $i\in \bT$ and $j\not\in \bT$, we use the Schur complement formula \eqref{e:Schur1}, i.e., $ \tilde{G}=-\tilde{G}{B}'\tGT,\quad  P=-P{B}'P^{(\bT)}$. 
By taking the difference of these two equations,
\begin{equation}\label{sumtwo2}
|\tilde{G}_{ij}-P_{ij}|\leq \frac{1}{\sqrt{d-1}} \sum_{\al \in\qq{1,\mu}}|\tG_{il_\al}||P^{(\bT)}_{\ta_\al j}-\tGT_{\ta_\al j}| + \frac{1}{\sqrt{d-1}}\sum_{\al\in\qq{1,\mu}}|\tG_{il_\al}-P_{il_\al }||P^{(\bT)}_{\ta_\al j}|.
\end{equation}
The first term in \eqref{sumtwo2} is bounded by
\begin{align*}
\sum_{\al \in\qq{1,\mu}}|\tG_{il_\al}||P^{(\bT)}_{\ta_\al j}-\tGT_{\ta_\al j}| 
&\lesssim \frac{\log N \varepsilon'}{(d-1)^{(\ell-\ell_i)/2}},
\end{align*}
where we used \eqref{sumtGtPi} and $|\tG_{il_\al}|\lesssim 1/(d-1)^{(\ell-\ell_i)/2}$.
  By  \eqref{e:lessshortdist}, $P^{(\bT)}_{\ta_\al i}$ are zero for all $\al\in \qq{1,\mu}$ except for at most $\OO_\fd(1)$ of them. Hence we can bound 
   the second term in \eqref{sumtwo2} by 
\begin{align*}
\frac{1}{\sqrt{d-1}}\sum_{\al\in\qq{1,\mu}}|\tG_{il_\al}-P_{il_\al}||P^{(\bT)}_{\ta_\al j}|
\lesssim  \frac{(\log N)^2\varepsilon'}{(d-1)^{(\ell-\ell_i)/2}},
\end{align*}
where we have used \eqref{e:Gammaxb}. Combining these bounds, we have proved that 
\begin{align*}
|\tilde{G}_{ij}-P_{ij}|\lesssim \frac{(\log N)^2\varepsilon'}{(d-1)^{(\ell-\ell_i)/2}}.
\end{align*}

\subsubsection{Proof of \eqref{e:GNNbound}}

For $i,j\not\in \bT$, by  \eqref{e:Schur} we have $  \tilde{G}= \tGT+\tGT{B}\tilde G{B}' \tGT,  \;{P}={P}^{(\bT)}+{P}^{(\bT)}{B}\tilde P{B}'{P}^{(\bT)}$. 
Taking their difference, we have 
\begin{align*}
 \tilde{G}-{P} =&\tGT-{P}^{(\bT)}+(\tGT-{P}^{(\bT)}){B}\tilde G {B}' \tGT\\
 +&\tilde{P}^{(\bT)}{B}(\tilde G-\tilde{P}) {B}' \tGT+{P}^{(\bT)}{B}{P}{B}' (\tGT-{P}^{(\bT)}).
\end{align*}
Notice that \eqref{e:stabilitytGT} implies that $|\tGT_{ij}-{P}^{(\bT)}_{ij}|\lesssim \varepsilon'$. Thus  we have
\begin{align}\begin{split}\label{GPxyout}
& |\tilde{G}_{ij}-{P}_{ij}|\lesssim \varepsilon'
 +\sum_{\al,\beta\in\qq{1,\mu}}|\tGT_{i\ta_\al}-{P}^{(\bT)}_{i\ta_\al}||\tG_{l_\al l_\beta}||\tGT_{\ta_\beta j}|\\
+&\sum_{\al,\beta\in\qq{1,\mu}}|{P}^{(\bT)}_{i\ta_\al}||\tG_{l_\al l_\beta}-{P}_{l_\al l_\beta}| |\tGT_{\ta_\beta j}|
+\sum_{\al,\beta\in\qq{1,\mu}}|{P}^{(\bT)}_{i\ta_\al}||{P}_{l_\al l_\beta}||\tGT_{\ta_\beta j}-{P}^{(\bT)}_{\ta_\beta j}|.
\end{split}
\end{align}
 From \eqref{e:greendistNJ}, we have the following estimates:
\begin{align*}
&\sum_{\al \in\qq{1,\mu}}|{P}^{(\bT)}_{i\ta_\al}|, \sum_{\beta \in\qq{1,\mu}}|\tGT_{\ta_\beta j}| \lesssim 1+\log N\varepsilon'+(d-1)^\ell\left((\varepsilon')^2+\varphi\right) \lesssim 1,\\
&\sum_{\al \in\qq{1,\mu}}|\tGT_{i\ta_\al}-P^{(\bT)}_{i\ta_\al}|, \sum_{\beta\in\qq{1,\mu}}|\tGT_{\ta_\beta j}-P^{(\bT)}_{\ta_\beta j}| \lesssim \log N\varepsilon'+(d-1)^\ell \left( (\varepsilon')^2+\varphi\right)\lesssim  \log N\varepsilon'.
\end{align*}
Therefore \eqref{GPxyout} can be  simplified  to
\begin{align*}
 |\tilde{G}_{ij}-{P}_{ij}|\lesssim  \log N \varepsilon'.
\end{align*}

\end{proof}

\subsection{Changing \texorpdfstring{$Q(\cG,z)$}{Q(G,z)} to \texorpdfstring{$Q(\tcG,z)$}{Q(tG,z)}}

\label{sec:changeQG}
%

In this section, we prove that we can replace $Q(\cG,z)$ by $Q(\tcG,z)$ up to a small error.
It follows from the general insensitivity of the quantity $Q$ to small changes of the graph. We recall the definition of $\Omega^+_o(z)$ from Definition
\ref{def:Omegao+}, and the set $F_2(\cG)$ from Definition \ref{def:F2}.

\begin{proposition}\label{l:IGchange}
For $z\in \bC^+$ with $\Im[z]\geq (\log N)^\fb/N$, $\cG \in \Omega_o^+(z)$ as in Definition \ref{def:Omegao+},  ${\bf S} \in F_2(\cG)$ and $\tcG=T_{\bf S}(\tG)$,
we have
\begin{align}\label{e:IGchange}
|\I(\cG,z)-\I(\tcG,z)|
\lesssim  \frac{(d-1)^{2\ell}{ (\Im[\md(z)]+\varepsilon'(z)+\varepsilon(z)/\sqrt{\kappa(z)+\eta(z)+\varepsilon(z)})}}{N\eta(z)}+\frac{1}{N^{1-\fc}}.
\end{align}%
\end{proposition}


\begin{proof}[Proof of Proposition \ref{l:IGchange}]
We denote the the indicator function $\chi_j(\cG^{(\bT\bW_{\bf S})})=1$ if the vertex $j$ is away from the vertices involved in the local resampling, i.e. $\dist_{\cG^{(\bT\bW_{\bf S})}}(j, \{a_1, b_1, c_1,\cdots, a_\mu, b_\mu, c_\mu\})>\fR/4$; otherwise, it is $0$. The number of vertices $j$ such that $\chi_j( \cG^{(\bT\bW_{\bf S})})=0$ is bounded by $\OO((d-1)^{\ell+\fR/4})\ll N^\fc$. 
Proposition \ref{l:IGchange} follows by proving
\begin{align}\begin{split}\label{GTtotGTconcentration}
&\phantom{{}={}}\Big | Q(\cG,z)-\frac{1}{Nd}\sum_{(i,j)\in \vec E(\cG^{(\bT\bW_{\bf S})})}\chi_j(\cG^{(\bT\bW_{\bf S})})G_{ii}^{(\bT\bW_{\bf S} j)}\Big |\\
&\lesssim   \frac{(d-1)^{2\ell}{ (\Im[\md(z)]+\varepsilon'+\varepsilon/\sqrt{\kappa+\eta+\varepsilon})}}{N\eta}+\frac{1}{N^{1-\fc}}, \\
&\phantom{{}={}}\Big |Q(\tcG,z)-\frac{1}{Nd}\sum_{(i,j)\in \vec E(\cG^{(\bT\bW_{\bf S})})}\chi_j( \cG^{(\bT\bW_{\bf S})})G_{ii}^{(\bT\bW_{\bf S} j)}\Big | \\
&\lesssim   \frac{(d-1)^{2\ell}{ (\Im[\md(z)]+\varepsilon'+\varepsilon/\sqrt{\kappa+\eta+\varepsilon})}}{N\eta}+\frac{1}{N^{1-\fc}},
\end{split}\end{align}
where $\vec E(\cG^{(\bT\bW_{\bf S})})$ is the set of directed edges of $\cG^{(\bT\bW_{\bf S})}$.

The proofs of the two estimates in \eqref{GTtotGTconcentration} are analogous, and we will only prove the last one.
We denote by $\vec{E}$ the set of oriented edges of $\tcG^{(\bT\bW_{\bf S})}$, and $\del_{\cG^{(j)}}(\bT\bW_{\bf S})$  the set of outer boundary vertices of $\bT\bW_{\bf S}$ in $\cG^{(j)} $. For vertex $j$ such that $\chi_j( \cG^{(\bT \bW_{\bf S})})=1$, it is far away from $\bT \bW_{\bf S}$.  For $\cG\in \Omega_o^+(z)$ and $\bfS\in F_2(\cG)$, \eqref{prop:tGweakstab} implies the Green's functions of $\tcG=T_\bfS(\cG)$ are bounded. Thus for $x,y\in \bT \bW_{\bf S}$, $|\tG^{(j)}_{xy}|=|\tG_{xy}-\tG_{xj}\tG_{jy}/\tG_{jj}|\lesssim 1$. Then,
by the Schur complement formula \eqref{e:Schur}, and noticing $|(\tG^{(j)}|_{\bT\bW_\bfS})_{xy}^{-1}|\lesssim 1$ for $x,y\in \del_{\cG^{(j)}}(\bT\bW_\bfS)$ we have 
\begin{align}\label{e:refk}
\begin{split}
 &\phantom{{}={}}\sum_{(i,j)\in \vec{E}}\chi_j(\cG^{(\bT\bW_{\bf S})})|\tG^{(\T\bW_{\bf S} j)}_{ii}-\tG^{(j)}_{ii}|\\
 &\lesssim \sum_{(i,j)\in \vec{E}}\chi_j(\cG^{(\bT\bW_{\bf S})})\sum_{x,y\in \del_{\cG^{(j)}}(\bT\bW_{\bf S})}|\tG^{(\T\bW_{\bf S} j)}_{ix}
 (\tG^{(j)}|_{\bT \bW_{\bf S}})^{-1}_{xy}
 \tG^{(\T\bW_{\bf S} j)}_{yi}|\\
&\lesssim\sum_{x,y\in \del_{\cG^{(j)}}(\bT\bW_{\bf S})}\sum_{(i,j)\in \vec{E}}\chi_j(\cG^{(\bT\bW_{\bf S})})|\tG^{(\T\bW_{\bf S} j)}_{ix}
 \tG^{(\T\bW_{\bf S} j)}_{yi}|\\
  &\lesssim\sum_{x,y\in \del_{\cG^{(j)}}(\bT\bW_{\bf S})}\Big (\sum_{(i,j)\in \vec{E}}\chi_j(\cG^{(\bT\bW_{\bf S})})|\tG^{(\T\bW_{\bf S} j)}_{ix}|^2\sum_{(i,j)\in \vec{E}}\chi_j(\cG^{(\bT\bW_{\bf S})})|\tG^{(\T\bW_{\bf S} j)}_{y i}|^2\Big)^{1/2}.
 \end{split}
\end{align}%
For the inner sum in the above expression, by   the resolvent identity 
\eqref{e:Schurixj}, 
$\tG^{(\T\bW_{\bf S} j)}_{ix} = \tG^{(\T\bW_{\bf S})}_{ix}-\tG^{(\T\bW_{\bf S})}_{ij}\tG^{(\T\bW_{\bf S})}_{jx}/\tG^{(\T\bW_{\bf S})}_{jj}$. 
By the definition of $\Omega_o^+(z)$ and \eqref{boundhGT}, we have the estimate $|\tG^{(\T\bW_{\bf S})}_{ij}|\lesssim 1$. 
Together with the fact that $\chi_j(\cG^{\bT\bW_{\bfS}}) |\tG^{(\T\bW_{\bf S})}_{jj}|\asymp \chi_j(\cG^{\bT\bW_{\bfS}}) $ 
and the Ward identity \eqref{e:Ward}, we can bound this sum by 
\begin{align*}
&\phantom{{}={}}\sum_{(i,j)\in \vec{E}}\chi_j(\cG^{(\bT\bW_{\bf S})})|\tG^{(\T\bW_{\bf S} j)}_{ix}|^2
=\sum_{(i,j)\in \vec{E}}\chi_j(\cG^{(\bT\bW_{\bf S})})\left|\tG^{(\T\bW_{\bf S})}_{ix}-\tG^{(\T\bW_{\bf S})}_{ij}\tG^{(\T\bW_{\bf S})}_{jx}/\tG^{(\T\bW_{\bf S})}_{jj}\right|^2\\
&\lesssim \sum_{(i,j)\in \vec{E}}|\tG^{(\T\bW_{\bf S})}_{ix}|^2+|\tG^{(\T\bW_{\bf S})}_{jx}|^2\lesssim\Im[\tG_{xx}^{(\bT\bW_{\bfS})}]/\eta\lesssim { (\Im[\md(z)]+\varepsilon'+\varepsilon/\sqrt{\kappa+\eta})}/\eta,
\end{align*}%
where the last inequality follows from the same argument used for \eqref{e:Imbound}.

The same argument also gives the relation $\sum_{(i,j)\in \vec{E}}\chi_j(\cG^{(\bT\bW_{\bf S})})|\tG^{(\T\bW_{\bf S} j)}_{y i}|^2\lesssim { (\Im[\md(z)]+\varepsilon'+\varepsilon/\sqrt{\kappa+\eta+\varepsilon})}/\eta$. Therefore, we have
\begin{align*}
\sum_{(i,j)\in \vec{E}}\chi_j(\cG^{(\bT\bW_{\bf S})})|\tG^{(\T\bW_{\bf S} j)}_{ii}-\tG^{(j)}_{ii}|\lesssim \frac{(d-1)^{2\ell}{ (\Im[\md(z)]+\varepsilon'+\varepsilon/\sqrt{\kappa+\eta+\varepsilon})}}{N\eta}, 
\end{align*}
where  the factor $(d-1)^{2\ell}$  comes from the summation of $x, y$ in \eqref{e:refk}. 
Thus 
\begin{align}\begin{split}\label{tGIGchange}
&\phantom{{}={}}\Big |Q(\tcG)-\frac{1}{Nd}\sum_{(i,j)\in \vec E(\cG^{(\bT\bW_{\bf S})})}\chi_j( \cG^{(\bT\bW_{\bf S})})G_{ii}^{(\bT\bW_{\bf S} j)}\Big |\\
&\leq \frac{1}{Nd}\sum_{(i,j)\not\in \vec E\text{ or }\atop \chi_j(\cG^{(\bT\bW_{\bf S})})=0}|\tG_{ii}^{(j)}+\tG_{jj}^{(i)}|
+\frac{1}{Nd} \sum_{(i,j)\in \vec{E}}\chi_j(\cG^{(\bT\bW_{\bf S})})|\tG^{(\T\bW_{\bf S} j)}_{ii}-\tG^{(j)}_{ii}|\\
&\lesssim   \frac{1}{N^{1-\fc}}+\frac{(d-1)^{2\ell}{ (\Im[\md(z)]+\varepsilon'+\varepsilon/\sqrt{\kappa+\eta+\varepsilon})}}{N\eta}.
\end{split}\end{align} 
where 
we have used  the fact that, by definition of $\Omega_o^+$,    $|\tG_{ii}^{(j)}|=|\tG_{ii}+\tG_{ij}\tG_{ji}/\tG_{jj}|\lesssim 1$, and the number of such terms is bounded by $\OO(N^\fc)$ to bound the first term on the right side. This proves the second bound in  \eqref{GTtotGTconcentration}.
\end{proof}

\subsection{Proof of Proposition \ref{p:Omega+}}
\label{sec:provekey1}
\begin{proof}[Proof of Proposition \ref{p:Omega+}]
The event $F_1(\cG)$ and $F_2(\cG)$ are constructed in Proposition \ref{p:newG} and in Definition \ref{def:F2} respectively, and they satisfy $\bP(F_1(\cG)\cap F_2(\cG))=1-\OO(N^{-\fd})$.  Proposition \ref{p:newG} states that for $\bfS \in F_1(\cG)$, $\tcG=T_\bfS(\cG)\in \bar \Omega^+$. 
The definition of spectral regular graphs $\Omega(z) \subset \bar\Omega$ as in Definition \ref{def:Omega} and Proposition \ref{prop:stabilityGT} implies that $\Omega(z)\subset\Omega_o^+(z)$.
Proposition \ref{l:IGchange} implies that for $\bfS \in F_2(\cG)$, $\tcG=T_\bfS(\cG)$ satisfies
\begin{align}\begin{split}\label{e:tta}
&\phantom{{}+{}}|\I(\tcG,z)-\msc(z)|\leq 
|\I(\cG,z)-\msc(z)|
+|\I(\cG,z)-\I(\tcG,z)|\\
&\lesssim  \frac{\varepsilon}{\sqrt{\kappa+\eta+\varepsilon}}+\frac{(d-1)^{2\ell}{(\Im[\md(z)]+\varepsilon'+\varepsilon/\sqrt{\kappa+\eta+\varepsilon})}}{N\eta(z)}+\frac{1}{N^{1-\fc}}\lesssim  \frac{\varepsilon}{\sqrt{\kappa+\eta+\varepsilon}}.
\end{split}\end{align}%
Thus
\eqref{e:weakQm} holds for the switched graph $\tcG$.
Thanks to Remark \ref{r:subgraph} and Proposition \ref{p:localization}, we can bound 
\begin{align*}\begin{split}
&\phantom{{}={}}|G_{ij}(\Ext(\cB_r(\{\bT, i,j\},\tcG), Q(\cG, z)),z)-G_{ij}(\Ext(\cB_r(\{i,j\},\tcG), Q(\tcG, z)),z)|\\
&\lesssim |G_{ij}(\Ext(\cB_r(\{\bT, i,j\},\tcG), Q(\cG, z)),z)-G_{ij}(\Ext(\cB_r(\{\bT, i,j\},\tcG),  Q(\tcG, z)),z)|\\
&+|G_{ij}(\Ext(\cB_r(\{\bT, i,j\},\tcG), Q(\tcG, z)),z)-G_{ij}(\Ext(\cB_r(\{i,j\},\tcG), Q(\tcG, z)),z)|\lesssim |Q(\cG,z)-Q(\tcG,z)|\\
&+\frac{1}{(d-1)^{r}}+\log N (\sqrt{\kappa+\eta}|Q(\tcG,z)-m_{sc}(z)|+|Q(\tcG,z)-m_{sc}(z)|^2)\lesssim \frac{1}{\log N}.
\end{split}\end{align*}%
Thus, Proposition \ref{prop:tGweakstab} then implies  
\begin{align*}
&\phantom{{}={}}|\tilde G_{ij}-G_{ij}(\Ext(\cB_r(\{i,j\},\tcG), Q(\tcG, z)),z)|
=|\tilde G_{ij}-G_{ij}(\Ext(\cB_r(\{\bT, i,j\},\tcG), Q(\cG, z)),z)|\\
&+|G_{ij}(\Ext(\cB_r(\{\bT, i,j\},\tcG), Q(\cG, z)),z)-G_{ij}(\Ext(\cB_r(\{i,j\},\tcG), Q(\tcG, z)),z)|
\lesssim 1/\log N,
\end{align*}%
which is the weak estimate
\eqref{e:weakrigid1} for the switched graph $\tcG$.
Finally, \eqref{e:stabilitytGT}  in Proposition \ref{prop:stabilitytGT} implies
\eqref{e:weakrigid2} for the switched graph $\tcG$. This finishes the proof of Proposition \ref{p:Omega+}.

\end{proof}

\section{Proof of Proposition \ref{p:Omega-}: Improved Green's Function Estimates}
\label{s:improveG}

In this Section, we prove Proposition \ref{p:Omega-} by a  concentration from local resampling: if the Green's functions of $\cG$ or $\cGT$ satisfy certain estimates, the Green's functions of the switched graphs, i.e. $\tcG$, satisfy improved estimates. 
In Section \ref{sec:concentration},  i.e. Proposition \ref{tGconcentration},  we will show that   weighted  averages 
of the Green's function of $\tcGT$ over the vertex boundary  of $\T$
concentrate under the local resampling. This concentration will be used in 
\eqref{e:ccopy1} to improve estimate for the difference $\tG_{oi}-P_{oi}$. This provides the key improvements on estimates of Green's functions    to implement continuity argument. 
These improvements were stated as  Proposition \ref{p:Omega-} and will be  proved  in Section \ref{sec:improved}.

\subsection{Concentration in the switched graph}
\label{sec:concentration}

We recall the notations and definitions  from  Section \ref{sec:switch} and Section \ref{sec:weak}. Our local resampling involves a fixed center vertex $o$,
and a radius $\ell$.
Given a $d$-regular graph $\cG$, we abbreviate $\cT=\cB_{\ell}(o,\cG)$ (which may not be a tree) and its vertex set $\bT$.
We enumerate edge boundary $\del_E \cT$ as $ \del_E \cT = \{e_1,e_2,\dots, e_\mu\}$, where $e_\al=\{l_\al, a_\al\}$ with $l_\al \in \T$ and $a_\al \in \qq{N} \setminus \T$. 
We recall that $\bI=\{\ta_1,\ta_2,\dots,\ta_\mu\}$ denotes the  set of outer boundary vertices in the switched graph $\tcG$, where $\ta_\al=a_\al$ if $\al\in \qq{1,\mu}\setminus \As_{\bf S}$ and $\ta_\al=c_\al$ if $\al\in \As_{\bf S}$. We also recall the definition of $Q(\cG,z)$ from \eqref{e:IG} and the function $Y_r(Q(\cG,z), z)$ from \eqref{def:Y}.

The result of this section is the following proposition, which asserts  that the average
of the Green's function of $\tcGT$ over the vertex boundary  of $\T$
concentrates under the local resampling.
This is where the condition that the edge boundary contains  more than $ \log N$ edges is
crucial.

\begin{proposition}\label{tGconcentration}
 Let $z\in\bC^+$ with $\Im[z]\geq (\log N)^\fb/N$ and $\cG \in \Omega_o^+(z)$ as in Definition
\ref{def:Omegao+}. Let 
\begin{align*}
P(z)=G(\Ext(\cB_r(\bT, \cG), Q(\cG, z)),z).
\end{align*}
For weights $w_1,w_2,\cdots, w_\mu$  given by either one of the two choices:
\begin{enumerate}
\item Fix any vertex $k\in \bT$. The weights $w_\al=P_{ol_\al}(z)P_{kl_\al }(z)$ for $\al\in \qq{1,\mu}$;
\item If $\cT$ is a truncated $d$-regular tree, fix any vertex $k\sim o$. The weights $w_\al=\bm1(\al\in \sA)/|\sA|$ where $\sA=\{\al\in \qq{1, \mu}$: path from $o$ to $\ell_\al$ does not pass through $k\}$. 
\end{enumerate}
Note that in both cases $\sum_{\al}|w_\al|\lesssim1$ (this follows from \eqref{e:sPPcc} for case (i)).

Then there exists an event $F_3(\cG) \subset F_2(\cG)$\index{$F_3(\cG)$} (as in Definition  \ref{def:F2}) with probability $\P_{\cG}(F_3(\cG))= 1-\OO(N^{-\fd})$
such that for any ${\bf S}\in F_3(\cG)$ with $\tcG=T_{\bf S}(\cG)$, the following holds
\begin{align}\begin{split}\label{eqn:tGconcentration}
&\phantom{{}={}}\left|\sum_{\al\in \qq{1,\mu}}w_\al(\tG^{(\T)}_{\ta_\al\ta_\al}(z)-G_{\ta_\al\ta_\al}(\Ext(\cB_r(\ta_\al,\tcG^{(\T )}), Q(\cG,z)),z))-{\sum_{\al\in\qq{1,\mu}}w_\al(Q(\cG,z)-Y_r(Q(\cG,z),z))}\right|\\
&\leq (\log N)\varepsilon'(z)\sqrt{\sum_\al |w_\al|^2}+\frac{(d-1)^{2\ell}{ (\Im[\md(z)]+\varepsilon'(z)+\varepsilon(z)/\sqrt{\kappa(z)+\eta(z)+\varepsilon(z)})}}{N\eta(z)}+\frac{1}{N^{1-\fc}}.
\end{split} \end{align}%
\end{proposition}

To prove Proposition~\ref{tGconcentration}, in Lemma~\ref{lem:econcentration},
we first show a similar statement for the unswitched graph $\cGT$
in which the problem becomes a concentration problem of independent random variables.
Then we prove Proposition~\ref{tGconcentration} by comparision,
using the estimates in Proposition~\ref{prop:stabilitytGT},
and the fact that the change from $\I(\tcG,z)$ to $\I(\cGT,z)$ is small (Proposition~\ref{l:IGchange}).
Proposition~\ref{prop:stabilitytGT}
is applicable since, by the definition of set $\Omega_o^+(z)$ in Definition \ref{def:Omegao+},
any graph $\cG\in \Omega_o^+(z)$ satisfies the assumptions in Proposition~\ref{prop:stabilitytGT}.

\subsubsection{Estimate for the unswitched graph}

The next lemma shows concentration of a certain average of the Green's function
in the unswitched graph.
Before stating it, we need to introduce some notations. Let
\begin{align*}
Z_\al=\chi_{b_\al}(\cGT)\left(G^{(\T b_\al)}_{c_\al c_\al}(z)-G_{c_\al c_\al}(\Ext(\cB_r(c_\al,\cal G^{(\T b_\al)}), Q(\cG,z)),z)\right),\quad \al\in \qq{1,\mu},
\end{align*}%
where the indicator function $\chi_j( \cG^{(\bT)})=1$ if the vertex $j$ is far away from the outer vertex boundary of $\bT$, i.e. $\dist_{\cG^{(\bT)}}(j, \{a_1,a_2,\cdots, a_\mu\})>{ \fR/4}$; otherwise, it is $0$. The number of vertices $j$ such that $\chi_j( \cG^{(\bT)})=0$ is bounded by $\OO((d-1)^{\ell+\fR/4})\ll N^\fc$. The deficit function of $\cG^{(\bT b_\al)}$ is given by our convention \eqref{e:removing}, explicitly given by $d-\deg_{\cG^{(\bT b_\al)}}(\cdot)$.

\begin{lemma} \label{lem:econcentration}
Under the assumptions of Proposition \ref{tGconcentration},
we define the set 
$F_3(\cG) \subset F_2(\cG)$ (as in Definition \ref{def:F2}) such that
\begin{align}\begin{split}
\label{econcentration}
  &\phantom{{}={}}\left| \sum_{\al\in \qq{1,\mu}}w_\al Z_\al-\sum_{\al\in\qq{1,\mu}}w_\al\left(Q(\cG,z)-Y_r(Q(\cG,z),z)\right)\right|\\
  &\lesssim (\log N)\varepsilon'(z)\sqrt{\sum_\al |w_\al|^2}+\frac{(d-1)^{2\ell}{ (\Im[\md(z)]+\varepsilon'(z)+\varepsilon(z)/\sqrt{\kappa(z)+\eta(z)+\varepsilon(z)})}}{N\eta(z)}+\frac{1}{N^{1-\fc}}.
 \end{split}\end{align}%
Then $\P_{\cG}(F_3(\cG))= 1-\OO(N^{-\fd})$.
\end{lemma}

\begin{proof}

Conditioned on the graph $\cGT$, the random sets $\vec S_1, \vec S_2,\dots, \vec S_\mu$ are independent
and identically distributed, and thus $Z_1,Z_2,\dots, Z_\mu$ are i.i.d random variables.
By the assumption that $\cG\in \Omega_o^+(z)$ as in Definition \ref{def:Omegao+}, using \eqref{e:Schurixj},
for any $\al\in \qq{1,\mu}$ , we have
\begin{align*}
|Z_\al|\lesssim \varepsilon'.
\end{align*}
The Azuma's inequality for independent random variables implies that
\begin{align}\label{azuma}
 \P_{\cG}\left(\left| \sum_{\al\in \qq{1,\mu}}w_\al (Z_\al-\E[Z_\al])\right|\geq t\varepsilon'\sqrt{\sum_\al |w_\al|^2}\right)\leq e^{-t^2/C}.
\end{align}
In the following, we still need to estimate $\E[Z_\al]$. Let $\vec E$ be the set of oriented edges of $\cGT$.
By definition, $\T$ is the radius $\ell$ neighborhood of the vertex $o$, and by the trivial bound it intersects at most
$\OO((d-1)^{\ell})$ edges.
Thus $Nd-\OO((d-1)^{\ell})\leq |\vec E|\leq Nd$.
If $\chi_j(\cG^{(\bT)})=1$,  then $j$ is far away from the outer vertex boundary of $\bT$, and $\cB_r(i,\cG^{(\bT)})$ has excess at most $\omega_d$ and zero deficit function. Thus we have $|G_{jj}(\Ext(\cB_r(i,\cGT),Q(\cG,z)))|\asymp 1$ from \eqref{e:boundPiicopy}.  By the property  \eqref{e:weakrigid2} of the set $\Omega_o^+$ and \eqref{e:Schurixj}, 
it follows that
\begin{align}\begin{split}\label{e:remove1}
&|G_{ii}^{(\T j)}-G_{ii}(\Ext(\cB_{r}(i,\cG^{(\T j)}),Q(\cG,z)))|\lesssim|G_{ii}^{(\T)}-G_{ii}(\Ext(\cB_{r}(i,\cG^{(\T)}),Q(\cG,z)))|\\
&+\left|\frac{G^{(\bT)}_{ij}G^{(\bT)}_{ij}}{G^{(\bT)}_{jj}}-\frac{G_{ij}(\Ext(\cB_{r}(i,\cG^{(\T)}),Q(\cG,z)))G_{ij}(\Ext(\cB_{r}(i,\cG^{(\T)}),Q(\cG,z)))}{G_{jj}(\Ext(\cB_{r}(i,\cG^{(\T)}),Q(\cG,z)))}\right|
\lesssim \varepsilon'.
\end{split}\end{align}%
It follows that
\begin{align}\begin{split}
\label{avIG}
\E[Z_\al]
&=\frac{1}{|\vec E|} \sum_{(i,j)\in\vec E}\chi_j(\cG^{(\bT)})\left(G_{ii}^{(\T j)}-G_{ii}(\Ext(\cB_{r}(i,\cG^{(\T j)}),Q(\cG,z)))\right)\\
&=\frac{1}{Nd}  \sum_{(i,j)\in\vec E}\chi_j(\cG^{(\bT)})\left(G_{ii}^{(\T j)}-G_{ii}(\Ext(\cB_{r}(i,\cG^{(\T j)}),Q(\cG,z)))\right)+\OO\left(\frac{(d-1)^{\ell}\varepsilon'}{Nd}\right), 
\end{split}
\end{align}%
where the error term comes from the change of normalization factors and we further bound these error terms by \eqref{e:remove1}. 
Similar to the proof of \eqref{GTtotGTconcentration}, we have
\begin{align}\label{e:avQQ}
\left|Q(\cG,z)-\frac{1}{Nd}  \sum_{(i,j)\in\vec E}\chi_j(\cG^{(\bT)})G_{ii}^{(\T j)}\right|\lesssim \frac{(d-1)^{2\ell}{ (\Im[\md(z)]+\varepsilon'+\varepsilon/\sqrt{\kappa+\eta+\varepsilon})}}{N\eta}+\frac{1}{N^{1-\fc}}.
\end{align}%
For those vertices $i$, whose radius-$\fR$ neighborhood in $\cGT$ is a $d$-regular tree truncated at level $\fR$, 
we have the equality $G_{ii}(\Ext(\cB_{r}(i,\cG^{(\T j)}), Q(\cG,z)),z)=Y_r(Q(\cG,z),z)$. Therefore
\begin{equation} \label{avemsc}
  \frac{1}{Nd} \sum_{(i,j)\in\vec E}\chi_j(\cG^{(\bT)})G_{ii}(\Ext(\cB_{r}(i,\cG^{(\T j)}), Q(\cG,z)))=Y_r(Q(\cG,z),z)+\OO\left(\frac{1}{N^{1-\fc}}\right)
  .
\end{equation}
Combining \eqref{avIG}, \eqref{e:avQQ}, \eqref{avemsc}, we get
\begin{align}\begin{split}\label{e:key1}
&\sum_{\al\in \qq{1,\mu}}w_\al\bE[Z_\al]=
\sum_{\al\in \qq{1,\mu}}w_\al(Q(\cG,z)-Y_r(Q(\cG,z),z))\\
&+\OO\left(\frac{(d-1)^{2\ell}{ (\Im[\md(z)]+\varepsilon'+\varepsilon/\sqrt{\kappa+\eta+\varepsilon})}}{N\eta}+\frac{1}{N^{1-\fc}}\right).
\end{split}\end{align}
By taking $t=\log N$ in \eqref{azuma}, using \eqref{e:key1} and taking a union bound over all possible choices of $\{w_\al\}_{\al\in\qq{1,\mu}}$, 
we concludes that \eqref{econcentration} holds with overwhelming probability,
\begin{align*}
\P_{\cG}(F_3(\cG))\geq \P_{\cG}(F_2(\cG))-(d+|\bT|)e^{-(\log N)^{2}}
=1-\OO(N^{-\fd}),
\end{align*}
where we have used \eqref{FGmeasure} to bound $\P_{\cG}(F_2(\cG))$. 
This completes the proof.
\end{proof}

\subsubsection{Adding of switched vertices}

We recall the index set $\sJ\subset\As_{\bf S}$ from Definition  \ref{def:sfJ}.
In this subsection, we show the following lemma, which states that $G_{c_\al c_\al}^{(\bT b_\al)}$ is close to $\tG_{c_\al c_\al}^{(\bT)}$.  We will later use 
$\tG_{c_\al c_\al}^{(\bT)}$  to replace $G_{c_\al c_\al}^{(\bT b_\al)}$ in Lemma \ref{lem:econcentration}, and complete the proof of Proposition \ref{tGconcentration}.

\begin{proposition} \label{lem:Gsw}
Under the assumptions of Proposition \ref{tGconcentration},
 for any $\al\in \sJ$, we have
\begin{align}\label{e:diff1360}
|\tGT_{c_\al c_\al}(z)-G^{(\T b_\al)}_{c_\al c_\al}(z)|
\lesssim (d-1)^\ell (\varphi(z))^2,
\end{align}
where $\varphi(z)$ is defined in \eqref{e:defphi}.
\end{proposition}

\begin{proof}[Proof of Proposition \ref{lem:Gsw}]
The claim \eqref{e:diff1360} follows from combining the next two estimates:
\begin{equation}\label{diff1361}
|G^{(\T \bW_{\bf S})}_{c_\al c_\al}-G^{(\T b_\al)}_{c_\al c_\al}|
\lesssim (d-1)^\ell \varphi^2,
\end{equation}
and
\begin{equation}\label{diff1362}
|\tGT_{c_\al c_\al}-\tG^{(\T \bW_{\bf S})}_{c_\al c_\al}|
\lesssim (d-1)^\ell \varphi^2.
\end{equation}

For \eqref{diff1361}, using the definition of the set $\sJ$, the same arguments as for \eqref{e:diffcellest}, for $\al\in \sJ$ it holds that $|G_{c_\al x}^{(\bT b_\al)}|, |G_{y c_\al}^{(\bT b_\al)}|\lesssim \varphi$, for any $x,y\in \bW_{\bf S}\setminus\{b_\al\}$. Moreover, by similar arguments for \eqref{e:Ginv} we have  the estimate $|(G^{(\bT b_\al)}|_{\bW_{\bf S}\setminus \{b_\al\}})^{-1}_{xy}|\lesssim \bm1_{x=y}+\varepsilon' \bm1_{x\simeq y}+\varphi$ for $x,y\in \bW_{\bf S}\setminus\{b_\al\}$. Combining them, and using the Schur complement formula \eqref{e:Schur1}, we have
\begin{align*}\begin{split}
&\phantom{{}={}}|G^{(\bT b_\al)}_{c_\al c_\al}-G^{(\bT\bW_{\bf S})}_{c_\al c_\al}|\leq \left|\sum_{x,y\in \bW_{\bf S}\setminus \{b_\al\}}G_{c_\al x}^{(\bT b_\al)}(G^{(\bT b_\al)}|_{\bW_{\bf S}\setminus \{b_\al\}})^{-1}_{xy}G^{(\bT b_\al)}_{yc_\al}\right|\\
&\lesssim  \varphi^2\sum_{x,y\in \bW_{\bf S}\setminus\{b_\al\}}|(G^{(\bT b_\al)}|_{\bW_{\bf S}\setminus \{b_\al\}})^{-1}_{xy}|\lesssim \left((d-1)^\ell+\varepsilon'(d-1)^{2\ell}\right)\varphi^2 \lesssim (d-1)^\ell \varphi^2,
\end{split}\end{align*}%
where we used that $(d-1)^\ell\varepsilon'\ll 1$.
This completes  the proof of \eqref{diff1361}.

For \eqref{diff1362}, 
we denote $\del_{\cG^{(\bT b_\al)}}(\bW_{\bf S}\setminus\{b_\al\})$ the set of outer boundary vertices of $\bW_{\bf S}\setminus\{b_\al\}$ in $\cG^{(\bT b_\al)}$. Since $\al\in \sJ$,  
\eqref{e:diffcellest} implies the estimates $|(\tG^{(\bT \bW_{\bf S})}B)_{c_\al x}|, |(B'\tG^{(\bT \bW_{\bf S})})_{y c_\al}|\lesssim \varphi$, for any $x,y\in \del_{\cG^{(\bT b_\al)}}(\bW_{\bf S}\setminus\{b_\al\})$. Moreover, for $\cG\in \Omega_o^{+}$ and $\bfS\in F_2(\cG)$, \eqref{e:stabilitytGT} with \eqref{e:Schurixj} implies that $|\tG_{xy}^{(\bT b_\al)}|=|\tG_{xy}^{(\bT)}-\tG_{xb_\al}^{(\bT)}\tG_{b_\al y}^{(\bT)}/\tG_{b_\al b_\al}^{(\bT)}|\lesssim \bm1_{x=y}+\varepsilon' $ for $x,y\in \bW_{\bf S}\setminus\{b_\al\}$. Combining them, and using the Schur complement formula \eqref{e:Schur1}, we have
\begin{align*}\begin{split}
&\phantom{{}={}}|\tG^{(\bT b_\al)}_{c_\al c_\al}-\tG^{(\bT\bW_{\bf S})}_{c_\al c_\al}|
\lesssim  \varphi^2\sum_{x,y\in \bW_{\bf S}\setminus\{b_\al\}}|\tG_{xy}^{(\bT b_\al)}|\lesssim \left((d-1)^\ell+\varepsilon'(d-1)^{2\ell}\right)\varphi^2 \lesssim (d-1)^\ell \varphi^2,
\end{split}\end{align*}%
where we used that $(d-1)^\ell\varepsilon'\ll 1$ from our choice of parameters \eqref{e:relation2}. This completes  the proof of \eqref{diff1362}.

\end{proof}

\subsubsection{Proof of Proposition~\ref{tGconcentration}}

Finally, using the previous lemmas and propositions, we can finish the proof Proposition~\ref{tGconcentration}.

\begin{proof}[Proof of Proposition~\ref{tGconcentration}]
For $\al \in \As_{\bfS}$, the radius-$\fR/4$ neighborhood of $ c_\al$ is a truncated $d$-regular tree with root degree $d-1$ in
both of the graphs $\cG^{(\T b_\al)}$,  $\tcG^{(\T)}$. Therefore, we have $\chi_{b_\al}( \cGT)=1$ and
\begin{align}\begin{split}\label{e:pcpc}
&G_{c_\al c_\al}(\Ext(\cB_r(c_\al,\cG^{(\T b_\al)}),Q(\cG,z)),z)=G_{c_\al c_\al}(\Ext(\cB_r(c_\al,\tcG^{(\T )}),\\
& Q(\cG,z)),z)=Y_r(Q(\cG,z),z).
\end{split}\end{align}
On the other hand, for any index $\al\in \qq{1,\mu}$, by  definition of  $\Omega_o^+(z)$  in Definition \ref{def:Omegao+}, and  using \eqref{e:Schurixj}
as in the proof of  \eqref{e:remove1}, we get
\begin{align}\begin{split}\label{eqn:somebounds1}
&\chi_{b_\al}( \cG^{(\bT)})\left|G^{(\T b_\al)}_{c_\al c_\al}-G_{c_\al c_\al}(\Ext(\cB_r(c_\al,\cG^{(\T b_\al)}),  Q(\cG,z)))\right|\lesssim \varepsilon'.
\end{split}\end{align}
It follows from combining \eqref{e:pcpc} and \eqref{eqn:somebounds1}, and noticing $\mu-|\As_\bfS|= \OO_\fd(1)$, we have
\begin{align}\label{hchange}
\begin{split}
&\phantom{{}={}}\sum_{\al\in \qq{1,\mu}}w_\al Z_\al=\sum_{\al\in \qq{1,\mu}}w_\al\chi_{b_\al}( \cGT)(G^{(\T b_\al)}_{c_\al c_\al}-G_{c_\al c_\al}(\Ext(\cB_r(c_\al,\cG^{(\T b_\al)}), Q(\cG,z))))\\
&=\sum_{\al\in \As_\bfS}w_\al(G^{(\T b_\al)}_{c_\al c_\al}-Y_r(Q(\cG,z),z))+\OO\left( \varepsilon'\max_\al|w_\al|\right).
\end{split}
\end{align}%
Using that for $\cG\in \Omega_o^+(z)$, the assumptions of Proposition \ref{prop:stabilitytGT} holds, 
we have that 
\eqref{e:stabilitytGT} implies for any index $\al\in \qq{1,\mu}$, 
\begin{align}\begin{split}\label{eqn:somebounds2}
\left|\tG^{(\T )}_{\ta_\al\ta_\al}-G_{\ta_\al \ta_\al}(\Ext(\cB_r(\ta_\al,\tcG^{(\T)}), Q(\tcG,z)))\right| \lesssim \varepsilon'.
\end{split}\end{align}
Combining \eqref{e:pcpc},  \eqref{eqn:somebounds2} and Proposition \ref{lem:Gsw} that  
\begin{align}\label{tGTchange}
\begin{split}
&\phantom{{}={}}\sum_{\al \in \qq{1,\mu}}w_\al(\tG^{(\T )}_{\ta_\al\ta_\al}-G_{\ta_\al\ta_\al}(\Ext(\cB_r(\ta_\al,\tcG^{(\T)}), Q(\cG,z))))\\
&=\sum_{\al\in \As_{\bfS}}w_\al(\tG^{(\T)}_{c_\al c_\al}-Y_r(Q(\cG,z),z))+\OO\left( \varepsilon'\max_\al|w_\al|\right)\\
&=\sum_{\al\in \As_\bfS}w_\al(G^{(\T b_\al)}_{c_\al c_\al}-Y_r(Q(\cG,z),z))+\OO\left( \varepsilon'\max_\al|w_\al|+(d-1)^\ell\varphi^2\right),
\end{split}
\end{align}
where we used that $\mu-|\As_\bfS|=\OO_\fd(1)$.
The claim \eqref{eqn:tGconcentration} follows from combining Lemma \ref{lem:econcentration}, estimates \eqref{hchange} and \eqref{tGTchange}.

\end{proof}

\subsection{Improved approximation in the switched graph}
\label{sec:improved}
In this section we prove the following proposition, which gives
improved estimates on the Green's function of the switched graph $\tcG$ centered at vertex $o$.
\begin{proposition}\label{improvetG}
Let $d\geq 3$, $\omega_d$ as in Definition \ref{d:defwd} and $z\in \bC^+$ with $\Im[z]\geq (\log N)^\fb/N$. 
Then for any $\cG\in \Omega_o^+(z)$ as in Definition \ref{def:Omegao+} and $\bfS\in F_3(\cG) \subset \sS(\cG)$ 
defined in Proposition \ref{tGconcentration}, $\tcG=T_{\bf S}\cG$ satisfies the following improved estimates centered at vertex $o$. 
For any vertex $i\in \qq{N}$ we define $\tP(z)=G(\Ext(\cB_r(o,i,\tcG), Q(\tcG, z)),z)$ and denote $\dist_{\tcG}(o,i)=\ell_i$, then  
 If vertex $i\in \bT$,  then
 \begin{align}
\begin{split}\label{G11bound}
 |\tG_{oi}(z)-\tP_{oi}(z)|
  &\lesssim
  \left(\frac{\log N }{(d-1)^{\ell/2}}+\frac{(\log N)^2}{(d-1)^{\ell-\ell_i/2}}\right)\varepsilon'(z) \\
&+\log N (\sqrt{\kappa(z)+\eta(z)}|Q(\tcG,z)-m_{sc}(z)|+|Q(\tcG,z)-m_{sc}(z)|^2),\quad \forall i\in \bT.
\end{split}
\\
\begin{split}\label{G11boundout}
  |\tG_{oi}(z)-\tP_{oi}(z)|
  &\lesssim
  \frac{(\log N)^2 \varepsilon'(z)}{(d-1)^{\ell/2}}\\
  &+\log N (\sqrt{\kappa(z)+\eta(z)}|Q(\tcG,z)-m_{sc}(z)|+|Q(\tcG,z)-m_{sc}(z)|^2), \quad \forall i\in \qq{N}\setminus \bT.
\end{split}\end{align}%
Moreover, if the vertex $o$ has radius-$\fR$ tree neighborhood in the graph $\tcG$,
then the following estimates hold:
\begin{align}
\label{G11treebound}
  \frac{1}{d}\sum_{i: o\sim i}|\tG_{oo}^{(i)}(z)-Y_\ell(Q(\tcG, z),z)| \lesssim \frac{(\log N) \varepsilon'(z)}{(d-1)^{\ell/2}}.
\end{align}
where $o \sim i$ means vertices $o$ and $i$ are adjacent. 
\end{proposition}

 We remark that the estimate \eqref{G11bound}  is better than \eqref{e:Gaxbound} in Proposition \ref{prop:tGweakstab}. If we take $i=j=o$ in \eqref{e:Gaxbound}, then $\ell_i=\ell_j=0$ and  the right side of \eqref{e:Gaxbound}  is $\varepsilon'(z)$. 
 On the other hand,  if we take $i=o$ in  \eqref{G11bound}, we have extra an $(d-1)^{ -\ell/2}$ factor in front of $\varepsilon'(z)$
 in  \eqref{G11bound}.

\begin{proof}[Proof of \eqref{G11bound} and \eqref{G11boundout}]  We first prove \eqref{G11boundout}, which is a simple consequence of the bound \eqref{e:GTNbound}.
Let $\cG_o=\cB_r(\{\bT, i\},\tcG)$. We remark that if $i\in \bT$ then $\cG_o=\cB_r(\bT, \tcG)$. The graph $\cG_o^{(\bT)}$ is obtained from $\cG_o$ by removing the vertex set $\bT$, and its deficit function is given by \eqref{e:removing}. We abbreviate
\begin{align*}
P=G(\Ext(\cG_o, Q(\cG, z))), \quad P^{(\bT)}=G(\Ext(\cG_o^{(\bT)}, Q(\cG, z))).
\end{align*}
 For  $i\in\qq{N}\setminus \bT$, the bound \eqref{e:GTNbound} gives 
\begin{align}\label{e:ccopy20}
|\tG_{oi}-P_{oi}|\lesssim\frac{(\log N)^2\varepsilon'}{(d-1)^{\ell/2}}.
\end{align}
Thanks to Remark \ref{r:subgraph}, we can bound 
\begin{align}\begin{split}\label{e:changeT}
&\phantom{{}={}}|P_{oi}-G_{oi}(\Ext(\cB_r(\{\bT, i\},\tcG), Q(\tcG, z)))|
\lesssim |Q(\cG,z)-Q(\tcG,z)|\\
&\lesssim \frac{(d-1)^{2\ell}{(\Im[\md(z)]+\varepsilon'+\varepsilon/\sqrt{\kappa+\eta+\varepsilon})}}{N\eta}+\frac{1}{N^{1-\fc}}, 
\end{split}\end{align}
where we have used  Proposition \ref{l:IGchange} in the last inequality.  We recall that $\tP(z)=G(\Ext(\cB_r(o,i,\tcG), Q(\tcG, z)),z)$, then
Proposition \ref{p:localization} implies that 
\begin{align}\label{replaceEr4}\begin{split}
&\phantom{{}={}}|G_{oi}(\Ext(\cB_r(\{\bT,i\},\tcG), Q(\tcG, z)))-\tilde P_{oi}|\lesssim 1/(d-1)^{r}\\
&+\log N (\sqrt{\kappa+\eta}|Q(\tcG,z)-m_{sc}(z)|+|Q(\tcG,z)-m_{sc}(z)|^2),
\end{split}\end{align}
where the factor $\log N$ comes  from the term $1+\diam(\cG)$ in \eqref{e:compatibility}.
Thus \eqref{e:changeT} and \eqref{replaceEr4} together imply
\begin{align}\begin{split}\label{e:replacePox}
|P_{oi}-\tP_{oi}|
&\lesssim \frac{1}{(d-1)^r}+\frac{(d-1)^{2\ell}(\Im[\md(z)]+\varepsilon'+\varepsilon/\sqrt{\kappa+\eta+\varepsilon})}{N\eta}\\
&+\log N (\sqrt{\kappa+\eta}|Q(\tcG,z)-m_{sc}(z)|+|Q(\tcG,z)-m_{sc}(z)|^2).
\end{split}\end{align}
The claim \eqref{G11boundout} follows from plugging \eqref{e:replacePox} into \eqref{e:ccopy20}.

We now  prove \eqref{G11bound}. 
By \eqref{Tterm},  for $i\in \bT$ we have
\begin{align}\begin{split}
  |\tG_{oi}-P_{oi}|& \lesssim \left|\sum_{\al\in \qq{1,\mu}}P_{o\ell_\al}P_{i\ell_\al}(\tilde G_{\tilde a_\al \tilde a_\al}^{(\bT)}-P^{(\bT)}_{\tilde a_\al \tilde a_\al})
\right| \\
&+\left|\sum_{\al\neq \beta\in \qq{1,\mu}}P_{o\ell_\al}P_{i\ell_\beta}(\tilde G_{\tilde a_\al \tilde a_\beta}^{(\bT)}-P^{(\bT)}_{\tilde a_\al \tilde a_\beta})
\right| +|((\tG-P){B}'(\tGT-P^{(\T)}){B}P)_{oi}|
 \\
 &\lesssim\Big |\sum_{\al\in \qq{1,\mu}}P_{o\ell_\al}P_{i\ell_\al}(\tilde G_{\tilde a_\al \tilde a_\al}^{(\bT)}-P^{(\bT)}_{\tilde a_\al \tilde a_\al})\Big |+\frac{(\log N)^2\varepsilon'}{(d-1)^{(2\ell-\ell_i)/2}}. \label{e:ccopy1}
 \end{split}
 \end{align}%
 For the second term on the righthand side of \eqref{e:ccopy1}, \eqref{e:off-d} implies that
 \begin{align}\label{e:bbt1}
 \left|\sum_{\al\neq \beta\in \qq{1,\mu}}P_{o\ell_\al}P_{i\ell_\beta}(\tilde G_{\tilde a_\al \tilde a_\beta}^{(\bT)}-P^{(\bT)}_{\tilde a_\al \tilde a_\beta})
\right|\lesssim\frac{(\log N)^2\varepsilon'}{(d-1)^{(2\ell-\ell_i)/2}}.
 \end{align}
 For the last term on the righthand side of \eqref{e:ccopy1}, \eqref{firstterm} and \eqref{e:Gammaxb} implies
 \begin{align}\begin{split}\label{e:bbt2}
 |((\tG-P){B}'(\tGT-P^{(\T)}){B}P)_{oi}|
 &\lesssim(\log N) (1+\ell_i)(d-1)^{(\ell-\ell_i)/2}\varepsilon'\Gamma_o\\
 &\lesssim \frac{(1+\ell_i)(\log N)^3}{(d-1)^{\ell_i/2}}(\varepsilon')^2.
\end{split} \end{align}
 By plugging \eqref{e:bbt1} and \eqref{e:bbt2} into \eqref{e:ccopy1}, we get
\begin{align}\begin{split}
  |\tG_{oi}-P_{oi}|& \lesssim\left |\sum_{\al\in \qq{1,\mu}}P_{o\ell_\al}P_{i\ell_\al}(\tilde G_{\tilde a_\al \tilde a_\al}^{(\bT)}-P^{(\bT)}_{\tilde a_\al \tilde a_\al})\right |+\frac{(\log N)^2\varepsilon'}{(d-1)^{(2\ell-\ell_i)/2}}. \label{e:ccopy2}
 \end{split}
 \end{align}
We will use the concentration inequality,  Proposition \ref{tGconcentration}, to bound the first term on the right hand side of \eqref{e:ccopy2}. Using Remark \ref{r:subgraph} and Proposition \ref{p:localization} we  first replace  $P_{\tilde a_\al \tilde a_\al}^{(\bT)}=G_{\tilde a_\al \tilde a_\al}(\Ext(\cG_o^{(\bT)}, Q(\cG, z)))$ in 
the first term on the right hand side of \eqref{e:ccopy2} by  $G_{\ta_\al\ta_\al}(\Ext(\cB_r(\ta_\al,\tcG^{(\T )}),  Q(\cG,z)))$  and an error term:
\begin{align}\label{replaceEr4.5}\begin{split}
&\phantom{{}={}}|P_{\tilde a_\al \tilde a_\al}^{(\bT)}-G_{\ta_\al\ta_\al}(\Ext(\cB_r(\ta_\al,\tcG^{(\T )}),  Q(\cG,z)))|\\
&\lesssim |G_{\tilde a_\al \tilde a_\al}(\Ext(\cG_o^{(\bT)},  Q(\tilde\cG, z)))-G_{\ta_\al\ta_\al}(\Ext(\cB_r(\ta_\al,\tcG^{(\T )}),  Q(\cG,z)))|+|Q(\cG,z)-Q(\tcG,z)|\\
&\lesssim 1/(d-1)^{r}+\log N (\sqrt{\kappa+\eta}|Q(\tcG,z)-m_{sc}(z)|+|Q(\tcG,z)-m_{sc}(z)|^2)+|Q(\cG,z)-Q(\tcG,z)|.
\end{split}\end{align}%
Hence we can  rewrite the first term in \eqref{e:ccopy2} as
\begin{align}\begin{split}\label{e:day1}
&\phantom{{}={}}\sum_{\al\in \qq{1,\mu}}P_{ol_\al}P_{il_\al}(\tilde G_{\tilde a_\al \tilde a_\al}^{(\bT)}-P^{(\bT)}_{\tilde a_\al \tilde a_\al})=\sum_{\al\in \qq{1,\mu}}P_{ol_\al}P_{il_\al}(\tilde G_{\tilde a_\al \tilde a_\al}^{(\bT)}-G_{\ta_\al\ta_\al}(\Ext(\cB_r(\ta_\al,\tcG^{(\T )}), Q(\cG,z)))\\
&+\OO\left(|Q(\cG,z)-Q(\tcG,z)|+\frac{1}{(d-1)^r}+\log N (\sqrt{\kappa+\eta}|Q(\tcG,z)-m_{sc}(z)|+|Q(\tcG,z)-m_{sc}(z)|^2)\right),
\end{split}\end{align}%
where we used that $\sum_\al |P_{o l_\al} P_{i l_\al}|\lesssim 1/(d-1)^{\ell_i/2}\lesssim 1$. We remark that for vertices $i\in \bT$, $\cB_r(\{\bT,i\}, \tcG)=\cB_r(\bT, \tcG)$. Moreover, Remark \ref{r:radiusR} says that the graphs $\cB_r(\bT, \tcG)$ and $\cB_r(\bT, \cG)$ are isomorphic, and thus $P_{ij}=G_{ij}(\Ext(\cB_r(\bT,\cG),Q(\cG,z)))$ for any $i,j\in \bT$.
The first term on the right side of \eqref{e:day1} is of the form given by Proposition \ref{tGconcentration}. We notice that $\sum_\al |P_{o l_\al} P_{i l_\al}|^2\lesssim 1/(d-1)^{\ell}$. (i) in Proposition \ref{tGconcentration} thus  implies
\begin{align}\begin{split}\label{e:day2}
&\phantom{{}={}}\left |\sum_{\al\in \qq{1,\mu}}P_{o l_\al}P_{i l_\al}(\tilde G_{\tilde a_\al \tilde a_\al}^{(\bT)}-G_{\ta_\al\ta_\al}(\Ext(\cB_r(\ta_\al,\tcG^{(\T )}),  Q(\cG,z))))\right |\\
&\lesssim \frac{|Q(\cG,z)-Y_r(Q(\cG,z))|}{(d-1)^{\ell_i/2}}+\frac{1}{(d-1)^r}+\frac{\log N \varepsilon'}{(d-1)^{\ell/2}}+\frac{(d-1)^{2\ell}(\Im[\md(z)]+\varepsilon'+\varepsilon/\sqrt{\kappa+\eta+\varepsilon})}{N\eta}.
\end{split}\end{align}%
Finally, thanks to \eqref{e:IGchange},  Remark \ref{r:subgraph} and Proposition \ref{p:fixpoint},  we have
\begin{align}\begin{split}\label{e:changeQ}
&\phantom{{}={}}|Q(\cG, z)-Y_r(Q(\cG, z))|=|Q(\tcG,z)-Y_r(Q(\tcG,z)|+\OO\left(|Q(\cG,z)-Q(\tcG,z)|\right)\\
&\lesssim \log N (\sqrt{\kappa+\eta}|Q(\tcG,z)-m_{sc}(z)|+|Q(\tcG,z)-m_{sc}(z)|^2)\\
&+\frac{(d-1)^{2\ell}{(\Im[\md(z)]+\varepsilon'+\varepsilon/\sqrt{\kappa+\eta+\varepsilon})}}{N\eta}+\frac{1}{N^{1-\fc}}.
\end{split}\end{align}
The claim \eqref{G11bound} follows from combining  \eqref{e:replacePox},\eqref{e:ccopy2}, \eqref{e:day1}, \eqref{e:day2} and \eqref{e:changeQ}.
\end{proof}

\begin{proof}[Proof of \eqref{G11treebound}]

If the radius-$\fR$ neighborhood of the vertex $o$ is a tree, then $\cT=\cB_\ell(o,\cG)=\cB_\ell(o,\tcG)$ is a truncated $d$-regular tree with depth $\ell$. We take 
\begin{align*}
P=G(\Ext(\cB_\ell(o,\tcG), Q(\cG, z)),z).
\end{align*}
For any vertex $i$ adjacent to $o$, let $P^{(i)}=G(\Ext(\cB_\ell(o,\tcG^{(i)}), Q(\cG, z)),z)$ .
Let $H$ be the normalized adjacency matrix of $\cT$, explicitly $P$ is given by \eqref{e:formula}
\begin{align}\begin{split}\label{e:Pexp}
&P
=\left(H-z-\frac{Q(\cG, z)}{d-1}\sum_{\al\in \qq{1,\mu}} e_{l_al_a}\right)^{-1}, \\
& P^{(i)}
=\left(H^{(i)}-z-\frac{Q(\cG, z)}{d-1}\sum_{\al\in \qq{1,\mu}} e_{l_al_a}\right)^{-1}.
\end{split}\end{align}
 Thanks to Remark \ref{r:subgraph}, $P_{oo}^{(i)}=Y_\ell(Q(\cG, z),z)=Y_\ell(Q(\tcG, z),z)+\OO(|Q(\cG, z)-Q(\tcG, z)|)$. 
The Schur Complement formula \eqref{e:Schur} gives that
\begin{align}
  \label{G1xSchur1tt}
 \tG^{(i)}|_\T&=(H^{(i)}-z-{B}'\tGT{B})^{-1}.
\end{align}
By taking difference of \eqref{G1xSchur1tt} and \eqref{e:Pexp} we have, 
\begin{align}
\begin{split}\label{G1x-P1x}
\tilde{G}_{oo}^{(i)}-P_{oo}^{(i)}
 &= \frac{1}{d-1}\sum_{\al\in\qq{1,\mu}}P^{(i)}_{o l_\al}P^{(i)}_{ol_\al}(\tGT_{\ta_\al\ta_\al}-Q(\cG,z))+ \frac{1}{d-1}\sum_{\al\neq \beta\in \qq{1,\mu}} P^{(i)}_{o l_\al}P^{(i)}_{o l_\beta}\tGT_{\ta_\al\ta_\beta} \\
 & +\frac{1}{d-1}\sum_{\al,\beta\in \qq{1,\mu}} (\tG^{(i)}_{ol_\al}-P^{(i)}_{o l_\al})P^{(i)}_{o l_\beta}(\tGT_{\ta_\al\ta_\beta}-Q(\cG, z)\delta_{\alpha\beta}).
\end{split}
\end{align}%
It follows from \eqref{e:stabilitytGT}, Propositions \ref{p:fixpoint} and \eqref{e:weakQm} from the definition of $\Omega_o^+(z)$ that 
\begin{align}\begin{split}\label{e:GQbb}
\left|\tGT_{\ta_\al\ta_\al}-Q(\cG, z)\right|
&\lesssim
\left|\tGT_{\ta_\al\ta_\al}-Y_r(Q(\cG,z),z)\right|+\left|Y_r(Q(\cG,z),z)-Q(\cG,z)\right|\\
&\lesssim \varepsilon'+\log N(\sqrt{\kappa+\eta}|Q(\cG,z)-m_{sc}|+|Q(\cG,z)-m_{sc}|^2)
\lesssim \varepsilon'.
\end{split}\end{align}
Using \eqref{e:GQbb} as input, by the same argument as for \eqref{e:ccopy2}, it gives us, 
\begin{align}\label{e:Goox}
\begin{split}
\tilde{G}_{oo}^{(i)}-P_{oo}^{(i)}
 &= \frac{1}{d-1}\sum_{\al\in\qq{1,\mu}}P^{(i)}_{o l_\al}P^{(i)}_{ol_\al}(\tGT_{\ta_\al\ta_\al}-Q(\cG,z)) +\OO\left(\frac{(\log N)^2\varepsilon'}{(d-1)^\ell}\right).
 \end{split}
\end{align}
Let $\sA$ be the subset of $\qq{1, \mu}$ given by  $\sA=\{\al\in \qq{1, \mu}$: path from $o$ to $\ell_\al$ does not pass through $i\}$. Then if $\al\in \qq{1,\mu}\setminus \sA$, $P_{o\ell_\al}^{(i)}=0$. For any $\al\in\sA$,  $P^{(i)}_{o l_\al}P^{(i)}_{ol_\al}\asymp|\msc(z)|^{2\ell+2}/(d-1)^\ell$ 
is independent of  $\al$. Moreover, in the graph $\tcGT$, $\tilde a_\al$ has radius $r$ tree neighborhood, then we have $G_{\ta_\al \ta_\al}(\Ext(\cB_r(\ta_\al, \tcGT), Q(\cG,z)),z)=Y_r(Q(\cG,z),z)$. Therefore, (ii) in Proposition \ref{tGconcentration} implies,
\begin{align}\begin{split}\label{e:olx}
\Big |\sum_{\al\in\sA}P^{(i)}_{o l_\al}P^{(i)}_{ol_\al}(\tGT_{\ta_\al\ta_\al}-Q(\cG,z))\Big |
&\asymp\Big | \frac{(\msc(z))^{2\ell+2}}{(d-1)^{\ell}}\sum_{\al\in\sA}(\tGT_{\ta_\al\ta_\al}-Q(\cG,z))\Big |\\
 &\lesssim \frac{\log N\varepsilon'}{\sqrt{|\sA|}}\asymp\frac{\log N\varepsilon'}{(d-1)^{\ell/2}}.
\end{split}\end{align}
Combining \eqref{e:Goox},  \eqref{e:olx}, and recalling $P_{oo}^{(i)}=Y_\ell(Q(\tcG, z))+\OO(|Q(\cG, z)-Q(\tcG, z)|)$, we get
\begin{align*}
\begin{split}
\left|\tilde{G}_{oo}^{(i)}-Y_\ell(Q(\tcG,z))\right|\lesssim \frac{\log N\varepsilon'}{(d-1)^{\ell/2}}.
 \end{split}
\end{align*}
The claim \eqref{G11treebound} follows.
\end{proof}

\begin{proof}[Proof of Proposition \ref{p:Omega-}]
The event $F_3(\cG)\subset F_2(\cG)$  is constructed in Proposition \ref{tGconcentration} and it satisfies $\bP(F_3(\cG))=1-\OO(N^{-\fd})$. As in \eqref{e:tta}, Proposition \ref{l:IGchange} implies that for $\bfS \in F_2(\cG)$, $\tcG=T_\bfS(\cG)$ satisfies \eqref{e:weakone}. The rest of Proposition \ref{p:Omega-} follows from Proposition \ref{improvetG}.
\end{proof}

\section{Eigenvalue Rigidity}

We recall the control parameters $\varepsilon(z), \varepsilon'(z)$ from Section \ref{sec:weak}. In this section we prove the following improved bound for  the Stieltjes transform $m_N(z)$ of the eigenvalues of $\cG$.  The bound of $m_N(z)-\md(z)$ in the following theorem is the square of the error bound in Theorem 
\ref{thm:mrmsc}. The improved bound implies that there are no eigenvalues (except for the trivial one) of $\cG$ outside the interval $[-2-N^{-\Omega(1)}, 2+N^{-\Omega(1)}]$.

\begin{theorem}[Improved Local Law] \label{thm:improvelocal} 
Fix $d\geq 3$, $\fc>0$, and recall the set of radius-$\fR$ tree like graphs $\bar\Omega\subset \GNd$ from Definition \ref{def:barOmega}.  For any large $\fC>0$, small $\fe>0$  and $N$ large enough, with probability $1-\OO(N^{-\fC})$ with respect to the uniform measure on $\bar\Omega$,  the Stieltjes transform of the eigenvalues of $\cG$ satisfies:  for any $z\in \bC^+$ such that $\Im[z]\geq (\log N)^\fb$ and $\varepsilon(z)\leq (\kappa(z)+\eta(z))/\log N$,
\begin{align} \begin{split}\label{e:improvelocal}
&\phantom{{}={}}\left|m_N(z) - \md(z)-\frac{1}{N(d/\sqrt{d-1}-z)}\right|\\
&\lesssim
   \frac{N^{\fe}}{\sqrt{\kappa(z)+\eta(z)}}\left((\varepsilon(z))^2+(\varepsilon(z))^{1/2}\left(\frac{{ \Im[\md(z)]+\varepsilon(z)/\sqrt{\kappa(z)+\eta(z)}}}{N\eta(z)}\right)^{1/2}\right).
\end{split}\end{align}%
\end{theorem}

In Section \ref{sec:fluctuation}, we prove a key proposition, 
which is reminiscent the Fluctuation Averaging Lemma for Wigner matrices \cite{MR2871147,MR3068390}. 
Using Proposition \ref{t:recursion} as input, we give the proof  of Theorem \ref{thm:improvelocal} in Section \ref{sec:proof}. Our main Theorem \ref{thm:evloc} on the rigidity of eigenvalues is a consequence of Theorem \ref{thm:mrmsc} and \ref{thm:improvelocal}.

\subsection{Fluctuation Averaging Proposition}\label{sec:fluctuation}
We recall the set $\Omega(z)$ of  spectral regular graphs from Definition \ref{def:Omega}.  A $d$-regular graph $\cG\in \Omega(z)$ then $\cG\in\bar\Omega$ as in Definition \ref{def:barOmega}, and its Green's function satisfies
 $|Q(\cG,z)-\msc(z)|\leq \varepsilon(z)/\sqrt{\kappa(z)+\eta(z)+\varepsilon(z)}$ and the relation
$\left|G_{ij}(z)-G_{ij}(\Ext(\cB_r(i,j,\cG),Q(\cG,z)),z)\right|
  \leq \varepsilon(z)$ for any two vertices $i,j\in \qq{N}$.

For any finite $d$-regular graph $\cG$, we recall the quantity $Q$ from \eqref{e:IG}
\begin{equation} \label{e:IGcopy}
Q(\cG,z)=\frac{1}{Nd}\sum_{i\in \qq N}\sum_{j: i\sim j} G_{ii}^{(j)}(\cG,z),
\end{equation}
$G^{(j)}(\cG,z)$ is the Green's function of the graph obtained from $\cG$ by removing the vertex $j$. The quantity $Q(\cG,z)$ satisfies the following self-consistent equation:
\begin{align}\label{e:self-con}
Q(\cG,z)\approx Y_\ell(Q(\cG,z),z),
\end{align}
where the function $Y_\ell$ is as defined in \eqref{def:Y}. The equation $w=Y_\ell(w,z)$ has a fix point at $w=m_{sc}(z)$.
In the proof of Proposition \ref{prop:bootstrap}, 
by averaging \eqref{e:improverigid2} over all $o\in \qq{N}$, we get with high probability
\begin{align}\label{e:epss}
Q(\cG,z)-Y_\ell(Q(\cG,z),z)=\OO\left(\frac{\log N \varepsilon'(z)}{(d-1)^{\ell/2}}\right).
\end{align}
In this section, we obtain the following estimate on the high moments of $Q(\cG,z)-Y_\ell(Q(\cG,z),z)-\delta_Q(z)$,  where $\delta_Q(z)$   is a deterministic quantity corresponding to the trivial eigenvalue at $d/\sqrt{d-1}$  defined in \eqref{e:dQ}.    
 This new estimate,   roughly  speaking, improves  the error in \eqref{e:epss} to  its square.

\begin{proposition}\label{t:recursion}
For any  $d\geq 3$ fixed, the Green's function of  a random $d$-regular graph satisfies that,  for any $z\in \bC^+$ with $\Im[z]\geq (\log N)^\fb/N$, and integer $p\geq 1$,
\begin{align}\begin{split}
\label{e:QY}&\phantom{{}={}}\bE[|Q(\cG,z)-Y_\ell(Q(\cG,z),z)-\delta_Q(z)|^{2p}{\bm1(\cG\in \Omega(z))}]^{1/2p}\\
&\lesssim (\varepsilon(z))^{1/2}\left(\frac{(d-1)^{2\ell}{(\Im[\md(z)]+\varepsilon'(z)+\varepsilon(z)/\sqrt{\kappa(z)+\eta(z)+\varepsilon(z)})}}{N\eta(z)}\right)^{1/2}+(d-1)^{2\ell} (\varepsilon'(z))^2,
\end{split}\end{align}%
where 
\begin{align}\label{e:dQ}
\delta_Q(z)\deq ((d-1)^{\ell+1}-1)(\msc(z))^{2\ell+2}\left(1+\frac{\msc(z)}{\sqrt{d-1}}\right)^2\frac{1}{N(d/\sqrt{d-1}-z)},
\end{align} 
is the correction corresponding to the trivial eigenvalue at $d/\sqrt{d-1}$. Moreover,  \eqref{e:QY} holds with the left side replaced by 
 $\bE[|m_N(z)-X_\ell(Q(\cG,z),z)-\delta_{m}(z)|^{2p}{\bm1(\cG\in \Omega(z))}]^{1/2p}$ where $X_\ell$ is defined in \eqref{def:Y} and 
\begin{align*}
\delta_{m}(z)\deq\frac{d}{d-1}(d(d-1)^\ell-1)(\md(z))^2(\msc(z))^{2\ell}\left(1+\frac{\msc(z)}{\sqrt{d-1}}\right)^2\frac{1}{N(d/\sqrt{d-1}-z)}.
\end{align*}
\end{proposition}

The proof of Proposition \ref{t:recursion} relies on the following lemma, which will be used to perform a discrete integration by part, when computing the expectation on the high moments of $Q(\cG,z)-Y_\ell(Q(\cG,z),z)-\delta_Q(z)$.
\begin{lemma}\label{l:exchange}
Fix $d\geq 3$. We recall the operator $T_\bfS$ from \eqref{e:Tdef1}. Let $\cG$ be a random $d$-regular graph  and $\bfS$ uniformly distributed over $\sS(\cG)$, then the graph pair $(\cG, T_{\bf S}(\cG))$ forms an exchangeable pair:
\begin{align*}
(\cG, T_{\bf S}(\cG))\stackrel{law}{=}(T_{\bf S}(\cG), \cG).
\end{align*}
\end{lemma}

\begin{proof}
Lemma 
The statement that $(\cG, T_{\bf S}(\cG))$ and $(T_{\bf S}(\cG), \cG)$ have the same law follows from
\begin{align}\label{e:exchange}
\bP(\cG=\cG_1, T_{\bf S}(\cG)=\cG_2)=\bP(T_{\bf S}(\cG)=\cG_1,\cG=\cG_2),
\end{align}
for any two $d$-regular graphs $\cG_1, \cG_2\in \GNd$. In the following we prove \eqref{e:exchange}. If $\cB_\ell(o, \cG_1)\neq \cB_\ell(o, \cG_2)$, then both sides of \eqref{e:exchange} are zero. There is nothing to prove. Otherwise, $\cB_\ell(o, \cG_1)=\cB_\ell(o, \cG_2)$, which implies that $|\sf S(\cG_1)|=|\sf S(\cG_2)|$. We recall $T(\bfS)$ from \eqref{e:Tdef2}. It is proven in Proposition \ref{e:tildeTT} that the operator $T$ is an involution, i.e. $T_{\bf S}(\cG_1)=\cG_2$ if and only if $T_{T(\bf S)}(\cG_2)=\cG_1$. This implies
\begin{align*}\begin{split}
&\phantom{{}={}}\bP(\cG=\cG_1, T_{\bf S}(\cG)=\cG_2)
=\frac{1}{|\GNd|}\frac{1}{|\sf S(\cG_1)|}\sum_{\bf S\in \sf S(\cG_1)}\bm 1(T_{\bf S}(\cG_1)=\cG_2)\\
&=\frac{1}{|\GNd|}\frac{1}{|\sf S(\cG_2)|}\sum_{\bf S\in \sf S(\cG_1)}\bm 1(T_{T(\bf S)}(\cG_2)=\cG_1)\\
&=\frac{1}{|\GNd|}\frac{1}{|\sf S(\cG_1)|}\sum_{\bf S\in \sf S(\cG_2)}\bm 1(T_{\bf S}(\cG_2)=\cG_1)
=\bP(\cG=\cG_2, T_{\bf S}(\cG)=\cG_1),
\end{split}\end{align*}
which finishes the proof of Lemma \ref{l:exchange}.
\end{proof}
\begin{remark}
Thanks to Proposition \ref{prop:reverse}, the local resampling of the edge boundaries of the radius-$\ell$ neighborhoods of $o$ is a reversible  Markov chain. Lemma \ref{l:exchange} is also an immediate consequence of the reversibility.
\end{remark}

\begin{proof}[Proof of Proposition \ref{t:recursion}]
For simplicity of notations, in this proof we will write $\varepsilon=\varepsilon(z)$, $\varepsilon'=\varepsilon'(z)$, $\Omega=\Omega(z)$, $Q=Q(\cG, z)$, $\tilde Q=Q(\tcG, z)$, $\delta_Q=\delta_Q(z)$, $\delta_m=\delta_m(z)$ and $Y_\ell(Q)=Y_\ell(Q(\cG,z),z)$. We will only prove the bound for $\bE[|Q-Y_\ell(Q)-\delta_Q|^{2p}{\bm1(\cG\in \Omega)}]^{1/2p}$, the bound for $\bE[|m_N(z)-X_\ell(Q)-\delta_m|^{2p}{\bm1(\cG\in \Omega)}]^{1/2p}$ can be proven in the same way. We denote that $Q'=Q-\delta_Q$ and $\tQ'=\tQ-\delta_Q$.
We introduce an indicator function $\chi_o( \cG)=1$ if the vertex $o$ has a radius $\fR$ tree neighborhood in $\cG$.  
\begin{align}\begin{split} \label{e:tt1}
&\phantom{{}={}}\bE[(Q'-Y_\ell(Q))^{p}\overline{(Q'-Y_\ell(Q))^p}{\bm1(\cG\in \Omega)}]\\
&=\frac{1}{Nd}\bE\left[\sum_{o,i}A_{oi}(G_{oo}^{(i)}-Y_\ell(Q)-\delta_Q)(Q'-Y_\ell(Q))^{p-1}\overline{(Q'-Y_\ell(Q))^p}{\bm1(\cG\in \Omega)}\right]\\
&=\frac{1}{Nd}\bE\left[\sum_{o,i}\chi_{o}(\cG)A_{oi}(G_{oo}^{(i)}-Y_\ell(Q)-\delta_Q)(Q'-Y_\ell(Q))^{p-1}\overline{(Q'-Y_\ell(Q))^p}{\bm1(\cG\in \Omega)}\right]\\
&+\OO\left(\frac{1}{N^{1-\fc}}\right)\bE\left[|Q'-Y_\ell(Q)|^{2p-1}{\bm1(\cG\in \Omega)}\right],
\end{split}\end{align}%
where we used that for $\cG\in \Omega$, $A_{oi}|G_{oo}^{(i)}|,|Y_\ell(Q)|, |\delta_Q|\lesssim 1$, and $\chi_o(\cG)=1$ except for $\OO(N^{\fc})$ vertices.

We recall the local resampling from Section \ref{sec:switch}. 
Given any $d$-regular graph $\cG$ and a fixed vertex $o$.  Take the resampling data $\bfS\in \sS(\cG)$, and denote $\tcG=T_{\bfS}(\cG)$. 
Using the definition  $\Delta(f(\cG))\deq f(\cG)-f(\tcG)$, 
we can rewrite the first term on the right side of \eqref{e:tt1} as
\begin{align}
\notag&\phantom{{}={}}\bE\left[\chi_{o}(\cG)A_{oi}(G_{oo}^{(i)}-Y_\ell(Q)-\delta_Q)(Q'-Y_\ell(Q))^{p-1}\overline{(Q'-Y_\ell(Q))^p}{\bm1(\cG\in \Omega)}\right]\\
&=\label{e:diff0}\bE\left[\chi_{o}(\cG)A_{oi}(-Y_\ell(Q)-\delta_Q)(Q'-Y_\ell(Q))^{p-1}\overline{(Q'-Y_\ell(Q))^p}{\bm1(\cG\in \Omega)}\right]\\
\label{e:diffI}&+\bE\left[\chi_o(\cG)A_{oi}\tG_{oo}^{(i)}(Q'-Y_\ell(Q))^{p-1}\overline{(Q'-Y_\ell(Q))^p}{\bm1(\cG\in \Omega)}\right]\\
\label{e:diffII}&+\bE\left[\chi_o(\cG)A_{oi}\Delta(G_{oo}^{(i)})(Q'-Y_\ell(Q))^{p-1}\overline{(Q'-Y_\ell(Q))^p}{\bm1(\cG\in \Omega)})\right].
\end{align}
We claim that  the following estimates for \eqref{e:diffI} and \eqref{e:diffII} hold: 
\begin{align}\begin{split}
&\eqref{e:diffI}
=\bE\left[\chi_o(\cG)A_{oi}(Y_\ell(Q)+\delta_Q)(Q'-Y_\ell(Q))^{p-1}\overline{(Q'-Y_\ell(Q))^p}{\bm1(\cG\in \Omega)}\right]\\
&+\OO((d-1)^{2\ell}(\varepsilon')^2)\frac{d}{N}\bE\left[|Q'-Y_\ell(Q)|^{2p}{\bm1(\cG\in \Omega)}\right]^{1-1/2p}
\label{e:diffIbound}
\end{split}
\\
\begin{split}
&\eqref{e:diffII}
\lesssim\varepsilon' N^{-2+\fc}\bE\left[|Q'-Y_\ell(Q)|^{2p}\bm1(\cG\in \Omega)\right]^{1-1/2p} \label{e:diffIIbound}\\
&+\sum_{q=1}^{2p-1} \frac{\varepsilon}{N}\left(\frac{(d-1)^{2\ell}{ (\Im[\md]+\varepsilon'+\varepsilon/\sqrt{\kappa+\eta+\varepsilon})}}{N\eta}+\frac{1}{N^{1-\fc}}\right)^q\bE\left[|Q'-Y_\ell(Q)|^{2p}{\bm1(\cG\in \Omega)}\right]^{1-(1+q)/2p}.
\end{split}
\end{align}%
We postpone their proofs to the next two sections. In the following we continue to complete the proof of \eqref{e:QY}.  Plugging \eqref{e:diffIbound} and \eqref{e:diffIIbound} into \eqref{e:tt1} and  noticing that the first term in \eqref{e:diffIbound} cancels with \eqref{e:diff0},  we get
\begin{align*}\begin{split}
&\bE[|Q'-Y_\ell(Q)|^{2p}{\bm1(\cG\in \Omega)}]
\lesssim (d-1)^{2\ell}(\varepsilon')^2\bE\left[|Q'-Y_\ell(Q)|^{2p}{\bm1(\cG\in \Omega)}\right]^{1-1/2p}\\
&+\sum_{q=1}^{2p-1} \frac{\varepsilon}{N}\left(\frac{(d-1)^{2\ell}{ (\Im[\md]+\varepsilon'+\varepsilon/\sqrt{\kappa+\eta+\varepsilon})}}{N\eta}+\frac{1}{N^{1-\fc}}\right)^q\bE\left[|Q'-Y_\ell(Q)|^{2p}{\bm1(\cG\in \Omega)}\right]^{1-(1+q)/2p}.
\end{split}\end{align*}%
The equation  \eqref{e:QY} follows from the above estimate by noticing that from our choice of parameters \eqref{e:defeps0} and \eqref{e:relation2} that $\varepsilon\gg (d-1)^{2\ell}(\Im[\md]+\varepsilon'+\varepsilon/\sqrt{\kappa+\eta+\varepsilon})/(N\eta)$
and $\varepsilon^{3/2}\gg N^{-1+\fc}$.
\end{proof}

\subsubsection{Proof of \eqref{e:diffIbound}}
In this section, we show that,  for $\cG\in \Omega$ with $\chi_o(\cG)=1$,
\begin{align}\label{e:top}
|\bE_\cG[A_{oi}(\tG_{oo}^{(i)}-Y_\ell(Q)-\delta_Q)]|\lesssim
  A_{oi}(d-1)^{2\ell}(\varepsilon')^2,
\end{align}
where the expectation is with respect to the randomness of the resampling data $\bfS$  defined in Definition \ref{def:enlarge}. 
If this is the case, then 
\begin{align*}\begin{split}
&\phantom{{}={}}\left|\bE\left[\chi_o(\cG)A_{oi}(\tG_{oo}^{(i)}-Y_\ell(Q)-\delta_Q)(Q'-Y_\ell(Q))^{p-1}\overline{(Q'-Y_\ell(Q))^p}{\bm1(\cG\in \Omega)}\right]\right|\\
&\lesssim\OO((d-1)^{2\ell}(\varepsilon')^2)\bE\left[A_{oi}|Q'-Y_\ell(Q)|^{2p-1}{\bm1(\cG\in \Omega)}\right]\\
& =\OO((d-1)^{2\ell}(\varepsilon')^2)\frac{1}{N}\bE\left[\sum_{i\in\qq{N}}A_{oi}|Q'-Y_\ell(Q)|^{2p-1}{\bm1(\cG\in \Omega)}\right]  \\
&\lesssim\OO((d-1)^{2\ell}(\varepsilon')^2)\frac{d}{N}\bE\left[|Q'-Y_\ell(Q)|^{2p}{\bm1(\cG\in \Omega)}\right]^{1-1/2p},
\end{split}\end{align*}
 where we have used the following facts: in the second line, we used that the index $i$ is  an arbitrary  index  from $\qq{N}$ and   thus 
 in the expectation we are allowed to  take the average  over $i$;  in the last line, we used that $\sum_{i\in \qq{N}}A_{oi}=d$ and   Jensen's inequality. Thus equation  \eqref{e:diffIbound} follows.

In the following we prove \eqref{e:top}.  If $A_{oi}=0$, \eqref{e:top} holds trivially. In the following we assume $A_{oi}=1$, i.e. $i$ is adjacent to $o$.  We recall the set $F_2(\cG)$ from Definition \ref{def:F2}. If $\chi_o(\cG)=1$, Remark \ref{r:radiusR} implies that  vertex $o$ has radius $\fR/4$ tree neighborhood in both $\cG$ and $\tcG$.  Denote by  
$P=G(\Ext(\cB_\ell(o,\cG), Q(\cG, z)))$ and thus 
$P_{oo}^{(i)}=Y_\ell(Q)$.   Using \eqref{G1x-P1x} we have for $\bfS\in F_2(\cG)$ that 
\begin{align}
\begin{split}\label{e:Goo-Poo}
\tilde{G}_{oo}^{(i)}-P_{oo}^{(i)}
 &= \frac{1}{d-1}\sum_{\al\in\qq{1,\mu}}P^{(i)}_{o l_\al}P^{(i)}_{ol_\al}(\tGT_{\ta_\al\ta_\al}-Q)\\
 &+ \frac{1}{d-1}\sum_{\al\neq \beta\in \qq{1,\mu}} P^{(i)}_{o l_\al}P^{(i)}_{o l_\beta}\tGT_{\ta_\al\ta_\beta} +\OO((d-1)^{2\ell}(\varepsilon')^2).
\end{split}
\end{align}
We notice that $|P^{(i)}_{o l_\al}|\lesssim |m_{sc}(z)|^{\ell+1}/(d-1)^{\ell/2}$ from \eqref{e:boundPijcopy} which is independent of $\al$ and the resampling data $\bfS$. As a consequence, 
the expectation of $\tilde{G}_{oo}^{(i)}-P_{oo}^{(i)}$ boils down to computing the expectations of $\tGT_{\ta_\al\ta_\al}-Q$ and  $\tGT_{\ta_\al\ta_\beta}$.
\begin{proposition}\label{c:expbound}
Under the assumptions of Proposition \ref{t:recursion},  $\cG\in \Omega$ and $\chi_o(\cG)=1$, the following estimates hold:
\begin{align}
\label{e:diag}&\left|\bE_\cG[(\tGT_{\ta_\al\ta_\al}-Q)\bm1(\bfS\in F_2(\cG))\bm]\right|\lesssim (d-1)^\ell(\varepsilon')^2,\\
\label{e:offdiag}&\left|\bE_\cG[\tGT_{\ta_\al\ta_\beta}\bm1(\bfS\in F_2(\cG))]-\left(1+\frac{\msc(z)}{\sqrt{d-1}}\right)^2\frac{1}{ N( d/\sqrt{d-1}-z)}\right|\lesssim (d-1)^\ell(\varepsilon')^2.
\end{align}
\end{proposition}

Using Proposition \ref{c:expbound}, we can wrap up the proof of \eqref{e:top}. We notice that for $\bfS\not\in F_2(\cG)$, we have the trivial bound $|\tG_{oo}^{(i)}-Y_\ell(Q)|\lesssim 1/\eta\lesssim N$. And this happens with small probability $\bP_\cG(\bfS\not\in F_2(\cG))\lesssim N^{-\fd}$. 
Moreover, from \eqref{e:Gtreemsccopy2} and \eqref{e:smalldiff}, we have $P^{(i)}_{ol_\al}= \msc(z)(-\msc(z)/\sqrt{d-1})^\ell+\OO(\ell |Q-\msc(z)|/(d-1)^{\ell/2})=\OO(1/(d-1)^{\ell/2})$. By taking expectation on both sides of \eqref{e:Goo-Poo}, and using Proposition \ref{c:expbound},
\begin{align*}\begin{split}
&\phantom{{}={}}\bE_\cG[\tG_{oo}^{(i)}-Y_\ell(Q)]
=\bE_\cG[\tG_{oo}^{(i)}-Y_\ell(Q)\bm1(\bfS\in F_2(\cG))]|+\OO(N^{-\fd+1})\\
 &=(d-1)^\ell((d-1)^{\ell+1}-1) P_{ol_\al}^{(i)}P_{ol_\beta}^{(i)}\bE_\cG[\tGT_{\ta_\al\ta_\beta}\bm1(\bfS\in F_2(\cG))] \\&+ \OO(|\bE_\cG[(\tGT_{\ta_\al\ta_\al}-Q)\bm1(\bfS\in F_2(\cG))]|+(d-1)^{2\ell}(\varepsilon')^2)\\
 &=((d-1)^{\ell+1}-1)(\msc(z))^{2\ell+2}\left(1+\frac{\msc(z)}{\sqrt{d-1}}\right)^2\frac{1}{N(d/\sqrt{d-1}-z)}+\OO(
  (d-1)^{2\ell}(\varepsilon')^2)\\
  &=\delta_Q+\OO(
  (d-1)^{2\ell}(\varepsilon')^2).
\end{split}\end{align*}
This completes the proof   of  \eqref{e:top}.

\begin{proof}[Proof of Proposition \ref{c:expbound}]
We recall the admissible set $\As_\bfS$ from \eqref{Wdef}. Since $\chi_o(\cG)=1$, the index $\al\in \As_\bfS$ if it holds $\dist_{\cGT}(a_\al, \{b_\al, c_\al\})\geq \fR/4$ and  $\dist_{\cG^{(\bT)}}(\{a_\al,b_\al,c_\al\}, \{a_\beta,b_\beta,c_\beta\})\geq \fR/4$ for all $\beta\in \qq{1,\mu}\setminus \{\al\}$. By a union bound, this holds with probability $1-\OO((d-1)^{\ell+\fR/4}/N)=1-\OO(N^{-1+\fc})$. We can rewrite \eqref{e:diag} as
\begin{align}\begin{split}\label{e:rep1}
&\phantom{{}={}}\left|\bE_{\cG}[(\tGT_{\ta_\al\ta_\al}-Q)\bm1(\bfS\in F_2(\cG))]\right|\\
&\leq \left|\bE_\cG[(\bm1(\al\in {\As_\bfS})+\bm1(\al\not\in {\As_\bfS}))(\tGT_{\ta_\al\ta_\al}-Q)\bm1(\bfS\in F_2(\cG))]\right|\\
&=
 \left|\bE_\cG[\bm1(\al\in {\As_\bfS})(\tGT_{c_\al c_\al}-Q)\bm1(\bfS\in F_2(\cG))]\right|+\OO(\varepsilon'/N^{1-\fc}).
\end{split}\end{align}
Next we show that if $\al\in \As_{\bfS}$ and $\bfS\in F_2(\cG)$, we can replace $\tG^{(\bT)}_{c_\al c_\al}$ in the above expression by $G_{c_\al c_\al}^{(b_\al)}$. In fact we show that both of them are close to $G^{(\bT\bW_{\bf S})}_{c_\al c_\al}$
with error $\OO\left((d-1)^{\ell}(\varepsilon')^2\right)$. For $\tG^{(\bT)}_{c_\al c_\al}$, the Schur complement formula \eqref{e:Schur1} implies 
\begin{align}\begin{split}\label{e:rp2}
\tGT_{c_\al c_\al}&=\tG^{(\bT\bW_{\bf S})}_{c_\al c_\al}-\sum_{x,y\in \bW_{\bf S}}|\tG^{(\bT)}_{c_\al x}(\tGT|_{\bW_\bfS})^{-1}_{xy}  \tG^{(\bT)}_{yc_\al}|\\
&=G^{(\bT\bW_{\bf S})}_{c_\al c_\al}+\OO\left((\varepsilon')^2\sum_{x,y\in \bW_{\bf S}}|(\tGT|_{\bW_\bfS})^{-1}_{xy}|\right)=G^{(\bT\bW_{\bf S})}_{c_\al c_\al}+\OO\left((d-1)^{\ell}(\varepsilon')^2\right),
\end{split}\end{align}
where we have used  \eqref{e:stabilitytGT} and  \eqref{e:Ginv} and $|(\tGT|_{\bW_\bfS})^{-1}_{xy}|\lesssim \bm1_{x=y}+\varepsilon'$.
And similarly 
\begin{align}\begin{split}\label{e:rp3}
G^{(b_\al)}_{c_\al c_\al}&=G^{(\bT\bW_{\bf S})}_{c_\al c_\al}-\sum_{x,y\in \bT\bW_{\bf S}\setminus\{b_\al\}}|G^{(b_\al)}_{c_\al x}(G^{(b_\al)}|_{\bT\bW_{\bf S}\setminus\{b_\al\}})^{-1}_{xy} G^{(b_\al)}_{yc_\al}|\\
&=G^{(\bT\bW_{\bf S})}_{c_\al c_\al}+\OO\left(\varepsilon^2\sum_{x,y\in \bT\bW_{\bf S}\setminus\{b_\al\}}|(G^{(b_\al)}|_{\bT\bW_{\bf S}\setminus\{b_\al\}})^{-1}_{xy}|\right)\\
&=G^{(\bT\bW_{\bf S})}_{c_\al c_\al}+\OO\left((d-1)^{\ell}(\varepsilon')^2\right),
\end{split}\end{align}
where, for $\cG\in \Omega$, $|(G^{(b_\al)}|_{\bT\bW_{\bf S}\setminus\{b_\al\}})^{-1}_{xy}|$ can be estimated as in Claims \ref{c:sumPTSbd} and \ref{c:GTinv}.
We can use \eqref{e:rp2} and \eqref{e:rp3} to  simplify \eqref{e:rep1} so that 
\begin{align}\begin{split}\label{e:rp4}
&\phantom{{}={}}\bE_{\cG}[\bm1(\al\in W_\bfS)(\tGT_{c_\al c_\al}-Q)\bm1(\bfS\in F_2(\cG))]\\
&=\bE_{\cG}[\bm1(\al\in W_\bfS)(G^{( b_\al)}_{c_\al c_\al}-Q)\bm1(\bfS\in F_2(\cG))]+\OO\left((d-1)^{\ell}(\varepsilon')^2\right)\\
&=\bE_\cG[G^{(b_{\al})}_{c_\al c_\al}-Q]+\OO\left((d-1)^{\ell}(\varepsilon')^2\right).
\end{split}\end{align}
We recall that $(b_\al, c_\al)$ is  uniformly randomly picked  from the directed edges of $\cGT$. Comparing with the definition \eqref{e:IGcopy} of $Q$ and noticing that  for $\cG\in \Omega$, $|G_{ii}^{(j)}|\lesssim 1$, we have
\begin{align}\label{e:rp5}
|\bE_\cG[G^{(b_{\al})}_{c_\al c_\al}-Q]|\lesssim (d-1)^\ell/N.
\end{align}
Thus  \eqref{e:diag} follows from combining estimates \eqref{e:rep1}, \eqref{e:rp4} and \eqref{e:rp5}.

In the following we prove \eqref{e:offdiag}. Similar to \eqref{e:rep1}, we split the expectation into two terms, one corresponds to $\al, \beta\in \As_\bfS$ and one corresponds to that at least one of $\al, \beta$ is not in $ \As_\bfS$, which happens with probability $\OO(N^{-1+\fc})$. Then we have 
\begin{align}\label{e:rp6}\begin{split}
\bE_\cG[\tGT_{\ta_\al\ta_\beta}\bm1(\bfS\in F_2(\cG))]
&= \bE_\cG[(\bm1(\al,\beta\in {\As_\bfS})+\bm1(\al\not\in {\As_\bfS} \text{ or } \beta\not\in {\As_\bfS}))\tGT_{\ta_\al\ta_\beta}\bm1(\bfS\in F_2(\cG))]\\
&=
\bE_\cG[\bm1(\al,\beta\in {\As_\bfS})\tGT_{c_\al c_\beta}\bm1(\bfS\in F_2(\cG))]
+ \OO(\varepsilon'/N^{-1+\fc}).
\end{split}\end{align}%
Similar to \eqref{e:rp2} and \eqref{e:rp3},  for $\al,\beta\in \As_\bfS$ and $\bfS\in F_2(\cG)$, we have 
\begin{align}\begin{split}\label{e:rp7}
\tGT_{c_\al c_\beta}
&=G^{(b_\al b_\beta)}_{c_\al c_\beta}+\OO\left((d-1)^{\ell}(\varepsilon')^2\right)
=G_{c_\al c_\beta}-\sum_{x,y\in\{b_\al, b_\beta\}}G_{c_\al x}(G|_{\{b_\al, b_\beta\}})^{-1}_{xy}G_{y c_\beta}+\OO\left((d-1)^{\ell}(\varepsilon')^2\right)\\
&=G_{c_\al c_\beta}-\frac{ X_r(Q) G_{c_\al b_\al}G_{b_\al c_\beta}+ X_r(Q) G_{c_\al b_\beta}G_{b_\beta c_\beta}-G_{c_\al b_\al}G_{b_\al b_\beta}G_{b_\beta c_\beta}}{X_r(Q)^2}+\OO((d-1)^\ell (\varepsilon')^2),
\end{split}\end{align}%
where $X_r(Q)=X_r(Q(\cG,z),z)=\md(z)+\OO(\varepsilon/\sqrt{\kappa+\eta+\varepsilon})$ is as defined in \eqref{def:Y}, and we used that $|G_{b_\al b_\beta}|, |G_{b_\al c_\beta}|, |G_{c_\al c_\beta}|\lesssim \varepsilon$ and $G_{b_\al b_\al}, G_{b_\beta b_\beta}=X_r(Q)+\OO(\varepsilon)$.

Since $\cG\in \Omega$, we have that 
\begin{align}\begin{split}\label{e:rp75}
&G_{c_\al b_\al}=G_{c_\al b_\al}(\Ext(\cB_r(c_\al, b_\al, \cG), Q))+\OO(\varepsilon)=-\frac{\md(z)\msc(z)}{\sqrt{d-1}}+\OO\left(\frac{\varepsilon}{\sqrt{\kappa+\eta+\varepsilon}}\right),\\
&G_{b_\beta c_\beta}=G_{b_\beta c_\beta}(\Ext(\cB_r(b_\beta, c_\beta, \cG),  Q))+\OO(\varepsilon)
=-\frac{\md(z)\msc(z)}{\sqrt{d-1}}+\OO\left(\frac{\varepsilon}{\sqrt{\kappa+\eta+\varepsilon}}\right).
\end{split}\end{align}
Since $\bm1$ is an eigenvector of $G$ with eigenvalue $(d/\sqrt{d-1}-z)^{-1}$, we have
$\sum_{i}G_{ij}=\sum_jG_{ij}=(d/\sqrt{d-1}-z)^{-1}$. We recall that $(b_\al, c_\al), (b_\beta, c_\beta)$ are uniform randomly picked directed edges from $\cGT$. We have, 
\begin{align}\label{e:rp8}
\bE_\cG[G_{c_\al c_\beta}],\bE_{\cG}[G_{b_\al c_\beta}], \bE_\cG[G_{c_\al b_\beta}], \bE_\cG[G_{b_\beta c_\beta}] =\frac{1}{N(d/\sqrt{d-1}-z)}+\OO( (d-1)^\ell/N).
\end{align}

By plugging \eqref{e:rp75} and \eqref{e:rp8} into \eqref{e:rp7} and taking expectation, we get
\begin{align}\label{e:rp9}
\bE_\cG[\bm1(\al,\beta\in {\As_\bfS})\tGT_{c_\al c_\beta}\bm1(\bfS\in F_2(\cG))]
=\left(1+\frac{\msc(z)}{\sqrt{d-1}}\right)^2\frac{1}{N(d/\sqrt{d-1}-z)}+\OO((d-1)^\ell(\varepsilon')^2).
\end{align}%
The claim \eqref{e:offdiag} follows from combining estimates \eqref{e:rp6}and \eqref{e:rp9}.
\end{proof}

\subsubsection{Proof of \eqref{e:diffIIbound}}

We split \eqref{e:diffII} into two terms
\begin{align}
\eqref{e:diffII}=
&\bE\left[\chi_o(\cG)A_{oi}\Delta(G_{oo}^{(i)})(Q-Y_\ell(Q))^{p-1}\overline{(Q-Y_\ell(Q))^p}{\bm1(\cG,\tcG\in \Omega)})\right]  \label{e:diffIII} \\
&+\bE\left[\chi_o(\cG)A_{oi}\Delta(G_{oo}^{(i)})(Q-Y_\ell(Q))^{p-1}\overline{(Q-Y_\ell(Q))^p}{\bm1(\cG\in\Omega,\tcG\not\in \Omega)})\right].\label{e:diffIIII}
\end{align}
We recall the set $F_2(\cG)$ from Section \ref{sec:stability}. If $\chi_o(\cG)=1$, Remark \ref{r:radiusR} implies vertex $o$ has radius $\fR/4$ tree neighborhood in both $\cG$ and $\tcG$. Proposition \ref{prop:tGweakstab} implies for $\cG\in \Omega$ and $\bfS\in F_2(\cG)$, $|\Delta(G_{oo}^{(i)})|\lesssim \varepsilon'$. For $\bfS\notin F_2(\cG)$, we have the trivial bound $|\Delta(G_{oo}^{(i)})|\lesssim 1/\eta\lesssim N$. And this happens with small probability $\bP_\cG(\bfS\not\in F_2(\cG))\lesssim N^{-\fd}$. We can further split \eqref{e:diffIIII} as
\begin{align}\begin{split}
\eqref{e:diffIIII}
&\lesssim
\bE\left[A_{oi}\varepsilon'|Q-Y_\ell(Q)|^{2p-1}{\bm1(\cG\in\Omega,\tcG\not\in \Omega)})\right]+
N^{-\fd+1}\bE\left[A_{oi}|Q-Y_\ell(Q)|^{2p-1}{\bm1(\cG\in\Omega)}\right]\\
&\lesssim \frac{\varepsilon'}{N}\bE\left[|Q-Y_\ell(Q)|^{2p-1}{\bm1(\cG\in\Omega,\tcG\not\in \Omega)})\right]+
N^{-\fd}\bE\left[|Q-Y_\ell(Q)|^{2p}{\bm1(\cG\in\Omega)}\right]^{1-1/2p},
\end{split}\label{e:hihi}
\end{align}%
where in the last line, we integrated out $A_{oi}$ which gives the factor $d/N$ and then used  Jensen's inequality. 

We recall the the set of radius $\fR$-tree like graphs $\bar\Omega$ from Definition \ref{def:barOmega}. We notice that for a graph $\cG\in \Omega\subset \bar\Omega$ as in Definition \ref{def:barOmega}, with $\chi_o(\cG)=1$, if $\cB_{2\fR}(\{b_\al, c_\al\},\cG)$ is a tree for all $\al\in \qq{1,\mu}$, and it holds $\dist_{\cG^{(\bT)}}(\{a_\al,b_\al,c_\al, a_\beta\}, \{b_\beta,c_\beta\})\geq 2\fR$ for all $\al\neq\beta\in \qq{1,\mu}$, then $\As_\bfS=\qq{1,\mu}$, $\chi_o(\tcG)=1$ and the excess of any radius $\fR$ neighborhood in $\tcG$ is at most $\omega_d$. Especially, in this case $\tcG\in \bar\Omega$ as in Definition \ref{def:barOmega}.
By a union bound, this event has probability at least $1-\OO(N^{2\ell+2\fR-1})=1-\OO(N^{-1+\fc})$,
which can be  rewritten as $\bP_\cG(\tcG\not\in\bar\Omega)\lesssim N^{-1+\fc}$. We can rewrite the first term on the righthand side of \eqref{e:hihi} as
\begin{align}\begin{split}\label{e:ffbb}
&\phantom{{}={}}\bE\left[|Q-Y_\ell(Q)|^{2p-1}{\bm1(\cG\in \Omega, \tcG\not\in \Omega)}\right]\\
&=\bE\left[|Q-Y_\ell(Q)|^{2p-1}({\bm1(\cG\in \Omega(z),\tcG\not\in \bar\Omega)}+{\bm1(\cG\in \Omega(z), \tcG\in \bar\Omega\setminus\Omega)})\right]\\
&\lesssim N^{-1+\fc}\bE\left[|Q-Y_\ell(Q)|^{2p-1}\bm1(\cG\in \Omega)\right]
+\bE\left[|Q-Y_\ell(Q)|^{2p}\bm1(\cG\in \Omega)\right]^{1-1/2p}\bE\left[\bm1( \tcG\in \bar\Omega\setminus\Omega)\right]^{1/2p}\\
&\lesssim \left(N^{-1+\fc}+N^{-(\fd-3)/2p}\right)\bE\left[|Q-Y_\ell(Q)|^{2p}\bm1(\cG\in \Omega)\right]^{1-1/2p},
\end{split}\end{align}%
where in the third line we used Holder's inequality; in the last line we used Jensen's inequality and $\bP(\bar\Omega\setminus \Omega)\lesssim N^{-\fd+3}$, which follows  from the proof of Proposition \ref{thm:mrQ}.
The estimates \eqref{e:hihi} and \eqref{e:ffbb} lead to the following upper bound for \eqref{e:diffIIII},
\begin{align}
\eqref{e:diffIIII}\lesssim \varepsilon' N^{-2+\fc}\bE\left[|Q-Y_\ell(Q)|^{2p}\bm1(\cG\in \Omega)\right]^{1-1/2p},
\label{e:upbb1}
\end{align}
provided we take that $\fd\geq 2p+3$.

Next we derive an upper bound for \eqref{e:diffIII}. For $\cG, \tcG\in \Omega$ and $\chi_o(\cG)=1$, we have $|\Delta(G_{oo}^{(i)})|\lesssim \varepsilon$. Moreover, from the discussion before \eqref{e:ffbb}, we have that $\bP_\cG(\chi_o(\tcG)=1)=1-\OO(N^{-1+\fc})$.
We can rewrite \eqref{e:diffIII} as
\begin{align}\begin{split}
\label{e:diff15}
&\phantom{{}={}}\bE\left[\chi_o(\cG)A_{oi}\Delta(G_{oo}^{(i)})(Q-Y_\ell(Q))^{p-1}\overline{(Q-Y_\ell(Q))^p}{\bm1(\cG,\tcG\in \Omega)})\right]\\
&=\bE\left[\chi_o(\cG)\chi_o(\tcG)A_{oi}\Delta(G_{oo}^{(i)})(Q-Y_\ell(Q))^{p-1}\overline{(Q-Y_\ell(Q))^p}{\bm1(\cG,\tcG\in \Omega)})\right]\\
&+\OO(\varepsilon N^{-1+\fc})\bE\left[A_{oi}|(Q-Y_\ell(Q)|^{2p-1}{\bm1(\cG\in \Omega)})\right].
\end{split}\end{align}
Thanks to Lemma \ref{l:exchange} that $(\cG,\tcG)$ forms an exchangeable pair, we can rewrite the first term on the righthand side of \eqref{e:diff15} as
\begin{align}\begin{split}
\label{e:diff2}
&\phantom{{}={}}\bE\left[\chi_o(\cG)\chi_o(\tcG)A_{oi}\Delta(G_{oo}^{(i)})(Q-Y_\ell(Q))^{p-1}\overline{(Q-Y_\ell(Q))^p}{\bm1(\cG,\tcG\in \Omega)})\right]\\
&=\frac{1}{2}\bE\left[\chi_o(\cG)\chi_o(\tcG)A_{oi}\Delta(G_{oo}^{(i)})\Delta\left((Q-Y_\ell(Q))^{p-1}\overline{(Q-Y_\ell(Q))^p}\right){\bm1(\cG,\tcG\in \Omega)}\right]\\
&\lesssim (\varepsilon/N)\bE \left[\left|\Delta\left((Q-Y_\ell(Q))^{p-1}\overline{(Q-Y_\ell(Q))^p}\right)\right|{\bm1(\cG,\tcG\in \Omega)}\right],
\end{split}\end{align}
where in the second line we used that $ \bE (\Delta f)g = \frac 1 2 \bE (\Delta f) \Delta g $, and in the last line we used that $|\Delta(G_{oo}^{(i)})|\lesssim \varepsilon$, and integrated out $A_{oi}$ which gives a factor $d/N$.

Remark \ref{r:subgraph} implies that 
$|(Q-Y_r(Q))-(\tQ-Y_r(\tQ))|\lesssim |Q-\tQ|$. 
For $\cG\in \Omega$ and $\bfS\in F_2(\cG)$,   we have by  Proposition \ref{l:IGchange} that 
\begin{align}
|Q-\tQ|\lesssim \frac{(d-1)^{2\ell}{ (\Im[\md(z)]+\varepsilon'+\varepsilon/\sqrt{\kappa+\eta+\varepsilon})}}{N\eta}+\frac{1}{N^{1-\fc}}.
\end{align}
For $\bfS\not\in F_2(\cG)$, we have the trivial bound $|Q-\tQ|\lesssim \varepsilon/\sqrt{\kappa+\eta+\varepsilon}$ and this happens with  probability $\bP_\cG(\bfS\not\in F_2(\cG))\lesssim N^{-\fd}$.
Thus we can bound \eqref{e:diff2} by 
\begin{align}\begin{split}\label{e:kk}
&\phantom{{}={}}\bE \left[ \left |\Delta\left(  (Q-Y_\ell(Q))^{p-1}\overline{(Q-Y_\ell(Q))^p}  \right) \right | {\bm1(\cG,\tcG\in \Omega)}  \right]\\
&\lesssim \sum_{r=1}^{2p-1} \bE\left[|Q-\tQ|^r(\bm1(\bfS \in F_2(\cG))+\bm1(\bfS \not\in F_2(\cG)))|Q-Y_\ell(Q)|^{2p-1-r}{\bm1(\cG,\tcG\in \Omega)}\right]\\
&\lesssim \sum_{r=1}^{2p-1} \left(\frac{(d-1)^{2\ell}{(\Im[\md(z)]+\varepsilon'+\varepsilon/\sqrt{\kappa+\eta+\varepsilon})}}{N\eta}+\frac{1}{N^{1-\fc}}\right)^r\bE\left[|Q-Y_\ell(Q)|^{2p-1-r}{\bm1(\cG\in \Omega)}\right]\\
&+\sum_{r=1}^{2p-1} \frac{1}{N^{\fd}}\left(\frac{\varepsilon}{\sqrt{\kappa+\eta+\varepsilon}}\right)^r\bE\left[|Q-Y_\ell(Q)|^{2p-1-r}{\bm1(\cG\in \Omega)}\right]\\
&\lesssim \sum_{r=1}^{2p-1} \left(\frac{(d-1)^{2\ell}{ (\Im[\md(z)]+\varepsilon'+\varepsilon/\sqrt{\kappa+\eta+\varepsilon})}}{N\eta}+\frac{1}{N^{1-\fc}}\right)^r\bE\left[|Q-Y_\ell(Q)|^{2p}{\bm1(\cG\in \Omega)}\right]^{1-(1+r)/2p},
\end{split}\end{align}%
where in the last inequality we assume that $\fd$ is large enough, i.e. $\fd\geq 2p$, and used Jensen's inequality.
The estimates \eqref{e:diff15}, \eqref{e:diff2} and \eqref{e:kk} lead to the following upper bound for \eqref{e:diffII},
\begin{align}\begin{split}
&\phantom{{}={}}\eqref{e:diffII}\lesssim \varepsilon N^{-2+\fc}\bE\left[|Q-Y_\ell(Q)|^{2p}\bm1(\cG\in \Omega)\right]^{1-1/2p}\\
&+\sum_{r=1}^{2p-1} \frac{\varepsilon}{N}\left(\frac{(d-1)^{2\ell}{ (\Im[\md(z)]+\varepsilon'+\varepsilon/\sqrt{\kappa+\eta+\varepsilon})}}{N\eta}+\frac{1}{N^{1-\fc}}\right)^r\\
&\times\bE\left[|Q-Y_\ell(Q)|^{2p}{\bm1(\cG\in \Omega)}\right]^{1-(1+r)/2p},
\end{split}\label{e:upbb2}
\end{align}
provided we take that $\fd\geq 2p$.
The claim \eqref{e:diffIIbound} follows from combining \eqref{e:upbb1} and \eqref{e:upbb2}.

\subsection{Proof of Theorem  \ref{thm:evloc} and \ref{thm:improvelocal}}\label{sec:proof}

\begin{proof}[Proof of Theorem \ref{thm:improvelocal}]

For any large $\fC>0$, small $\fe>0$, by taking $p\geq \fC/\fe$ in Proposition \ref{t:recursion}, Markov's inequality implies
\begin{align}\label{744}
\begin{split}
&\phantom{{}={}}|Q(\cG,z)-Y_\ell(Q(\cG,z),z)-\delta_Q|\bm1(\cG\in \Omega(z))\\
&\lesssim
N^{\fe}\left(\varepsilon^2+\varepsilon^{1/2}\left(\frac{{ \Im[\md(z)]+\varepsilon/\sqrt{\kappa+\eta}}}{N\eta}\right)^{1/2}\right).
\end{split}\end{align}
with probability $1-\OO(N^{-\fC})$.  Similar inequality holds if we replace the left side  by 
$|m_N(z)-X_\ell(Q(\cG,z),z)-\delta_m|\bm1(\cG\in \Omega(z))$.
Thanks to Proposition \ref{p:recurbound}, we can expand $Y_\ell(Q(\cG,z),z)$ around $m_{sc}(z)$ and get
\begin{align*}\begin{split}
&\delta_Q+\OO(N^{\fe})\left(\varepsilon^2+\varepsilon^{1/2}\left(\frac{{ \Im[\md(z)]+\varepsilon/\sqrt{\kappa+\eta}}}{N\eta}\right)^{1/2}\right)=(1-(m_{sc}(z))^{2\ell+2})(Q(\cG,z)-m_{sc}(z))\\
&-\msc^{2\ell+2}(z)\md(z)\left(\frac{1-\msc^{2\ell+2}(z)}{d-1}+\frac{d-2}{d-1}\frac{1-\msc^{2\ell+2}(z)}{1-\msc^2(z)}\right)(Q(\cG,z)-m_{sc}(z))^2.
\end{split}\end{align*}%
We take $\ell\in\qq{\fa \log_{d-1}\log N, 2\fa\log_{d-1}\log N}$ such that $|1+(\msc(z))^2+(\msc(z))^4+\cdots+(\msc(z))^{2\ell}|\gtrsim 1$.
A stability analysis of the previous equation, similar to the arguments given after  \eqref{e:expand},  implies that on $\Omega(z)$ with probability $1-\OO(N^{-\fC})$,
\begin{align*}
\left|Q(\cG,z)-\msc(z)-\frac{\delta_Q}{1-(m_{sc}(z))^{2\ell+2}}\right|\lesssim \frac{N^{\fe}}{\sqrt{\kappa+\eta}}\left(\varepsilon^2+\varepsilon^{1/2}\left(\frac{{ \Im[\md(z)]+\varepsilon/\sqrt{\kappa+\eta}}}{N\eta}\right)^{1/2}\right).
\end{align*}%
Rewriting $m_N(z)-\md(z)=(m_N(z)-X_\ell(Q(\cG,z),z))+(X_\ell(Q(\cG,z),z)-\md(z))$ and using 
Proposition \ref{p:recurbound} to estimate $X_\ell(Q(\cG,z),z)-\md(z)$,  we have that 
\begin{align}\begin{split}\label{e:plug}
&\phantom{{}={}}m_N(z)-\md(z)=(m_N(z)-X_\ell(Q(\cG,z),z))+(X_\ell(Q(\cG,z),z)-\md(z))\\
&= (m_N(z)-X_\ell(Q(\cG,z),z))+\frac{d}{d-1}(\md(z))^2(m_{sc}(z))^{2\ell}(Q(\cG,z)-m_{sc}(z))
+\OO\left(\ell|Q(\cG,z)-m_{sc}(z)|^2\right)\\
&=\delta_m+\frac{d}{d-1}\frac{(\md(z))^2(m_{sc}(z))^{2\ell}}{1-(m_{sc}(z))^{2\ell+2}}\delta_Q+\frac{\OO(N^{\fe})}{\sqrt{\kappa+\eta}}\left(\varepsilon^2+\varepsilon^{1/2}\left(\frac{{ \Im[\md(z)]+\varepsilon/\sqrt{\kappa+\eta}}}{N\eta}\right)^{1/2}\right)
\end{split}\end{align}%
holds with probability $1-\OO(N^{-\fC})$ on $\Omega(z)$.
Noticing that for $z=d/\sqrt{d-1}$, $\msc(z)=-1/\sqrt{d-1}$ and $|\msc(z)+1/\sqrt{d-1}|\asymp |z-d/\sqrt{d-1}|$. Straightforward computation implies that 
\begin{align}\label{e:cancel}
\delta_m+\frac{d}{d-1}\frac{(\md(z))^2(m_{sc}(z))^{2\ell}}{1-(m_{sc}(z))^{2\ell+2}}\delta_Q
=\frac{1}{N(d/\sqrt{d-1}-z)}+\OO\left(\frac{(d-1)^\ell}{N\sqrt{\kappa+\eta}}\right).
\end{align}
Theorem \ref{thm:improvelocal} follows from plugging \eqref{e:cancel} into \eqref{e:plug}, and  recalling that $\bP(\bar\Omega\setminus \Omega(z))\lesssim N^{-\fd+3}$ from the proof of Proposition \ref{e:Q-msc}.

  \end{proof}

\begin{proof}[Proof of Theorems \ref{thm:evloc} and \ref{thm:extremeev}]

We take $z=2+\kappa+\ri\eta$, or $z=2-\kappa-\ri\eta$ with $\eta=(d-1)^{3r/2}/N$ and $ N^{2\fe}/(d-1)^{r/3}\leq\kappa\leq d$. Then $\varepsilon(z)\asymp (\log N)^{8\fa}/(d-1)^{r}$ and $\Im[\md(z)]\asymp \eta/\sqrt{\kappa+\eta}\ll 1/{N\eta}$.
In  this regime, Theorem  \ref{thm:improvelocal} implies that
\begin{align*} \begin{split}
\left|m_N(z) - \md(z)-\frac{1}{N(d/\sqrt{d-1}-z)}\right|
&\lesssim
   \frac{N^{\fe}}{\sqrt{\kappa}}\left(\varepsilon^2+\varepsilon^{1/2}\left(\frac{1}{N\sqrt{\kappa}}+\frac{\varepsilon}{N\eta\sqrt{\kappa}}\right)^{1/2}\right)\\
  & \lesssim\frac{N^{\fe}\varepsilon}{(N\eta)^{1/2}\kappa^{3/4}}\lesssim \frac{(\log N)^{8\fa}}{N^{\fe/2}}\frac{1}{(d-1)^{3r/2}}\ll \frac{1}{N\eta}.
\end{split}\end{align*}%
with probability $1-\OO(N^{-\fC})$ with respect to the uniform measure on $\bar\Omega$.
Using again, that $\Im[\md(z)]\asymp \eta/\sqrt{\kappa+\eta}\ll 1/{N\eta}$, 
 we have in this regime 
\begin{align*}
\Im\left[m_N(z) -\frac{1}{N(d/\sqrt{d-1}-z)}\right]\ll \frac{1}{N\eta}.
\end{align*}
A standard argument, as in \cite[Section 11]{MR3699468}, implies there are no eigenvalues on the interval $[\Re[z]-\eta, \Re[z]+\eta]$. 
By a union bound, we conclude that with probability $1-\OO(N^{-\fC})$ with respect to the uniform measure on $\bar\Omega'$, it holds that 
\begin{align}\label{e:extremebb}
\la_2,|\la_N|\leq 2+ N^{2\fe}/(d-1)^{r/3}. 
\end{align}
Since $\fe$ is arbitrarily small and $(d-1)^{r/3} \ge N^{\fr/3}$, this proves  Theorem \ref{thm:extremeev}. In fact, we can take $\fc<1$, and $\tau=\fc/32$, $\fa,\fe$ sufficiently small, then we can take the exponent $\Omega(1)$ in Theorem \ref{thm:extremeev} to be $0.01$.

Theorem \ref{thm:evloc} follows from the bound \eqref{e:extremebb} of the extremal eigenvalues with the local law, Theorems \ref{thm:mrmsc} and \ref{thm:improvelocal} by a standard argument, see \cite[Section 11]{MR3699468}.
\end{proof}

\appendix

\section{Properties of the Green's functions}
\label{app:Green}

Throughout this paper, we repeatedly use some (well-known) identities for Green's functions,
which we collect in this appendix.

\subsection{Resolvent identity}

The following well-known identity is referred as resolvent identity:
for two invertible matrices $A$ and $B$ of the same size, we have
\begin{equation} \label{e:resolv}
  A^{-1} - B^{-1} = A^{-1}(B-A)B^{-1}=B^{-1}(B-A)A^{-1}.
\end{equation}

\subsection{Schur complement formula}

Given an $N\times N$ matrix $M$ and an index set $\T \subset \qq{N}$, recall that we denote by
$M|_\T$ the $\T \times \T$-matrix obtained by restricting $M$ to $\T$,
and that by $M^{(\T)} = M|_{\bT^c}$, where $\bT^c=\qq{N}\setminus \T$, we denote the matrix obtained by removing
the rows and columns corresponding to indices in $\T$.
Thus, for any $\T \subset \qq{N}$,
any symmetric matrix $H$ can be written (up to rearrangement of indices) in the block form
\begin{equation*}
  H = \begin{bmatrix} A& B'\\ B &D  \end{bmatrix},
\end{equation*}
with $A=H|_{\T}$ and $D=H^{(\T)}$.
The Schur complement formula asserts that, for any $z\in \bC^+$,
\begin{equation} \label{e:Schur}
 G=(H-z)^{-1}= \begin{bmatrix}
   (A-B'\GT B)^{-1} & -(A-B'\GT B)^{-1}B'\GT \\
   -\GT B(A-B'\GT B)^{-1} & \GT+\GT B(A-B'\GT B)^{-1}B'\GT 
 \end{bmatrix},
\end{equation}
where $\GT=(D-z)^{-1}$.
Throughout the paper, we often use the following special cases of \eqref{e:Schur}:
\begin{align} \begin{split}\label{e:Schur1}
  G|_{\T} &= (A-B'\GT B)^{-1},\\
  G|_{\T^c}-G^{(\T)}&=G|_{\T^c\T}(G|_{\T})^{-1}G|_{\T\T^c},\\
  G|_{\T\T^c}&=-G|_{\T}B'\GT,
  \end{split}
\end{align}
as well as the special case
\begin{equation} \label{e:Schurixj}
G_{ij}^{(k)} = G_{ij}-\frac{G_{ik}G_{kj}}{G_{kk}}.
\end{equation}

\subsection{Ward identity}

For any symmetric $N\times N$ matrix $H$, its Green's function $G(z)=(H-z)^{-1}$ satisfies
the Ward identity
\begin{equation} \label{e:Ward}
  \sum_{j=1}^{N} |G_{ij}(z)|^2= \frac{\im G_{jj}(z)}{\eta},
\end{equation}
where $\eta=\Im [z]$. This identity provides a bound for the sum $\sum_{j=1}^{N} |G_{ij}(z)|^2$
in terms of the diagonal entries of the Green's function.
For an explanation why this algebraic identity has the interpretation of a Ward, see e.g.\ \cite[p.147]{MR3204347}.

\subsection{Covering map}
Given a graph $\cG$ (which might be infinite), with vertex set $\bG$ and Green's function $G$.
For any vertex $i$, the vector $(G_{i1}, G_{i2}, G_{i3}, \dots)\in \ell^2(\bG)$ is uniquely determined by the following relations:
\begin{align}\label{relation}
\begin{split}
1+zG_{ii}=\frac{1}{\sqrt{d-1}}\sum_{k:i\sim k}G_{ik}\\
zG_{ij}=\frac{1}{\sqrt{d-1}}\sum_{k:j\sim k}G_{ik}
\end{split}
\end{align}
where $i\sim k$ denotes that $i$ and $k$ are adjacent in $\cG$, i.e., that $A_{ik}=A_{ki}=1$.

\begin{lemma}\label{l:covering}
Given a covering $\pi:\tcG \to \cG$ of graphs,
denote the Green's function of $\tcG$ by $\tG$ and that of $\cG$ by $G$. Then
for all vertices $i,j$ in $\cG$, the Green's functions obey
\begin{align}\label{defg-bis}
G_{ij}=\sum_{y:\pi(y)=j}\tilde{G}_{xy},
\end{align}
where $x$ is any vertex of $\tcG$ with $\pi(x)=i$, provided the righthand side of \eqref{defg-bis} converges.
\end{lemma}

\begin{proof}
See \cite[Lemma B.1]{MR3962004}.
\end{proof}

\begin{proposition}\label{p:wd}
For any $d\geq 3$, $\n_d$  defined in \ref{e:boundPij} satisfies $\n_d\geq 1$.
\end{proposition}
\begin{proof}
By Definition \ref{d:defwd}, $\n_d\geq 1$ if 
for any connected graph $\cG$ with vertex set $\bG$ and deficit function $g$ such that  $\sum_{v\in \bG}g(v)+\text{excess}(\cG)\le 1$, \eqref{e:boundPij} and \eqref{e:boundPii} hold. For any $z\in \bC^+$, from \eqref{e:mdstieltjes}, we have
\begin{align*}
|m_{sc}(z)|\asymp |m_d(z)|\asymp 1,\quad \Im[m_{sc}(z)]\asymp\Im[m_{d}(z)]\geq0.
\end{align*}
If $\sum_{v\in \bG}g(v)+\text{excess}(\cG)=0$, then $\TE(\cG)$ is the infinite $d$-regular tree, \eqref{e:boundPij} and \eqref{e:boundPii} follow from \eqref{e:Gtreemkm}. If $\sum_{v\in \bG}g(v)=1$ and $\text{excess}(\cG)=0$, then $\TE(\cG)$ is the infinite $(d-1)$-ary tree, \eqref{e:boundPij} and \eqref{e:boundPii} follow from \eqref{e:Gtreemsc}. If $\sum_{v\in \bG}g(v)=0$ and $\text{excess}(\cG)=1$, then $\TE(\cG)$ is an infinite $d$-regular graph with exact one cycle. The infinite $d$-regular tree $\cX$ is the universal covering of $\cG$, $\pi: \cX\mapsto \cG$. The Green's function of $\cX$ is computed in \eqref{e:Gtreemkm}. Lemma \ref{l:covering} implies for any vertices $i,j\in \cG$, let $x$ be any vertex of $\cX$ with $\pi(x)=i$, we have
\begin{align}\label{defg-bis1}
G_{ij}=m_{d}(z)\sum_{y:\pi(y)=j}\left(-\frac{\msc(z)}{\sqrt{d-1}}\right)^{\dist_{\cX}(x,y)}.
\end{align}
Each vertex $y$ with $\pi(y)=j$ corresponds to a non-backtracking path from $i$ to $j$, and $\dist_{\cX}(x,y)$ is the length of the non-backtracking path. Since $\cG$ contains exact one cycle, we denote it as $\cal C$. The lengths of non-backtracking paths from $i$ to $j$ form arithmetic progressions with a common difference of $|\cal C|\geq 3$. Therefore, the righthand side \eqref{defg-bis1} converges, \eqref{e:boundPij} and \eqref{e:boundPii} follow.

\end{proof}

\section{Proof of Estimates in Section \ref{sc:ext}}
\label{a:ext}
\begin{proof}[Proof of Proposition \ref{p:recurbound}]
The proofs of \eqref{e:Xrecurbound} and \eqref{e:recurbound} are the same, we will only prove \eqref{e:recurbound}.
We denote $\cG=\cB_\ell(o, \cY)$, which is the truncated $(d-1)$-ary tree, and its vertex set $\bG$. Let $\mathbb I^\del=\sum_{v:\dist_\cG(o,v)=\ell} e_{vv} $. The normalized adjacency matrices of $\Ext(\cG, m_{sc}(z))$ and $\Ext(\cG, \Delta(z))$ are given by
\begin{align*}
H-\msc(z)\mathbb I^\del,\quad
H-\Delta(z)\mathbb I^\del,
\end{align*}
where $H$ is the normalized adjacency matrix of $\cG$. We denote $P$ the Green's function of $\Ext(\cG, m_{sc}(z))$, 
\begin{align*}
P=G(\Ext(\cG,m_{sc}(z),z)
=\left(H-z-\msc(z)\mathbb I^\del\right)^{-1}.
\end{align*}
The entries of $P$ are explicitly given in \eqref{e:Gtreemsc}.
We can compute the Green's function of $\Ext(\cG, \Delta(z))$ by a perturbation argument,
\begin{align}\begin{split}\label{e:expansionGP}
&\phantom{{}={}}G(\Ext(\cG,\Delta(z)), z)
=\left(H-z-\Delta(z)\mathbb I^\del\right)^{-1}=\left(H-z-\msc(z)\mathbb I^\del-(\Delta(z)-\msc(z))\mathbb I^\del\right)^{-1}\\
&=P+\sum_{k\geq 1}(\Delta(z)-m_{sc}(z))^kP\left(\mathbb I^\del P\right)^k.
\end{split}\end{align}
For $k\geq 2$ we will show the following two relations for $P$,
\begin{align}
   \label{e:Pboundary} &\left(P\mathbb I^\del P\mathbb I^\del P\right)_{oo}=\msc^{2\ell+2}\md\left(\frac{1-\msc^{2\ell+2}}{d-1}+\frac{d-2}{d-1}\frac{1-\msc^{2\ell+2}}{1-\msc^2}\right) ,\\
    \label{e:Ptotalsum}&\left(P(\mathbb I^\del P)^{k}\right)_{oo}\lesssim (C\ell)^{k-1}.
\end{align}
The claim \eqref{e:recurbound} follows from plugging \eqref{e:Pboundary} and \eqref{e:Ptotalsum} into \eqref{e:expansionGP}.

We recall that for $l,l'\in \bG$, ${\rm anc}(l,l')$ is the distance from the common ancestor of the vertices $l,l'$ to the root $o$ in $\bG$.
The relation \eqref{e:Pboundary} follows from explicit computation using \eqref{e:Gtreemsc} and \eqref{e:Gtreemsc2}
\begin{align*}
&\phantom{{}={}}\sum_{l, l'}P_{ol}P_{ll'}P_{l'o}=\sum_{l,l'}\msc^2 \left(-\frac{\msc}{\sqrt{d-1}}\right)^{2\ell}\md\left(1-\left(-\frac{\msc}{\sqrt{d-1}}\right)^{2+2{\rm anc}(l, l')}\right)\left(-\frac{\msc}{\sqrt{d-1}}\right)^{\dist_{\cG}(l,l')}\\
&= \frac{\msc^{2\ell+2}}{(d-1)^{2\ell}}\sum_{l,l'} \md\left(1-\left(\frac{\msc}{\sqrt{d-1}}\right)^{2+2{\rm anc}(l, l')}\right)\left(-\frac{\msc}{\sqrt{d-1}}\right)^{\dist_{\cG}(l,l')}\\
&=\msc^{2\ell+2}\md\left(1-\left(\frac{\msc}{\sqrt{d-1}}\right)^{2+2\ell}
+\sum_{r=1}^\ell \left(1-\left(\frac{\msc}{\sqrt{d-1}}\right)^{2+2(\ell-r)}\right)\left(\frac{\msc}{\sqrt{d-1}}\right)^{2r}(d-2)(d-1)^{r-1}
\right)\\
&=\msc^{2\ell+2}\md\left(\frac{1-\msc^{2\ell+2}}{d-1}+\frac{d-2}{d-1}\frac{1-\msc^{2\ell+2}}{1-\msc^2}\right),
\end{align*}
where the summation in the first line is over $l,l'$ such that $\dist_{\cG}(l,o)=\dist_{\cG}(l',o)=\ell$; in the second to last line we used that for a given $l$, there are $(d-2)(d-1)^{r-1}$ values of $l'$ such that $\dist_{\cG}(l, l')=2r, {\rm anc}(l, l')=\ell-r$ for $1\leq r\leq \ell$. When $l=l'$, it holds $\dist_{\cG}(l, l')=0, {\rm anc}(l, l')=\ell$.

For each $i,j\in \bH$, let $|P|_{ij}\deq |P_{ij}|$.
From the above computation and notice $|m_{sc}(z)|\leq 1$, for any $l'$ such that $\dist_{\cG}(l',o)=\ell$, we have
\begin{align*}
\sum_l |P_{ll'}|\lesssim \ell.
\end{align*}
where the summation in the first term is over $l$ such that $\dist_{\cG}(l,o)=\ell$. Then it follows 
\begin{align*}
\left(P(\mathbb I^\del P)^{k}\right)_{oo}
\lesssim \sum_{l_1, l_2,\cdots, l_k} |P_{ol_1}| |P_{l_1 l_2}|\cdots |P_{l_k o}|
\lesssim (C\ell)^{k-1}.
\end{align*}

%
\end{proof}

\begin{proof}[Proof of Proposition \ref{p:fixpoint}]
We denote $P=G(\Ext(\cG, m_{sc}(z)),z)$ and use the same notations as in the proof of Proposition \ref{p:recurbound}. Similar to the expansion \eqref{e:expansionGP},  we have
\begin{align}\begin{split}\label{e:expansionGP15}
&\phantom{{}={}}G_{oo}(\Ext(\cG, \Delta(z)),z)
=P_{oo}+\sum_{k\geq 1}(\Delta(z)-m_{sc}(z))^k\left(P\left(\sum_{x\in \bG}\frac{d-g(x)-\deg_{\cG}(x)}{d-1}e_{xx}P\right)^k\right)_{oo},
\end{split}\end{align}%
For $k=1$, we have
\begin{align}\begin{split}\label{e:ooterm}
&\phantom{{}={}}\left(P\sum_{x\in \bG}\frac{d-g(x)-\deg_{\cG}(x)}{d-1}e_{xx}P\right)_{oo}\\
&=\sum_{x\in \bG}\frac{d-g(x)-\deg_{\cG}(x)}{d-1}P_{ox}P_{xo}\\
&=\sum_{x\in \bG}\frac{d-g(x)-\deg_{\cG}(x)}{d-1}\frac{(m_{sc}(z))^{2\dist_{\cG}(o,x)+2}}{(d-1)^{\dist_{\cG}(o,x)}}\\
&=1+\sum_{x\in \bG}(d-g(x)-\deg_{\cG}(x))\frac{(m_{sc}(z))^{2\dist_{\cG}(o,x)+2}-1}{(d-1)^{\dist_{\cG}(o,x)+1}},
\end{split}\end{align}
where we used \eqref{e:Gtreemsccopy2} in the second equality; For the last equality, since $\cG$ has no cycles, $g(v)=1$ if $v=o$ otherwise $g(v)=0$, it is a combinatorial fact that
\begin{align*}
\sum_{x\in \bG}\frac{d-g(x)-\deg_{\cG}(x)}{(d-1)^{\dist_{\cG}(o,x)+1}}=1.
\end{align*}
More generally for any $k\geq 1$, we have the following upper bound,
\begin{align*}\begin{split}
&\phantom{{}={}}\left|\left(P\left(\sum_{x\in \bG}\frac{d-g(x)-\deg_{\cG}(x)}{d-1}e_{xx}P\right)^k\right)_{oo}\right|\\
&\lesssim \sum_{x_1\in \bG: \deg_\cG(x_1)<d-g(x_1)}\cdots\sum_{x_{k}\in \bG: \deg_\cG(x_{k})<d-g(x_{k})}|P_{ox_1}||P_{x_1x_2}|\cdots |P_{x_{k}o}|.
\end{split}\end{align*}
Since $P_{ij}$ satisfies $|P_{ij}|\lesssim (|m_{sc}(z)|/\sqrt{d-1})^{\dist_{\cG}(i,j)}$, we have
\begin{align}\begin{split}\label{e:sumPix2}
\sum_{x\in \bG: \deg_\cG(x)<d-g(x)}|P_{ix}||P_{jx}|
&\lesssim 
\sum_{x\in \bG: \deg_\cG(x)<d-g(x)}\left(\frac{|m_{sc}(z)|}{\sqrt{d-1}}\right)^{\dist_{\cG}(i,x)+\dist_{\cG}(j, x)}\\
&\lesssim 
(1+\dist_{\cG}(i,j))\left(\frac{|m_{sc}(z)|}{\sqrt{d-1}}\right)^{\dist_{\cG}(i,j)}.
\end{split}\end{align}
We can repeatedly use the above estimate \eqref{e:sumPix2} to sequentially sum over vertices $x_1,x_2,\cdots, x_{k}$, 
\begin{align}\begin{split}\label{e:PBBPbound}
&\phantom{{}={}}\sum_{x_1\in \bG: \deg_\cG(x_1)<d-g(x_1)}\cdots\sum_{x_{k-1}\in \bG: \deg_\cG(x_{k})<d-g(x_{k})}|P_{ox_1}||P_{x_1x_2}|\cdots |P_{x_{k}o}|\\
&\lesssim \sum_{x_2\in \bG: \deg_\cG(x_2)<d-g(x_2)}\cdots\sum_{x_{k}\in \bG: \deg_\cG(x_{k})<d-g(x_{k})}(1+\diam(\cG))|P_{ox_2}||P_{x_2x_3}|\cdots |P_{x_{k}o}|\\
&\lesssim\cdots
\lesssim  (1+\diam(\cG))^{k-1}|P_{oo}|
\lesssim
 (1+\diam(\cG))^{k-1}.
\end{split}\end{align}
By plugging \eqref{e:ooterm} and \eqref{e:PBBPbound} into \eqref{e:expansionGP15}, and notice $P_{oo}=m_{sc}(z)$, we get
\begin{align*}\begin{split}
&\phantom{{}={}}G_{oo}(\Ext(\cG,\Delta(z)),z)-\Delta(z)\\
&=(m_{sc}(z)-\Delta(z))
+(\Delta(z)-m_{sc}(z))\sum_{x\in \bG}\frac{(d-g(x)-\deg_\cG(x))P_{ox}P_{xo}}{d-1}
+\OO\left((1+\diam(\cG))|\Delta(z)-m_{sc}(z)|^2\right)
\\
&=\sum_{x\in \bG}(d-g(x)-\deg_{\cG}(x))\frac{(m_{sc}(z))^{2\dist_{\cG}(o,x)+2}-1}{(d-1)^{\dist_{\cG}(o,x)+1}}(\Delta(z)-m_{sc}(z))+\OO\left((1+\diam(\cG))|\Delta(z)-m_{sc}(z)|^2\right)\\
&=\OO\left((1+\diam(\cG)\left(\sqrt{\kappa+\eta}|\Delta(z)-m_{sc}(z)|+|\Delta(z)-m_{sc}(z)|^2\right)\right),
\end{split}\end{align*}%
where in the last line we used 
\begin{align*}
|(m_{sc}(z))^{2\dist_{\cG}(o,x)+2}-1|
\leq (1+\dist_{\cG}(o,x))|(m_{sc}(z))^2-1|\lesssim (1+\diam(\cG))\sqrt{\kappa+\eta}.
\end{align*}
\end{proof}

\begin{proof}[Proof of Proposition \ref{p:subgraph}]
We denote $P=G(\Ext(\cG, m_{sc}(z)),z)$ and use the same notations as in the proof of Proposition \ref{p:recurbound}. 
Thanks to the expansion \eqref{e:expansionGP15},  we have
\begin{align}\begin{split}\label{e:expansionGP2}
&\phantom{{}={}}G_{ij}(\Ext(\cG, \Delta(z)),z)
=P_{ij}+\sum_{k\geq 1}(\Delta(z)-m_{sc}(z))^k\left(P\left(\sum_{x\in \bG}\frac{d-g(x)-\deg_{\cG}(x)}{d-1}e_{xx}P\right)^k\right)_{ij},
\end{split}\end{align}%
By the same argument as in \eqref{e:PBBPbound}, we have the following upper bound,
\begin{align}\begin{split}\label{e:PBBPboundcopy}
\left|\left(P\left(\sum_{x\in \bG}\frac{d-g(x)-\deg_{\cG}(x)}{d-1}e_{xx}P\right)^k\right)_{ij}\right|\lesssim
 (1+\diam(\cG))^{k-1}(1+\dist_{\cG}(i,j))\left(\frac{|m_{sc}(z)|}{\sqrt{d-1}}\right)^{\dist_\cG(i,j)}.
\end{split}\end{align}%
The claim \eqref{e:smalldiff} follows by plugging \eqref{e:PBBPboundcopy} into \eqref{e:expansionGP2}, and using the estimate $\diam(\cG)|\Delta(z)-m_{sc}(z)|\ll1$. Since the righthand side of \eqref{e:smalldiff} is much smaller than $(|m_{sc}(z)|/\sqrt{d-1})^{\dist_{\cG}(i,j)}$, it follows that 
\begin{align}\begin{split}\label{e:PDelta}
&G_{ij}(\Ext(\cG,z,\Delta(z)),z)\lesssim \left(\frac{|m_{sc}|}{\sqrt{d-1}}\right)^{\dist_\cG(i,j)}, 
\\ &G_{ii}(\Ext(\cG,\Delta(z)),z)=G_{ii}(\Ext(\cG,m_{sc}(z)),z)+\oo(1). 
\end{split}\end{align}
The estimates \eqref{e:boundPijcopy} follows \eqref{e:PDelta}. If the sum of deficit function $\sum_{v\in \bG}g(v)$  plus the excess of $\cG$ is at most $\n_d$, then $|G_{ii}(\Ext(\cG,m_{sc}),z)|\asymp 1$, thus \eqref{e:boundPiicopy} holds.  The upper bounds \eqref{e:sumPixcopy} and \eqref{e:sumPix2copy} are the consequences of \eqref{e:boundPijcopy}.
\end{proof}

\begin{proof}[Proof of Proposition \ref{p:localization}]
If vertices $i,j$ are in different connected components of $\cH$, then $\dist_{\cG}(i,j)\geq 2r$ and $G_{ij}(\Ext(\cH,\Delta(z)),z)=0$. Thanks to Proposition \ref{p:subgraph}, we have
\begin{align*}
|G_{ij}(\Ext(\cG,\Delta(z)),z)|\lesssim \left(\frac{|m_{sc}(z)|}{\sqrt{d-1}}\right)^{2r},
\end{align*}
and the claim \eqref{e:compatibility} follows. 
In the following we assume $\cH$ is connected, and denote its vertex set by $\bH$. We can  also assume that in $\cG$, vertices  in $ \bH$ are not extensible. If this is not the case, we can add a weight $(d-g(v)-\deg_\cG(v))\Delta(z)/\sqrt{d-1}$ to any extensible vertex $v\in \bH$ to both $\cH$ and $\cG$. This does not affect their extension. Under this assumption, the extensible vertices of $\bH$ are $\{v\in\bH: \deg_{\cH}(v)<\deg_{\cG}(v)\}$.

We denote $\cG'$ the same graph as $\cG$ but with a different deficit function $g'$: $g'(v)=g(v)$ for $v\in \bH$ and $g'(v)=0$ for $v\in \bG\setminus \bH$. The claim \eqref{e:compatibility} follows from the next two estimates:
\begin{align}\label{e:GG'}
|G_{ij}(\Ext(\cG,\Delta(z)),z)-G_{ij}(\Ext(\cG', \Delta(z)),z)|\lesssim  \left(\frac{|m_{sc}(z)|}{\sqrt{d-1}}\right)^{2r},
\end{align}
and 
\begin{align}\begin{split}\label{e:G'H}
&\phantom{{}={}}\left|G_{ij}(\Ext(\cG',\Delta(z)),z)-G_{ij}(\Ext(\cH, \Delta(z)),z)\right|\\
&\lesssim (1+\diam(\cG))\left(\sqrt{\kappa+\eta}|\Delta(z)-m_{sc}(z)|+|\Delta(z)-m_{sc}(z)|^2\right)
+\left(\frac{|m_{sc}(z)|}{\sqrt{d-1}}\right)^{2r}.
\end{split}\end{align}
For estimate \eqref{e:GG'}, we use \eqref{e:formula} and the resolvent identity \eqref{e:resolv},
\begin{align*}
&\phantom{{}={}}G_{ij}(\Ext(\cG,\Delta(z)))-G_{ij}(\Ext(\cG', \Delta(z)))\\
&=\frac{\Delta(z)}{d-1}\sum_{v\in \bG\setminus \bH}g(v)G_{iv}(\Ext(\cG,\Delta(z)))G_{vj}(\Ext(\cG',\Delta(z))).
\end{align*}
Since both graphs $\cG$ and $\cG'$ satisfy the assumptions in Proposition \ref{p:subgraph}, \eqref{e:boundPijcopy} implies $|G_{iv}(
\Ext(\cG,\Delta(z)))|$, $|G_{vj}(\Ext(\cG', \Delta(z)))|\lesssim (|m_{sc}|/\sqrt{d-1})^r$, and the above expression is bounded by $\OO((|m_{sc}|/\sqrt{d-1})^{2r})$.

For the proof of \eqref{e:G'H}, We further enumerate $\del_E \bH$ as $\{e_1,e_2,\cdots,e_\mu\}$, i.e., $e_\al=\{l_\al, a_\al\}$ are the edges with one vertex $l_\al \in \bH$ and one $a_\al\in \bG\setminus\bH$. 
We denote the normalized adjacency matrices of $\Ext(\cG',\Delta(z))$ and $\Ext(\cH, \Delta(z))$ in the block  matrix forms as
\begin{align*}
 \left[\begin{array}{cc}
        H & B'\\
        B & D
       \end{array}
\right],\quad 
H-\frac{\Delta(z)}{d-1}\sum_{\al\in \qq{1,\mu}}e_{l_\al l_\al}
\end{align*}
where $H$ is the normalized adjacency matrix of $\cH$.
$B$ corresponds to the boundary edges $e_1, e_2,\cdots, e_\mu$. $D$ is the normalized adjacency matrix of $\Ext(\cG, \Delta(z))|_{\bG\setminus \bH}$. It follows from Proposition \ref{p:deficitbound}, 
the graph $\cG^{(\bH)}$ breaks down into many connected components, all components satisfy the assumptions in Proposition \ref{p:subgraph}, which implies $|(D-z)^{-1}_{xy}|\lesssim 1$; there exists a subset $\sA\subset \qq{1,\mu}$ with $|\sA|\lesssim \omega_d$, such that if $\al\not\in \sA$ the component of $\cG^{(\bH)}$ containing $a_\al$ doesn't contain $a_\beta$ for any $\beta\neq \al$. It gives that $(D-z)^{-1}_{a_\al a_\beta}=0$. Moreover, the component containing $a_\al$ satisfies the assumption in  Proposition \ref{p:fixpoint} with root vertex $a_\al$, and we have
\begin{align*}
(D-z)^{-1}_{a_\al a_\al}-\Delta(z)=\OO\left((1+\diam(\cG)\left(\sqrt{\kappa+\eta}|\Delta(z)-m_{sc}(z)|+|\Delta(z)-m_{sc}(z)|^2\right)\right).
\end{align*}

By the Schur complement formula \eqref{e:Schur1}, we have
\begin{align}\label{e:schur1}
G(\Ext(\cH,\Delta(z)), z)
=\left(H-z-\frac{\Delta(z)}{d-1}\sum_{\al\in \qq{1,\mu}}e_{l_\al l_\al}\right)^{-1},
\end{align}
and 
\begin{align}\begin{split}\label{e:schur2}
G(\Ext(\cG',\Delta(z)),z)|_{\bH}
&=\left(H-z-B'(D-z)^{-1}B\right)^{-1}\\
&=\left(H-z-\frac{1}{d-1}\sum_{\al,\beta\in \qq{1,\mu}}(D-z)^{-1}_{a_\al a_\beta}e_{l_\al l_\beta}\right)^{-1}.
\end{split}\end{align}
For simplification of notations, we write $P\deq G(\Ext(\cH,\Delta(z)),z)$ and $P'\deq G(\Ext(\cG',\Delta(z)),z)$.
We take the difference of \eqref{e:schur1} and \eqref{e:schur2}, and use the resolvent identity \eqref{e:resolv}
\begin{align*}\begin{split}
&\phantom{{}={}}P_{ij}-P'_{ij}
=-\frac{1}{d-1}\sum_{\al,\beta\in \qq{1,\mu}}P_{il_\al}\left((D-z)^{-1}_{a_\al a_\beta}-\delta_{\al\beta}\Delta(z)\right)P'_{l_\beta j}\\
&=-\frac{1}{d-1}\sum_{\al\neq \beta\in \sA}P_{il_\al}P'_{l_\beta j}(D-z)_{a_\al b_\beta}^{-1}-\frac{1}{d-1}\sum_{\al\not\in \sA}P_{il_\al}P'_{l_\al j}\left((D-z)^{-1}_{a_\al a_\al}-\Delta(z)\right)\\
&\lesssim \sum_{\al\neq\beta\in \sA}|P_{il_\al}||P'_{l_\beta j}|+\sum_{\al\not\in \sA}|P_{il_\al}||P'_{l_\al j}|(1+\diam(\cG))\left(\sqrt{\kappa+\eta}|\Delta(z)-m_{sc}(z)|+|\Delta(z)-m_{sc}(z)|^2\right)\\
&\lesssim\left(\frac{|m_{sc}(z)|}{\sqrt{d-1}}\right)^{2r}+(1+\diam(\cG))\left(\sqrt{\kappa+\eta}|\Delta(z)-m_{sc}(z)|+|\Delta(z)-m_{sc}(z)|^2\right),
\end{split}\end{align*}%
where in the last inequality we used \eqref{e:sumPixcopy} and \eqref{e:sumPix2copy} in Proposition \ref{p:subgraph}. This finishes the proof of \eqref{e:G'H} and Proposition \ref{p:localization}.

\end{proof}

\printindex

\bibliography{all}
\bibliographystyle{abbrv}

\end{document}

%% file: TE.pspdftex
\begin{picture}(0,0)%
\includegraphics{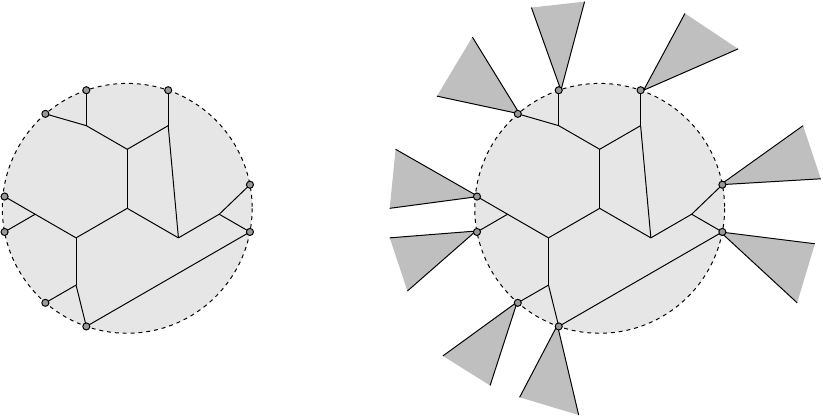}%
\end{picture}%
\setlength{\unitlength}{4144sp}%
\begingroup\makeatletter\ifx\SetFigFont\undefined%
\gdef\SetFigFont#1#2#3#4#5{%
  \reset@font\fontsize{#1}{#2pt}%
  \fontfamily{#3}\fontseries{#4}\fontshape{#5}%
  \selectfont}%
\fi\endgroup%
\begin{picture}(6267,3174)(1666,-5698)
\put(6841,-5461){\makebox(0,0)[b]{\smash{{\SetFigFont{11}{13.2}{\familydefault}{\mddefault}{\updefault}{\color[rgb]{0,0,0}$\TE(\cG_0)$}%
}}}}
\put(2611,-5461){\makebox(0,0)[b]{\smash{{\SetFigFont{11}{13.2}{\familydefault}{\mddefault}{\updefault}{\color[rgb]{0,0,0}$\cG_0$}%
}}}}
\end{picture}%

%% file: switch.pspdftex
\begin{picture}(0,0)%
\includegraphics{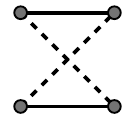}%
\end{picture}%
\setlength{\unitlength}{3947sp}%
\begingroup\makeatletter\ifx\SetFigFont\undefined%
\gdef\SetFigFont#1#2#3#4#5{%
  \reset@font\fontsize{#1}{#2pt}%
  \fontfamily{#3}\fontseries{#4}\fontshape{#5}%
  \selectfont}%
\fi\endgroup%
\begin{picture}(1080,995)(-1964,3315)
\put(-1949,4139){\makebox(0,0)[rb]{\smash{{\SetFigFont{11}{13.2}{\familydefault}{\mddefault}{\updefault}{\color[rgb]{0,0,0}$v_1$}%
}}}}
\put(-1949,3389){\makebox(0,0)[rb]{\smash{{\SetFigFont{11}{13.2}{\familydefault}{\mddefault}{\updefault}{\color[rgb]{0,0,0}$v_3$}%
}}}}
\put(-899,3389){\makebox(0,0)[lb]{\smash{{\SetFigFont{11}{13.2}{\familydefault}{\mddefault}{\updefault}{\color[rgb]{0,0,0}$v_4$}%
}}}}
\put(-899,4139){\makebox(0,0)[lb]{\smash{{\SetFigFont{11}{13.2}{\familydefault}{\mddefault}{\updefault}{\color[rgb]{0,0,0}$v_2$}%
}}}}
\end{picture}%